\documentclass[12pt,a4paper]{article}

\usepackage{latexsym}
\usepackage{amssymb}

\newtheorem{thm}{Theorem}[section]
\newtheorem{prop}{Proposition}[section]
\newtheorem{lem}{Lemma}[section]
\newtheorem{cor}{Corollary}[section]
\newtheorem{rmr}{Remark}[section]

\begin{document}
{



\begin{center}
{\Large\bf Generalized resolvents of symmetric and isometric operators: the Shtraus approach }
\end{center}

\begin{center}
{\large S.M. Zagorodnyuk }
\end{center}

\tableofcontents

\section{Introduction.}

This exposition paper is devoted to the theory of Abram Vilgelmovich Shtraus and his disciples and followers. This theory
studies the so-called generalized resolvents of symmetric and isometric operators in a Hilbert space.

The object, which now is called a generalized resolvent of a symmetric operator, appeared for the first time
in papers of Neumark and Krein in 40--th of the previous century.
At that time it was already known, due to the result of Carleman, that
each densely defined symmetric operator has a spectral function (the term "spectral function"$~$ is due to Neumark).
Neumark showed that each spectral function of a (densely defined) symmetric operator is generated by
the orthogonal spectral function of a self-adjoint extension of the given symmetric operator in a possibly larger Hilbert space
(Neumark's dilation theorem)~\cite{cit_12000_Neumark}.
In particular, this fact allowed him to describe all solutions of the Hamburger moment problem
in terms of spectral functions of the operator defined by the Jacobi matrix~\cite[pp. 303-305]{cit_12000_Neumark}.
However, as was stated by Neumark himself in his paper~\cite[p. 285]{cit_14000_Neumark},
despite of its theoretical generality this result did not give practical tools for finding spectral functions in various concrete
cases.

For operators with the deficiency index $(1,1)$ Neumark proposed a description of the generalized resolvents which was
convenient for practical applications.
As one such an application Nevanlinna's formula for all solutions of the Hamburger moment problem was derived~\cite[pp. 292-294]{cit_14000_Neumark}.
In 1944, independently from Neumark, another description of the generalized resolvents for operators with the deficiency index $(1,1)$
was given by Krein~\cite{cit_14000_Krein}.
Krein also noticed a possibility of an application of his results to a study of the moment problem and to the
problem of the interpolation of bounded analytic functions.
In 1946 appeared the famous paper of Krein which gave an analytic description of the generalized resolvents of a symmetric operator
with arbitrary equal finite defect numbers. Also there appeared a description of all $\Pi$-resolvents in the case that the operator
is non-negative
(a $\Pi$-resolvent --- a generalized resolvent generated by a non-negative self-adjoint extension of the given operator).
This paper of Krein became a starting point for a series of papers devoted to constructions of Krein-type formulas for various
classes of operators and relations in a Hilbert space, in the Pontriagin space, in the Krein space. Among other mathematicians
we may mention the following authors here (in the alphabetical order): 
Behrndt, Derkach, de Snoo, Hassi, Krein, Kreusler, Langer, Malamud, Mogilevskii, Ovcharenko, Sorjonen, see, e.g.:
\cite{cit_1400_Krein_Langer}, 
\cite{cit_1500_Krein_Langer}, 
\cite{cit_2500_Langer_Sorjonen}, 
\cite{cit_1000_Krein_Ovcharenko}, 
\cite{cit_9000_Derkach_Malamud}, 
\cite{cit_9000_Malamud}, 
\cite{cit_9500_Derkach_Malamud}, 
\cite{cit_8500_Derkach}, 
\cite{cit_30000_Malamud_Mogilevskii}, 
\cite{cit_1000_Mogilevskii}, 
\cite{cit_1000_Behrndt_Kreusler}, 
\cite{cit_9000_Derkach_Hassi_Malamud_de_Snoo}, 
and references therein.

Starting from the paper of Krein~\cite{cit_14000_Krein} the theory of the resolvent matrix has been developed intensively, see., e.g.,
papers of Krein and Saakyan~\cite{cit_2000_Krein_Saakyan}, 
Krein and Ovcharenko\cite{cit_1500_Krein_Ovcharenko}, 
Langer and Textorius~\cite{cit_1000_Langer_Textorius}, 
Derkach~\cite{cit_8000_Derkach} 
and papers cited therein.
Krein-type formulas are close to the linear fractional transformation and this fact suggested to apply them to interpolation problems,
where such form appears by the description of solutions.

On the other hand, in 1954 году, in his remarkable work~\cite{cit_900_Shtraus}, Shtraus obtained another analytic
description of the generalized resolvents of a (densely defined) symmetric operator with arbitrary defect numbers.
Characteristic features of Shtraus's formula, and of Shtraus-type formulas in general, are the following:
\begin{itemize}
\item[1)] An analytic function-parameter is {\it bounded};

\item[2)] The minimal number of all parameters in the formula: an analytic function-parameter and the given operator.
\end{itemize}
These features allow to solve various matrix interpolation problems, in the non-degenerate and degenerate cases.
Moreover, these features give possibility to solve not only one-dimensional, but also two-dimensional interpolation problems
and to obtain analytic descriptions of these problems.

The Shtraus formula set the beginning of a series of papers devoted to the derivation of such type formulas
for another classes of operators and their applications. The content of papers devoted to isometric and symmetric
operators in a Hilbert space is used in preparation of this survey.
Formal references and list of papers will be given below.
Besides that, we should mention papers for the generalized resolvents of operators
related to the conjugation, see., e.g., papers of Kalinina~\cite{cit_500_Kalinina},
Makarova~\cite{cit_100_Makarova}, \cite{cit_150_Makarova}, and references therein.
We should also mention papers on the generalized resolvents in spaces with an indefenite metric, see, e.g., papers
of Gluhov, Nikonov~\cite{cit_50_Gluhov_Nikonov}, \cite{cit_5_Nikonov},
\cite{cit_3_Nikonov}, Etkin~\cite{cit_10_Etkin}, and papers cited therein.

Notice, that Shtraus-type formulas can be used to obtain Krein-type formulas, see e.g.~\cite{cit_1000_Aleksandrov_Ilmushkin}.

Thus, in the theory of generalized resolvents one can see two directions: the derivation of Krein-type formulas and и
the derivation of Shtraus-type formulas.
The Krein formulas are partially described in books of Akhiezer and Glazman~\cite{cit_10_Akhiezer_Glazman},
M.L.~Gorbachuk and V.I.~Gorbachuk~\cite{cit_20_Gorbachuk}.
The Shtraus-type formulas can be found only in papers.
We think that it would be convenient to have an exposition of the corresponding theory with full proofs and necessary corrections.

We omit references in the following text using only known facts from classical books.
Formal references to all appeared results will be given afterwards, in the section "Formal credits".

We hope that this exposition paper will be useful to all those who want to study this remarkable theory.

\newpage


\noindent
\textbf{Notations.}

\begin{tabular}{p{3cm}p{10cm}}
$\mathbb{R}$ & the set of all real numbers \\
$\mathbb{C}$ & the set of all complex numbers \\
$\mathbb{N}$ & the set of all positive integers \\
$\mathbb{Z}$ & the set of all integers \\
$\mathbb{Z_+}$ & the set of all non-negative integers \\
$\mathbb{C}_+$  & $\{ z\in \mathbb{C}:\ \mathop{\rm Im}\nolimits z > 0 \}$ \\
$\mathbb{C}_-$  & $\{ z\in \mathbb{C}:\ \mathop{\rm Im}\nolimits z < 0 \}$ \\
$\mathbb{D}$  & $\{ z\in \mathbb{C}:\ |z|<1 \}$ \\
$\mathbb{T}$  & $\{ z\in \mathbb{C}:\ |z|=1 \}$ \\
$\mathbb{D}_e$  & $\{ z\in \mathbb{C}:\ |z|>1 \}$ \\
$\mathbb{T}_e$  & $\{ z\in \mathbb{C}:\ |z|\not=1 \}$ \\
$\mathbb{R}_e$  & $\{ z\in \mathbb{C}:\ \mathop{\rm Im}\nolimits z\not= 0 \}$ \\
$\mathbb{R}^d$  & the real Euclidean $d$-dimentional space, $d\in \mathbb{N}$ \\
$k\in\overline{0,\rho}$ & means that $k\in \mathbb{Z}_+$, $k\leq\rho$, if $\rho <\infty$;
or $k\in \mathbb{Z}_+$, if $\rho = \infty$ \\
$\mathbb{C}_{K\times N}$ & the set of all complex matrices of size $(K\times N)$, $K,N\in \mathbb{N}$ \\
$\mathbb{C}_{N\times N}^\geq$ & the set of all non-negative matrices from  $\mathbb{C}_{N\times N}$,
$N\in \mathbb{N}$ \\
$I_N$ & the identity matrix of size $(N\times N)$, $N\in \mathbb{N}$ \\
$C^T$ & the transpose of the matrix  $C\in (K\times N)$, $K,N\in \mathbb{N}$ \\
$C^*$ & the adjoint of the matrix  $C\in (K\times N)$, $K,N\in \mathbb{N}$ \\
$\Pi_{\lambda}$ & the half-plane $\mathbb{C}_+$ or $\mathbb{C}_-$, which contains a point $\lambda\in \mathbb{R}_e$ \\
$\Pi_{\lambda}^\varepsilon$ & $\{ z\in\Pi_{\lambda}:\ \varepsilon < |\arg z| < \pi - \varepsilon \},\qquad
0 < \varepsilon < \frac{\pi}{2}$, $\lambda\in \mathbb{R}_e$ \\
$\mathfrak{B}(M)$ & the set of all Borel subsets of a set $M$, which belongs $\mathbb{C}$ or $\mathbb{R}^d$, $d\in \mathbb{N}$ \\
$\mathbb{P}$ & the set of all scalar algebraic polynomials with complex coefficients 
\end{tabular}

\begin{tabular}{p{3cm}p{10cm}}
$\mathop{\rm card}\nolimits (M)$ & the quantity of elements of a set $M$ with a finite number of elements \\
\end{tabular}

\begin{tabular}{p{3cm}p{10cm}}
$\mathcal{S}(D;N,N')$ & a class of all analytic in a domain $D\subseteq \mathbb{C}$
operator-valued functions $F(z)$, which values are linear non-expanding operators mapping the whole
$N$ into $N'$, where $N$ and $N'$ are some Hilbert spaces. \\
\end{tabular}

All Hilbert spaces will be assumed to be separable, and operators in them are supposed to be linear.
If $H$ is a Hilbert space then

\begin{tabular}{p{3cm}p{10cm}}
$(\cdot,\cdot)_H$ & the scalar product in $H$ \\
$\| \cdot \|_H$ & the norm in $H$ \\
$\overline{M}$ & the closure of a set $M\subseteq H$ in the norm of $H$ \\
$\mathop{\rm Lin}\nolimits M$ & the linear span of a set $M\subseteq H$ \\
$\mathop{\rm span}\nolimits M$ & the closed linear span of a set $M\subseteq H$ \\
$E_H$ & the identity operator in $H$, i.e. $E_H x = x$, $\forall x\in H$ \\
$O_H$ & the null operator in $H$, i.e. $0_H x = 0$, $\forall x\in H$ \\
$o_H$ & the operator in $H$ with $D(o_H)=\{ 0 \}$, $o_H 0 = 0$ \\
$P_{H_1} = P^H_{H_1}$ & the operator of the orthogonal projection on a subspace $H_1$ in $H$ \\
\end{tabular}

Indices may be omitted in obvious cases.
If $A$ is a linear operator in $H$, then

\begin{tabular}{p{3cm}p{10cm}}
$D(A)$ & the domain of  $A$ \\
$R(A)$ & the range of $A$ \\
$\mathop{\rm Ker}\nolimits A$ & $\{ x\in H:\ Ax=0 \}$ (the kernel of $A$) \\
$A^*$ & the adjoint, if it exists \\
$A^{-1}$ & the inverse, if it exists \\
$\overline{A}$ & the closure, if $A$ admits the closure \\
$A|_M$ & the restriction of $A$ to the set $M\subseteq H$ \\
$\rho_r(A)$ & a set of all points of the regular type of $A$ \\
$\| A \|$ & the norm of $A$, if it is bounded \\
$w.-\lim$ & the limit in the sense of the weak operator topology \\
$s.-\lim$ & the limit in the sense of the strong operator topology \\
$u.-\lim$ & the limit in the sense of the uniform operator topology \\
\end{tabular}

If $A$ is a closed isometric operator then

\begin{tabular}{p{3cm}p{10cm}}
$M_\zeta = M_\zeta(A)$ & $(E_H - \zeta A) D(A)$, where $\zeta\in \mathbb{C}$  \\
$N_\zeta = N_\zeta(A)$ & $H\ominus M_\zeta$, where $\zeta\in \mathbb{C}$ \\
$M_\infty = M_\infty(A)$ & $R(A)$  \\
$N_\infty = N_\infty(A)$ & $H\ominus R(A)$  \\
$\mathcal{R}_\zeta = \mathcal{R}_\zeta(V)$ & $(E_H-\zeta V)^{-1}$,
$\zeta\in \mathbb{C}\backslash \mathbb{T}$ \\
\end{tabular}

If $A$ is a closed symmetric operator (not necessarily densely defined) then

\begin{tabular}{p{3cm}p{10cm}}
$\mathcal{M}_z = \mathcal{M}_z(A)$ & $(A - zE_H) D(A)$, where $z\in \mathbb{C}$  \\
$\mathcal{N}_z = \mathcal{N}_z(A)$ & $H\ominus \mathcal{M}_z$, where $z\in \mathbb{C}$ \\
$R_z = R_z(A)$ & $(A-zE_H)^{-1}$, $z\in \mathbb{C}\backslash \mathbb{R}$ \\
\end{tabular}


\section{Generalized resolvents of isometric and densely-defined symmetric operators.}\label{chapter1}

\subsection{Properties of a generalized resolvent of an isometric operator.}\label{section1_1}

Consider an arbitrary closed isometric operator $V$ in a Hilbert space $H$.
As it is well known, for $V$ there always exists an unitary extension $U$, which acts
in a Hilbert space  $\widetilde H\supseteq H$.
Define an operator-valued function $\mathbf{R}_\zeta$ in the following way:
$$ \mathbf{R}_\zeta = \mathbf{R}_\zeta(V) = \mathbf{R}_{u;\zeta}(V) =
P^{\widetilde{H}}_H \left( E_{\widetilde{H}} - \zeta U \right)^{-1}|_H,\qquad
\zeta\in \mathbb{T}_e. $$
The function $\mathbf{R}_\zeta$ is said to be  {\bf the generalized resolvent}
of an isometric operator $V$ (corresponding to the extension $U$).

Let $\{ F_t \}_{t\in [0,2\pi]}$ be a left-continuous orthogonal resolution of the unity of the operator $U$.
The operator-valued function
$$ \mathbf{F}_t = P^{\widetilde{H}}_H F_t,\qquad t\in [0,2\pi], $$
is said to be a (left-continuous) {\bf spectral function} of the isometric operator $V$
(corresponding to the extension $U$).
Let $F(\delta)$, $\delta\in\mathfrak{B}(\mathbb{T})$, be the orthogonal spectral measure of the unitary operator $U$.
The following function
$$ \mathbf{F}(\delta) = P^{\widetilde{H}}_H F(\delta),\qquad \delta\in \mathfrak{B}(\mathbb{T}), $$
is said to be {\bf a spectral measure} of the isometric operator $V$
(corresponding to the extension $U$).
Of course, the spectral function and the spectral measure are related by the following equality:
$$ \mathbf{F}(\delta_t) = \mathbf{F}_t,\qquad \delta_t=\{ z=e^{i\varphi}:\ 0\leq \varphi < t \},\quad
t\in [0,2\pi], $$
what follows from the analogous property of orthogonal spectral measures.
Moreover, generalized resolvents and spectral functions (measures) are connected by the following relation:
\begin{equation}
\label{f1_1_p1_1}
(\mathbf R_z h,g)_H = \int_{\mathbb{T}} \frac{1}{1-z\zeta} d(\mathbf{F}(\cdot) h,g)_H =
\int_0^{2\pi} \frac{1}{1-ze^{it}} d(\mathbf{F}_t h,g)_H,\quad \forall h,g\in H,
\end{equation}
which follows directly from their definitions.
This relation allows to talk about the one-to-one correspondence
between generalized resolvents and spectral measures, in the accordance with the well-known inversion formula for such integrals.

The generalized resolvents, as it is not surprizing from their definition, have much of the properties of the Fredholm resolvent
$( E_{\widetilde{H}} - \zeta U )^{-1}$ of the unitary operator $U$.
The aim of this subsection is to investigate these properties and to find out whether they are characteristic.

\begin{thm}
\label{t1_1_p1_1}
Let $V$ be a closed isometric operator in a Hilbert space $H$, and
$\mathbf{R}_\zeta$ be an arbitrary generalized resolvent of $V$. The following relations hold:

\begin{itemize}

\item[1)] $( z \mathbf{R}_z - \zeta \mathbf{R}_\zeta ) f = (z-\zeta) \mathbf{R}_z \mathbf{R}_\zeta f$,
for arbitrary $z,\zeta\in \mathbb{T}_e$, $f\in M_\zeta(V)$;

\item[2)] $\mathbf{R}_0 = E_H$;

\item[3)] For an arbitrary $h\in H$ the following inequalities hold:
$$ \mathop{\rm Re}\nolimits (\mathbf{R}_\zeta h,h)_H \geq \frac{1}{2} \| h \|_H^2,\qquad \zeta\in \mathbb{D}; $$
$$ \mathop{\rm Re}\nolimits (\mathbf{R}_\zeta h,h)_H \leq \frac{1}{2} \| h \|_H^2,\qquad \zeta\in \mathbb{D}_e; $$

\item[4)] $\mathbf{R}_\zeta$ is an analytic operator-valued function of a parameter $\zeta$
in $\mathbb{T}_e$;

\item[5)] For an arbitrary $\zeta\in \mathbb{T}_e\backslash\{ 0 \}$ it holds:
$$ \mathbf{R}_\zeta^* = E_H - \mathbf{R}_{\frac{1}{ \overline{\zeta} }}. $$
\end{itemize}
\end{thm}
{\bf Proof. } Let $U$ be a unitary operator in a Hilbert space
$\widetilde H\supseteq H$, which corresponds to $\mathbf{R}_\zeta(V)$.
For the Fredholm resolvent
$\mathcal{R}_\zeta = \mathcal{R}_\zeta(U) =( E_{\widetilde{H}} - \zeta U )^{-1}$ of the unitary operator $U$
the following relation holds
$$ z \mathcal{R}_z - \zeta \mathcal{R}_\zeta = (z-\zeta) \mathcal{R}_z \mathcal{R}_\zeta,\quad
z,\zeta\in \mathbb{T}_e, $$
which follows easily from the spectral decomposition of the unitary operator.
By applying the operator $P^{\widetilde H}_H$ to the latter relation, we get
$$ z \mathbf{R}_z - \zeta \mathbf{R}_\zeta = (z-\zeta) \mathbf{R}_z \mathcal{R}_\zeta,\quad
z,\zeta\in \mathbb{T}_e. $$
Taking into account that for an arbitrary element $f\in M_\zeta$, $f=( E_H - \zeta V ) g_V$, $g_V\in D(V)$,
holds
$$ \mathcal{R}_\zeta f = ( E_{\widetilde{H}} - \zeta U )^{-1} ( E_H - \zeta V ) g_V =
g_V = P^{\widetilde H}_H g_V $$
$$ = P^{\widetilde H}_H ( E_{\widetilde{H}} - \zeta U )^{-1} ( E_{\widetilde{H}} - \zeta V ) g_V
= P^{\widetilde H}_H ( E_{\widetilde{H}} - \zeta U )^{-1} f =
\mathbf{R}_\zeta f,   $$
we conclude that the first property of $\mathbf{R}_z$ is true.

\noindent
The second property follows directly from the definition of a generalized resolvent.

\noindent
Denote by $\{ F_t \}_{t\in[0,2\pi]}$ the left-continuous orthogonal resolution of the unity of the operator $U$.
For an arbitrary $h\in H$ holds
$$ \mathop{\rm Re}\nolimits (\mathbf{R}_\zeta h,h)_H =
\mathop{\rm Re}\nolimits (\mathcal{R}_\zeta h,h)_{\widetilde H} =
\mathop{\rm Re}\nolimits \int_0^{2\pi} \frac{1}{1-\zeta e^{it}} d(F_t h,h)_{\widetilde H} $$
$$ = \int_0^{2\pi} \frac{ 2-\overline{\zeta}e^{-it} - \zeta e^{it} }{2|1-\zeta e^{it}|^2}
d(F_t h,h)_{\widetilde H},\qquad
\zeta\in \mathbb{T}_e. $$
On the other hand, we may write
$$ (1-|\zeta|^2) (\mathcal{R}_\zeta^* h,\mathcal{R}_\zeta^* h)_{\widetilde H} + (h,h)_H =
\int_0^{2\pi} \frac{ 1-|\zeta|^2 }{|1-\zeta e^{it}|^2} d(F_t h,h)_{\widetilde H} + (h,h)_H $$
$$ = \int_0^{2\pi} \frac{ 2-\overline{\zeta}e^{-it} - \zeta e^{it} }{|1-\zeta e^{it}|^2}
d(F_t h,h)_{\widetilde H},\qquad
\zeta\in \mathbb{T}_e. $$
Consequently, we have
\begin{equation}
\label{f1_2_p1_1}
\mathop{\rm Re}\nolimits (\mathbf{R}_\zeta h,h)_H =
\frac{1}{2}
\left(
(1-|\zeta|^2) \| \mathcal{R}_\zeta^* h \|_{\widetilde H}^2 + \| h \|_H^2
\right),\qquad \zeta\in \mathbb{T}_e,
\end{equation}
and the third property of the generalized resolvent of an isometric operator follows.

\noindent
The fourth property follows directly from the same property of the resolvent of the unitary operator $U$.

\noindent
Using one more time the spectral resolution of the unitary operator $U$, we easily check that
$$ \mathcal{R}_\zeta^* + \mathcal{R}_{ \frac{1}{ \overline{\zeta} } } = E_{\widetilde H},\qquad
\zeta\in \mathbb{T}_e\backslash\{ 0 \}. $$
For arbitrary elements $f,g$ from $H$ we may write
$$ (\mathbf{R}_\zeta f,g)_H = (\mathcal{R}_\zeta f,g)_{\widetilde H}
= (f,\mathcal{R}_\zeta^* g)_{\widetilde H} $$
$$ = \left( f,
\left( E_{\widetilde H} - \mathcal{R}_{ \frac{1}{ \overline{\zeta} } } \right) g
\right)_{\widetilde H}
=
\left( f,
\left( E_{H} - \mathbf{R}_{ \frac{1}{ \overline{\zeta} } } \right) g
\right)_{H},\qquad
\zeta\in \mathbb{T}_e\backslash\{ 0 \}. $$
From the latter relation it follows the fifth property of the generalized resolvent.
$\Box$

\begin{thm}
\label{t1_2_p1_1}
Let an operator-valued function  $R_\zeta$ in a Hilbert space $H$ be given, which depends on a complex parameter
$\zeta\in \mathbb{T}_e$ and which values are linear operators defined on the whole $H$.
This function is a generalized resolvent of a
closed isometric operator in $H$ if an only if the following conditions are satisfied:

\begin{itemize}

\item[1)] There exists a number $\zeta_0\in \mathbb{D}\backslash\{ 0 \}$ and a subspace
$L\subseteq H$ such that
$$ ( \zeta R_\zeta - \zeta_0 R_{\zeta_0} ) f
= (\zeta-\zeta_0) R_\zeta R_{\zeta_0} f, $$
for arbitrary $\zeta\in \mathbb{T}_e$ and $f\in L$;

\item[2)] The operator $R_0$ is bounded and $R_0 h = h$, for all $h\in H\ominus \overline{ R_{\zeta_0}L }$;

\item[3)] For an arbitrary $h\in H$ the following inequality holds:
$$ \mathop{\rm Re}\nolimits (R_\zeta h,h)_H \geq \frac{1}{2} \| h \|_H^2,\qquad \zeta\in \mathbb{D}; $$

\item[4)] For an arbitrary $h\in H$
$R_\zeta h$ is an analytic vector-valued function of a parameter $\zeta$
in $\mathbb{D}$;

\item[5)] For an arbitrary $\zeta\in \mathbb{D}\backslash\{ 0 \}$ holds:
$$ R_\zeta^* = E_H - R_{\frac{1}{ \overline{\zeta} }}. $$
\end{itemize}
\end{thm}
{\bf Proof. }
The necessity of the properties 1)-5) follows immediately from the previous theorem.

Suppose that properties 1)-5) are true.
At first we check that properties~1),2) imply $R_0 = E_H$.
In fact, for an arbitrary element $h\in H$, $h=h_1 + h_2$, $h_1\in \overline{ R_{\zeta_0}L }$,
$h_2\in H\ominus \overline{ R_{\zeta_0}L }$, using~2) we may write
$$ R_0 h = R_0 h_1 + R_0 h_2 = R_0 h_1 + h_2. $$
There exists a sequence $h_{1,n}$ ($n\in \mathbb{N}$) of elements of $R_{\zeta_0}L$, tending to
$h_1$ as $n\rightarrow\infty$. Since $h_{1,n} = R_{\zeta_0} f_{1,n}$, $f_{1,n}\in L$
($n\in \mathbb{N}$), by condition~1) with $\zeta =0$ and $f=f_{1,n}$ we get
$$ -\zeta_0 R_{\zeta_0} f_{1,n} = -\zeta_0 R_0 R_{\zeta_0} f_{1,n},\qquad n\in \mathbb{N}. $$
Dividing by $-\zeta_0$ and passing to the limit as $n\rightarrow\infty$ we obtain that
$R_0 h_1 = h_1$. Here we used the fact that $R_0$ is continuous. Therefore we get $R_0 h = h$.

By property~4) the following function
$$ F(\zeta) := (R_\zeta h,h)_H - \frac{1}{2} \| h \|_H^2,\qquad \zeta\in \mathbb{D}, $$
is analytic and
$$ \mathop{\rm Im}\nolimits F(0) = \mathop{\rm Im}\nolimits \frac{1}{2} \| h \|_H^2 = 0. $$
From property~3) it follows that
$$ \mathop{\rm Re}\nolimits F(\zeta) =
\mathop{\rm Re}\nolimits (R_\zeta h,h)_H - \frac{1}{2} \| h \|_H^2 \geq 0,\qquad
\zeta\in \mathbb{D}. $$
This means that the function $F(\zeta)$ belongs to the Carath\'eodory class of all analytic in $\mathbb{D}$
functions satisfying the condition $\mathop{\rm Re}\nolimits F(\zeta) \geq 0$. Functions of this class  admit
the Riesz-Herglotz integral representation.
Using this representation we get
$$ (R_\zeta h,h)_H =
F(\zeta) + \frac{1}{2} \| h \|_H^2 =
\int_0^{2\pi} \frac{1}{1-\zeta e^{it}} d\sigma(t;h,h),\qquad \zeta\in \mathbb{D}, $$
where $\sigma(t;h,h)$ is a left-continuous non-decreasing function on the interval $[0,2\pi]$,
$\sigma(0;h,h)=0$. The function $\sigma(t;h,h)$ with such normalization is defined uniquely by the Riesz-Herglotz integral representation.
Since
$$ \int_0^{2\pi} \frac{1}{1-\zeta e^{it}} d\sigma(t;h,h) =
(R_\zeta h,h) = (h,R_\zeta^* h) = (h,h) - (h,R_{\frac{1}{ \overline{\zeta} }} h),\quad
\zeta\in \mathbb{D}\backslash\{ 0 \}, $$
where we used property~5), we get
$$ ( R_{\frac{1}{ \overline{\zeta} }} h, h) = (h,h) -
\int_0^{2\pi} \frac{1}{1- \overline{\zeta} e^{-it}} d\sigma(t;h,h) =
\int_0^{2\pi} \frac{1}{1- \frac{1}{\overline{\zeta}} e^{it}} d\sigma(t;h,h). $$
Therefore
$$ (R_\zeta h,h)_H =
\int_0^{2\pi} \frac{1}{1-\zeta e^{it}} d\sigma(t;h,h),\qquad \zeta\in \mathbb{T}_e,\ h\in H. $$
For arbitrary $h,g$ from $H$ we set
$$ \sigma(t;h,g) :=
\frac{1}{4} \sigma(t;h+g,h+g) - \frac{1}{4} \sigma(t;h-g,h-g)
+ \frac{i}{4} \sigma(t;h+ig,h+ig) $$
$$ - \frac{i}{4} \sigma(t;h-g,h-g). $$
Then
\begin{equation}
\label{f1_3_p1_1}
\int_0^{2\pi} \frac{1}{1-\zeta e^{it}} d\sigma(t;h,g) = (R_\zeta h,g)_H,\qquad
\zeta\in \mathbb{T}_e,\ h,g\in H.
\end{equation}
The function $\sigma(t;h,g)$ is a left-continuous complex-valued function of bounded variation with the normalization $\sigma(0;h,g)=0$.
If for the function $(R_\zeta h,g)_H$ there would exist a representation of a type~(\ref{f1_3_p1_1}) with another
left-continuous complex-valued function of bounded variation $\widetilde\sigma(t;h,g)$  with the normalization
$\widetilde\sigma(0;h,g)=0$, then
\begin{equation}
\label{f1_4_p1_1}
\int_0^{2\pi} \frac{1}{1-\zeta e^{it}} d\mu(t;h,g) = 0,\qquad \zeta\in \mathbb{T}_e,
\end{equation}
where $\mu(t;h,g) := \widetilde\sigma(t;h,g) - \sigma(t;h,g)$.
The function $\mu$ is defined uniquely by its trigonometric moments $c_k = \int_0^{2\pi} e^{ikt} d\mu(t;h,g)$, $k\in \mathbb{Z}$.
For $|\zeta|<1$, writing the expression under the integral sign in~(\ref{f1_4_p1_1}) as a sum of the geometric progression
and integrating we get $c_k = 0$, $k\in \mathbb{Z}_+$.

\noindent
For $|\zeta|>1$, the expression under the integral sign we can write as
$$ \frac{1}{1-\zeta e^{it}} = 1 - \frac{1}{1 - \frac{1}{\zeta} e^{-it}} = 1-
\sum_{k=0}^\infty \frac{1}{\zeta^k} e^{-ikt}, $$
and integrating we obtain that $c_{-k} = 0$, $k\in \mathbb{N}$.

\noindent
Thus, the representation~(\ref{f1_3_p1_1}) for $(R_\zeta h,g)_H$ defines the function $\sigma(t;h,g)$
with the above-mentioned properties uniquely.

\noindent
For an arbitrary $\zeta\in \mathbb{T}_e$ and $g,h\in H$, taking into account property~5), we may write
$$ \int_0^{2\pi} \frac{d\sigma(t;g,h)}{1-\zeta e^{it}} =
(R_\zeta g,h) = (g,R_\zeta^* h) = (g,h-R_{ \frac{1}{\overline{\zeta}} } h) $$
$$ = \overline{(h,g)} - \overline{(R_{ \frac{1}{\overline{\zeta}} } h,g)}
= \int_0^{2\pi} d \overline{\sigma(t;h,g)} -
\int_0^{2\pi} \frac{ d\overline{\sigma(t;h,g)} }{1- \frac{1}{\zeta} e^{-it}} $$
$$ = \int_0^{2\pi} \frac{ d\overline{\sigma(t;h,g)} }{1-\zeta e^{it}}. $$
By the uniqueness of the representation~(\ref{f1_3_p1_1}) it follows that
$$ \overline{\sigma(t;g,h)} = \sigma(t;h,g),\qquad h,g\in H,\ t\in [0,2\pi]. $$
From representation~(\ref{f1_3_p1_1}) and the linearity of the resolvent it follows that for arbitrary
 $\alpha_1,\alpha_2\in \mathbb{C}$ and  $h_1,h_2,g\in H$ holds
$$ \sigma(t;\alpha_1 h_1 + \alpha_2 h_2,g) = \alpha_1\sigma(t;h_1,g) + \alpha_2
\sigma(t;h_2,g),\qquad t\in [0,2\pi]. $$
Moreover, since the function $\sigma(t;h,h)$ is non-decreasing, the following estimate is true:
\begin{equation}
\label{f1_5_p1_1}
\sigma(t;h,h) \leq \sigma(2\pi;h,h) = \int_0^{2\pi} d\sigma(\tau;h,h) = (h,h)_H,\quad
t\in [0,2\pi],\ h\in H.
\end{equation}
This means that $\sigma(t;h,g)$ for each fixed $t$ from the interval $[0,2\pi]$
is a bounded bilinear functional in $H$, with the norm less or equal to $1$.
Consequently, it admits the following representation
$$ \sigma(t;h,g) = (E_t h,g)_H,\qquad t\in [0,2\pi], $$
where $\{ E_t \}_{t\in [0,2\pi]}$ is an operator-valued function of a parameter $t$,
which values are linear non-expanding operators, defined on the whole $H$ .
Let us study the properties of this operator-valued function.
Since
$$ (E_t h,g) = \sigma(t;h,g) = \overline{ \sigma(t;g,h) } = \overline{ (E_t g,h) } = (h,E_t g), $$
for arbitrary $h,g\in H$, the operators $E_t$ are all self-adjoint.
The function $(E_t h,h) = \sigma(t;h,h)$ is non-decreasing of $t$ on the segment $[0,2\pi]$,
for an arbitrary $h\in H$.
Since $\sigma(0;h,g)=0$, then $E_0 = E_H$, and from~(\ref{f1_5_p1_1}) it is seen that
$E_{2\pi} = E_H$.
Taking into account that for arbitrary $t\in (0,2\pi]$ and $h,g\in H$ holds
$$ \lim_{s\to t-0} (E_s h,g)_H = \lim_{s\to t-0} \sigma(s;h,g) = \sigma(t;h,g) = (E_t h,g)_H, $$
we conclude that $E_t$ is left-continuous on the segment $[0,2\pi]$ in the weak operator topology sense.
Also we have
$$ \| E_t h - E_s h \|^2_H = ((E_t-E_s)^2 h,h)_H \leq ((E_t-E_s) h,h)_H \rightarrow 0, $$
as $s\rightarrow t-0$. Here we used the fact that for a self-adjoint operator $(E_t-E_s)$ with the spectrum
in $[0,1]$ holds the inequality from the latter relation. This fact follows directly from
the spectral resolution of the self-adjoint operator $(E_t-E_s)$.
Therefore the function $E_t$ is left-continuous on the segment $[0,2\pi]$ in the strong operator topology sense, as well.

Consequently, $\{ E_t \}_{t\in[0,2\pi]}$ is a generalized resolution of the identity.
By the well-known Neumark's dilation theorem there exists an orthogonal resolution of unity $\{ \widetilde E_t \}_{t\in[0,2\pi]}$ in
a Hilbert space $\widetilde H\supseteq H$, such that
$$ E_t h = P^{\widetilde H}_H \widetilde E_t h,\qquad t\in [0,2\pi],\ h\in H. $$
From~(\ref{f1_3_p1_1}) and the definition of $E_t$ it follows that
\begin{equation}
\label{f1_6_p1_1}
R_\zeta = \int_0^{2\pi} \frac{1}{1-\zeta e^{it}} d E_t,\qquad
\zeta\in \mathbb{T}_e,\ h,g\in H.
\end{equation}
Here the convergence of the integral in the right-hand side is understood in the sense of the strong convergence of the integral sums.
The existence of the integral follows from the corresponding property of the orthogonal resolution of the identity
$\widetilde E_t$.

For arbitrary $\zeta\in \mathbb{T}_e$: $\zeta\not= \zeta_0$;
$f\in L$ and $g\in H$, we may write:
$$ \left( (\zeta R_\zeta - \zeta_0 R_{\zeta_0}) f,g \right)_H =
\zeta \int_0^{2\pi} \frac{d(E_t f,g)_H}{ 1 - \zeta e^{it} }
-
\zeta_0 \int_0^{2\pi} \frac{d(E_t f,g)_H}{ 1 - \zeta_0 e^{it} } $$
$$ = (\zeta - \zeta_0)
\int_0^{2\pi} \frac{1}{ 1 - \zeta e^{it} } \frac{1}{ 1 - \zeta_0 e^{it} } d(E_t f,g)_H $$
$$ = (\zeta - \zeta_0)
\int_0^{2\pi} \frac{1}{ 1 - \zeta e^{it} } d\int_0^t \frac{1}{ 1 - \zeta_0 e^{is} } d(E_s f,g)_H;
$$
$$ \left( (\zeta - \zeta_0) R_\zeta R_{\zeta_0} f,g \right)_H =
(\zeta - \zeta_0)
\int_0^{2\pi} \frac{1}{ 1 - \zeta e^{it} } d\left( E_t R_{\zeta_0} f,g \right)_H. $$
By the first condition of the theorem we conclude that
$$ \int_0^{2\pi} \frac{1}{ 1 - \zeta e^{it} } d\int_0^t \frac{1}{ 1 - \zeta_0 e^{is} } d(E_s f,g)_H
=
\int_0^{2\pi} \frac{1}{ 1 - \zeta e^{it} } d\left( E_t R_{\zeta_0} f,g \right)_H. $$
By the continuity from the left of the function $E_t$ and the initial condition $E_0 = 0$, the last relation inply that
\begin{equation}
\label{f1_7_p1_1}
\int_0^t \frac{1}{ 1 - \zeta_0 e^{is} } d(E_s f,g)_H =
\left( E_t R_{\zeta_0} f,g \right)_H,\qquad t\in [0,2\pi],\ f\in L,\ g\in H.
\end{equation}
Therefore
\begin{equation}
\label{f1_8_p1_1}
\int_0^t \frac{1}{ 1 - \zeta_0 e^{is} } dE_s f = E_t R_{\zeta_0} f,\qquad t\in [0,2\pi],\ f\in L.
\end{equation}
From condition~3) of the theorem it follows that the equality $R_\zeta h = 0$, for some $h\in H$ and $\zeta\in \mathbb{D}$,
implies the equality $h=0$. Thus, the operator $R_\zeta$ is invertible for all $\zeta\in \mathbb{D}$.

Define the following operator
$$ V g = \frac{1}{\zeta_0} \left( E_H - R_{\zeta_0}^{-1} \right) g,\qquad g\in R_{\zeta_0} L, $$
with the domain $D(V)=R_{\zeta_0} L$.
For an arbitrary element $g\in R_{\zeta_0} L$, $g=R_{\zeta_0}f$, $f\in L$, and an arbitrary $h\in H$,
we may write:
$$ (Vg,h)_H = \left( \frac{1}{\zeta_0} (R_{\zeta_0} - E_H) f,h \right)_H =
\frac{1}{\zeta_0}
\left( (R_{\zeta_0}f,h)_H  - (f,h)_H \right)_H $$
$$ =
\frac{1}{\zeta_0}
\left(
\int_0^{2\pi} \frac{1}{ 1 - \zeta_0 e^{it} } d( E_t f,h)_H
-
\int_0^{2\pi} d( E_t f,h)_H
\right) $$
$$ =
\int_0^{2\pi} \frac{e^{it}}{ 1 - \zeta_0 e^{it} } d( E_t f,h)_H
=
\int_0^{2\pi} e^{it} d\int_0^t \frac{1}{ 1 - \zeta_0 e^{is} } d( E_s f,h)_H $$
\begin{equation}
\label{f1_8_1_p1_1}
=
\int_0^{2\pi} e^{it} d( E_t g,h)_H,
\end{equation}
where we have used~(\ref{f1_7_p1_1}) in the last equality.
Therefore
$$ Vg = \int_0^{2\pi} e^{it} d E_t g,\qquad g\in D(V). $$
Using the latter equality we shall check that the operator $V$ is isometric. In fact, for arbitrary $g_1,g_2\in D(V)$ holds
\begin{equation}
\label{f1_9_p1_1}
(Vg_1, Vg_2)_H = \int_0^{2\pi} e^{it} d (E_t g_1, Vg_2)_H =
\int_0^{2\pi} e^{it} d (g_1, E_t V g_2)_H.
\end{equation}
On the other hand, for $g_2 = R_{\zeta_0} f_2$, $f_2\in L$, using the definition of the operator $V$ and
relation~(\ref{f1_8_p1_1}) we write
$$ E_t V g_2 = \frac{1}{\zeta_0} E_t (R_{\zeta_0} - E_t) f_2 =
\frac{1}{\zeta_0} \left(
E_t R_{\zeta_0} f_2 - E_t f_2
\right)
$$
$$ =
\frac{1}{\zeta_0} \left(
\int_0^t \frac{1}{ 1 - \zeta_0 e^{is} } dE_s f_2
-
\int_0^t dE_s f_2
\right) =
\int_0^t \frac{e^{is}}{ 1 - \zeta_0 e^{is} } dE_s f_2. $$
It follows that
$$ (g_1, E_t V g_2)_H = \overline{ (E_t V g_2, g_1)_H }
= \overline{
\int_0^t \frac{e^{is}}{ 1 - \zeta_0 e^{is} } d (E_s f_2, g_1)_H
} $$
$$ =
\overline{
\int_0^t e^{is} d\int_0^s \frac{1}{ 1 - \zeta_0 e^{ir} } d (E_r f_2, g_1)_H
}
=
\overline{
\int_0^t e^{is} d(E_s g_2, g_1)_H
} $$
$$ =
\int_0^t e^{-is} (E_s g_1, g_2)_H, $$
where we used~(\ref{f1_7_p1_1}) to obtain the last equality.
Substituting the obtained relation in~(\ref{f1_9_p1_1}), we get
$(Vg_1, Vg_2)_H = (g_1,g_2)_H$.
Thus, the operator $V$ is isometric.

\noindent
From the definition of the operator $V$ it follows that
$$ (E_H - \zeta_0 V) g = R_{\zeta_0}^{-1} g,\qquad g\in D(V)=R_{\zeta_0}L, $$
and therefore
$$ (E_H - \zeta_0 V) D(V) = L. $$
Consequently, the bounded operator $(E_H - \zeta_0 V)^{-1}$ is defined on the subspace $L$,
and therefore it is closed. Hence, $V$ is closwd, as well.

Consider a unitary operator $U$ in $\widetilde H$, which corresponds to the resolution
единицы $\{ \widetilde E_t \}_{t\in[0,2\pi]}$:
$$ U = \int_0^{2\pi} e^{it} d\widetilde E_t. $$
Let us show that $U\supseteq V$.
For arbitrary elements $g\in D(V)$ and $h\in H$ holds
$$ (Vg,h)_H = \int_0^{2\pi} e^{it} d(E_t g,h)_H =
\int_0^{2\pi} e^{it} d(\widetilde E_t g,h)_{\widetilde H} =
(U g,h)_{\widetilde H} = (P^{\widetilde H}_H U g,h)_H, $$
where we used relation~(\ref{f1_8_1_p1_1}). Therefore
\begin{equation}
\label{f1_10_p1_1}
Vg = P^{\widetilde H}_H U g,\qquad g\in D(V).
\end{equation}
Moreover, we have
$$ \| g \|_H = \| Vg \|_H = \| P^{\widetilde H}_H U g \|_H \leq
\| U g \|_{\widetilde H} = \| g \|_H. $$
Therefore $\| P^{\widetilde H}_H U g \|_H = \| U g \|_{\widetilde H}$, and
$$ P^{\widetilde H}_H U g = U g,\qquad g\in D(V). $$
Comparing the latter equality with~(\ref{f1_10_p1_1}) we conclude that $U\supseteq V$.

From~(\ref{f1_3_p1_1}) it follows that
$$ (R_\zeta h,g)_H =
\int_0^{2\pi} \frac{1}{1-\zeta e^{it}} d(E_t h,g)_H
=
\int_0^{2\pi} \frac{1}{1-\zeta e^{it}} d(\widetilde E_t h,g)_{\widetilde H} $$
$$ =
( (E_{\widetilde H} - \zeta U)^{-1} h,g )_{\widetilde H}
= ( P^{\widetilde H}_H (E_{\widetilde H} - \zeta U)^{-1} h,g )_H,\qquad
\zeta\in \mathbb{T}_e,\ h,g\in H. $$
Consequently,  $R_\zeta$, $\zeta\in \mathbb{T}_e$, is a generalized resolvent of a closed isometric
operator $V$.
$\Box$

Let us give conditions for the given operator-valued function $R_\zeta$ to be a generalized resolvent
of a concrete prescribed closed isometric operator.

\begin{thm}
\label{t1_3_p1_1}
Let an operator-valued function  $R_\zeta$ in a Hilbert space $H$ be given, which depends on a complex parameter
$\zeta\in \mathbb{T}_e$ and which values are linear operators defined on the whole $H$.
Let  a closed isometric operator $V$ in $H$ be given.
The function $R_\zeta$, $\zeta\in \mathbb{T}_e$, is a generalized resolvent of the operator $V$
if an only if the following conditions hold:

\begin{itemize}

\item[1)] For all $\zeta\in \mathbb{T}_e$ and for all $g\in D(V)$ the following equality holds:
$$ R_\zeta (E_H - \zeta V) g = g; $$

\item[2)] The operator $R_0$ is bounded and $R_0 h = h$, for all $h\in H\ominus D(V)$;

\item[3)] For an arbitrary $h\in H$ the following inequality holds:
$$ \mathop{\rm Re}\nolimits (R_\zeta h,h)_H \geq \frac{1}{2} \| h \|_H^2,\qquad \zeta\in \mathbb{D}; $$

\item[4)] For an arbitrary $h\in H$
$R_\zeta h$ is an analytic vector-valued function of a parameter $\zeta$
in $\mathbb{D}$;

\item[5)] For an arbitrary  $\zeta\in \mathbb{D}\backslash\{ 0 \}$ the following equality is true:
$$ R_\zeta^* = E_H - R_{\frac{1}{ \overline{\zeta} }}. $$
\end{itemize}
\end{thm}
{\bf Proof. }
The necessity of properties~3)-5) follows from the previous theorem.
Let $R_\zeta$ be a generalized resolvent of $V$, corresponding to a unitary operator $U\supseteq V$ in
a Hilbert space $\widetilde H\supseteq H$. Then
$$ R_\zeta (E_H - \zeta V) g = P^{\widetilde{H}}_H \left( E_{\widetilde{H}} - \zeta U \right)^{-1}
 (E_H - \zeta V) g $$
$$ =
P^{\widetilde{H}}_H \left( E_{\widetilde{H}} - \zeta U \right)^{-1}
 (E_H - \zeta U) g = g,\qquad g\in D(V). $$
Therefore condition~1) holds, as well. The second property follows from Theorem~\ref{t1_1_p1_1},
property~2).

Suppose that conditions 1)-5) from the statement of the theorem are satisfied.
Let us check that for the function $R_\zeta$ conditions 1),2) from the previous theorem are true.
Choose an arbitrary $\zeta_0\in \mathbb{D}\backslash\{ 0 \}$. Let
$L := (E_H - \zeta_0 V) D(V)$. By the first condition of the theorem we may write
$$ R_{\zeta_0} (E_H-\zeta_0 V) g = (E_H-\zeta_0 V)^{-1} (E_H-\zeta_0 V) g,\qquad g\in D(V), $$
and therefore
\begin{equation}
\label{f1_10_1_p1_1}
R_{\zeta_0} f = (E_H-\zeta_0 V)^{-1} f\in D(V),\qquad f\in L;
\end{equation}
$$ f = (E_H-\zeta_0 V) R_{\zeta_0} f,\qquad f\in L. $$
\begin{equation}
\label{f1_11_p1_1}
R_\zeta f = R_\zeta (E_H-\zeta_0 V) R_{\zeta_0} f,\qquad f\in L,\ \zeta\in \mathbb{T}_e.
\end{equation}
Let us use condition~1) of the theorem with the vector $g=R_{\zeta_0} f\in D(V)$:
$$ R_{\zeta_0} f = R_\zeta (E_H - \zeta V) R_{\zeta_0} f,\qquad f\in L. $$
Assuming that $\zeta\not=0$, we divide relation~(\ref{f1_11_p1_1}) by $\zeta_0$ and subtract the last relation divided by $\zeta$ from it.
Multiplying the obtained relation by $\zeta \zeta_0$,
we get the required equality from condition~1) of Theorem~\ref{t1_1_p1_1}.
In the case $\zeta=0$, the required equality takes the following form:
$R_{\zeta_0} f = R_0 R_{\zeta_0} f$, $f\in L$,
and it follows directly from the first condition of the theorem with $\zeta=0$, $g=R_{\zeta_0} f\in D(V)$.

\noindent
By relation~(\ref{f1_10_1_p1_1}) and the definition of $L$ it follows that
\begin{equation}
\label{f1_12_p1_1}
R_{\zeta_0} L = (E_H-\zeta_0 V)^{-1} L = D(V).
\end{equation}
Consequently, the second condition of the theorem~\ref{t1_1_p1_1} follows from the second condition of this theorem.
Observe that conditions 3)-5) of Theorem~\ref{t1_1_p1_1} coincide with the corresponding conditions
of this theorem. Thus, applying Theorem~\ref{t1_1_p1_1} we conclude that
the function $R_\zeta$, $\zeta\in \mathbb{T}_e$, is a generalized resolvent of a closed isometric operator in $H$.
Moreover, in the proof of Theorem~\ref{t1_1_p1_1}
such an operator  $V_1$ was constructed explicitly. Recall that its domain was
$R_{\zeta_0} L$, and
$$ V_1 g = \frac{1}{\zeta_0} (E_H - R_{\zeta_0}^{-1}) g,\qquad g\in R_{\zeta_0} L. $$
Let us check that $V=V_1$.
By~(\ref{f1_12_p1_1}) we see that the domains of the operators $V$ and $V_1$ coincide.
For an arbitrary $g\in D(V)=D(V_1)$, using condition~1) we write
$$ V_1 g = \frac{1}{\zeta_0} (E_H - R_{\zeta_0}^{-1}) g =
\frac{1}{\zeta_0} (E_H - (E_H - \zeta_0 V)) g = Vg. $$
Thus, $V_1=V$, and therefore the function $R_\zeta$, $\zeta\in \mathbb{T}_e$, is a generalized resolvent of the operator
$V$.
$\Box$

In conclusion of this subsection, we establish some propositions which will be used in what follows.
\begin{prop}
\label{p1_1_p1_1}
Let $F(\lambda)$ be an operator-valued analytic function in a domain $D\subseteq \mathbb{C}$,
which values are linear bounded operators in a Hilbert space $H$, defined
on the whole $H$.
Suppose that at each point of the domain $D$ there exists the inverse $F^{-1}(\lambda)$,
defined on the whole $H$. Assume that for each $\lambda\in D$ there exists
an open neighborhood of this point in which
$\| F^{-1}(\lambda) \|$ is bounded, as a function of $\lambda$.
Then the operator-valued function $F^{-1}(\lambda)$ is analytic in $D$, as well.
\end{prop}
{\bf Proof. } Consider an arbitrary point $\mu\in D$. We may write:
\begin{equation}
\label{f1_14_p1_1}
F^{-1}(\lambda) - F^{-1}(\mu) = - F^{-1}(\lambda) (F(\lambda)-F(\mu)) F^{-1}(\mu),\qquad
\lambda\in D,
\end{equation}
and therefore 
$$ \left\| F^{-1}(\lambda) - F^{-1}(\mu) \right\| \leq
\left\| F^{-1}(\lambda) \right\| \left\|
F(\lambda)-F(\mu) \right\|
\left\| F^{-1}(\mu) \right\| \rightarrow 0, $$
as $\lambda\rightarrow \mu$.
Thus, $F^{-1}(\lambda)$ is continuous in the domain $D$ in the uniform operator topology.

\noindent
From equality~(\ref{f1_14_p1_1}) it follows that
$$ \left\| \frac{ 1 }{\lambda - \mu} (F^{-1}(\lambda) - F^{-1}(\mu)) +
F^{-1}(\mu) F'(\mu) F^{-1}(\mu) \right\| $$
$$ \leq
\left\| -F^{-1}(\lambda) \left( \frac{1}{\lambda - \mu} (F(\lambda) - F(\mu)) \right)
+ F^{-1}(\mu) F'(\mu) \right\| \left\| F^{-1}(\mu) \right\| $$
$$ =
\left\| -F^{-1}(\lambda) \left( \frac{1}{\lambda - \mu} (F(\lambda) - F(\mu))
- F'(\mu)
\right)
+ (F^{-1}(\mu) - F^{-1}(\lambda)) F'(\mu) \right\| $$
$$ * \left\| F^{-1}(\mu) \right\| $$
$$ \leq \left[ \left\| F^{-1}(\lambda) \right\|
\left\|
\frac{1}{\lambda - \mu} (F(\lambda) - F(\mu))
- F'(\mu)
\right\|
+
\left\|
F^{-1}(\mu) - F^{-1}(\lambda)
\right\|
\left\|
F'(\mu)
\right\| \right] $$
$$ * \left\| F^{-1}(\mu) \right\| \rightarrow 0, $$
as $\lambda\rightarrow\mu$.
$\Box$

\begin{prop}
\label{p1_2_p1_1}
Let $W$ be a linear non-expanding operator in a Hilbert space $H$, and $f$ be an
arbitrary element from $D(W)$ such that:
$$ \| Wf \|_H = \| f \|_H. $$
Then
$$ ( Wf, Wg )_H = ( f, g )_H,\qquad \forall g\in D(W). $$
\end{prop}
{\bf Proof. }
Since $W$ is non-expanding, we may write
$$ \| \zeta f+g \|_H^2 - \| W ( \zeta f+g ) \|_H^2 \geq 0,\qquad g\in D(W),\ \zeta\in \mathbb{C}. $$
Then
$$ (\zeta f + g, \zeta f + g) - (\zeta Wf + Wg, \zeta Wf + Wg) =
|\zeta|^2 \| f \|^2 + 2 \mathop{\rm Re}\nolimits (\zeta (f,g))  $$
$$ + \| g \|^2 - |\zeta|^2 \| Wf \|^2 - 2 \mathop{\rm Re}\nolimits (\zeta (Wf,Wg))  - \| Wg \|^2 $$
$$ = \| g \|^2 - \| Wg \|^2 + 2 \mathop{\rm Re}\nolimits (\zeta [(f,g)-(Wf,Wg)]). $$
Thus, we obtain the inequality
$$ \| g \|^2 - \| Wg \|^2 + 2 \mathop{\rm Re}\nolimits (\zeta [(f,g)-(Wf,Wg)]) \geq 0,\quad
g\in D(W),\ \zeta\in \mathbb{C}. $$
Set $\zeta = -k \overline{((f,g) - (Wf,Wg))}$, $k\in \mathbb{N}$.
Then we get
$$ \| g \|^2 - \| Wg \|^2 - 2k |(f,g) - (Wf,Wg)|^2 \geq 0,\qquad k\in \mathbb{N}. $$
If $k$ tends to the infinity we get $( Wf, Wg )_H = ( f, g )_H$.
$\Box$

\begin{prop}
\label{p1_3_p1_1}
Let $W$ be a linear non-expanding operator in a Hilbert space $H$, and
$V$ be a closed isometric operator in $H$. Then the following conditions are equivalent:
\begin{itemize}
\item[(i)] $W\supseteq V$;

\item[(ii)] $W = V\oplus T$, where $T$ is a linear non-expanding operator in $H$, such that
$D(T) \subseteq N_0(V)$, $R(T) \subseteq N_\infty(V)$.
\end{itemize}
\end{prop}
{\bf Proof. }
(i)$\Rightarrow$(ii). Let $W\supseteq V$.
For an arbitrary element $h\in D(V)$ holds
$W h = Vh$, $\| Wh \| = \| Vh \| = \| h \|$.
Applying the previous proposition we obtain that
$$ (Wh,Wg) = (h,g),\qquad h\in D(V),\ g\in D(W). $$
Denote $M = D(W)\cap N_0(V) \subseteq N_0(V)$. Then
$$ (Wh,Wg) = 0,\qquad h\in D(V),\ g\in M. $$
Therefore $WM \subseteq N_\infty(V)$.
Set $T = W|_M$. The operator $T$ is a linear non-expanding operator and
$D(T) \subseteq N_0(V)$, $R(T) \subseteq N_\infty(V)$.

\noindent
Choose an arbitrary element $f\in D(W)$. It can be decomposed as a sum:
$f = f_1 + f_2$, where $f_1\in M_0(V) = D(V)$, $f_2\in N_0(V)$. Since $f$ and $f_1$
belong to $D(W)$, then $f_2$ belongs, as well. Therefore $f_2\in M$.
Thus, we obtain: $D(W) = D(V) \oplus D(T)$, where the manifold $D(T)$ is not supposed
to be closed.

\noindent
The implication (ii)$\Rightarrow$(i) is obvious.
$\Box$

\subsection{Chumakin's formula for the generalized resolvents of an isometric operator.}\label{section1_2}

The construction of the generalized resolvents using the definition of the generalized resolvent is hard.
It requires a construction of extensions of the given isometric operator in larger spaces.
Moreover, it may happen that different extensions produce the same generalized resolvent.
A convenient description of generalized resolvents in terms of an analytic class
of operator-valued functions is provided by the following theorem.

\begin{thm}
\label{t2_1_p1_1}
An arbitrary generalized resolvent $\mathbf{R}_\zeta$ of a closed isometric operator $V$,
acting in a Hilbert space $H$, has the following form:
\begin{equation}
\label{f2_1_p1_1}
\mathbf{R}_\zeta = \left[ E_H - \zeta (V\oplus F_\zeta) \right]^{-1},\qquad \zeta\in \mathbb{D},
\end{equation}
where $F_\zeta$ is a function from $\mathcal{S}(\mathbb{D};N_0(V),N_\infty(V))$.
Conversely, an arbitrary function $F_\zeta\in\mathcal{S}(\mathbb{D};N_0(V),N_\infty(V))$
defines by relation~(\ref{f2_1_p1_1}) a generalized resolvent $\mathbf{R}_\zeta$ of the operator $V$.
Moreover, to different functions from $\mathcal{S}(\mathbb{D};N_0(V),N_\infty(V))$ there correspond
different generalized resolvents of the operator $V$.
\end{thm}
{\bf Proof. }
Let $\mathbf{R}_\zeta$ be an arbitrary generalized resolvent of the operator $V$ from the statement of the theorem,
which corresponds to a unitary extension $U\supseteq V$ in a Hilbert space
$\widetilde H\supseteq H$.
For an arbitrary $\zeta\in \mathbb{T}_e$ and $h\in M_\zeta(V)$, $h=(E_H - \zeta V) f$, $f\in D(V)$,
we may write:
$$ \mathbf{R}_\zeta h = P^{\widetilde H}_H (E_{\widetilde H} - \zeta U)^{-1} h =
P^{\widetilde H}_H (E_{\widetilde H} - \zeta U)^{-1} (E_{\widetilde H} - \zeta U) f $$
$$ = f = (E_H - \zeta V)^{-1} h; $$
$$ \mathbf{R}_\zeta \supseteq (E_H - \zeta V)^{-1}. $$
By property 3) of Theorem~\ref{t1_1_p1_1}, the operator $\mathbf{R}_\zeta$ is invertible for $\zeta\in \mathbb{D}$.
Then
$$ \mathbf{R}_\zeta^{-1} \supseteq (E_H - \zeta V),\qquad \zeta\in \mathbb{D}. $$
Set
$$ T_\zeta = \frac{1}{\zeta} \left(
E_H - \mathbf{R}_\zeta^{-1}
\right),\qquad \zeta\in \mathbb{D}\backslash\{ 0 \}. $$
The operator $T_\zeta$ is an extension of the operator $V$.
Suppose that a vector $g\in H$ is orthogonal to the domain of the operator $T_\zeta$, i.e.
$$ 0 = (\mathbf{R}_\zeta h,g)_H,\qquad \forall h\in H. $$
Then
$\mathbf{R}_\zeta^* g = 0$.
Using property~3) of Theorem~\ref{t1_1_p1_1} we write
$$ 0 = (g,\mathbf{R}_\zeta^* g)_H = (\mathbf{R}_\zeta g,g)_H
= \mathop{\rm Re}\nolimits (\mathbf{R}_\zeta g,g)_H \geq
\frac{1}{2} (g,g)_H \geq 0. $$
Therefore $g=0$. Thus, we see that $\overline{ \mathbf{R}_\zeta H } = H$.

Using property~3) with arbitrary $h\in H$ we may write
$$ \| h \|_H^2 = (h,h)_H \leq 2 \mathop{\rm Re}\nolimits (\mathbf{R}_\zeta h,h)_H \leq
2 \left|
(\mathbf{R}_\zeta h,h)_H
\right|
\leq 2 \| \mathbf{R}_\zeta h \|_H \| h \|_H; $$
$$ \| h \|_H \leq 2 \| \mathbf{R}_\zeta h \|_H,\qquad \zeta\in \mathbb{D}. $$
Consequently, the operator $\mathbf{R}_\zeta^{-1}$ is bounded and
\begin{equation}
\label{f2_2_p1_1}
\| \mathbf{R}_\zeta^{-1} \| \leq 2,\qquad \zeta\in \mathbb{D}.
\end{equation}
Since $\mathbf{R}_\zeta$ is closed, the operator $\mathbf{R}_\zeta^{-1}$ is closed, as well. Hence, it is defined
on the whole space $H$.
Thus, we conclude that the operator $T_\zeta$ ($\zeta\in \mathbb{D}\backslash\{ 0 \}$) is defined on the whole
$H$.
Set
$$ B_\zeta = E_H - \mathbf{R}_\zeta^{-1},\qquad \zeta\in \mathbb{D}, $$
and we see that for $\zeta\in \mathbb{D}\backslash\{ 0 \}$ we have $B_\zeta = \zeta T_\zeta$.
Using this fact and the equality $\mathbf{R}_0 = E_H$, we obtain that $B_\zeta$ is defined on the whole
space $H$, for all $\zeta\in \mathbb{D}$.
For arbitrary $\zeta\in \mathbb{D}$ and $f\in D(B_\zeta)$, $f=\mathbf{R}_\zeta h$, $h\in H$, we write
$$ \| B_\zeta f \|^2_H = ( \mathbf{R}_\zeta h - h, \mathbf{R}_\zeta h - h )_H =
\| \mathbf{R}_\zeta h \|_H^2  - 2\mathop{\rm Re}\nolimits (\mathbf{R}_\zeta h,h)_H +
\| h \|_H^2. $$
Using one more time property~3) of Theorem~\ref{t1_1_p1_1} we conclude that
$$ \| B_\zeta f \|_H \leq \| f \|_H, $$
and therefore $\| B_\zeta \| \leq 1$.
Taking into account inequality~(\ref{f2_2_p1_1}) and using Proposition~\ref{p1_1_p1_1} we conclude that the function
$\mathbf{R}_\zeta^{-1}$ is analytic in $\mathbb{D}$.
Therefore the function $B_\zeta$ is also analytic in $\mathbb{D}$, and $B_0 = 0$.
By the Schwarz lemma for operator-valued functions the following inequality holds:
$\| B_\zeta \| \leq \zeta$, $\zeta\in \mathbb{D}$.
Therefore $\| T_\zeta \| \leq 1$, $\zeta\in \mathbb{D}\backslash\{ 0 \}$.

The function $B_\zeta$, as an analytic operator-valued function, can be expanded into the Maclaurin series:
$$ B_\zeta = \sum_{k=1}^\infty C_k \zeta^k,\qquad \zeta\in \mathbb{D}, $$
where $C_k$ are some linear bounded operators defined on the whole $H$.
Dividing this representation by $\zeta$ for $\zeta\not=0$, we see that the function $T_\zeta$ coincides with
a Maclaurin series. Defining $T_\zeta$ at zero by continuity, we get an analytic function
in $\mathbb{D}$.

Applying to $T_\zeta$ Proposition~\ref{p1_3_p1_1}, we conclude that
$$ T_\zeta = V \oplus F_\zeta,\qquad \zeta\in \mathbb{D}, $$
where $F_\zeta$ is a linear non-expanding operator, defined on the whole $N_0(V)$, with values in
$N_\infty(V)$. Since $T_\zeta$ is analytic in $\mathbb{D}$, then $F_\zeta$ is analytic, as well.
Therefore $F_\zeta\in\mathcal{S}(\mathbb{D};N_0(V),N_\infty(V))$.
By the definitions of the functions $T_\zeta$ and $B_\zeta$, from the latter relation we easily get
formula~(\ref{f2_1_p1_1}) with $\zeta\in \mathbb{D}\backslash\{ 0 \}$. For $\zeta = 0$
the validity of~(\ref{f2_1_p1_1}) is obvious.

Let us check the second statement of the theorem.
Consider an arbitrary function $F_\zeta\in\mathcal{S}(\mathbb{D};N_0(V),N_\infty(V))$.
Set
$$ T_\zeta = V\oplus F_\zeta,\qquad \zeta\in \mathbb{D}. $$
The function $T_\zeta$ is an analytic. Observe that there exists the inverse
$(E_H - \zeta T_\zeta)^{-1}$, which is bounded and defined on the whole $H$, since
$\| \zeta T_\zeta \| \leq |\zeta| < 1$. Moreover, we have
$$ \| (E_H - \zeta T_\zeta) h \|_H \geq \left|
\| h \|_H - |\zeta| \| T_\zeta h \|_H
\right| \geq (1-|\zeta|) \| h \|_H,\quad h\in H,\ \zeta\in \mathbb{D}. $$
Therefore
$$ \| (E_H - \zeta T_\zeta)^{-1} \| \leq \frac{1}{1-|\zeta|},\qquad \zeta\in \mathbb{D}. $$
By Proposition~\ref{p1_1_p1_1}, the function $(E_H - \zeta T_\zeta)^{-1}$ is analytic in
$\mathbb{D}$.
Set
$$ R_\zeta = (E_H - \zeta T_\zeta)^{-1},\qquad \zeta\in \mathbb{D}, $$
and
$$ R_\zeta = E_H - R^*_{\frac{1}{ \overline{\zeta} }},\qquad \zeta\in \mathbb{D}_e. $$
For the function $R_\zeta$ conditions~4) and~5) of Theorem~\ref{t1_3_p1_1} hold.
Let us check the validity of the rest of conditions of this Theorem.
Condition~1) for $\zeta\in \mathbb{D}$ follows from the definition of the function $R_\zeta$.
Choose arbitrary $\zeta\in \mathbb{D}_e$ and $g\in D(V)$. Then
$$ R_\zeta = E_H - R^*_{\frac{1}{ \overline{\zeta} }} =
E_H - \left( \left(E_H - \frac{1}{ \overline{\zeta} } T_{ \frac{1}{ \overline{\zeta} } } \right)^{-1} \right)^*
= E_H - \left( E_H - \frac{1}{\zeta} T^*_{\frac{1}{ \overline{\zeta} }} \right)^{-1} $$
$$ = -\frac{1}{\zeta} T^*_{\frac{1}{ \overline{\zeta} }}
\left( E_H - \frac{1}{\zeta} T^*_{\frac{1}{ \overline{\zeta} }} \right)^{-1}; $$
$$ R_\zeta (E_H - \zeta V) g = -\zeta R_\zeta \left( E_H - \frac{1}{\zeta} V^{-1} \right) Vg $$
$$ =
T^*_{\frac{1}{ \overline{\zeta} }}
\left( E_H - \frac{1}{\zeta} T^*_{\frac{1}{ \overline{\zeta} }} \right)^{-1}
\left( E_H - \frac{1}{\zeta} V^{-1} \right) Vg = g, $$
since $T^*_{\frac{1}{ \overline{\zeta} }} \supseteq V^{-1}$.
Thus, condition~1) of Theorem~\ref{t1_3_p1_1} is true.
Since $R_0 = E_H$, then condition~2) of Theorem~\ref{t1_3_p1_1} is true, as well.

\noindent
Choose arbitrary $h\in H$ and $\zeta\in \mathbb{D}$. Notice that
$$ \mathop{\rm Re}\nolimits (R_\zeta h,h) =
\frac{1}{2} \left(
(R_\zeta h,h) + (h,R_\zeta h)
\right) $$
$$ =
\frac{1}{2} \left\{
((E_H - \zeta T_\zeta)^{-1} h,h) + (h,(E_H - \zeta T_\zeta)^{-1} h)
\right\}. $$
Set $f = (E_H - \zeta T_\zeta)^{-1} h$. Then
$$ \mathop{\rm Re}\nolimits (R_\zeta h,h) =
\frac{1}{2} \left\{
( f, (E_H - \zeta T_\zeta)f ) + ((E_H - \zeta T_\zeta) f, f)
\right\} $$
$$ =
\frac{1}{2} \left\{
(f,f) - (f, \zeta T_\zeta f) + (f,f) - (\zeta T_\zeta f, f)
\right\}. $$
On the other hand, the following equality is true:
$$ \frac{1}{2} (h,h) = \frac{1}{2} ((E_H - \zeta T_\zeta)f, (E_H - \zeta T_\zeta)f) $$
$$ =
\frac{1}{2} \left\{
(f,f) - (f, \zeta T_\zeta f) - (\zeta T_\zeta f, f) + |\zeta|^2 \| T_\zeta f \|^2
\right\}. $$
Comparing last relations and taking into account the inequality $|\zeta|^2 \| T_\zeta f \|^2 \leq \| f \|^2$,
we get:
$\mathop{\rm Re}\nolimits (R_\zeta h,h) \geq \frac{1}{2} (h,h)$.
Consequently, condition~3) of Theorem~\ref{t1_3_p1_1} is true, as well. Therefore the function $R_\zeta$
is a generalized resolvent of the operator $V$.

Let us check the last assertion of the theorem.
Consider an arbitrary generalized resolvent $\mathbf{R}_\zeta$ of the operator $V$.
Suppose that $\mathbf{R}_\zeta$ admits two representations of a type~(\ref{f2_1_p1_1})
with some functions $F_\zeta$ and $\widetilde F_\zeta$. Then
$$ F_\zeta = \widetilde F_\zeta = \frac{1}{\zeta} (E_H - \mathbf{R}_\zeta^{-1})|_{N_0(V)},\qquad
\zeta\in \mathbb{D}\backslash\{ 0 \}. $$
By continuity these functions coincide at zero, as well.
$\Box$

Formula~(\ref{f2_1_p1_1}) is said to be the Chumakin formula for the generalized resolvents of a closed isometric operator.
This formula contains the minimal number of parameters: the operator and the function-parameter $F_\zeta$.
Due to this fact, it will be not hard to use the Chumakin formula by solving interpolation problems.

\subsection{Inin's formula for the generalized resolvents of an isometric operator.}\label{section1_3}

Consider a closed isometric operator $V$ in a Hilbert space $H$.
For the operator $V$ we shall obtain another description of the generalized resolvents --- Inin's formula,
see formula~(\ref{f1_7}) below. It has a less transparent structure, but instead of this it is more general
and coincides with Chumakin's formula in the case $z_0=0$. Inin's formula turns out to be very useful
in the investigation of the generalized resolvents of isometric operators with gaps in their spectrum
(by a gap we mean an open arc of the unit circle $\mathbb{D}$, which consists
of points of the regular type of an isometric operator).

An important role in the sequel will be played by the following operator:
\begin{equation}
\label{f1_2}
V_z = (V-\overline{z}E_H) (E_H - zV)^{-1},\qquad z\in \mathbb{D}.
\end{equation}
Notice that $D(V_z) = M_z$ and $R(V_z) = M_{\frac{1}{ \overline{z} }}$.
It is readily checked that the operator $V_z$ is isometric and
\begin{equation}
\label{f1_3}
V = (V_z + \overline{z}E_H) (E_H + zV_z)^{-1} = \left( V_z \right)_{-z}.
\end{equation}
Moreover, if $V$ is unitary, then $V_z$ is unitary, as well, and vice versa (this follows from~(\ref{f1_3})).

Let $\widehat V_z\supseteq V_z$ be a unitary extension of the operator $V_z$, acting in
a Hilbert space $\widehat H\supseteq H$.
We can define the following operator
\begin{equation}
\label{f1_4}
\widehat V = (\widehat V_z + \overline{z} E_{\widehat H}) (E_{\widehat H} + z\widehat V_z)^{-1},
\end{equation}
which will be a unitary extension of the operator $V$ in $\widehat H$.
Formula~(\ref{f1_4}) establishes a one-to-one correspondence between all unitary extensions $\widehat V_z$
of the operator $V_z$ in a Hilbert space $\widehat H\supseteq H$, and
all unitary extensions $\widehat V$ of the operator $V$ in $\widehat H$.

Fix an arbitrary point $z_0\in \mathbb{D}$.
Consider an arbitrary linear non-expanding operator $C$ with the domain
$D(C) = N_{z_0}$ and the range $R(C)\subseteq N_{\frac{1}{ \overline{z_0} }}$.
Denote
\begin{equation}
\label{f1_5}
V^+_{z_0;C} = V_{z_0} \oplus C;
\end{equation}
\begin{equation}
\label{f1_6}
V_{C}=V_{C;z_0} =  (V^+_{z_0;C} + \overline{z_0} E_{H}) (E_{H} + z_0 V^+_{z_0;C})^{-1}.
\end{equation}
If $z_0\not=0$, we may also write:
\begin{equation}
\label{f1_6_1}
V_{C}=V_{C;z_0} =  \frac{1}{z_0} E_H + \frac{|z_0|^2-1}{z_0} (E_{H} + z_0 V^+_{z_0;C})^{-1};
\end{equation}
\begin{equation}
\label{f1_6_2}
V^+_{z_0;C} =  -\frac{1}{z_0} E_H + \frac{ 1 - |z_0|^2}{z_0} (E_{H} - z_0 V_{C;z_0})^{-1}.
\end{equation}
The operator $V_C$ is said to be {\bf an orthogonal extension} of the closed isometric operator $V$,
defined by the operator $C$.

\begin{thm}
\label{t4_1_p1_1}
Let $V$ be a closed isometric operator in a Hilbert space $H$.
Fix an arbitrary point $z_0\in \mathbb{D}$.
An arbitrary generalized resolvent $\mathbf{R}_\zeta$ of the operator $V$ has the following representation:
\begin{equation}
\label{f1_7}
\mathbf R_{\zeta} = \left[
E - \zeta V_{C(\zeta;z_0)}
\right]^{-1},\qquad
\zeta\in \mathbb{D},
\end{equation}
where $C(\zeta) = C(\zeta;z_0)$ is a functionfrom $\mathcal{S}(\mathbb{D}; N_{z_0},N_{\frac{1}{ \overline{z_0} }})$.
Conversely, an arbitraryfunction $C(\zeta;z_0)\in \mathcal{S}(\mathbb{D};N_{z_0},N_{\frac{1}{ \overline{z_0} }})$
generates by relation~(\ref{f1_7}) a generalized resolvent
$\mathbf{R}_\zeta$ of the operator $V$.
To different functions from $\mathcal{S}(\mathbb{D}; N_{z_0}, N_{\frac{1}{ \overline{z_0} }})$ there correspond
different generalized resolvents of the operator $V$.
\end{thm}
{\bf Proof. }
Let $V$ be a closed isometric operator in a Hilbert space $H$, and
$z_0\in \mathbb{D}\backslash\{ 0 \}$ be a fixed point.
Consider the following linear fractional transformation:
\begin{equation}
\label{f2_1}
t = t(u) = \frac{u-\overline{z_0}}{1-z_0 u},
\end{equation}
which maps the unit circle $\mathbb{T}$ on $\mathbb{T}$, and $\mathbb{D}$ on $\mathbb{D}$.

\noindent
Let $\widehat V_{z_0}$ be an arbitrary unitary extension of the operator $V_{z_0}$, acting in
a Hilbert space $\widehat H\supseteq H$, and $\widehat V$ be the corresponding unitary extension
of the operator $V$, defined by relation~(\ref{f1_4}).
Choose an arbitrary number $u\in \mathbb{T}_e\backslash\{ 0, \overline{z_0}, \frac{1}{z_0} \}$.
Then
$t = t(u)\in \mathbb{T}_e\backslash\{ 0, -\overline{z_0}, -\frac{1}{z_0} \}$.
Moreover, the following conditions are equivalent:
\begin{equation}
\label{f2_1_1}
u\in \mathbb{T}_e\backslash\{ 0, \overline{z_0}, \frac{1}{z_0} \}
\Leftrightarrow
t \in \mathbb{T}_e\backslash\{ 0, -\overline{z_0}, -\frac{1}{z_0} \}.
\end{equation}
We may write:
$$ ( \widehat V_{z_0} - tE_{\widehat H} )^{-1} =
\left(
(\widehat V-\overline{z_0}E_{\widehat H} ) (E_{\widehat H} - z_0 \widehat V)^{-1}
 \right. $$
$$ \left.
- \frac{u-\overline{z_0}}{1-z_0 u} (E_{\widehat H} - z_0\widehat V) (E_{\widehat H} - z_0\widehat V)^{-1}
\right)^{-1}
$$
$$ =
\left(
\frac{1-|z_0|^2}{1-z_0 u}
(\widehat V - uE_{\widehat H})
(E_{\widehat H} - z_0\widehat V)^{-1}
\right)^{-1}.
$$
$$ = \frac{1-z_0 u}{1-|z_0|^2} (E_{\widehat H} - z_0\widehat V)
(\widehat V - uE_{\widehat H})^{-1} $$
$$ = - \frac{z_0(1-z_0 u)}{1-|z_0|^2} E_{\widehat H} +
\frac{(1-z_0 u)^2}{1-|z_0|^2} (\widehat V - uE_{\widehat H})^{-1}. $$
Therefore
$$ -\frac{1}{t} (E_{\widehat H} - \frac{1}{t}\widehat V_{z_0})^{-1}
=
- \frac{z_0(1-z_0 u)}{1-|z_0|^2} E_{\widehat H} -
\frac{(1-z_0 u)^2}{u(1-|z_0|^2)} (E_{\widehat H} - \frac{1}{u}\widehat V)^{-1}; $$
$$ (E_{\widehat H} - \frac{1}{u}\widehat V)^{-1} =
-\frac{z_0 u}{1-z_0 u} E_{\widehat H} + \frac{u(1-|z_0|^2)}{(1-z_0 u)^2 t}
(E_{\widehat H} - \frac{1}{t}\widehat V_{z_0})^{-1} $$
$$ = -\frac{z_0 u}{1-z_0 u} E_{\widehat H} + \frac{u(1-|z_0|^2)}{(1-z_0 u)(u-\overline{z_0})}
(E_{\widehat H} - \frac{1}{t}\widehat V_{z_0})^{-1}. $$
Set $\widetilde u = \frac{1}{u}$, $\widetilde t = \frac{1}{t}$. Observe that
$\widetilde u\in \mathbb{T}_e\backslash\{ 0, \frac{1}{\overline{z_0}}, z_0 \}$,
$\widetilde t\in \mathbb{T}_e\backslash\{ 0, -z_0, -\frac{1}{\overline{z_0}} \}$.
Moreover, we have
\begin{equation}
\label{f2_1_2}
\widetilde u\in \mathbb{T}_e\backslash\{ 0, \frac{1}{\overline{z_0}}, z_0 \}
\Leftrightarrow
\widetilde t\in \mathbb{T}_e\backslash\{ 0, -z_0, -\frac{1}{\overline{z_0}} \},
\end{equation}
and the latter conditions are equivalent to the conditions from relation~(\ref{f2_1_1}).
Then
$$ (E_{\widehat H} - \widetilde{u}\widehat V)^{-1} =
-\frac{z_0}{\widetilde{u}-z_0} E_{\widehat H} + \frac{\widetilde{u} (1-|z_0|^2)}{ (\widetilde{u} - z_0)
(1 - \overline{z_0} \widetilde{u}) }
(E_{\widehat H} - \widetilde{t}\widehat V_{z_0})^{-1}. $$
Applying the projection operator $P^{\widehat H}_H$ to the both sides of the last relation
we come to the following equality:
\begin{equation}
\label{f2_2}
\mathbf{R}_{\widetilde u}(V) = -\frac{z_0}{\widetilde{u}-z_0} E_{H}
+ \frac{\widetilde{u} (1-|z_0|^2)}{ (\widetilde{u} - z_0)(1 - \overline{z_0} \widetilde{u}) }
\mathbf{R}_{\frac{ \widetilde{u} - z_0 }{ 1 - \overline{z_0} \widetilde{u} }}(V_{z_0}),\quad
\widetilde u\in \mathbb{T}_e\backslash\{ 0, \frac{1}{\overline{z_0}}, z_0 \},
\end{equation}
where $\mathbf{R}_{\widetilde u}(V)$, $\mathbf{R}_{\widetilde t}(V_{z_0})$,
are the generalized resolvents of the operators $V$ and $V_{z_0}$, respectively.

\noindent
Since $\mathbf{R}_{\widetilde u}(V)$ is analytic in $\mathbb{T}_e$, it is uniquely defined by the generalized resolvent
$\mathbf{R}_{\widetilde{t}}(V_{z_0})$, according to relation~(\ref{f2_2}).
The same relation~(\ref{f2_2}) uniquely defines the generalized resolvent
$\mathbf{R}_{\widetilde{t}}(V_{z_0})$ by the generalized resolvent $\mathbf{R}_{\widetilde u}(V)$.

\noindent
Thus, relation~(\ref{f2_2})
establishes a one-to-one correspondence between all generalized resolvents of the operator $V_{z_0}$,
and all generalized resolvents of the operator $V$.

Let us apply Chumakin's formula~(\ref{f2_1_p1_1}) to the operator $V_{z_0}$:
\begin{equation}
\label{f2_3}
\mathbf{R}_{\widetilde t}(V_{z_0}) =
\left[
E_H - \widetilde t ( V_{z_0} \oplus F(\widetilde t) )
\right]^{-1},\qquad
\widetilde t\in \mathbb{D},
\end{equation}
where $F(\widetilde t)$ is a function of the class
$\mathcal{S}(\mathbb{D};N_{z_0},N_{\frac{1}{ \overline{z_0} }})$.

Consider relation~(\ref{f2_2}) for points $\widetilde u\in \mathbb{D}\backslash\{ 0,z_0 \}$,
what is equivalent to the condition $\widetilde t\in \mathbb{D}\backslash\{ 0, -z_0 \}$.
In this restricted case relation~(\ref{f2_2}) also establishes
a one-to-one correspondence between generalized resolvents.
Using~(\ref{f2_2}),(\ref{f2_3}) we get
$$ \mathbf{R}_{\widetilde u}(V) = -\frac{z_0}{\widetilde{u}-z_0} E_{H} $$
\begin{equation}
\label{f2_4}
+ \frac{\widetilde{u} (1-|z_0|^2)}{ (\widetilde{u} - z_0)(1 - \overline{z_0} \widetilde{u}) }
\left[
E_H - \frac{ \widetilde{u} - z_0 }{ 1 - \overline{z_0} \widetilde{u} }
\left( V_{z_0} \oplus
F\left(
\frac{ \widetilde{u} - z_0 }{ 1 - \overline{z_0} \widetilde{u} }
\right)
\right)
\right]^{-1},\
\widetilde u\in \mathbb{D}\backslash\{ 0,z_0 \},
\end{equation}
where $F(\widetilde t)\in\mathcal{S}(\mathbb{D};N_{z_0},N_{\frac{1}{ \overline{z_0} }})$.

\noindent
Relation~(\ref{f2_4}) establishes a one-to-one correspondence between all functions
$F(\widetilde t)$ from $\mathcal{S}(\mathbb{D};N_{z_0},N_{\frac{1}{ \overline{z_0} }})$, and
all generalized resolvents of the operator $V$.

\noindent
Set $C(\widetilde u) = F(
\frac{ \widetilde{u} - z_0 }{ 1 - \overline{z_0} \widetilde{u} })$, $u\in \mathbb{D}$.
Notice that $C(\widetilde u)\in\mathcal{S}(N_{z_0};N_{\frac{1}{ \overline{z_0} }})$.
We may write:
$$ E_H - \frac{ \widetilde{u} - z_0 }{ 1 - \overline{z_0} \widetilde{u} }
\left( V_{z_0} \oplus
F\left(
\frac{ \widetilde{u} - z_0 }{ 1 - \overline{z_0} \widetilde{u} }
\right)
\right) =
E_H - \frac{ \widetilde{u} - z_0 }{ 1 - \overline{z_0} \widetilde{u} }
V^+_{z_0;C(\widetilde u)} $$
$$ = (E_H - z_0 V_{C(\widetilde u);z_0}) (E_H - z_0 V_{C(\widetilde u);z_0})^{-1} $$
$$ - \frac{ \widetilde{u} - z_0 }{ 1 - \overline{z_0} \widetilde{u} }
(V_{C(\widetilde u);z_0} - \overline{z_0} E_H) (E_H - z_0 V_{C(\widetilde u);z_0})^{-1} $$
$$ = \frac{1-|z_0|^2}{ 1 - \overline{z_0}\widetilde u} (E_H - \widetilde u V_{C(\widetilde u);z_0})
(E_H - z_0 V_{C(\widetilde u);z_0})^{-1} $$
By substitution of the last relation into~(\ref{f2_4}), after elementary calculations we get
\begin{equation}
\label{f2_5}
\mathbf{R}_{\widetilde u}(V) =
(E_H - \widetilde u V_{C(\widetilde u);z_0})^{-1},\qquad
\widetilde u\in \mathbb{D}\backslash\{ 0,z_0 \}.
\end{equation}
Of course, in the case $\widetilde u = 0$ relation~(\ref{f2_5}) is true, as well.
It remains to check the validity of relation~(\ref{f2_5}) in the case $\widetilde u = z_0$.

\noindent
By Chumakin's formula for $\mathbf{R}_{\widetilde u}(V)$ we see that
$(\mathbf{R}_{\widetilde u}(V))^{-1}$ is an analytic operator-valued function in $\mathbb{D}$.
Using~(\ref{f2_5}) we may write:
\begin{equation}
\label{f2_5_1}
(\mathbf{R}_{z_0}(V))^{-1} = u.-\lim_{\widetilde u\to z_0} (\mathbf{R}_{\widetilde u}(V))^{-1}
= E_H - u.-\lim_{\widetilde u\to z_0} \widetilde u V_{C(\widetilde u);z_0},
\end{equation}
where the limits are understood in the sense of the uniform operator topology.

\noindent
The operator-valued function
$V^+_{z_0;C(\widetilde u)} = V_{z_0}\oplus C(\widetilde u)$ is analytic in $\mathbb{D}$,
and its values are non-expanding operators in $H$. Then
$$ \| (E_H + z_0 V^+_{z_0;C(\widetilde u)}) h \| \geq
| \| h \| - |z_0| \| V^+_{z_0;C(\widetilde u)}) h \| | \geq
(1 - |z_0|) \| h \|,\qquad h\in H; $$
\begin{equation}
\label{f2_6}
\| (E_H + z_0 V^+_{z_0;C(\widetilde u)})^{-1} \| \leq \frac{1}{1 - |z_0|},\qquad
\widetilde u\in \mathbb{D}.
\end{equation}
By Proposition~\ref{p1_1_p1_1} we obtain that the operator-valued function
$(E_H + z_0 V^+_{z_0;C(\widetilde u)})^{-1}$ is analytic in $\mathbb{D}$.
Therefore
$V_{C(\widetilde u);z_0} = (V^+_{z_0;C(\widetilde u)} + \overline{z_0} E_H)
(E_H + z_0 V^+_{z_0;C(\widetilde u)})^{-1}$
is analytic in $\mathbb{D}$, as well.
Passing to the limit in relation~(\ref{f2_5_1}) we get
$$ (\mathbf{R}_{z_0}(V))^{-1} =
E_H - z_0 V_{C(z_0);z_0}. $$
Thus, relation~(\ref{f2_5}) holds in the case $\widetilde u = z_0$, as well.
$\Box$

\subsection{The generalized resolvents of a symmetric operator. A connection with the generalized resolvents of the Cayley transformation
of the operator.}\label{section1_4}

Consider an arbitrary closed symmetric operator $A$ in a Hilbert space $H$.
The domain of $A$ is not supposed to be necessarily dense in $H$.
It is well known that for the operator $A$ there always exists a self-adjoint extension $\widetilde A$, acting
in a Hilbert space  $\widetilde H\supseteq H$.
Define an operator-valued function $\mathbf{R}_\lambda$ in the following way:
$$ \mathcal{\mathbf{R}}_\lambda = \mathbf{R}_\lambda(A) = \mathbf{R}_{s;\lambda}(A) =
P^{\widetilde{H}}_H \left( \widetilde A - \lambda E_{\widetilde{H}} \right)^{-1}|_H,\qquad
\lambda\in \mathbb{R}_e. $$
The function $\mathbf{R}_\lambda$ is said to be {\bf a generalized resolvent}
of the symmetric operator $A$ (corresponding to the extension $\widetilde A$).
The additional index $s$ is necessary in those case, where there can appear a muddle with the generalized resolvent of
an isometric operator.

Let $\{ \widetilde E_t \}_{t\in \mathbb{R}}$ be a left-continuous orthogonal resolution
of the identity of the operator $\widetilde A$. The operator-valued function
$$ \mathbf{E}_t = P^{\widetilde{H}}_H \widetilde E_t,\qquad t\in \mathbb{R}, $$
is said to be a (left-continuous) {\bf spectral function} of a symmetric operator $A$
(corresponding to the extension $\widetilde A$).
Let $\widetilde E(\delta)$, $\delta\in\mathfrak{B}(\mathbb{R})$, be the orthogonal spectral measure
of the self-adjoint operator $\widetilde A$.
The function
$$ \mathbf{E}(\delta) = P^{\widetilde{H}}_H \widetilde E(\delta),\qquad \delta\in \mathfrak{B}(\mathbb{R}), $$
is said to be {\bf a spectral measure} of a symmetric operator $A$
(corresponding to the extension $\widetilde A$).
The spectral function and the spectral measure of the operator $A$ are connected by the following relation:
$$ \mathbf{E}([0,t)) = \mathbf{E}_t,\qquad t\in \mathbb{R}, $$
what follows from the analogous property of the orthogonal spectral measures.
Moreover, generalized resolvents and spectral functions (measures) of the operator $A$ связаны are connected by the following equality:
\begin{equation}
\label{f3_1_p1_1}
(\mathbf R_\lambda h,g)_H = \int_{\mathbb{R}} \frac{1}{t-\lambda} d(\mathbf{E}(\cdot) h,g)_H =
\int_{\mathbb{R}} \frac{1}{t - \lambda} d(\mathbf{E}_t h,g)_H,\quad \forall h,g\in H,
\end{equation}
which follows directly from their definitions. By the Stieltjes-Perron inversion formula this means that
between generalized resolvents and spectral measures there exists a one-to-one correspondence.

\noindent
For arbitrary elements $h,g\in H$ we may write:
$$ (\mathbf R_\lambda h,g)_H = (R_\lambda(\widetilde A) h,g)_{\widetilde H} =
(h, R_\lambda^*(\widetilde A) g)_{\widetilde H} =
(h, R_{\overline \lambda}(\widetilde A) g)_{\widetilde H} $$
$$ = (h, \mathbf R_\lambda^* g)_H,\qquad \lambda\in \mathbb{R}_e, $$
where $R_\lambda(\widetilde A) = (\widetilde A - \lambda E_{\widetilde H})^{-1}$ is
the resolvent of the self-adjoint operator $\widetilde A$ in a Hilbert space $\widetilde H$,
corresponding to $\mathbf{R}_{\lambda}$.
Therefore
\begin{equation}
\label{f3_1_1_p1_1}
\mathbf R_{s;\lambda}^*(A) = \mathbf R_{s;\overline{\lambda}}(A),\qquad \lambda\in \mathbb{R}_e.
\end{equation}
Choose and fix an arbitrary number $z\in \mathbb{R}_e$.
Consider the Cayley transformation of the operator $A$:
\begin{equation}
\label{f3_1_2_p1_1}
U_z = U_z(A) = (A - \overline{z} E_H)(A-zE_H)^{-1} = E_H + (z-\overline{z}) (A-zE_H)^{-1}.
\end{equation}
The operator $U_z$ is closed and $D(U_z)=\mathcal{M}_z$, $R(U_z)=\mathcal{M}_{\overline{z}}$.
It is readily checked that $U_z$ is isometric.
Since $U_z - E_H = (z-\overline{z}) (A-zE_H)^{-1}$, the operator $U_z$ has no non-zero fixed points.
The operator $A$ is expressed by $U_z$ in the following way:
$$ A = (zU_z - \overline{z} E_H) (U_z - E_H)^{-1} = zE_H + (z-\overline{z})(U_z - E_H)^{-1}. $$
Let $\widehat A\supseteq A$ be a self-adjoint extension of the operator $A$, acting in
a Hilbert space $\widehat H\supseteq H$. Then the following operator:
\begin{equation}
\label{f3_2_p1_1}
W_z = (\widehat A - \overline{z} E_{\widehat H})(\widehat A-z E_{\widehat H})^{-1}
= E_{\widehat H} + (z-\overline{z}) (\widehat A-z E_{\widehat H})^{-1},
\end{equation}
is a unitary extension of the operator $U_z$, acting in $\widehat H$, and having no
non-zero fixed points.

Conversely, if there is a unitary extension $W_z$ of the operator $U_z$, acting in
a Hilbert space $\widehat H$ and having no
non-zero fixed points, then the operator
\begin{equation}
\label{f3_3_p1_1}
\widehat A = (z W_z - \overline{z} E_{\widehat H})
(W_z - E_{\widehat H})^{-1} = z E_{\widehat H} + (z-\overline{z})(W_z - E_{\widehat H})^{-1},
\end{equation}
will be a self-adjoint extension of the operator $A$, acting in $\widehat H$.
Thus, between self-adjoint extensions $\widehat A$ of the operator $A$ and
unitary extensions $W_z$ of the operator $U_z$, acting in $\widehat H$ and having no
non-zero fixed points, there exists a one-to-one correspondence according to formulas~(\ref{f3_1_p1_1}),(\ref{f3_2_p1_1}).

\begin{thm}
\label{t3_1_p1_1}
Let $A$ be a closed symmetric operator in a Hilbert space $H$,
which domain is not necessarily dense in $H$.
Let $z\in \mathbb{R}_e$ be an arbitrary fixed point, and $U_z$ be the Cayley transformation of $A$.
The following equality:
\begin{equation}
\label{f3_4_p1_1}
\mathbf{R}_{u;\frac{\lambda - z}{\lambda - \overline{z}}} (U_z)
= \frac{\lambda - \overline{z}}{z-\overline{z}} E_H +
\frac{(\lambda - \overline{z})(\lambda - z)}{z-\overline{z}} \mathbf{R}_{s;\lambda} (A),\qquad
\lambda\in \mathbb{R}_e\backslash\{ z,\overline{z} \},
\end{equation}
establishes a one-to-one correspondence between all generalized resolvents $\mathbf{R}_{s;\lambda}(A)$
of the operator $A$
and those generalized resolvents $\mathbf{R}_{u;\zeta} (U_z)$ of the closed isometric
operator $U_z$ which are generated by extensions of $U_z$ without non-zero fixed points.

In the case $\overline{ D(A) } = H$, equality~(\ref{f3_4_p1_1})
establishes a one-to-one correspondence between all generalized resolvents $\mathbf{R}_{s;\lambda}(A)$
of the operator $A$
and all generalized resolvents $\mathbf{R}_{u;\zeta} (U_z)$ of the operator $U_z$.
\end{thm}
{\bf Proof. }
Choose and fix a point $z\in \mathbb{R}_e$.
Let $\mathbf{R}_{s;\lambda}(A)$ be an arbitrary generalized resolventof a closed symmetric
operator $A$. It is generated by a self-adjoint operator $\widehat A\supseteq A$
in a Hilbert space $\widehat H\supseteq H$.
Consider the Cayleytransformation $W_z$ of the operator $\widehat A$,
given by~(\ref{f3_2_p1_1}). As it was mentioned above, the operator $W_z$
is a unitary extension of the operator $U_z$, acting in $\widehat H$, and
having no non-zero fixed elements.
Consider the following
linear fractional transformation:
$$ \zeta = \frac{ \lambda - z }{ \lambda - \overline{z} },\quad
\lambda = \frac{ \overline{z} \zeta - z }{ \zeta - 1 }. $$
Supposing that
$\lambda\in \mathbb{R}_e\backslash\{ z,\overline{z} \}$, what is equivalent to the inclusion
$\zeta\in \mathbb{T}_e\backslash\{ 0 \}$,
we write:
$$ \left( E_{\widehat H} - \zeta W_z \right)^{-1} =
\left( E_{\widehat H} - \frac{\lambda - z}{ \lambda - \overline{z}
} \left( E_{\widehat H} + (z-\overline{z}) (\widehat A-z
E_{\widehat H})^{-1} \right) \right)^{-1} $$
$$ = \frac{\lambda - \overline{z}}{z-\overline{z}} \left(
(\widehat A - \lambda E_{\widehat H}) (\widehat A - z E_{\widehat
H})^{-1} \right)^{-1} $$
$$ = \frac{\lambda - \overline{z}}{z-\overline{z}}
(\widehat A - z E_{\widehat H})
(\widehat A - \lambda E_{\widehat H})^{-1}  $$
$$ = \frac{\lambda - \overline{z}}{z-\overline{z}} E_{\widehat H} +
\frac{(\lambda - \overline{z})(\lambda - z)}{z-\overline{z}}
(\widehat A - \lambda E_{\widehat H})^{-1}. $$
Applying the projection operator $P^{\widehat H}_H$ to the first and to the last parts of this relation we get relation~(\ref{f3_4_p1_1}).
The function
$\mathbf{R}_{u;\zeta} (U_z)$ is a generalized resolvent of $U_z$ of the required class.

Conversely, suppose that $\mathbf{R}_{u;\zeta} (U_z)$ is a generalized resolvent of $U_z$,
generated by a unitary extension $W_z\supseteq U_z$, acting in a Hilbert space
$\widehat H\supseteq H$ and having no non-zero fixed elements.
Define the operator $\widehat A$ by equality~(\ref{f3_3_p1_1}). As it was said above, the operator
$\widehat A$ is a unitary extension of $A$ in $\widehat H$. The operator $W_z$ is its
Cayley's transformation. Repeating for the operator $\widehat A$ considerations at the beginning of the proof
we shall come to relation~(\ref{f3_4_p1_1}).

The bijectivity oof the correspondence~(\ref{f3_4_p1_1}) is obvious. In fact, relation~(\ref{f3_4_p1_1})
connects all values of two generalized resolvents except two points where they are defined by the continuity.

Consider the case $\overline{ D(A) } = H$.
Let $\mathbf{R}_{u;\zeta} (U_z)$ be an arbitrary generalized resolvent of the operator $U_z$,
which is generated by a unitary extension $\widetilde W_z\supseteq U_z$, acting
in a Hilbert space $\widetilde H\supseteq H$.
Let $h\in \widetilde H$ be a fixed point of the operator $\widetilde W_z$:
$\widetilde W_z h = h$.
For arbitrary element $g\in\widetilde H$ we may write:
$$ (h,\widetilde W_z g)_{\widetilde H} = (\widetilde W_z h,\widetilde W_z g)_{\widetilde H} =
(h,g)_{\widetilde H}; $$
$$ (h, (\widetilde W_z - E_{\widetilde H}) g )_{\widetilde H} = 0,\qquad g\in\widetilde H. $$
In particular, this implies that $h\perp (U_z - E_H) D(U_z) = D(A)$.
Therefore $h\in\widetilde H\ominus H$.
Denote by $M$ a set of all fixed elements of the operator $\widetilde W_z$.
The set $M$ is a subspace in $\widetilde H$ which is orthogonal to $H$.
Thus, we have $\widetilde H = H\oplus M\oplus H_1$, where $H_1$ is a subspace in $\widetilde H$.

\noindent
Consider an operator $W_z = \widetilde W_z|_{H\oplus H_1}$. The operator $W_z$ has no non-zero fixed points
and it is a unitary extension of $U_z$, if it is considered as an operator in
$H\oplus H_1$.
Since for an arbitrary $f\in H$ we have:
$$ (E_{\widetilde H} - \zeta \widetilde W_z)^{-1} f =
(E_{H\oplus H_1} - \zeta W_z)^{-1} f, $$
then the operator $W_z$ also generates the generalized resolvent $\mathbf{R}_{u;\zeta} (U_z)$.

\noindent
Thus, the set of those generalized resolvents of the operator $U_z$, which are generated by unitary
extensions without non-zero fixed elements, coincides with the set of all generalized resolvents of the operator
$U_z$.
$\Box$

\subsection{Shtraus's formula for the generalized resolvents of a symmetric densely defined operator.}\label{section1_5}

Consider a closed symmetric operator $A$ in a Hilbert space $H$,
assuming that its domain is dense in $H$,
$\overline{D(A)} = H$. Fix an arbitrary point $z\in \mathbb{R}_e$, and consider the Cayley transformation $U_z$
of the operator $A$ from~(\ref{f3_1_2_p1_1}).

\noindent
Let $F$ be an arbitrary linear bounded operator with the domain
$D(F) = \mathcal{N}_z(A)$ and range $R(F) \subseteq \mathcal{N}_{\overline{z}} (A)$.
Set
$$ U_{z;F} = U_{z;F}(A) = U_z \oplus F. $$
Thus, the operator $U_{z;F}(A)$ is a linear bounded operator defined on
the whole $H$. This operator has no non-zero fixed elements. In fact, let for an element
$h\in H$, $h = h_1 + h_2$, $h_1\in \mathcal{M}_z(A)$, $h_2\in \mathcal{N}_z(A)$, we have
$$ 0 = (U_{z;F}(A) - E_H) h = (U_z - E_H) h_1 + Fh_2 - h_2. $$
The first summand in the right-hand side of the last equality belongs to $D(A)$, the second summand
belongs to $\mathcal{N}_{\overline{z}} (A)$, and the third belongs to $\mathcal{N}_z(A)$.
But in the case $\overline{D(A)} = H$ the linear manifolds $D(A)$, $\mathcal{N}_{\overline{z}} (A)$ and
$\mathcal{N}_z(A)$ are linearly independent, and therefore we get an iquality $h_1=h_2=0$.
Thus, we get $h=0$.

Define the following linear operator in $H$:
$$ A_F = A_{F,z} =
(z U_{z;F} - \overline{z} E_H) (U_{z;F} - E_H)^{-1} = zE_H + (z-\overline{z})(U_{z;F} - E_H)^{-1}. $$
Notice that the operator $A_{F,z}$ is an extension of the operator $A$.
The operator $A_F=A_{F,z}$ is said to be {\bf quasi-self-adjoint extension of a symmetric operator $A$,
defined by the operator $F$}.
We outline the following property:
\begin{equation}
\label{f5_1_p1_1}
A_{F,z}^* = \overline{z} E_H + (\overline{z} - z)(U_{z;F}^* - E_H)^{-1} =
\overline{z} E_H + (\overline{z} - z)(U_z^{-1}\oplus F^* - E_H)^{-1} =
A_{F^*,\overline{z}}.
\end{equation}

\begin{thm}
\label{t5_1_p1_1}
Let $A$ be a closed symmetric operator in a Hilbert space $H$ with the dense domain: $\overline{D(A)} = H$.
Fix an arbitrary point $\lambda_0\in \mathbb{R}_e$.
An arbitrary generalized resolvent $\mathbf{R}_{s;\lambda}$ of the operator $A$ has the following form:
\begin{equation}
\label{f5_2_p1_1}
\mathbf R_{s;\lambda} = \left\{ \begin{array}{cc}
\left( A_{F(\lambda)} - \lambda E_H \right)^{-1}, & \lambda\in \Pi_{\lambda_0}\\
\left( A_{F^*(\overline{\lambda})} - \lambda E_H \right)^{-1}, &
\overline{\lambda}\in \Pi_{\lambda_0}
\end{array}
\right.,
\end{equation}
where $F(\lambda)$ is a function from $\mathcal{S}(\Pi_{\lambda_0}; \mathcal{N}_{\lambda_0},
\mathcal{N}_{\overline{\lambda_0}})$.
Conversely, an arbitrary function $F(\lambda)\in \mathcal{S}(\Pi_{\lambda_0}; \mathcal{N}_{\lambda_0},
\mathcal{N}_{\overline{\lambda_0}})$
defines by relation~(\ref{f5_2_p1_1}) a generalized resolvent
$\mathbf{R}_{s;\lambda}$ of the operator $A$.
To different functions from
$\mathcal{S}(\Pi_{\lambda_0}; \mathcal{N}_{\lambda_0},
\mathcal{N}_{\overline{\lambda_0}})$
there correspond different generalized resolvents of the operator $A$. Here $\Pi_{\lambda_0}$ is that half-plane of
$\mathbb{C}_+$ and $\mathbb{C}_-$, which contains the point $\lambda_0$.
\end{thm}
{\bf Proof. }
Let $A$ be a closed symmetric operator in a Hilbert space $H$, $\overline{D(A)} = H$, and
$\lambda_0\in \mathbb{R}_e$ be a fixed point.
Consider an arbitrary  generalized resolvent $\mathbf{R}_{s;\lambda} (A)$ of the operator $A$.
Let us use Theorem~\ref{t3_1_p1_1} with $z=\lambda_0$.
From equality~(\ref{f3_4_p1_1}) we express $\mathbf{R}_{s;\lambda} (A)$:
$$
\mathbf{R}_{s;\lambda} (A)
=
\frac{\lambda_0-\overline{\lambda_0}}{(\lambda - \overline{\lambda_0})(\lambda - \lambda_0)} \left(
\mathbf{R}_{u;\frac{\lambda - \lambda_0}{\lambda - \overline{\lambda_0}}} (U_{\lambda_0})
- \frac{\lambda - \overline{\lambda_0}}{\lambda_0-\overline{\lambda_0}} E_H
\right),\ \lambda\in \mathbb{R}_e\backslash\{ \lambda_0,\overline{\lambda_0} \}.
$$
Observe that the linear fractional transformation $\zeta =
\frac{\lambda - \lambda_0}{\lambda - \overline{\lambda_0}}$ maps the half-plane $\Pi_{\lambda_0}$
on $\mathbb{D}$ (and $\mathbb{R}$ on $\mathbb{T}$).
Thus, we can restrict the last relation and consider it only for
$\lambda\in \Pi_{\lambda_0}\backslash\{ \lambda_0 \}$. In this case we may apply
to the generalized resolvent $\mathbf{R}_{u;\frac{\lambda - \lambda_0}{\lambda - \overline{\lambda_0}}} (U_{\lambda_0})$
Chumakin's formula~(\ref{f2_1_p1_1}). Then
$$
\mathbf{R}_{s;\lambda} (A)
=
\frac{\lambda_0-\overline{\lambda_0}}{(\lambda - \overline{\lambda_0})(\lambda - \lambda_0)} \left(
\left[ E_H - \frac{\lambda - \lambda_0}{\lambda - \overline{\lambda_0}}
(U_{\lambda_0}\oplus \Phi_{\frac{\lambda - \lambda_0}{\lambda - \overline{\lambda_0}}}) \right]^{-1}
\right. $$
\begin{equation}
\label{f5_2_1_p1_1}
\left. - \frac{\lambda - \overline{\lambda_0}}{\lambda_0-\overline{\lambda_0}} E_H
\right),\qquad \lambda\in \Pi_{\lambda_0}\backslash\{ \lambda_0 \},
\end{equation}
where $\Phi_\zeta$ is a function from $\mathcal{S}(\mathbb{D};\mathcal{N}_{\lambda_0}(A),
\mathcal{N}_{\overline{\lambda_0}}(A))$.
Denote
$$ F(\lambda) = \Phi_{\frac{\lambda - \lambda_0}{\lambda - \overline{\lambda_0}}},\qquad \lambda\in\Pi_{\lambda_0}. $$
Notice that $F(\lambda)\in \mathcal{S}(\Pi_{\lambda_0};\mathcal{N}_{\lambda_0}(A),
\mathcal{N}_{\overline{\lambda_0}}(A))$.
Then
$$
\mathbf{R}_{s;\lambda} (A)
=
\frac{\lambda_0-\overline{\lambda_0}}{(\lambda - \overline{\lambda_0})(\lambda - \lambda_0)} \left(
\left[ E_H - \frac{\lambda - \lambda_0}{\lambda - \overline{\lambda_0}}
(U_{\lambda_0}\oplus F(\lambda)) \right]^{-1}
\right. $$
$$ \left. - \frac{\lambda - \overline{\lambda_0}}{\lambda_0-\overline{\lambda_0}}
\left[ E_H - \frac{\lambda - \lambda_0}{\lambda - \overline{\lambda_0}}
(U_{\lambda_0}\oplus F(\lambda)) \right]
\left[ E_H - \frac{\lambda - \lambda_0}{\lambda - \overline{\lambda_0}}
(U_{\lambda_0}\oplus F(\lambda)) \right]^{-1}
\right) $$
$$ = \frac{1}{\lambda - \overline{\lambda_0}}
\left(
- E_H + (U_{\lambda_0}\oplus F(\lambda))
\right)
\left[ E_H - \frac{\lambda - \lambda_0}{\lambda - \overline{\lambda_0}}
(U_{\lambda_0}\oplus F(\lambda)) \right]^{-1}, $$
\begin{equation}
\label{f5_3_p1_1}
\lambda\in \Pi_{\lambda_0}\backslash\{ \lambda_0 \}.
\end{equation}
In the half-plane $\Pi_{\lambda_0}$ the following inequality holds:
$$ \left| \frac{\lambda - \lambda_0}{\lambda - \overline{\lambda_0}} \right| < 1,\qquad
\lambda\in\Pi_{\lambda_0}. $$
Consider the open neighborhood of the point $\lambda_0$:
$$ D(\lambda_0) = \left\{ z\in \Pi_{\lambda_0}:\ |z-\lambda_0| < \frac{1}{2} \mathop{\rm Im}\nolimits \lambda_0
\right\}. $$
In the closed neighborhood $\overline{ D(\lambda_0) }$ the function
$\left| \frac{\lambda - \lambda_0}{\lambda - \overline{\lambda_0}} \right|$
is continuous and attains its maximal value $q < 1$. Therefore
$$ \left| \frac{\lambda - \lambda_0}{\lambda - \overline{\lambda_0}} \right| \leq q < 1,\qquad
\lambda\in D(\lambda_0). $$
For an arbitrary element $h\in H$ we may write:
$$ \left\| \left[ E_H - \frac{\lambda - \lambda_0}{\lambda - \overline{\lambda_0}}
(U_{\lambda_0}\oplus F(\lambda)) \right] h \right\| \geq
\left|
\| h \| - \left| \frac{\lambda - \lambda_0}{\lambda - \overline{\lambda_0}} \right|
\| (U_{\lambda_0}\oplus F(\lambda)) h \|
\right| $$
$$ \geq
(1-q) \| h \|, $$
and therefore
$$ \left\| \left[ E_H - \frac{\lambda - \lambda_0}{\lambda - \overline{\lambda_0}}
(U_{\lambda_0}\oplus F(\lambda)) \right]^{-1}
\right\| \leq \frac{1}{1-q},\qquad \lambda\in D(\lambda_0). $$
Applying Proposition~\ref{p1_1_p1_1} we conclude that the function
$\left[ E_H - \frac{\lambda - \lambda_0}{\lambda - \overline{\lambda_0}}
(U_{\lambda_0}\oplus F(\lambda)) \right]^{-1}$ is analytic in $D(\lambda_0)$.
Passing to the limit in relation~(\ref{f5_3_p1_1}) as $\lambda\rightarrow\lambda_0$
we obtain:
$$
\mathbf{R}_{s;\lambda_0} (A)
= \frac{1}{\lambda_0 - \overline{\lambda_0}}
\left(
- E_H + (U_{\lambda_0}\oplus F(\lambda_0))
\right). $$
Thus, relation~(\ref{f5_3_p1_1}) holds for $\lambda = \lambda_0$, as well.

\noindent
In our notations, maid at the beginning of this subsection, we have
$U_{\lambda_0}\oplus F(\lambda) = U_{\lambda_0; F(\lambda)}$.
According to our above remarks $U_{\lambda_0; F(\lambda)}$ has no non-zero fixed elements.
Therefore there exists the inverse
$$ \mathbf{R}_{s;\lambda}^{-1} (A) =
(\lambda - \overline{\lambda_0})
\left[ E_H - \frac{\lambda - \lambda_0}{\lambda - \overline{\lambda_0}}
U_{\lambda_0; F(\lambda)} \right]
\left(
U_{\lambda_0; F(\lambda)}
- E_H
\right)^{-1} $$
$$  =
\left[
\lambda E_H - \overline{\lambda_0} E_H - \lambda U_{\lambda_0; F(\lambda)}
+ \lambda_0 U_{\lambda_0; F(\lambda)}
\right]
\left(
U_{\lambda_0; F(\lambda)}
- E_H
\right)^{-1} $$
$$ = -\lambda E_H +
\left[
\lambda_0 U_{\lambda_0; F(\lambda)} - \overline{\lambda_0} E_H
\right]
\left(
U_{\lambda_0; F(\lambda)}
- E_H
\right)^{-1} $$
$$ = -\lambda E_H + A_{F(\lambda);\lambda_0},\qquad \lambda\in \Pi_{\lambda_0}. $$
From this relation we conclude that~(\ref{f5_2_p1_1}) holds for $\lambda\in\Pi_{\lambda_0}$.

Suppose now that $\lambda$ is such that $\overline{\lambda}\in\Pi_{\lambda_0}$.
Applying the proved part of the formula for $\overline{\lambda}$ we get
$$ \mathbf R_{s;\overline{\lambda}} =
\left( A_{F(\overline{\lambda})} - \overline{\lambda} E_H \right)^{-1}. $$
Using property~(\ref{f3_1_1_p1_1}) of the generalized resolvent and using relation~(\ref{f5_1_p1_1}) we write
$$ \mathbf R_{s;\lambda} =
\mathbf R_{s;\overline{\lambda}}^* =
\left( A_{F(\overline{\lambda})}^* - \lambda E_H \right)^{-1}
= \left( A_{F^*(\overline{\lambda})} - \lambda E_H \right)^{-1}. $$
Consequently, relation~(\ref{f5_2_p1_1}) for $\overline{\lambda}\in\Pi_{\lambda_0}$
is true, as well.

Conversely, let  an operator-valued function
$F(\lambda)\in \mathcal{S}(\Pi_{\lambda_0}; \mathcal{N}_{\lambda_0},
\mathcal{N}_{\overline{\lambda_0}})$ be given. Set
$$ \Phi(\zeta) = F\left(
\frac{\overline{\lambda_0}\zeta - \lambda_0}{\zeta - 1}
\right),\qquad \zeta\in \mathbb{D}. $$
The function $\Phi(\zeta)$ belongs to $\mathcal{S}(\mathbb{D};\mathcal{N}_{\lambda_0}(A),
\mathcal{N}_{\overline{\lambda_0}}(A))$.
By Chumakin's formula~(\ref{f2_1_p1_1}), to this function there corresponds a generalized resolvent
$\mathbf{R}_{u;\zeta} (U_{\lambda_0})$.
Define a generalized resolvent $\mathbf{R}_{s;\lambda} (A)$ of the operator $A$ by relation~(\ref{f5_2_1_p1_1}).
Repeating for the generalized resolvent $\mathbf{R}_{s;\lambda} (A)$
considerations after~(\ref{f5_2_1_p1_1}) we come to relation~(\ref{f5_2_p1_1}).
Therefore the function $F(\lambda)$ generates by relation~(\ref{f5_2_p1_1}) a generalized resolvent of the operator $A$.

Suppose that two operator-valued functions $F_1(\lambda)$ and $F_2(\lambda)$, which belong to
$\mathcal{S}(\Pi_{\lambda_0}; \mathcal{N}_{\lambda_0}, \mathcal{N}_{\overline{\lambda_0}})$,
generate by relation~(\ref{f5_2_p1_1}) the same generalized resolvent. In this case
we have
$$ \left( A_{F_1(\lambda)} - \lambda E_H \right)^{-1} =
\left( A_{F_2(\lambda)} - \lambda E_H \right)^{-1},\qquad \lambda\in \Pi_{\lambda_0}, $$
and therefore
$$ A_{F_1(\lambda)} = A_{F_2(\lambda)},\qquad \lambda\in \Pi_{\lambda_0}. $$
From the last relation it follows that $U_{\lambda_0;F_1(\lambda)} = U_{\lambda_0;F_2(\lambda)}$,
and we get the equality $F_1(\lambda) = F_2(\lambda)$.
$\Box$

Formula~(\ref{f5_2_p1_1}) is said to be the Shtraus formula for the generalized resolvents of a symmetric operator with
a dense domain.


\section{Generalized resolvents of non-densely defined symmetric operators.}\label{chapter2}

\subsection{The forbidden operator.}
Throughout this subsection we shall consider a closed symmetric operator $A$ in a Hilbert
space $H$, which domain can be non-dense in $H$.


\begin{prop}
\label{p1_1_p2_1}
Let $A$ be a closed symmetric operator in a Hilbert space $H$ and $z\in \mathbb{R}_e$. Let
elements $\psi\in \mathcal{N}_z(A)$ and $\varphi\in \mathcal{N}_{\overline{z}}(A)$ be such that
$\varphi - \psi\in D(A)$, i.e. they are comparable by modulus $D(A)$.
Then $\| \varphi \|_H = \| \psi \|_H$.
\end{prop}
{\bf Proof. }
Since $D(A) = (U_z - E_H) \mathcal{M}_z(A)$, then there exists an element $g\in \mathcal{M}_z(A)$ such that
$$ \varphi - \psi = U_z g - g. $$
Then
\begin{equation}
\label{f1_1_p2_1}
\varphi + U_z g = \psi +  g.
\end{equation}
Using the orthogonality of summands in the left and right sides we conclude that
$$ \| \psi \|_H^2 + \| g \|_H^2 = \| \varphi \|_H^2 + \| U_z g \|_H^2 = \| \varphi \|_H^2 + \| g \|_H^2, $$
and the required equality follows. $\Box$

\begin{cor}
\label{c1_1_p2_1}
Let $A$ be a closed symmetric operator in a Hilbert space $H$.
Then $D(A)\cap \mathcal{N}_z(A) = \{ 0 \}$, $\forall z\in \mathbb{R}_e$.
\end{cor}
{\bf Proof. }
If $\psi\in D(A)\cap \mathcal{N}_z(A)$, then applying Proposition~\ref{p1_1_p2_1}
to elements $\psi$ and $\varphi = 0$, we get $\psi = 0$.
$\Box$

Let $z\in \mathbb{R}_e$ be an arbitrary number.
Consider the following operator:
$$ X_z \psi = X_z(A) \psi
= \varphi,\qquad \psi\in \mathcal{N}_z(A)\cap \left( \mathcal{N}_{\overline{z}}(A) \dotplus D(A)
\right), $$
where $\varphi\in \mathcal{N}_{\overline{z}}(A)$: $\psi - \varphi\in D(A)$.

\noindent
This definition will be correct if an element $\varphi\in \mathcal{N}_{\overline{z}}$
such that $\psi - \varphi\in D(A)$,
is unique (its existence is obvious, since
$\psi\in \left( \mathcal{N}_{\overline{z}}(A) \dotplus D(A) \right)$).
Let $\varphi_1\in \mathcal{N}_{\overline{z}}$: $\psi - \varphi_1\in D(A)$.
Then $\varphi - \varphi_1\in D(A)$. Since $\varphi - \varphi_1\in \mathcal{N}_{\overline{z}}$,
then by applying Corollary~\ref{c1_1_p2_1} we get $\varphi = \varphi_1$.

From Proposition~\ref{p1_1_p2_1} it follows that the operator $X_z$ is isometric.
Observe that $D(A) = \mathcal{N}_z(A)\cap \left( \mathcal{N}_{\overline{z}}(A) \dotplus D(A)
\right)$ and $R(A) = \mathcal{N}_{\overline z}(A)\cap \left( \mathcal{N}_z(A) \dotplus D(A)
\right)$. The operator $X_z=X_z(A)$ is said to be \textbf{forbidden with respect to the symmetric operator $A$}.

Now we shall obtain another representation for the operator $X_z$. Before this we shall prove some auxilliary
results:
\begin{prop}
\label{p1_2_p2_1}
Let $A$ be a closed symmetric operator in a Hilbert space $H$ and $z\in \mathbb{R}_e$.
The following two conditions are equivalent:
\begin{itemize}
\item[$\mathrm{(i)}$] Elements $\psi\in \mathcal{N}_z(A)$ and $\varphi\in \mathcal{N}_{\overline{z}}(A)$
are comparable by modulus $D(A)$;

\item[$\mathrm{(ii)}$] Elements $\psi\in \mathcal{N}_z(A)$ and $\varphi\in \mathcal{N}_{\overline{z}}(A)$
admit the following representation:
\begin{equation}
\label{f1_2_p2_1}
\psi = P^H_{\mathcal{N}_z(A)} h,\quad \varphi = P^H_{\mathcal{N}_{\overline{z}}(A)} h,
\end{equation}
where $h\in H$ is such that
\begin{equation}
\label{f1_3_p2_1}
P^H_{\mathcal{M}_{\overline{z}}(A)} h = U_z P^H_{\mathcal{M}_z(A)} h.
\end{equation}

\end{itemize}

If these conditions are satisfied, the element $h$ is defined uniquely and admits the following representation:
\begin{equation}
\label{f1_4_p2_1}
h = \frac{1}{z - \overline{z}}
\left(
A(\psi - \varphi) - \overline{z}\psi + z\varphi
\right).
\end{equation}
\end{prop}
{\bf Proof. } $\mathrm{(i)} \Rightarrow \mathrm{(ii)}$.
As at the beginning of the proof of Proposition~\ref{p1_1_p2_1}, we derive formula~(\ref{f1_1_p2_1}),
where $g\in \mathcal{M}_z(A)$.
Set $h = g + \psi$.
Applying
$P^H_{\mathcal{N}_{z}(A)}$ to this equality we obtain the first equality in~(\ref{f1_2_p2_1}).
From~(\ref{f1_1_p2_1}) it follows that
$$ \varphi + U_z (h-\psi) = h. $$
The summands in the left-hand side belong to orthogonal subspaces. Applying the projection operators
$P^H_{\mathcal{N}_{\overline{z}}(A)}$ and $P^H_{\mathcal{M}_{\overline{z}}(A)}$, we obtain the following equalities:
$$ \varphi = P^H_{\mathcal{N}_{\overline{z}}(A)} h,\quad
U_z(h-\psi) = P^H_{\mathcal{M}_{\overline{z}}(A)} h, $$
and therefore the second equality in~(\ref{f1_2_p2_1}) is satisfied. Using the latter equality and
the equality $P^H_{\mathcal{M}_z} h = P^H_{\mathcal{M}_z} (g+\psi) =
P^H_{\mathcal{M}_z} g = g$,
we may write
$$ P^H_{\mathcal{M}_{\overline{z}}(A)} h = U_z g = U_z P^H_{\mathcal{M}_z(A)} h. $$
Consequently, relation~(\ref{f1_3_p2_1}) is true, as well.

\noindent
$\mathrm{(ii)} \Rightarrow \mathrm{(i)}$.
Set $g = P^H_{\mathcal{M}_z} h$. Then
$$ h = P^H_{\mathcal{M}_z} h + P^H_{\mathcal{N}_z} h = g + \psi, $$
$$ h = P^H_{\mathcal{M}_{\overline z}} h + P^H_{\mathcal{N}_{\overline z}} h =
U_z g + \varphi, $$
and subtracting last equalities we obtain that $\psi - \varphi = U_z g - g\in D(A)$.
Moreover, we have:
$$ h = g+\psi = (U_z - E_H)^{-1} (\psi - \varphi) + \psi =
\frac{1}{z-\overline{z}} (A - zE_H) (\psi - \varphi) + \psi, $$
and relation~(\ref{f1_4_p2_1}) follows.
$\Box$

\begin{prop}
\label{p1_3_p2_1}
Let $A$  be a closed symmetric operator in a Hilbert space $H$, and
$z\in \mathbb{R}_e$.
Then
$$ H\ominus \overline{D(A)} = \left\{
h\in H:\ P^H_{\mathcal{M}_{\overline z}(A)} h = U_z P^H_{\mathcal{M}_z(A)} h
\right\}. $$
\end{prop}
{\bf Proof. }
For arbitrary elements $g\in \mathcal{M}_z(A)$ and $h\in H$ we may write:
$$ (U_z g, h)_H = (U_z g, P^H_{ \mathcal{M}_{ \overline{z} }(A) } h)_H; $$
$$ (g, h)_H = (g, P^H_{ \mathcal{M}_z(A) } h)_H = (U_z g, U_z P^H_{ \mathcal{M}_z(A) } h)_H, $$
and therefore
$$ ((U_z - E_H)g,h)_H = \left(
U_z g, \left( P^H_{ \mathcal{M}_{ \overline{z} } } - U_z P^H_{ \mathcal{M}_z } \right) h
\right)_H. $$
If $h\perp D(A)$, then the left-hand side of the last equality is equal to zero, since $(U_z - E_H)g\in D(A)$.
Then $\left( P^H_{ \mathcal{M}_{ \overline{z} } } - U_z P^H_{ \mathcal{M}_z } \right) h = 0$.

Conversely, if $h\in H$ is such that
$\left( P^H_{ \mathcal{M}_{ \overline{z} } } - U_z P^H_{ \mathcal{M}_z } \right) h = 0$, then
we get $((U_z - E_H)g,h)_H =0$, $g\in \mathcal{M}_z(A)$, i.e. $h\perp D(A)$.
$\Box$

\begin{cor}
\label{c1_2_p2_1}
Let $A$  be a closed symmetric operator in a Hilbert space $H$.
Then $\left( H\ominus\overline{D(A)} \right)\cap \mathcal{M}_z(A) = \{ 0 \}$, $\forall z\in \mathbb{R}_e$.
\end{cor}
{\bf Proof. }
Let $h\in \left( H\ominus\overline{D(A)} \right)\cap \mathcal{M}_z(A)$, $z\in \mathbb{R}_e$.
By Proposition~\ref{p1_3_p2_1} the following equality holds:
$$ P^H_{\mathcal{M}_{\overline z}(A)} h = U_z P^H_{\mathcal{M}_z(A)} h. $$
By Proposition~\ref{p1_2_p2_1} we obtain that elements
$\psi := P^H_{\mathcal{N}_z(A)} h = 0$ and $\varphi := P^H_{\mathcal{N}_{\overline{z}}(A)} h$
are comparable by modulo $D(A)$. By Proposition~\ref{p1_1_p2_1} we get $\varphi = 0$.
From~(\ref{f1_4_p2_1}) it follows that $h=0$.
$\Box$

Let us check that
\begin{equation}
\label{f1_5_p2_1}
D(X_z) = P^H_{\mathcal{N}_z(A)} (H\ominus \overline{D(A)}),
\end{equation}
\begin{equation}
\label{f1_6_p2_1}
X_z P^H_{\mathcal{N}_z(A)} h = P^H_{\mathcal{N}_{\overline z}(A)} h,\qquad h\in H\ominus \overline{D(A)}.
\end{equation}
In fact, if $\psi\in D(X_z)$, then $\psi\in \mathcal{N}_z$ and there exists
$\varphi\in \mathcal{N}_{\overline{z}}$ such that $\psi - \varphi\in D(A)$.
By Propositions~\ref{p1_2_p2_1} and~\ref{p1_3_p2_1} we see that
there exists $h\in H\ominus \overline{D(A)}$:
$\psi = P^H_{\mathcal{N}_z(A)} h$.
Therefore $D(X_z) \subseteq P^H_{\mathcal{N}_z(A)} (H\ominus \overline{D(A)})$.

Conversely, let $\psi\in P^H_{\mathcal{N}_z(A)} (H\ominus \overline{D(A)})$,
$\psi = P^H_{\mathcal{N}_z(A)} h$, $h\in H\ominus \overline{D(A)}$.
By Propositions~\ref{p1_2_p2_1} and~\ref{p1_3_p2_1} we obtain that $\psi$
and $\varphi := P^H_{\mathcal{N}_{\overline{z}}(A)} h$ are comparable by modulus $D(A)$.
Consequently, we have $\psi\in D(X_z)$ and therefore $P^H_{\mathcal{N}_z(A)} (H\ominus \overline{D(A)})\subseteq D(X_z)$.
Thus, equality~(\ref{f1_5_p2_1}) is true.

For an arbitrary element $h\in H\ominus \overline{D(A)}$ we set
$\psi = P^H_{\mathcal{N}_z(A)} h$, $\varphi = P^H_{\mathcal{N}_{\overline{z}}(A)} h$.
By Propositions~\ref{p1_2_p2_1} and~\ref{p1_3_p2_1}, we have $\psi - \varphi\in D(A)$.
Therefore $\psi\in D(X_z)$ and  $X_z \psi = \varphi$, and equality~(\ref{f1_6_p2_1}) follows.

\subsection{Admissible operators.}
We continue considering of a closed symmetric operator $A$ in a Hilbert space $H$.
Fix an arbitrary point $z\in \mathbb{R}_e$.
As before, $X_z$ will denote the forbidden operator with respect to $A$.
Let $V$ be an arbitrary operator with the domain
$D(V)\subseteq \mathcal{N}_z(A)$ and the range $R(V)\subseteq \mathcal{N}_{\overline z}(A)$.
The operator $V$ is said to be {\bf $z$-admissible (or admissible) with respect to symmetric operator $A$},
if the operator $V-X_z$ is invertible.
In other words, $V$ is admissible with respect to the operator $A$, if the equality $V\psi = X_z\psi$
can happen only if $\psi = 0$.
Using the definition of the forbidden operator $X_z$, the latter definition can be formulated in the following way:
$V$ is admissible with respect to the operator $A$, if $\psi\in D(V)$, $V\psi -\psi \in D(A)$,
implies $\psi = 0$.

If the domain of the operator $A$ is dense in $H$, then by~(\ref{f1_5_p2_1}) we get
$D(X_z) = \{ 0 \}$. Consequently, {\it in the case $\overline{ D(A) } = H$, an arbitrary operator
with the domain
$D(V)\subseteq \mathcal{N}_z(A)$ and the range $R(V)\subseteq \mathcal{N}_{\overline z}(A)$
is admissible with respect to $A$. }

Let $B$ be a symmetric extension of the operator $A$ in the space $H$.
Then its Cayley's transformation
\begin{equation}
\label{f2_0_p2_1}
W_z = (B - \overline{z} E_H)(B-zE_H)^{-1} = E_H + (z-\overline{z}) (B-zE_H)^{-1},
\end{equation}
is an isometric extension of the operator $U_z(A)$ in $H$, which has no non-zero
fixed elements. Moreover, we have:
\begin{equation}
\label{f2_0_1_p2_1}
B = (z W_z - \overline{z} E_H) (W_z - E_H)^{-1} = zE_H + (z-\overline{z})(W_z - E_H)^{-1}.
\end{equation}
Conversely, for an arbitrary isometric operator $W_z\supseteq U_z(A)$, having no non-zero fixed elements,
by~(\ref{f2_0_1_p2_1}) one can define a symmetric operator
$B\supseteq A$.
Thus, formula~(\ref{f2_0_p2_1}) establishes a one-to-one correspondence between all symmetric
extensions $B\supseteq A$ in $H$, and all isometric
extensions $W_z\supseteq U_z(A)$  in $H$, having no non-zero fixed elements.

By Proposition~\ref{p1_3_p1_1}, all such extensions
$W_z$ have the following form:
\begin{equation}
\label{f2_0_2_p2_1}
W_z = U_z(A)\oplus T,
\end{equation}
where $T$ is a isometric operator with the domain $D(T)\subseteq \mathcal{N}_z(A)$
and the range $R(T)\subseteq \mathcal{N}_{\overline{z}}(A)$.

Conversely, if we have an arbitrary isometric operator $T$ with the domain $D(T)\subseteq \mathcal{N}_z(A)$
and the range $R(T)\subseteq \mathcal{N}_{\overline{z}}(A)$, then by~(\ref{f2_0_2_p2_1})
we define an isometric extension $W_z$ of the operator $U_z(A)$.
Question: what additional property should have the operator $T$ which garantees that the operator $W_z$ has no non-zero fixed elements?
An answer on this question is provided by the following theorem.

\begin{thm}
\label{t2_1_p2_1}
Let $A$ be a closed symmetric operator in a Hilbert space $H$, $z\in \mathbb{R}$
be a fixed point, and $U_z$ be the Cayley transformation of the operator $A$.
Let $V$ be an arbitrary operator with the domain
$D(V)\subseteq \mathcal{N}_z(A)$ and the range $R(V)\subseteq \mathcal{N}_{\overline z}(A)$.
The operator $U_z\oplus V$ having no non-zero fixed elements if and only if the operator
$V$ is $z$-admissible with respect to $A$.
\end{thm}
{\bf Proof. }
{\it Necessity. } To the contrary, suppose that
$U_z\oplus V$ has no non-zero fixed elements but $V$ is not admissible with respect to $A$.
This means that there exists a non-zero element $\psi\in D(V)\cap D(X_z)$ such that $V\psi = X_z\psi$,
where $X_z$ denotes the forbidden operator.
By~(\ref{f1_5_p2_1}),(\ref{f1_6_p2_1}) there exists an element $h\in H\ominus \overline{D(A)}$:
$\psi = P^H_{\mathcal{N}_z(A)} h$ and
$$ V P^H_{\mathcal{N}_z(A)} h = X_z P^H_{\mathcal{N}_z(A)} h = P^H_{\mathcal{N}_{\overline z}(A)} h. $$
By Proposition~\ref{p1_3_p2_1} we obtain the following equality:
$$ P^H_{\mathcal{M}_{\overline z}(A)} h = U_z P^H_{\mathcal{M}_z(A)} h. $$
Then
$$ (U_z\oplus V) h =  U_z P^H_{\mathcal{M}_z(A)} h + V P^H_{\mathcal{N}_z(A)} h =
P^H_{\mathcal{M}_{\overline z}(A)} h + P^H_{\mathcal{N}_{\overline z}(A)} h = h, $$
i.e. $h$ is a fixed element of the operator $U_z\oplus V$.
Since we assumed that $U_z\oplus V$ has no non-zero fixed elements, then
$h=0$, and therefore $\psi = 0$. We obtained a contradiction.

\noindent
{\it Sufficiency. } Suppose to the contrary that
$V$ is admissible with respect to $A$, but there exists a non-zero element $h\in H$ such that
$$ (U_z\oplus V) h = U_z P^H_{\mathcal{M}_z(A)} h + V P^H_{\mathcal{N}_z(A)} h = h =
P^H_{\mathcal{M}_{\overline z}(A)} h + P^H_{\mathcal{N}_{\overline z}(A)} h. $$
Then the following equalities hold:
$$ U_z P^H_{\mathcal{M}_z(A)} h = P^H_{\mathcal{M}_{\overline z}(A)} h, $$
$$ V P^H_{\mathcal{N}_z(A)} h = P^H_{\mathcal{N}_{\overline z}(A)} h. $$
By Proposition~\ref{p1_3_p2_1} and the first of these equalities we get
$h\in H\ominus \overline{D(A)}$. By relations~(\ref{f1_5_p2_1}),(\ref{f1_6_p2_1})
we see that $\psi := P^H_{\mathcal{N}_z(A)} h$ belongs to $D(X_z)$, and
$$ X_z P^H_{\mathcal{N}_z(A)} h = P^H_{\mathcal{N}_{\overline z}(A)} h = V P^H_{\mathcal{N}_z(A)} h, $$
i.e. $V\psi = X_z\psi$. Since $V$ is admissible with respect to $A$, then $\psi = 0$.
By the definition of the operator $X_z$, the element $\varphi := X_z\psi$ is comparable with $\psi$
by modulus $D(A)$. By Proposition~\ref{p1_1_p2_1} we conclude that $\varphi = 0$.
By Proposition~\ref{p1_2_p2_1} and formula~(\ref{f1_4_p2_1}) we get $h=0$.
The obtained contradiction completes the proof.
$\Box$

\begin{rmr}
\label{r2_1_p2_1}
The last Theorem shows that formulas~(\ref{f2_0_p2_1}),(\ref{f2_0_1_p2_1}) and
(\ref{f2_0_2_p2_1}) establish a one-to-one correspondence between all symmetric
extensions $B$ of the operator $A$ in $H$, and all isometric
operators $T$  with the domain $D(T)\subseteq \mathcal{N}_z(A)$
and the range $R(T)\subseteq \mathcal{N}_{\overline{z}}(A)$, which are admissible with respect to $A$.
\end{rmr}
A symmetric operator $A$ is said to be {\it maximal} if it has no proper
(i.e. different from $A$) symmetric extensions in $H$. By the above remark
the maximality depends on the number of isometric operators $T$ with $D(T)\subseteq \mathcal{N}_z(A)$,
$R(T)\subseteq \mathcal{N}_{\overline{z}}(A)$, which are admissible with respect to $A$.

Consider arbitrary subspaces of the same dimension: $N_z^0\in \mathcal{N}_z(A)$
and $N_{\overline{z}}^0\in \mathcal{N}_{\overline{z}}(A)$. Let us check that there always exists
an isometric operator $V$, which maps $N_z^0$ on $N_{\overline{z}}^0$,
and which is admissible with respect to the operator $A$.
In order to do that we shall need the following simple proposition.
\begin{prop}
\label{p2_1_p2_1}
Let $H_1$ and $H_2$ are subspaces of the same dimension in a Hilbert space $H$.
Then there exists an isometric operator $V$ with the domain $D(V)=H_1$ and
the range $R(V)=H_2$, which has no non-zero fixed elements.
\end{prop}
{\bf Proof. }
Choose orthonormal bases $\{ f_k \}_{k=0}^d$ and $\{ g_k \}_{k=0}^d$, $d\leq\infty$, in
$H_1$ and $H_2$, respectively. The operator
$$ U \sum_{k=0}^d \alpha_k f_k = \sum_{k=0}^d \alpha_k g_k,\qquad \alpha_k\in \mathbb{C}, $$
maps isometrically $H_1$ on $H_2$. Denote
$$ H_0 = \{ h\in H_1\cap H_2:\ Uh = h \}, $$
i.e. $H_0$ is a set of the fixed elements of the operator $U$.
Notice that $H_0$ is a subspace in $H_1$.
Choose and fix a number $\alpha\in (0,2\pi)$.
Set
$$ V h = e^{i\alpha} P^{H_1}_{H_0} h + U P^{H_1}_{H_1\ominus H_0} h,\qquad h\in H_1. $$
Suppose that $g\in H_1$ is a fixed element $V$, i.e.
$$ Vg = e^{i\alpha} P^{H_1}_{H_0} g + U P^{H_1}_{H_1\ominus H_0} g = g =
P^{H_1}_{H_0} g + P^{H_1}_{H_1\ominus H_0} g. $$
Since
$$ (U P^{H_1}_{H_1\ominus H_0} g, h ) = (U P^{H_1}_{H_1\ominus H_0} g, U h ) =
(P^{H_1}_{H_1\ominus H_0} g, h ) = 0,\qquad \forall h\in H_0, $$
then $U P^{H_1}_{H_1\ominus H_0} g\perp H_0$, and equating summands in the previous equality we get
$(e^{i\alpha} - 1) P^{H_1}_{H_0} g=0$ and $U P^{H_1}_{H_1\ominus H_0} g = P^{H_1}_{H_1\ominus H_0} g$.
Therefore $P^{H_1}_{H_1\ominus H_0} g\in H_0$, and $P^{H_1}_{H_1\ominus H_0} g = 0$.
Thus, we get $g=0$.
$\Box$

Consider the following linear manifold:
$$ M = \left\{
\psi\in D(X_z)\cap N_z^0:\ X_z\psi\in N_{\overline{z}}^0
\right\}. $$
Set
$$ X_{z;0} = X_z|_M. $$
Denote
$$ N_z' = \overline{M} \subseteq N_z^0,\quad
N_{\overline{z}}' = \overline{R(X_{z;0})} \subseteq N_{\overline{z}}^0. $$
Since $X_z$ is isometric, we get
$$ \dim N_z' = \dim N_{\overline{z}}'. $$
Consider an arbitrary isometric operator $W$, which maps the subspace $N_z^0$
on the subspace $N_{\overline{z}}^0$ of the same dimention.
Set
$$ K := W N_z' \subseteq N_{\overline{z}}^0. $$
Since $W$ is isometric we get
\begin{equation}
\label{f2_1_p2_1}
\dim K = \dim N_z' = \dim N_{\overline{z}}',\quad
\dim (N_{\overline{z}}^0\ominus K) = \dim(N_z^0\ominus N_z').
\end{equation}

We shall say that {\bf an isometric operator $S$ corrects $X_{z;0}$ in $N_{\overline{z}}^0$}, if
the following three conditions are satisfied:

\begin{itemize}
\item[1)] $D(S) = N_{\overline{z}}'$, $R(S) \subseteq N_{\overline{z}}^0$;

\item[2)] $\dim (N_{\overline{z}}^0\ominus R(S)) = \dim(N_z^0\ominus N_z')$;

\item[3)] $S$ has no non-zero fixed elements in $R(X_{z;0})$.
\end{itemize}

By Proposition~\ref{p2_1_p2_1} and the first equality in~(\ref{f2_1_p2_1}), there exists an isometric operator, which maps
$N_{\overline{z}}'$ on $K$, and has no non-zero fixed elements.
By the second equality in~(\ref{f2_1_p2_1}) condition~2) for this operator is satisfied, as well.
Therefore {\it an isometric operator correcting $X_{z;0}$ in $N_{\overline{z}}^0$
always exists}.

\begin{thm}
\label{t2_2_p2_1}
Let $A$ be a closed symmetric operator in a Hilbert space $H$, $z\in \mathbb{R}$ be a fixed point, $N_z^0\in \mathcal{N}_z(A)$
and $N_{\overline{z}}^0\in \mathcal{N}_{\overline{z}}(A)$ be subspaces of the same dimension.

An arbitrary admissible with respect to $A$ isometric operator $V$, which maps
$N_z^0$ on $N_{\overline{z}}^0$, has the following form:
\begin{equation}
\label{f2_2_p2_1}
V = S \overline{X_{z;0}} \oplus T,
\end{equation}
where $S$ is an isometric operator, which corrects $X_{z;0}$ in $N_{\overline{z}}^0$,
and $T$ is an isometric operator with the domain $D(T) = N_z^0\ominus N_z'$
and the range $R(T) = N_{\overline{z}}^0\ominus R(S)$.

Conversely, an arbitrary isometric operator $S$, which corrects $X_{z;0}$ in $N_{\overline{z}}^0$,
and an arbitrary isometric operator $T$ with the domain $D(T) = N_z^0\ominus N_z'$
and the range $R(T) = N_{\overline{z}}^0\ominus R(S)$, generate by relation~(\ref{f2_2_p2_1})
an admissible with respect to $A$ isometric operator $V$, which maps
$N_z^0$ on $N_{\overline{z}}^0$.
\end{thm}
{\bf Proof. }
Let us check the first assertion of the Theorem. Let $V$  be an admissible with respect to $A$ isometric operator, which
maps $N_z^0$ on $N_{\overline{z}}^0$.
Denote
$$ Q = V|_{N_z'},\qquad T = V|_{N_z^0\ominus N_z'}, $$
and
$$ S = Q (\overline{X_{z;0}})^{-1}:\ N_{\overline{z}}'\rightarrow V N_z'. $$
Then
$$ V = S \overline{X_{z;0}} \oplus T. $$
The operator $S$ is isometric. Let us check that it corrects $X_{z;0}$ in $N_{\overline{z}}^0$.
Since
$$ \dim (N_{\overline{z}}^0\ominus R(S)) =
\dim (N_{\overline{z}}^0\ominus (V N_z')) =
\dim (V (N_z^0\ominus N_z')) =
\dim(N_z^0\ominus N_z'), $$
then the second condition from the definition of the correcting operator is satisfied.
Let $\varphi\in R(X_{z;0})$ is a fixed element of the operator $S$:
$S\varphi = Q (\overline{X_{z;0}})^{-1} \varphi = \varphi$.
Denote $\psi = (\overline{X_{z;0}})^{-1} \varphi = X_{z;0}^{-1} \varphi \in M$. Then
$V\psi = Q\psi = X_{z;0} \psi$.
Since $V$  is admissible with respect to $A$ isometric operator then
$\psi = 0$, and $\varphi = \overline{X_{z;0}} \psi = 0$.
Thus, the operator $S$ corrects $X_{z;0}$ in $N_{\overline{z}}^0$.

Conversely, let $S$ be an arbitrary isometric operator which corrects $X_{z;0}$ in $N_{\overline{z}}^0$,
and $T$ be an arbitrary isometric operator with the domain $D(T) = N_z^0\ominus N_z'$
and the range $R(T) = N_{\overline{z}}^0\ominus R(S)$. Define the operator $V$ by relation~(\ref{f2_2_p2_1}).
The operator $V$ is isometric and maps $N_z^0$ on $N_{\overline{z}}^0$.
Let us check that $V$ is admissible with respect to $A$.
Consider an arbitrary element $\psi\in D(V)\cap D(X_z)$ such that
$V\psi = X_z\psi$.

At first, we consider the case $\psi\in N_z'$. In this case we have $S\overline{X_{z;0}}\psi = X_z\psi$.
By the definition of the closure of an operator and by the continuity of $X_z$, we can assert that there exists a sequence
$\psi_n\in D(X_{z;0})$, $n\in \mathbb{N}$, such that
$\psi_n\rightarrow\psi$, and $X_z\psi_n = X_{z;0}\psi_n\rightarrow \overline{X_{z;0}}\psi = X_z\psi$,
as $n\rightarrow\infty$.
Therefore $S X_z\psi = X_z\psi$. By the definition of the correcting operator, $S$ can not have
non-zero fixed elements in $R(X_{z;0})$. Two cases are posiible:

\noindent
1) $X_z\psi\in R(X_{z;0})$ и $X_z\psi = 0$. By the invertibility of $X_z$ this means that $\psi = 0$.

\noindent
2) $X_z\psi\notin R(X_{z;0})$. This means that $\psi\notin D(X_{z;0})$.
Since $\psi\in D(X_z)\cap N_z'$,  then it is possible only in the case $X_z\psi\notin N_{\overline{z}}^0$.
But this is impossible since $X_z\psi = V\psi \in N_{\overline{z}}^0$. Thus, this case should be excluded.

Consider the case $\psi\notin N_z'$. In this case $\psi\notin D(X_{z;0})$ and
$\psi\in N_z^0 = D(V)$. By the definition of $X_{z;0}$, it is possible only if $X_z\psi\notin N_{\overline z}^0$.
But then the equality $V\psi = X_z\psi$ will be impossible, since
$V\psi\in N_{\overline{z}}^0$. Consequently, the case $\psi\notin N_z'$ is impossible.
Thus, we obtain that the operator $V$ is admissible with respect to $A$.
$\Box$

\begin{cor}
\label{c2_1_p2_1}
Let $A$ be a closed symmetric operator in a Hilbert space $H$.
The following assertions hold:
\begin{itemize}
\item[(i)]
The operator $A$ is maximal if and only if (at least) one of its
defect numbers is equal to zero.
\item[(ii)]
If $A$ is maximal then $\overline{D(A)} = H$.
\end{itemize}
\end{cor}
{\bf Proof. }
(i): {\it Necessity. } Suppose to the contrary that $A$ is maximal but
its both defect numbers are non-zero. Choose and fix an arbitrary number
$z\in \mathbb{R}_e$.
Let $N_z^0\subseteq \mathcal{N}_z(A)$ and
$N_{\overline{z}}^0\subseteq \mathcal{N}_{\overline{z}}(A)$ be some one-dimensional subspaces.
By Theorem~\ref{t2_2_p2_1}, there exists an operator $V$,
which is admissible with respect to $A$ and which maps $N_z^0$ on $N_{\overline{z}}^0$.
By Remark~\ref{r2_1_p2_1} it follows that to the operator $V$ there corresponds a symmetric extension
$B$ of the operator $A$, according to relations~(\ref{f2_0_p2_1})-(\ref{f2_0_2_p2_1}).
This extension is proper since to $A$ by~(\ref{f2_0_p2_1})-(\ref{f2_0_2_p2_1}) there corresponds the null isometric operator $T$.
We obtained a contradiction.

\noindent
(ii): {\it Sufficiency. } Fix an arbitrary number $z\in \mathbb{R}_e$.
If one of the defect numbers is equal to zero, then the Cayley transformation
$U_z(A)$ does not have proper extensions. Consequently, taking into account considerations after~(\ref{f2_0_1_p2_1}),
the operator $A$ has no proper extensions.

\noindent
(ii). Let $A$ be maximal. From the proved assertion~(i) it follows that
one of the defect numbers is equal to zero, i.e. $\mathcal{N}_z(A) = \{ 0 \}$, for some
$z\in \mathbb{R}_e$. By Corollary~\ref{c1_2_p2_1} we may write:
$$ H\ominus\overline{D(A)} = \left( H\ominus\overline{D(A)} \right)\cap \mathcal{M}_z(A)
= \{ 0 \}, $$
и значит $\overline{D(A)} = H$.
$\Box$

\subsection{Properties of dissipative and accumulative operators.}

A linear operator $A$ in a Hilbert space $H$ is said to be {\bf dissipativeдиссипативным (accumulative)},
if $\mathop{\rm Im}\nolimits (A h,h)_H \geq 0$
(respectively $\mathop{\rm Im}\nolimits (A h,h)_H \leq 0$) for all elements $h$ from $D(A)$.
A dissipative (accumulative) operator $A$ is said to be {\bf maximal}, if it has no proper
(i.e. different from $A$) dissipative (respectively accumulative) extensions in $H$.

We establish some properties of dissipative and accumulative operators, which will be used later.

\begin{thm}
\label{t3_1_p2_1}
Let $A$ be a closed dissipative (accumulative) operator in a Hilbert space $H$.
The following assertions are true:

\begin{itemize}
\item[1)] Points of $\mathbb{C}_-$ (respectively $\mathbb{C}_+$) are points of the regular type
of the operator $A$. The following inequality holds:
\begin{equation}
\label{f3_4_p2_1}
\left\| ( A - \xi E_H )^{-1} \right\| \leq \frac{1}{|\mathop{\rm Im}\nolimits \xi|},\qquad
\xi\in \mathbb{C}_-\ (\mathbb{C}_+).
\end{equation}
\item[2)] Choose and fix a number $z$ from $\mathbb{C}_-$ (respectively from $\mathbb{C}_-$).
Let $B$ be a dissipative (respectively accumulative) extension of the operator $A$ in the space $H$.
Consider the following operator
\begin{equation}
\label{f3_5_p2_1}
W_z = (B - \overline{z} E_H)(B-zE_H)^{-1} = E_H + (z-\overline{z}) (B-zE_H)^{-1}.
\end{equation}
The operator $W_z$ is a non-expanding extension of the operator $U_z(A) = (A - \overline{z} E_H)(A-zE_H)^{-1}$ в $H$,
having no non-zero fixed elements. Moreover, we have:
\begin{equation}
\label{f3_6_p2_1}
B = (z W_z - \overline{z} E_H) (W_z - E_H)^{-1} = zE_H + (z-\overline{z})(W_z - E_H)^{-1}.
\end{equation}
Conversely, for an arbitrary non-expanding operator $W_z\supseteq U_z(A)$, having no non-zero fixed elements,
by~(\ref{f3_6_p2_1}) one defines a dissipative (respectively accumulative) operator
$B\supseteq A$.
Formula~(\ref{f3_5_p2_1}) establishes a one-to-one correspondence between all dissipative (respectively accumulative)
extensions $B\supseteq A$ in $H$, and all non-expanding
extensions $W_z\supseteq U_z(A)$  в $H$, having no non-zero fixed elements.

\item[3)] If $A$ is maximal, then $(A-zE_H) D(A) = H$, for all points $z$ from $\mathbb{C}_-$
(respectively $\mathbb{C}_+$). Conversely, if there exists a point $z_0$
from $\mathbb{C}_-$ (respectively $\mathbb{C}_+$) such that $(A-z_0 E_H) D(A) = H$, then
the operator $A$ is maximal;

\item[4)] If $A$ is maximal, then $\overline{D(A)} = H$.
\end{itemize}
\end{thm}
{\bf Proof. }
1) Consider the case of a dissipative operator $A$. Choose an arbitrary point $z$ from $\mathbb{C}_-$,
$z=x+iy$, $x\in \mathbb{R}$, $y<0$.
For an arbitrary element $f\in D(A)$ we may write:
$$ \| (A - zE_H) f \|_H^2 = ((A-xE_H)f - iyf,(A-xE_H)f - iyf)_H $$
$$ =
\| (A-xE_H)f \|^2_H - 2y \mathop{\rm Im}\nolimits (Af,f)_H + y^2 \| f \|_H^2
\geq y^2 \| f \|_H^2. $$
Therefore the operator $A - zE_H$ has a bounded inverse and inequality~(\ref{f3_4_p2_1}) holds.

\noindent
In the case of an accumulative operator $A$, the operator $-A$ is dissipative and we may apply to it the proved part of the
assertion.

\noindent
2) Let $z$ from $\mathbb{C}_-$ (respectively from $\mathbb{C}_+$) and
$B$ be a dissipative (respectively accumulative) extension of the operator $A$ in the space $H$.
By the proved assertion~1), the operator $W_z$ is correctly defined and bounded.
For an arbitrary element $g\in D(W_z) = (B-zE_H)D(B)$, $g=(B-zE_H)f$, $f\in D(B)$, we may write:
$$ \| W_z g \|_H^2 = (W_z g, W_z g)_H = ((B-\overline{z}E_H)f, (B-\overline{z}E_H)f)_H $$
$$ = \| Bf \|_H^2 - z(Bf,f)_H - \overline{z}(f,Bf)_H + |z|^2 \| f \|_H^2; $$
$$ \| g \|^2_H = ((B-zE_H)f, (B-zE_H)f)_H $$
$$ = \| Bf \|_H^2 - \overline{z}(Bf,f)_H - z(f,Bf)_H + |z|^2 \| f \|_H^2, $$
and therefore
$$ \| W_z g \|_H^2 - \| g \|^2_H = (\overline{z} - z)((Bf,f)_H - (f,Bf)_H) =
4 (\mathop{\rm Im}\nolimits z) \mathop{\rm Im}\nolimits (Bf,f)_H \leq 0. $$
Since $W_z - E_H = (z-\overline{z}) (B-zE_H)^{-1}$, then $W_z$ has no non-zero fixed elements.
Formula~(\ref{f3_6_p2_1}) follows directly from~(\ref{f3_5_p2_1}).

Conversely, let $W_z\supseteq U_z(A)$ be a non-expanding operator having no non-zero
fixed elements. Define an operator $B$ by equality~(\ref{f3_6_p2_1}).
For an arbitrary element $f\in D(B) = (W_z - E_H)D(W_z)$, $f=(W_z-E_H)g$, $g\in D(W_z)$, we write:
$$ (Bf,f)_H = ((zW_z - \overline{z}E_H)g, (W_z - E_H)g )_H = z\| W_z g \|_H^2 - z(W_z g,g)_H $$
$$ - \overline{z} (g,W_z g)_H + \overline{z}\| g \|_H^2; $$
$$ \mathop{\rm Im}\nolimits (Bf,f)_H = \mathop{\rm Im}\nolimits (z) \| W_z g \|_H^2 +
\mathop{\rm Im}\nolimits (\overline{z}) \| g \|_H^2
= \mathop{\rm Im}\nolimits (z) \left( \| W_z g \|_H^2 - \| g \|_H^2 \right). $$
Consequently, the operator $B$ is dissipative (respectively accumulative).
Since $W_z\supseteq U_z(A)$, then $B\supseteq A$.
The correspondence, established by formulas~(\ref{f3_5_p2_1}),(\ref{f3_6_p2_1}), is obviously one-to-one.

\noindent
3) Suppose to the contrary that the operator $A$ is maximal but there exists a number $z$ from $\mathbb{C}_-$
(respectively $\mathbb{C}_-$) such that $(A-z E_H) D(A)\not= H$. Consider the following operator
$U_z = U_z(A) = (A - \overline{z} E_H)(A-zE_H)^{-1}$
By the proved assertion~2), the operator
$U_z$ is non-expanding and has no non-zero fixed elements.
Moreover, $U_z$ is closed, since $A$ is closed by assumption. Therefore $U_z$ is defined on
a subspace $H_1 = (A-z E_H) D(A)\not= H$. Set $H_2 = H\ominus H_1$, $\dim H_2\geq 1$.
Consider the following operator
$$ W_z := U_z\oplus 0_{H_2}. $$
The operator $W_z$ is a non-expanding extension of the operator $U_z$. Let us check that the operator $W_z$
has no non-zero fixed elements. Suppose to the contrary that there exists a non-zero
element $h\in H$ such that $W_z h = h$. Then $h\notin D(U_z)$, since
$U_z$ has no non-zero fixed elements. Let $h = h_1 + h_2$, $h_1\in H_1$, $h_2\in H_2$,
and $h_2\not= 0$. Then
$$ W_z h = (U_z\oplus 0_{H_2}) (h_1+h_2) = U_z h_1 = h = h_1 + h_2. $$
Therefore
$$ \| h_1 \|_H^2 \geq \| U_z h_1 \|_H^2 = \| h_1 \|_H^2 + \| h_2 \|_H^2, $$
and we get $h_2 = 0$. The obtained contradiction shows that $W_z$
has no non-zero fixed elements.
By the proved assertion~2) to the operator $W_z\not= U_z$ there corresponds
a dissipative (respectively accumulative) extension
$B\not= A$ of the operator $A$. This contradicts to the maximality of the operator $A$, and the first part of
assertion~3) is proved.

Suppose that we have a closed dissipative (respectively accumulative) operator $A$
in a Hilbert space $H$, and there exists a point $z_0$
from $\mathbb{C}_-$ (respectively $\mathbb{C}_+$) such that $(A-z_0 E_H) D(A) = H$.
Suppose to the contrary that $A$ is not maximal. Let
$B$ be a proper dissipative (respectively accumulative) extension of the operator $A$.
By the proved assertion~2), for the operator $B$ there corresponds a non-expanding operator $W_{z_0}$
from~(\ref{f3_5_p2_1}) with $z=z_0$, having no non-zero fixed elements. Moreover, we have
$D(W_{z_0}) = (B-z_0 E_H) D(B)\supseteq (A-z_0 E_H) D(A) = H$. Therefore $D(W_{z_0}) = H = D(U_{z_0})$,
and we get $W_{z_0} = U_{z_0}$, $B=A$. We obtained a contradiction, since $B$ is a proper extension.

\noindent
4) Consider an arbitrary element $h\in H$: $h\perp D(B)$. Choose and fix an arbitrary
number $z$ from $\mathbb{C}_-$ (respectively $\mathbb{C}_+$). By the proved assertion~3),
the following equality holds: $(A - zE_H) D(A) = H$. Therefore there exists an element $f\in D(A)$ such that
$(A-zE_H) f = h$. Then
$$ 0 = ( (A - z E_H) f, f)_H = (Af,f)_H - z\| f \|_H^2; $$
$$ \mathop{\rm Im}\nolimits (Af,f)_H = \mathop{\rm Im}\nolimits (z) \| f \|_H^2.  $$
Since the left-hand side of the last equality is non-negative (respectively non-positive),
and the right-hand side is non-positive (respectively non-negative), then $f = 0$.
$\Box$

\begin{rmr}
\label{r3_1_p2_1}
Consider a closed symmetric operator $A$ in a Hilbert space $H$.
Choose and fix a number $z$ from $\mathbb{C}_-$ ($\mathbb{C}_-$).
Let $B$ be a dissipative (respectively accumulative) extension of the operator $A$ in $H$.
By assertion~2) of Theorem~\ref{t3_1_p2_1} for it there corresponds a non-expanding operator
$W_z$ from~(\ref{f3_5_p2_1}), having no non-zero fixed elements.
By Proposition~\ref{p1_3_p1_1},                                        
such extensions $W_z\supseteq U_z(A)$ have the following form:
\begin{equation}
\label{f3_7_p2_1}
W_z = U_z(A)\oplus T,
\end{equation}
where $T$ is a linear nonexpanding operator with the domain $D(T)\subseteq \mathcal{N}_z(A)$
and the range $R(T)\subseteq \mathcal{N}_{\overline{z}}(A)$.
By Theorem~\ref{t2_1_p2_1} it follows that $T$ is admissible with respect to the operator $A$.

Conversely, if we have an arbitrary non-expanding operator $T$,
with the domain $D(T)\subseteq \mathcal{N}_z(A)$,
the range $R(T)\subseteq \mathcal{N}_{\overline{z}}(A)$,
admissible with respect to $A$,
then by~(\ref{f3_7_p2_1})
one defines a non-expanding extension $W_z$ of the operator $U_z(A)$.
By Theorem~\ref{t2_1_p2_1} this extension has no non-zero fixed elements.
Thus, there is a one-to-one correspondence between admissible with respect to
$A$ non-expanding operators $T$ and non-expanding operators $W_z\supseteq U_z(A)$ without
non-zero fixed elements.

Using the proved theoremу we conclude that formulas~(\ref{f3_5_p2_1})-(\ref{f3_7_p2_1}) establish
a one-to-one correspondence between admissible with respect to
$A$ non-expanding operators $T$,
$D(T)\subseteq \mathcal{N}_z(A)$, $R(T)\subseteq \mathcal{N}_{\overline{z}}(A)$,
and dissipative (respectively accumulative)
extensions of the closed symmetric operator $A$.
\end{rmr}

\subsection{Generalized Neumann's formulas.}

The classical Neumann's formulas describe all symmetric extensions of a given closed symmetric
operator $A$ in a Hilbert space $H$ in the case $D(A)=H$.
The following theorem provides a description such extensions without the assumption $D(A)=H$,
and it describes dissipative and accumulative extensions of $A$, as well.

\begin{thm}
\label{t4_1_p2_1}
Let $A$ be a closed symmetric operator in a Hilbert space $H$, and
$z$ from $\mathbb{C}_-$ ($\mathbb{C}_+$) be a fixed number.
The following formulas
\begin{equation}
\label{f4_1_p2_1}
D(B) = D(A) \dotplus (T-E_H) D(T),
\end{equation}
\begin{equation}
\label{f4_2_p2_1}
B ( f + T\psi - \psi ) = Af + zT\psi - \overline{z}\psi,\qquad f\in D(A),\ \psi\in D(T),
\end{equation}
establish a one-to-one correspondence between all admissible with respect to $A$ isometric
operators $T$, $D(T)\subseteq \mathcal{N}_z(A)$, $R(T)\subseteq \mathcal{N}_{\overline{z}}(A)$,
and all symmetric extensions $B$ of the operator $A$.
Moreover, we have
\begin{equation}
\label{f4_3_p2_1}
D(T) = \mathcal{N}_z(A)\cap R(B-zE_H),
\end{equation}
\begin{equation}
\label{f4_4_p2_1}
T\subseteq (B-\overline{z}E_H)(B-zE_H)^{-1}.
\end{equation}
A symmetric operator $B$ is: closed / closed and maximal / self-adjoint,
if and only if
respectively: $D(T)$ is: a subspace / $D(T)=\mathcal{N}_z(A)$ or $R(T)=\mathcal{N}_{\overline z}(A)$ /
$D(T)=\mathcal{N}_z(A)$ and $R(T)=\mathcal{N}_{\overline z}(A)$.

Formulas~(\ref{f4_1_p2_1}),(\ref{f4_2_p2_1})
define a one-to-one correspondence between all admissible with respect to $A$ non-expanding
operators $T$, $D(T)\subseteq \mathcal{N}_z(A)$, $R(T)\subseteq \mathcal{N}_{\overline{z}}(A)$,
and all dissipative (respectively accumulative) extensions $B$ of the operator $A$.
Also relations~(\ref{f4_3_p2_1}),(\ref{f4_4_p2_1}) hold.
A dissipative (respectively accumulative) operator $B$: closed / closed and maximal,
if and only if
respectively: $D(T)$ is: a subspace / $D(T)=\mathcal{N}_z(A)$.
\end{thm}
{\bf Proof. }
Let $A$ be a closed symmetric operator in a Hilbert space $H$, and
$z$ be from $\mathbb{C}_-$ ($\mathbb{C}_+$).
By Remark~\ref{r2_1_p2_1}, formulas~(\ref{f2_0_p2_1}),(\ref{f2_0_1_p2_1}) and
(\ref{f2_0_2_p2_1}) establish a one-to-one correspondence between all symmetric
extensions $B$ of the operator $A$ in $H$, and all isometric
operators $T$  with the domain $D(T)\subseteq \mathcal{N}_z(A)$
and the range $R(T)\subseteq \mathcal{N}_{\overline{z}}(A)$, which are admissible with respect to $A$.
Let us check that formulas~(\ref{f2_0_1_p2_1}),(\ref{f2_0_2_p2_1}) define the same operator $B$,
as formulas~(\ref{f4_1_p2_1}),(\ref{f4_2_p2_1}).
Let the operator $B$ be defined by formulas~(\ref{f2_0_1_p2_1}),(\ref{f2_0_2_p2_1}).
The domain of $B$ is
$$ D(B) = (W_z - E_H) D(W_z) = (U_z\oplus T - E_H)(D(U_z) + D(T)) $$
$$ = (U_z - E_H)D(U_z) + (T-E_H)D(T). $$
If $f\in (U_z - E_H)D(U_z) \cap (T-E_H)D(T)$, then
$f = (U_z - E_H)g$, $g\in D(U_z)$, and $f = (T-E_H) h$, $h\in D(T)$.
Then
$$ (U_z\oplus T - E_H) (g - h) = f - f = 0, $$
i.e. $g-h$ is a fixed element of the operator $W_z$. Since $W_z$ has no non-zero fixed
elements, we get $g=h$. Since $g\perp h$, then $g=0$ and $f=0$. Thus, for the operator
$B$ holds~(\ref{f4_1_p2_1}).

\noindent
Consider arbitrary elements $f\in D(A)$ and $\psi\in D(V)$. Then
$$ B(f + V\psi - \psi) = Bf + B(V-E_H)\psi = Af + B(W_z - E_H)\psi $$
$$ = Af + (zW_z - \overline{z}E_H)\psi = Af + zV\psi - \overline{z}\psi, $$
and therefore formula~(\ref{f4_2_p2_1}) holds, as well.

\noindent
Since
$$ \mathcal{N}_z(A)\cap (B-zE)D(B) = \mathcal{N}_z(A)\cap D(W_z) =
\mathcal{N}_z(A)\cap (\mathcal{M}_z(A) + D(T)) $$
$$ = \mathcal{N}_z(A)\cap D(T) = D(T), $$
then relation~(\ref{f4_3_p2_1}) holds.
Relation~(\ref{f4_4_p2_1}) means that $T\subseteq W_z$.

If the operator $B$ is closed, then $W_z$ is closed, as well. Since $W_z$ is bounded, it is defined
on a subspace. Therefore $D(W_z) = R(B-zE_H)$ is closed, $D(T) = \mathcal{N}_z(A)\cap R(B-zE_H)$ is closed
and it is a subspace.

Conversely, if $D(T)$ is a subspace then the direct sum $\mathcal{M}_z(A)\oplus D(T) = D(W_z)$
is closed. Therefore the operators $W_z$ and $B$ are closed.

By Corollary~\ref{c2_1_p2_1} a closed symmetric operator $B$ in a Hilbert space $H$
is maximal if and only if (at least) one of its defect numbers is equal to zero.
This is equivqlent to the following condition: $D(T)=\mathcal{N}_z(A)$ or
$R(T)=\mathcal{N}_{\overline z}(A)$.

If $B$ if self-adjoint then its Cayley's transformation $W_z$ is unitary, i.e.
$D(T)=\mathcal{N}_z(A)$ and $R(T)=\mathcal{N}_{\overline z}(A)$. Conversely, the last conditions imply that
$W_z$ is unitary, and therefore $B$ is self-adjoint.

By Remark~\ref{r3_1_p2_1} formulas~(\ref{f3_5_p2_1})-(\ref{f3_7_p2_1})
establish a one-to-one correspondence between admissible with respect to $A$ non-expanding operators $T$,
$D(T)\subseteq \mathcal{N}_z(A)$, $R(T)\subseteq \mathcal{N}_{\overline{z}}(A)$,
and dissipative (respectively accumulative)
extensions of a closed symmetric operator $A$.
Formulas~(\ref{f4_1_p2_1}),(\ref{f4_2_p2_1}) define the same operator $B$ as formulas~(\ref{f3_5_p2_1})-(\ref{f3_7_p2_1}).
It is checked similar as it was done above for the case of symmetric extensions.
The proof of relations~(\ref{f4_3_p2_1}),(\ref{f4_4_p2_1}), and the proof of
equivalence: ($B$ is closed) $\Leftrightarrow$ ($B$ is a subspace) are the same.

Let $B$ be closed and maximal.
By assertion~3) of Theorem~\ref{t3_1_p2_1} we obtain that $R(B-zE_H) = (B-zE_H)D(B) =H$.
From~(\ref{f4_3_p2_1}) it follows that $D(T) = \mathcal{N}_z(A)$.

Conversely, let $D(T) = \mathcal{N}_z(A)$. Then
$$ (B-zE_H)D(B) = D(W_z) = \mathcal{M}_z(A)\oplus D(T) = \mathcal{M}_z(A)\oplus \mathcal{N}_z(A) = H. $$
By assertion~3) of Theorem the operator $B$ is maximal.
$\Box$

Let $A$ be a closed symmetric operator in a Hilbert space $H$, and
$z\in\mathbb{R}_e$  be a  fixed number.
Let $T$ be an admissible with respect to $A$ non-expanding
operator with $D(T)\subseteq \mathcal{N}_z(A)$, $R(T)\subseteq \mathcal{N}_{\overline{z}}(A)$.
An extension $B$ of the operator $A$ which corresponds to $T$ in formulas~(\ref{f4_1_p2_1}),(\ref{f4_2_p2_1})
we denote by $A_T=A_{T,z}$ and call {\bf a quasi-self-adjoint extension of a symmetric operator $A$
defined by the operator $T$}. This definition agrees with the above given definition for the case
of a densely defined symmetric operator $A$.

\subsection{Extensions of a symmetric operator with an exit out of the space.}

Let $A$ be a closed symmetric operator in a Hilbert space $H$.
In order to construct a generalized resolvent of $A$, according to its definition,
one needs to construct self-adjoint extensions $\widetilde A$ of $A$ in larger Hilbert spaces $\widetilde H\supseteq H$.

Choose and fix a Hilbert space $\widetilde H\supseteq H$. Decompose it as a direct sum:
\begin{equation}
\label{f5_1_p2_1}
\widetilde H = H\oplus H_e,
\end{equation}
where $H_e = \widetilde H\ominus H$.
The operator $A$ can be identified with the operator $A\oplus o_{H_e}$ in $\widetilde H$,
where $o_{H_e}$ is an operator in $H_e$ with $D(o_{H_e})=\{ 0 \}$, $o_{H_e} 0 = 0$.

Keeping in mind this important case, we shall now assume that in a Hilbert space
$\widetilde H$ of a type~(\ref{f5_1_p2_1}), where $H$, $H_e$ are some Hilbert spaces, there is given a
symmetric operator $\mathcal{A}$ of the following form:
\begin{equation}
\label{f5_2_p2_1}
\mathcal{A} = A\oplus A_e,
\end{equation}
where $A$, $A_e$ are symmetric operators in Hilbert spaces $H$,$H_e$, respectively.
Let us investigate a possibility for the construction of symmetric and self-adjoint extensions of the operator $\mathcal{A}$
and their properties. In particular, we can use the generalized Neumann formulas.

Fix an arbitrary number $z\in \mathbb{R}_e$.
Since $D(\mathcal{A}) = D(A)\oplus D(A_e)$, then it can be directly verified that
\begin{equation}
\label{f5_3_p2_1}
\mathcal{M}_z(\mathcal{A}) = \mathcal{M}_z(A) \oplus \mathcal{M}_z(A_e),
\end{equation}
\begin{equation}
\label{f5_4_p2_1}
\mathcal{N}_z(\mathcal{A}) = \mathcal{N}_z(A) \oplus \mathcal{N}_z(A_e),
\end{equation}
\begin{equation}
\label{f5_5_p2_1}
X_z(\mathcal{A}) = X_z(A) \oplus X_z(A_e),
\end{equation}
where $X_z(\cdot)$ denotes the forbidden operator.

Let $T$ be an arbitrary non-expanding operator in $\widetilde H$ with the domain
 $D(T) = \mathcal{N}_z(\mathcal{A})$ and the range $R(T) \subseteq \mathcal{N}_{\overline z}(\mathcal{A})$.
With respect to the decomposition~(\ref{f5_4_p2_1}) (for $z$ and $\overline{z}$) the operator $T$ has the following
block representation:
\begin{equation}
\label{f5_6_p2_1}
T = \left(
\begin{array}{cc} T_{11} & T_{12}\\
T_{21} & T_{22}\end{array}
\right),
\end{equation}
where $T_{11} = P_{\mathcal{N}_{\overline{z}}(A)} T P_{\mathcal{N}_z(A)}$,
$T_{12} = P_{\mathcal{N}_{\overline{z}}(A)} T P_{\mathcal{N}_z(A_e)}$,
$T_{21} = P_{\mathcal{N}_{\overline{z}}(A_e)} T P_{\mathcal{N}_z(A)}$,
$T_{22} = P_{\mathcal{N}_{\overline{z}}(A_e)} T P_{\mathcal{N}_z(A_e)}$.

\noindent
By the definition of the admissible operator the operator $T$ is admissible with respect to operator
$\mathcal{A}$ if and only if relation
$$ \psi\in D(X_z(\mathcal{A})),\quad X_z(\mathcal{A})\psi = T\psi, $$
implies $\psi = 0$. Using decompositions~(\ref{f5_5_p2_1}),(\ref{f5_6_p2_1}) we see that the latter condition
is equivalent to the following one: relation
$$ \psi_1\in D(X_z(A)),\ \psi_2\in D(X_z(A_e)), $$
\begin{equation}
\label{f5_7_p2_1}
T_{11}\psi_1 + T_{12}\psi_2 = X_z(A)\psi_1,\
T_{21}\psi_1 + T_{22}\psi_2 = X_z(A_e)\psi_2,
\end{equation}
implies $\psi_1 = \psi_2 = 0$.

\begin{thm}
\label{t5_0_p2_1}
Let $\mathcal{A}$ be a closed symmetric operator of the form~(\ref{f5_2_p2_1}),
where $A$, $A_e$ are symmetric operators in some Hilbert spaces $H$ and $H_e$, respectively.
Let $z\in \mathbb{R}_e$ be a fixed number and
$T$ be a non-expanding operator in $\widetilde H = H\oplus H_e$,
with the domain $D(T) = \mathcal{N}_z(\mathcal{A})$ and
the range $R(T) \subseteq \mathcal{N}_{\overline z}(\mathcal{A})$.
The following assertions are true:
\begin{itemize}
\item[1)]
If $T$ is $z$-admissible with respect to the operator $\mathcal{A}$, then the operators $T_{11}$ and $T_{22}$
are $z$-admissible  with respect to operators $A$ and $A_e$, respectively.

\item[2)]
If $\overline{D(A)}=H$ and $T_{22}$ is $z$-admissible with respect to the operator $A_e$, then
the operator $T$ is $z$-admissible with respect to the operator $\mathcal{A}$.
\end{itemize}
Here operators $T_{11}$ and $T_{22}$
are the operators from the block representation~(\ref{f5_6_p2_1}) for $T$.
\end{thm}
{\bf Proof. }
Let us check assertion~1).
Let $\psi_1$ be an element from $D(X_z(A))$ such that
$T_{11} \psi_1 = X_z(A) \psi_1$. Since $T$ is non-expanding and the forbidden operator
is isometric, we may write:
$$ \| \psi_1 \|_H^2 \geq \| T\psi_1 \|^2 = \| T_{11} \psi_1 \|_H^2 + \| T_{21} \psi_1 \|_H^2 $$
$$ = \| X_z(A) \psi_1 \|_H^2 + \| T_{21} \psi_1 \|_H^2 =
\| \psi_1 \|_H^2 + \| T_{21} \psi_1 \|_H^2. $$
Therefore $T_{21} \psi_1 = 0$. Set $\psi_2 = 0$. We see that condition~(\ref{f5_7_p2_1})
is satisfied. By the condition of the theorem $T$ is admissible with respect to $A$, and therefore $\psi_1=0$.
This means that the operator $T_{11}$ is admissible with respect to $A$.

\noindent
In a similar way, if $\psi_2$ is an element from $D(X_z(A_e))$ such that
$T_{22} \psi_2 = X_z(A_e) \psi_2$, then
$$ \| \psi_2 \|_H^2 \geq \| T\psi_2 \|^2 = \| T_{12} \psi_2 \|_H^2 + \| T_{22} \psi_2 \|_H^2 $$
$$ = \| T_{12} \psi_2 \|_H^2 + \| X_z(A_e) \psi_2 \|_H^2 =
\| T_{12} \psi_2 \|_H^2 + \| \psi_2 \|_H^2. $$
Therefore $T_{12} \psi_2 = 0$. Set $\psi_1 = 0$. We conclude that condition~(\ref{f5_7_p2_1})
holds. Consequently, we get $\psi_2 = 0$.

Let us check assertion~2). Let $\psi_1\in D(X_z(A))$, $\psi_2\in D(X_z(A_e))$, and relation~(\ref{f5_7_p2_1}) holds.
Since $\overline{D(A)}=H$, then by~(\ref{f1_5_p2_1})
the domain of the operator $X_z(A)$ consists of the null element and
$\psi_1 = 0$. The second equality in~(\ref{f5_7_p2_1}) may be written as
$$ T_{22}\psi_2 = X_z(A_e)\psi_2. $$
Since $T_{22}$ is admissible with respect to $A_e$, we get $\psi_2 = 0$.
Consequently, $T$ is admissible with respect to $A$.
$\Box$

Let us continue our considerations started before the formulation of the last theorem.
Let $\psi_1\in \mathcal{N}_z(A)$ and $\psi_2\in D(X_z(A_e))$ be such elements (for example, null elements) that
\begin{equation}
\label{f5_8_p2_1}
T_{21}\psi_1 + T_{22}\psi_2 = X_z(A_e)\psi_2.
\end{equation}
Then
\begin{equation}
\label{f5_9_p2_1}
\| T_{11}\psi_1 + T_{12}\psi_2 \|_H \leq \| \psi_1 \|_H.
\end{equation}
In fact, since $T$ is non-expanding, we have
$$ \| T(\psi_1 + \psi_2) \|_H^2 = \| T_{11}\psi_1 + T_{12}\psi_2 \|_H^2 +
\| T_{21}\psi_1 + T_{22}\psi_2 \|_H^2 $$
$$ \leq \| \psi_1 \|_H^2 + \| \psi_2 \|_H^2. $$
If~(\ref{f5_8_p2_1}) holds, then since $X_z(A_e)$ is isometric we get~(\ref{f5_9_p2_1}).
Moreover, if an equality $\| T(\psi_1 + \psi_2) \|_H = \| \psi_1 + \psi_2 \|_H$ holds,
then in~(\ref{f5_9_p2_1}) we have an equality, and vice versa.

Define the following operator:
$$ \Phi \psi_1 = \Phi(z;\mathcal{A},T) \psi_1 = T_{11}\psi_1 + T_{12}\psi_2, $$
where $\psi_1\in \mathcal{N}_z(A)$ is such an element, for which there exists $\psi_2\in D(X_z(A_e))$:
$$ T_{21}\psi_1 + T_{22}\psi_2 = X_z(A_e)\psi_2. $$
Let us check that this definition is correct. If for $\psi_1\in \mathcal{N}_z(A)$ there exists another element
$\widetilde\psi_2\in D(X_z(A_e))$:
$$ T_{21}\psi_1 + T_{22}\widetilde\psi_2 = X_z(A_e)\widetilde\psi_2, $$
then subtracting last relations we obtain that
$$ T_{22}(\psi_2 - \widetilde\psi_2) = X_z(A_e)(\psi_2 - \widetilde\psi_2). $$
Set $\psi_1 = 0$. We see that equality~(\ref{f5_8_p2_1}) holds, and therefore
inequality~(\ref{f5_9_p2_1}) is true and it takes the following form:
$$ \| T_{12} (\psi_2 - \widetilde\psi_2) \|_H \leq 0, $$
i.e. $T_{12} \psi_2 = T_{12} \widetilde\psi_2$.
Thus, the definition of the operator $\Phi$ does not depend on the choice of $\psi_2$ and it is correct.
By~(\ref{f5_9_p2_1}) the operator $\Phi$ is non-expanding.
Moreover, {\it the operator $\Phi$ is isometric, if $T$ is isometric} (what follows from the considerations before the definition
of $\Phi$).

\begin{rmr}
\label{r5_1_p2_1}
If it is additionally known that the operator  $T_{22}$ is admissible with respect to the operator $A_e$,
then it exists $(X_z(A_e) - T_{22})^{-1}$.
Consider an arbitrary element $\psi_1\in D(\Phi)$. By the definition of $\Phi$,
for the corresponding to
$\psi_1$ the element $\psi_2$ eqality~(\ref{f5_8_p2_1}) holds, and therefore
$\psi_2 = (X_z(A_e) - T_{22})^{-1} T_{21}\psi_1$.
Substitute this expression in the definition of the operator $\Phi$ to obtain the following representation:
\begin{equation}
\label{f5_10_p2_1}
\Phi \psi_1 = T_{11} \psi_1 + T_{12} (X_z(A_e) - T_{22})^{-1} T_{21} \psi_1,\qquad \psi_1\in D(\Phi).
\end{equation}
\end{rmr}

\begin{thm}
\label{t5_1_p2_1}
Let $\mathcal{A}$ be a closed symmetric operator of the form~(\ref{f5_2_p2_1}),
where $A$, $A_e$ are symmetric operators in some Hilbert spaces $H$ and $H_e$, respectively.
Let $z\in \mathbb{R}_e$ be a fixed number and
$T$ be a non-expanding operator in $\widetilde H = H\oplus H_e$,
with the domain $D(T) = \mathcal{N}_z(\mathcal{A})$ and
the range $R(T) \subseteq \mathcal{N}_{\overline z}(\mathcal{A})$.
The operator $T$ is $z$-admissible with respect to $\mathcal{A}$ вif and only if
the operators $\Phi(z;\mathcal{A},T)$ and $T_{22}$
are $z$-admissible with respect to the operators $A$ and $A_e$, respectively.
Here the operators  $T_{22}$
is an operator from the block representation~(\ref{f5_6_p2_1}) for $T$.
\end{thm}
{\bf Proof. }
{\it Necessity.} Let an operator $T$ be admissible with respect to $\mathcal{A}$.
By Theorem~\ref{t5_0_p2_1} the operator $T_{22}$
is admissible with respect to $A_e$.
Consider an element $\psi_1 \in D(\Phi(z;\mathcal{A},T))$ such that
$\Phi(z;\mathcal{A},T) \psi_1 = X_{z}(A) \psi_1$. By the definition of the operator $\Phi(z;\mathcal{A},T)$,
the element $\psi_1$ belongs to $\mathcal{N}_z(A)$ and there exists $\psi_2\in D(X_z(A_e))$:
$T_{21}\psi_1 + T_{22}\psi_2 = X_z(A_e)\psi_2$.
Moreover, we have $\Phi(z;\mathcal{A},T) \psi_1 = T_{11}\psi_1 + T_{12}\psi_2$. Thus, conditions~(\ref{f5_7_p2_1})
are true. Since $T$ is admissible with respect to $\mathcal{A}$,
then these inequalities imply $\psi_1=\psi_2=0$. Therefore the operator $\Phi$ is admissible with respect to $A$.

{\it Sufficiency.} Suppose that $\Phi(z;\mathcal{A},T)$ and $T_{22}$
are admissible with respect to $A$ and $A_e$, respectively.
Consider some elements $\psi_1\in D(X_z(A))$, $\psi_2\in D(X_z(A_e))$, such that~(\ref{f5_7_p2_1}) holds.
This relation means that $\psi_1\in D(\Phi)$ and
$\Phi\psi_1 = X_z(A)\psi_1$. Since $\Phi$ is admissible with respect to $A$, then $\psi_1 = 0$.
From the second equality in~(\ref{f5_7_p2_1}) it follows that
$T_{22}\psi_2 = X_z(A_e)\psi_2$. Since $T_{22}$ is admissible with respect to
$A_e$, we get $\psi_2 = 0$.
Consequently, the operator  $T$ is admissible with respect to the operator $\mathcal{A}$.
$\Box$

Now we shall obtain a more explicit expression for the domain of $\Phi$ and for its action.

\begin{thm}
\label{t5_2_p2_1}
Let $\mathcal{A}$ be a closed symmetric operator of the form~(\ref{f5_2_p2_1}),
where $A$, $A_e$ are symmetric operators in some Hilbert spaces $H$ and $H_e$, respectively.
Let $z\in \mathbb{R}_e$ be a fixed number and
$T$ be a non-expanding operator in  $\widetilde H = H\oplus H_e$,
with the domain $D(T) = \mathcal{N}_z(\mathcal{A})$, and
with the range $R(T) \subseteq \mathcal{N}_{\overline z}(\mathcal{A})$.
Then
\begin{equation}
\label{f5_11_p2_1}
D(\Phi(z;\mathcal{A},T)) = P^{\widetilde H}_{\mathcal{N}_z(A)} \widetilde\Theta_z,
\end{equation}
\begin{equation}
\label{f5_12_p2_1}
\Phi(z;\mathcal{A},T) P^{\widetilde H}_{\mathcal{N}_z(A)} h = P^{\widetilde H}_{\mathcal{N}_{\overline z}(A)}
(U_z(\mathcal{A})\oplus T) h,\qquad h\in \widetilde\Theta_z,
\end{equation}
where
\begin{equation}
\label{f5_14_p2_1}
\widetilde\Theta_z = \left\{
h\in\widetilde H:\ P^{\widetilde H}_{H_e} (U_z(\mathcal{A})\oplus T) h = P^{\widetilde H}_{H_e} h
\right\}.
\end{equation}
\end{thm}
{\bf Proof. }
We shall need the next lemma on the structure of the set $\widetilde\Theta_z$.
\begin{lem}
\label{l5_1_p2_1}
In conditions of Theorem~\ref{t5_2_p2_1}, the set $\widetilde\Theta_z$ consists of elements
$h\in \widetilde H$,
\begin{equation}
\label{f5_15_p2_1}
h = g_1 + \psi_1 + g_2 + \psi_2,
\end{equation}
where $g_1\in \mathcal{M}_z(A)$, $\psi_1\in \mathcal{N}_z(A)$, $g_2\in \mathcal{M}_z(A_e)$,
$\psi_2\in \mathcal{N}_z(A_e)$, such that
\begin{equation}
\label{f5_16_p2_1}
\psi_2\in D(X_z(A_e)),\ X_z(A_e)\psi_2 = T_{21}\psi_1 + T_{22}\psi_2,
\end{equation}
\begin{equation}
\label{f5_17_p2_1}
g_2 = \frac{1}{\overline{z}-z} (A_e - zE_{H_e}) (X_z(A_e) - E_{H_2}) \psi_2.
\end{equation}
\end{lem}
{\it Proof of Lemma.} Consider an arbitrary element $h\in \widetilde\Theta_z$.
As each element of $\widetilde H$, the element $h$ has a (unique) decomposition~(\ref{f5_15_p2_1}).
Since
$$ P^{\widetilde H}_{H_e} h = g_2 + \psi_2, $$
$$ P^{\widetilde H}_{H_e} (U_z(\mathcal{A})\oplus T) h =
P^{\widetilde H}_{H_e} \left( U_z(\mathcal{A}) (g_1 + g_2)  + T (\psi_1 + \psi_2) \right) $$
$$ = U_z(A_e) g_2 + T_{21}\psi_1 + T_{22}\psi_2, $$
then by the definition of the set $\widetilde\Theta_z$ we obtain:
\begin{equation}
\label{f5_18_p2_1}
U_z(A_e) g_2 + T_{21}\psi_1 + T_{22}\psi_2 = g_2 + \psi_2.
\end{equation}
Then
$$ T_{21}\psi_1 + T_{22}\psi_2 - \psi_2 = (E_{H_e} - U_z(A_e)) g_2 = (\overline{z} - z)
(A_e - zE_{H_e})^{-1} g_2 \in D(A_e). $$
Since $T_{21}\psi_1 + T_{22}\psi_2\in \mathcal{N}_{\overline{z}}(A_e)$, then
$\psi_2\in D(X_z(A_e))$ и $X_z(A_e) \psi_2 = T_{21}\psi_1 + T_{22}\psi_2$. From the latter relations
it follows~(\ref{f5_16_p2_1}),(\ref{f5_17_p2_1}).

Conversely, if for an element $h\in H$ of form~(\ref{f5_16_p2_1}) holds~(\ref{f5_16_p2_1}) and
(\ref{f5_17_p2_1}), then
$$ (E_{H_e} - U_z(A_e)) g_2 = (\overline{z} - z) (A_e - zE_{H_e})^{-1} g_2 = (X_z(A_e) - E_{H_e}) \psi_2 $$
$$ = T_{21}\psi_1 + T_{22}\psi_2 - \psi_2. $$
Therefore equality~(\ref{f5_18_p2_1}) holds, which means that $h\in D(\widetilde\Theta_z)$.
$\Box$ (end of the proof of Lemma).

By the proved Lemma it follows that if $\psi_1\in P^{\widetilde H}_{\mathcal{N}_z(A)} \widetilde\Theta_z$, then
$\psi_1\in \mathcal{N}_z(A)$, and there exists element $\psi_2\in D(X_z(A_e)$ such that
$X_z(A_e)\psi_2 = T_{21}\psi_1 + T_{22}\psi_2$.

Conversely, if $\psi_1\in \mathcal{N}_z(A)$, and there exists an element $\psi_2\in D(X_z(A_e))$ such that
$X_z(A_e)\psi_2 = T_{21}\psi_1 + T_{22}\psi_2$, then $X_z(A_e)\psi_2 - \psi_2\in D(A_e)$.
Define $g_2$ by~(\ref{f5_17_p2_1}), and as $g_1$ we take an arbitrary element from
$\mathcal{M}_z(A)$. Then $h := \psi_1 + g_1 + \psi_2 + g_2\in \widetilde\Theta_z$,
and therefore $\psi_1\in P^{\widetilde H}_{\mathcal{N}_z(A)} \widetilde\Theta_z$.

Comparing this with the definition of the operator $\Phi$ we conclude that relation~(\ref{f5_11_p2_1})
holds.

Choose an arbitrary element $h\in P^{\widetilde H}_{\mathcal{N}_z(A)} \widetilde\Theta_z$. By the proved
Lemma, it has representations~(\ref{f5_15_p2_1}), which elements satisfy
relations~(\ref{f5_16_p2_1}),(\ref{f5_17_p2_1}).
Then
$$ P^{\widetilde H}_{\mathcal{N}_{\overline z}(A)}
(U_z(\mathcal{A})\oplus T) h = P^{\widetilde H}_{\mathcal{N}_{\overline z}(A)}
T (\psi_1 + \psi_2) = T_{11}\psi_1 + T_{12}\psi_2; $$
$$ \Phi(z;\mathcal{A},T) P^{\widetilde H}_{\mathcal{N}_z(A)} h =
\Phi(z;\mathcal{A},T) \psi_1 = T_{11}\psi_1 + T_{12}\psi_2, $$
and the required equality~(\ref{f5_12_p2_1}) follows.
$\Box$

\begin{cor}
\label{c5_1_p2_1}
In conditions of Theorem~\ref{t5_2_p2_1} consider the following operator
$W_z = U_z(A)\oplus \Phi(z;\mathcal{A},T)$. Then
\begin{equation}
\label{f5_19_p2_1}
D(W_z) = P^{\widetilde H}_{H} \widetilde\Theta_z,
\end{equation}
\begin{equation}
\label{f5_20_p2_1}
W_z P^{\widetilde H}_{H} h = P^{\widetilde H}_{H}
(U_z(\mathcal{A})\oplus T) h,\qquad h\in \widetilde\Theta_z.
\end{equation}
\end{cor}
{\it Proof. }
By Lemma~\ref{l5_1_p2_1} we get $\mathcal{M}_z(A)\subseteq \widetilde\Theta_z$. Therefore
$$ P^{\widetilde H}_{H} \widetilde\Theta_z = \mathcal{M}_z(A) +
P^{\widetilde H}_{\mathcal{N}_z(A)} \widetilde\Theta_z = \mathcal{M}_z(A) + D(\Phi) = D(W_z).  $$
Let $h$ be an arbitrary element from $\widetilde\Theta_z$, having representation~(\ref{f5_15_p2_1}) which elements
satisfy, according to Lemma~\ref{l5_1_p2_1}, relations~(\ref{f5_16_p2_1}),(\ref{f5_17_p2_1}).
Using~(\ref{f5_12_p2_1}), we may write
$$ W_z P^{\widetilde H}_{H} h = (U_z(A)\oplus \Phi(z;\mathcal{A},T)) (g_1 + \psi_1) =
U_z(A) g_1 + P^{\widetilde H}_{\mathcal{N}_z(A)} h $$
$$ = P^{\widetilde H}_{\mathcal{M}_{\overline z}(A)} (U_z(\mathcal{A})\oplus T) h
+ P^{\widetilde H}_{\mathcal{N}_{\overline z}(A)}
(U_z(\mathcal{A})\oplus T) h =
P^{\widetilde H}_{H} (U_z(\mathcal{A})\oplus T) h.  $$
Therefore equality~(\ref{f5_20_p2_1}) holds.
$\Box$

Consider our considerations interrupted by formulations of the last Theorem and its Corollary.
In the sequel we assume that the operator $T$ is admissible with respect to the symmetric operator
$\mathcal{A}$.
In this case, by Theorem~\ref{t5_1_p2_1} and Remark~\ref{r5_1_p2_1}
the non-expanding operator $\Phi(z;\mathcal{A},T)$ is admissible with respect to $A$ and has a representation~(\ref{f5_10_p2_1}).
By the generalized Neumann formulas (Theorem~\ref{t4_1_p2_1}), for the operator $T$ it corresponds
a quasi-self-adjoint extension $\mathcal{A}_T$ of the symmetric operator $\mathcal{A}$.

Define the operator $\mathfrak{B}$ on the manifold $D(\mathcal{A}_T) \cap H$ in $H$ in the following way:
\begin{equation}
\label{f5_22_p2_1}
\mathfrak{B} h = \mathfrak{B}(z;\mathcal{A},T) h
= P^{\widetilde H}_H \mathcal{A}_T h,\qquad h\in D(\mathcal{A}_T) \cap H.
\end{equation}

\begin{thm}
\label{t5_3_p2_1}
Let $\mathcal{A}$ be a closed symmetric operator of the form~(\ref{f5_2_p2_1}),
where $A$, $A_e$ are symmetric operators in some Hilbert spaces $H$ and $H_e$, respectively.
Let $z\in \mathbb{R}_e$ be a fixed number and
$T$ be a non-expanding operator in  $\widetilde H = H\oplus H_e$,
with $D(T) = \mathcal{N}_z(\mathcal{A})$, $R(T) \subseteq \mathcal{N}_{\overline z}(\mathcal{A})$,
which is admissible with respect to $\mathcal{A}$.
The operator $\mathfrak{B}$, defined on the manifold $D(\mathcal{A}_T) \cap H$ in $H$ by equality~(\ref{f5_22_p2_1}), admits
the following representation:
$$ \mathfrak{B}(z;\mathcal{A},T) = A_{\Phi(z;\mathcal{A},T)}. $$
\end{thm}
{\bf Proof. } Consider a set $\widetilde \Theta_z$, defined by~(\ref{f5_14_p2_1}).
This definition we can rewrite in another form:
\begin{equation}
\label{f5_23_p2_1}
\widetilde\Theta_z = \left\{
h\in\widetilde H:\ \left(
(U_z(\mathcal{A})\oplus T) - E_{\widetilde H}
\right) h\in H
\right\}.
\end{equation}
Therefore
$$ \left(
(U_z(\mathcal{A})\oplus T) - E_{\widetilde H}
\right) \widetilde\Theta_z \subseteq H. $$
Since $R \left( (U_z(\mathcal{A})\oplus T) - E_{\widetilde H}
\right) = D(\mathcal{A}_T)$, then
$$ \left(
(U_z(\mathcal{A})\oplus T) - E_{\widetilde H}
\right) \widetilde\Theta_z \subseteq H\cap D(\mathcal{A}_T). $$
Conversely, an arbitrary element $g\in H\cap D(\mathcal{A}_T)$ has the following form
$g = \left( (U_z(\mathcal{A})\oplus T) - E_{\widetilde H}
\right) f$, $f\in \widetilde H$. By~(\ref{f5_23_p2_1}) the element $f$ belongs to
the set $\widetilde\Theta_z$. Therefore
\begin{equation}
\label{f5_24_p2_1}
\left(
(U_z(\mathcal{A})\oplus T) - E_{\widetilde H}
\right) \widetilde\Theta_z = H\cap D(\mathcal{A}_T) = D(\mathfrak{B}),
\end{equation}
where the operator $\mathfrak{B}$ is defined by equality~(\ref{f5_22_p2_1}).

\noindent
Consider the operator $W_z$ from the Corollary~\ref{c5_1_p2_1}. By~(\ref{f5_20_p2_1}) we
may write
\begin{equation}
\label{f5_25_p2_1}
(W_z - E_H) P^{\widetilde H}_{H} h = P^{\widetilde H}_{H}
\left( (U_z(\mathcal{A})\oplus T) - E_{\widetilde H}
\right) h
=
\left( (U_z(\mathcal{A})\oplus T) - E_{\widetilde H}
\right) h,\quad
h\in \widetilde\Theta_z.
\end{equation}
Therefore
$$ (W_z - E_H) P^{\widetilde H}_{H} \widetilde\Theta_z = D(\mathfrak{B}). $$
Taking into account relation~(\ref{f5_19_p2_1}) we get
\begin{equation}
\label{f5_26_p2_1}
(W_z - E_H)D(W_z) = D(\mathfrak{B}).
\end{equation}
Choose an arbitrary element $u\in D(\mathfrak{B})$. By~(\ref{f5_24_p2_1}) the element $u$ has the following form:
$u = \left(
(U_z(\mathcal{A})\oplus T) - E_{\widetilde H}
\right) h$, $h\in\widetilde\Theta_z$. Then
$$ \mathfrak{B} u = P^{\widetilde H}_H \mathcal{A}_T u $$
$$ =
P^{\widetilde H}_H
\left( z (U_z(\mathcal{A})\oplus T) - \overline{z} E_{\widetilde H}
\right)
\left( (U_z(\mathcal{A})\oplus T) - E_{\widetilde H}
\right)^{-1}
\left(
(U_z(\mathcal{A})\oplus T) - E_{\widetilde H}
\right) h $$
$$ =
P^{\widetilde H}_H
\left( z (U_z(\mathcal{A})\oplus T) - \overline{z} E_{\widetilde H}
\right) h
=
z P^{\widetilde H}_H (U_z(\mathcal{A})\oplus T) h - \overline{z} P^{\widetilde H}_H h $$
$$ = z W_z P^{\widetilde H}_H h - \overline{z} P^{\widetilde H}_H h
= ( z W_z - \overline{z} E_H ) P^{\widetilde H}_H h , $$
where we used~(\ref{f5_20_p2_1}).
Since the operator $\Phi$ is admissible with respect to the operator $A$, then $W_z$ has no non-zero
fixed elements.
By formula~(\ref{f5_25_p2_1}) we get
$$ \mathfrak{B} u =
( z W_z - \overline{z} E_H ) (W_z - E_H)^{-1}
\left( (U_z(\mathcal{A})\oplus T) - E_{\widetilde H}
\right) h $$
$$ = ( z W_z - \overline{z} E_H ) (W_z - E_H)^{-1} u
= A_{\Phi(z;\mathcal{A},T)} u,\qquad u\in D(\mathfrak{B}). $$
Equality~(\ref{f5_26_p2_1}) shows that $D(\mathfrak{B}) = D(A_{\Phi(z;\mathcal{A},T)})$.
Therefore $\mathfrak{B} = A_{\Phi(z;\mathcal{A},T)}$.
$\Box$

\subsection{The operator-valued function $\mathfrak{B}_\lambda$ generated by an extension of a symmetric operator
with an exit out of the space. }

Consider a closed symmetric operator $A$ in a Hilbert space $H$.
Let $\widetilde A$ be a self-adjoint extension of $A$, acting in
a Hilbert space $\widetilde H\supseteq H$.
Denote
\begin{equation}
\label{f6_3_p2_1}
\widetilde{\mathfrak{L}}_\lambda =
\left\{
h\in D(\widetilde A):\ (\widetilde A - \lambda E_{\widetilde H})h \in H
\right\},\qquad \lambda\in \mathbb{C};
\end{equation}
\begin{equation}
\label{f6_4_p2_1}
\mathfrak{L}_\lambda = P^{\widetilde H}_H   \widetilde{\mathfrak{L}}_\lambda,\qquad \lambda\in \mathbb{C};
\end{equation}
\begin{equation}
\label{f6_5_p2_1}
\mathfrak{L}_\infty = D(\widetilde A)\cap H.
\end{equation}
Notice that the sets $\widetilde{\mathfrak{L}}_\lambda$, $\mathfrak{L}_\lambda$ and $\mathfrak{L}_\infty$ contain
$D(A)$.
Since for non-real $\lambda$ for the self-adjoint operator $\widetilde A$ holds
$(\widetilde A - \lambda E_{\widetilde H}) D(\widetilde A) = \widetilde H$, then
\begin{equation}
\label{f6_9_p2_1}
(\widetilde A - \lambda E_{\widetilde H}) \widetilde{\mathfrak{L}}_\lambda = H,\qquad \lambda\in \mathbb{R}_e.
\end{equation}
Observe that
\begin{equation}
\label{f6_10_p2_1}
\widetilde{\mathfrak{L}}_\lambda \cap (\widetilde H\ominus H) = \{ 0 \},\qquad \lambda\in \mathbb{R}_e.
\end{equation}
In fact, if $h\in \widetilde{\mathfrak{L}}_\lambda \cap (\widetilde H\ominus H)$, then
$( (\widetilde A - \lambda E_{\widetilde H})h, h)_{\widetilde H} = 0$.
Therefore $0 = \mathop{\rm Im}\nolimits (\widetilde A h,h)_{\widetilde H} = (\mathop{\rm Im}\nolimits \lambda)
(h,h)_{\widetilde H}$, и $h=0$.

Define an operator-valued function $\mathfrak{B}_\lambda = \mathfrak{B}_\lambda (A,\widetilde A)$
for $\lambda\in \mathbb{R}_e$, which values are operators in $H$ with the domain:
\begin{equation}
\label{f6_10_0_p2_1}
D(\mathfrak{B}_\lambda) = \mathfrak{L}_\lambda,
\end{equation}
and
\begin{equation}
\label{f6_10_1_p2_1}
\mathfrak{B}_\lambda P^{\widetilde H}_H h = P^{\widetilde H}_H \widetilde A h,\qquad
h\in \widetilde{\mathfrak{L}}_\lambda.
\end{equation}
Let us check that such definition of the operator $\mathfrak{B}_\lambda$ is correct.
If an element $g\in \mathfrak{L}_\lambda$ admits two representations:
$g = P^{\widetilde H}_H h_1 = P^{\widetilde H}_H h_2$, $h_1,h_2\in \widetilde{\mathfrak{L}}_\lambda$,
then
$$ h_1 - h_2 = P^{\widetilde H}_H h_1 - P^{\widetilde H}_H h_2 +
P^{\widetilde H}_{\widetilde H\ominus H} h_1 - P^{\widetilde H}_{\widetilde H\ominus H} h_2 $$
$$ = P^{\widetilde H}_{\widetilde H\ominus H} h_1 - P^{\widetilde H}_{\widetilde H\ominus H} h_2
= P^{\widetilde H}_{\widetilde H\ominus H} (h_1 - h_2). $$
Therefore $h_1 - h_2\in \widetilde{\mathfrak{L}}_\lambda \cap (\widetilde H\ominus H)$.
By~(\ref{f6_10_p2_1}) we obtain that $h_1 = h_2$.

\noindent
The operator $\mathfrak{B}_\lambda$ for each $\lambda\in \mathbb{R}_e$ is an extension of $A$.
Define an operator $\mathfrak{B}_\infty = \mathfrak{B}_\infty(A,\widetilde A)$ in $H$ in the following way:
\begin{equation}
\label{f6_10_2_p2_1}
\mathfrak{B}_\infty = P^{\widetilde H}_H \widetilde A|_{\mathfrak{L}_\infty}.
\end{equation}
The operator $\mathfrak{B}_\infty$ is an extension of $A$, as well. Let us check that $\mathfrak{B}_\infty$
is symmetric. In fact, for arbitrary elements
$f,g\in \mathfrak{L}_\infty$ we may write
$$ (\mathfrak{B}_\infty f,g)_H = ( P^{\widetilde H}_H \widetilde A f,g)_H =
(\widetilde A f,g)_{\widetilde H} = (f, \widetilde A g)_{\widetilde H} $$
$$ = (f, P^{\widetilde H}_H \widetilde A g)_H = (f, \mathfrak{B}_\infty g)_H. $$
We emphasize that the operator $\mathfrak{B}_\infty$ {\it is not necessarily closed}.
The operator-valued function $\mathfrak{B}_\lambda$ and the operator $\mathfrak{B}_\infty$ will play an important role
in a description of the generalized resolvents of $A$.

Notice that
\begin{equation}
\label{f6_11_p2_1}
(\widetilde A - \lambda E_{\widetilde H}) h = (\mathfrak{B}_\lambda - \lambda E_H) P^{\widetilde H}_H h,\qquad
h\in \widetilde{\mathfrak{L}}_\lambda,\ \lambda\in \mathbb{R}_e.
\end{equation}
In fact, by the definitions of the sets $\widetilde{\mathfrak{L}}_\lambda$ and $\mathfrak{B}_\lambda$,
for an arbitrary element $h$ from $\widetilde{\mathfrak{L}}_\lambda$ we may write:
$$ (\widetilde A - \lambda E_{\widetilde H}) h = P^{\widetilde H}_H (\widetilde A - \lambda E_{\widetilde H}) h
= P^{\widetilde H}_H \widetilde A h - \lambda P^{\widetilde H}_H h
$$
$$ = \mathfrak{B}_\lambda P^{\widetilde H}_H h - \lambda P^{\widetilde H}_H h =
(\mathfrak{B}_\lambda - \lambda E_H) P^{\widetilde H}_H h. $$
Let us check that {\it for $\lambda$ from $\mathbb{C}_-$ ($\mathbb{C}+$)
the operator $\mathfrak{B}_\lambda$ is a maximal closed
dissipative (respectively accumulative) extension of $A$.}
In fact, using relation~(\ref{f6_11_p2_1}) for
an arbitrary element $h$ from $\widetilde{\mathfrak{L}}_\lambda$
we may write:
$$ ( \mathfrak{B}_\lambda P^{\widetilde H}_H h, P^{\widetilde H}_H h )_H =
( \mathfrak{B}_\lambda P^{\widetilde H}_H h, h )_H = ( \widetilde
A h - \lambda ( E_{\widetilde H} - P^{\widetilde H}_H ) h, h )_H
$$
$$ = (\widetilde A h,h)_H - \lambda \| P^{\widetilde H}_{\widetilde H\ominus H} h \|_H^2. $$
Then
$$ \mathop{\rm Im}\nolimits ( \mathfrak{B}_\lambda P^{\widetilde H}_H h, P^{\widetilde H}_H h )_H =
- \| P^{\widetilde H}_{\widetilde H\ominus H} h \|_H^2 \mathop{\rm Im}\nolimits \lambda, $$
and therefore $\mathfrak{B}_\lambda$ is a dissipative (respectively accumulative)
extension of $A$.
By~(\ref{f6_9_p2_1}),(\ref{f6_11_p2_1}) we get
$R(\mathfrak{B}_\lambda - \lambda E_H) = (\mathfrak{B}_\lambda - \lambda E_H) P^{\widetilde H}_H
\widetilde{\mathfrak{L}}_\lambda = (\widetilde A - \lambda E_{\widetilde H}) \widetilde{\mathfrak{L}}_\lambda
= H$.
Using formula~(\ref{f4_3_p2_1}) and the generalized Neumann formulas for the operator
$\mathfrak{B}_\lambda$ it corresponds a non-expanding operator $T$ with $D(T) = \mathcal{N}_z(A)$.
Therefore $\mathfrak{B}_\lambda$ is closed and maximal.

\noindent
By the fourth assertion of Theorem~\ref{t3_1_p2_1} {\it the operator $\mathfrak{B}_\lambda$ is densely defined,
$\lambda\in \mathbb{R}_e$}.

Consider a generalized resolvent of $A$ which corresponds to the self-adjoint extension $\widetilde A$:
$$ \mathbf{R}_\lambda = P^{\widetilde H}_H (\widetilde A - \lambda E_{\widetilde H})^{-1},\qquad
\lambda\in \mathbb{R}_e. $$
Let us check that
\begin{equation}
\label{f6_12_p2_1}
\mathbf{R}_\lambda = (\mathfrak{B}_\lambda - \lambda E_H)^{-1},\qquad \lambda\in \mathbb{R}_e.
\end{equation}
Consider an arbitrary element $g\in H$ and denote
$h = (\widetilde A - \lambda E_{\widetilde H})^{-1} g$, where $\lambda\in \mathbb{R}_e$. Since
$h\in D(\widetilde A)$ and $(\widetilde A - \lambda E_{\widetilde H}) h\in H$, then
$h\in \widetilde{\mathfrak{L}}_\lambda$. Moreover, we have $P^{\widetilde H}_H h = \mathbf{R}_\lambda g$.
By~(\ref{f6_11_p2_1}) we may write:
$$ g = (\widetilde A - \lambda E_{\widetilde H}) h =
(\mathfrak{B}_\lambda - \lambda E_H) P^{\widetilde H}_H h =
(\mathfrak{B}_\lambda - \lambda E_H) \mathbf{R}_\lambda g. $$
By the properties of maximal dissipative and accumulative operators (see~Theorem~\ref{t3_1_p2_1}),
there exists the bounded inverse $(\mathfrak{B}_\lambda - \lambda E_H)^{-1}$ defined
on the whole $H$. Applying this operator to the latter equality we get the required relation~(\ref{f6_12_p2_1}).

By~(\ref{f6_12_p2_1}) it follows that {\it the generalized resolvent $\mathbf{R}_\lambda$
of the symmetric operator $A$ is invertible for all $\lambda\in \mathbb{R}_e$}, and
the operator $\mathbf{R}_\lambda^{-1}$ is densely defined.
The operator $\mathfrak{B}_\lambda$ admits another definition by the generalized resolvent:
\begin{equation}
\label{f6_12_1_p2_1}
\mathfrak{B}_\lambda = \mathbf{R}_\lambda^{-1} + \lambda E_H,\qquad \lambda\in \mathbb{R}_e.
\end{equation}
Therefore
$$ \mathfrak{B}_\lambda^* = (\mathbf{R}_\lambda^*)^{-1} + \overline{\lambda} E_H =
\mathbf{R}_{\overline{\lambda}}^{-1} + \overline{\lambda} E_H $$
\begin{equation}
\label{f6_14_p2_1}
= \mathfrak{B}_{\overline{\lambda}},\qquad \lambda\in \mathbb{R}_e.
\end{equation}
where we used the property~(\ref{f3_1_1_p1_1}) of the generalized resolvent. 

Choose and fix a number $\lambda_0\in \mathbb{R}_e$.
By Theorem~\ref{t4_1_p2_1}, for the operator $\mathfrak{B}_\lambda$
for $\lambda\in \Pi_{\lambda_0}$ it corresponds an admissible with respect to $A$ non-expanding operator,
which we denote by $\mathfrak{F}(\lambda)$, сwith the domain $\mathcal{N}_{\lambda_0}(A)$ and
the range in $\mathcal{N}_{\overline{\lambda_0}}(A)$. Namely, we set
\begin{equation}
\label{f6_14_1_p2_1}
\mathfrak{F}(\lambda) = \mathfrak{F}(\lambda; \lambda_0, A, \widetilde A) =
\left. (\mathfrak{B}_\lambda - \overline{\lambda_0}E_H)
(\mathfrak{B}_\lambda - \lambda_0 E_H)^{-1} \right|_{\mathcal{N}_{\lambda_0}(A)},\qquad \lambda\in \Pi_{\lambda_0}.
\end{equation}
Notice that
\begin{equation}
\label{f6_14_2_p2_1}
\mathfrak{B}_\lambda = A_{\mathfrak{F}(\lambda),\lambda_0},\qquad \lambda\in \Pi_{\lambda_0}.
\end{equation}
i.e. $\mathfrak{B}_\lambda$ is a quasi-self-adjoint extension of the symmetric operator $A$,
defined by $\mathfrak{F}(\lambda)$.

By the generalized Neumann formulas for the operator $\mathfrak{B}_\infty$  there corresponds an
admissible with respect to $A$ isometric operator
$\Phi_\infty = \Phi_\infty(\lambda_0; A,\widetilde A)$ with the domain
$D(\Phi_\infty)\subseteq N_{\lambda_0}(A)$ and
the range $R(\Phi_\infty)\subseteq N_{ \overline{\lambda_0} }(A)$.
This operator will be used later.


\begin{prop}
\label{p6_1_p2_1}
Let $A$ be a closed symmetric operator in a Hilbert space $H$, $z\in \mathbb{R}_e$ be an arbitrary
fixed number, and $T$ be a linear non-expanding operator with
$D(T) = \mathcal{N}_z(A)$ and $R(T)\subseteq \mathcal{N}_{\overline{z}}(A)$, which is $z$-admissible
with respect to $A$. Then the operator $T^*$ is $\overline{z}$-admissible with respect to $A$ and
$$ (A_{T,z})^* = A_{T^*,\overline{z}}. $$
\end{prop}
{\bf Proof. } By Theorem~\ref{t4_1_p2_1} the operator $A_{T,z}$ is maximal dissipative or accumulative,
since $D(T) = \mathcal{N}_z(A)$.
By Theorem~\ref{t3_1_p2_1} we conclude that $A_{T,z}$ is densely defined and therefore there exists its adjoint.
By~(\ref{f3_6_p2_1}) we may write:
$$ ( A_{T,z} )^* = \left( zE_H + (z-\overline{z})(U_z(A)\oplus T - E_H)^{-1} \right)^* $$
$$ = \overline{z} E_H + (\overline{z}-z)(U_{\overline{z}}(A)\oplus T^* - E_H)^{-1}. $$
Consequently, the operator $U_{\overline{z}}(A)\oplus T^*$ has no non-zero fixed elements.
By Theorem~\ref{t2_1_p2_1} the operator $T^*$ is admissible with respect to $A$, and the last
equality gives the required formula.
$\Box$

By the proved Proposition and formula~(\ref{f6_14_p2_1}) we may write:
$$ \mathfrak{B}_{\overline{\lambda}} = \mathfrak{B}_\lambda^* =
( A_{\mathfrak{F}(\lambda),\lambda_0} )^* =
A_{\mathfrak{F}^*(\lambda),\overline{\lambda_0}},\qquad
\lambda\in \Pi_{\lambda_0}. $$
Let us check that the operator-valued function $\mathfrak{F}(\lambda)$ is an analytic function
of $\lambda$ in a half-plane $\Pi_{\lambda_0}$.
By~(\ref{f6_12_p2_1}) we may write:
$$ \mathfrak{B}_\lambda - \lambda_0 E_H =
(\mathfrak{B}_\lambda - \lambda_0 E_H) \mathbf{R}_\lambda (\mathfrak{B}_\lambda - \lambda E_H) $$
$$ = \left( (\mathfrak{B}_\lambda - \lambda E_H) + (\lambda - \lambda_0) E_H \right)
\mathbf{R}_\lambda (\mathfrak{B}_\lambda - \lambda E_H) $$
$$ = (E_H + (\lambda - \lambda_0) \mathbf{R}_\lambda ) (\mathfrak{B}_\lambda - \lambda E_H),\qquad
\lambda\in \mathbb{R}_e. $$
By the maximality of the dissipative or accumulative operator
$\mathfrak{B}_\lambda$ the operator $(\mathfrak{B}_\lambda - \lambda E_H)^{-1}$
exists and it is defined on the whole $H$. If we additionally assume that $\lambda\in\Pi_{\lambda_0}$,
then the operator $(\mathfrak{B}_\lambda - \lambda_0 E_H)^{-1}$ exists and is defined on the whole $H$, as well.
Then the operator
$$ E_H + (\lambda - \lambda_0) \mathbf{R}_\lambda =
(\mathfrak{B}_\lambda - \lambda_0 E_H) (\mathfrak{B}_\lambda - \lambda E_H)^{-1}, $$
has a bounded inverse:
$$ ( E_H + (\lambda - \lambda_0) \mathbf{R}_\lambda )^{-1} =
(\mathfrak{B}_\lambda - \lambda E_H) (\mathfrak{B}_\lambda - \lambda_0 E_H)^{-1}
= E_H + (\lambda_0 - \lambda) (\mathfrak{B}_\lambda - \lambda_0 E_H)^{-1}, $$
defined on the whole $H$, for $\lambda\in\Pi_{\lambda_0}$.
By property~(\ref{f3_4_p2_1}) we get
$$ \left\| ( E_H + (\lambda - \lambda_0) \mathbf{R}_\lambda )^{-1} \right\| \leq
1 + \frac{|\lambda_0 - \lambda|}{|\mathop{\rm Im}\nolimits \lambda_0|},\qquad \lambda\in\Pi_{\lambda_0}. $$
By Proposition~\ref{p1_1_p1_1} 
the function $( E_H + (\lambda - \lambda_0) \mathbf{R}_\lambda )^{-1}$ is analytic in the halfplane
$\Pi_{\lambda_0}$.
Using relation~(\ref{f6_12_p2_1}) we may write that
$$ (\mathfrak{B}_\lambda - \overline{\lambda_0}E_H)
(\mathfrak{B}_\lambda - \lambda_0 E_H)^{-1} =
E_H + (\lambda_0 - \overline{\lambda_0}) (\mathfrak{B}_\lambda - \lambda_0 E_H)^{-1} $$
$$ = E_H + (\lambda_0 - \overline{\lambda_0}) (\mathfrak{B}_\lambda - \lambda E_H)^{-1}
(\mathfrak{B}_\lambda - \lambda E_H)(\mathfrak{B}_\lambda - \lambda_0 E_H)^{-1} $$
$$ = E_H + (\lambda_0 - \overline{\lambda_0}) \mathbf{R}_\lambda
( E_H + (\lambda - \lambda_0) \mathbf{R}_\lambda )^{-1}, $$
is an analytic function in $\Pi_{\lambda_0}$.
Thenthe function $\mathfrak{F}(\lambda)$ is analytic in
$\Pi_{\lambda_0}$, as well.

From our considerations we conclude that the operator-valued function $\mathfrak{B}_\lambda$
admits the following representation:
\begin{equation}
\label{f6_15_p2_1}
\mathfrak{B}_\lambda = \left\{ \begin{array}{cc}
A_{\mathfrak{F}(\lambda),\lambda_0}, & \lambda\in \Pi_{\lambda_0}\\
A_{\mathfrak{F}^*(\overline{\lambda}),\overline{\lambda_0}}, &
\overline{\lambda}\in \Pi_{\lambda_0}\end{array}\right.,
\end{equation}
where $\mathfrak{F}(\lambda)$ is an analytic function of $\lambda$
in the half-plane $\Pi_{\lambda_0}$, which values are non-expanding operators
with $D(\mathfrak{F}(\lambda)) = \mathcal{N}_{\lambda_0}(A)$ and
$R(\mathfrak{F}(\lambda)) \subseteq \mathcal{N}_{\overline{\lambda_0}}(A)$, admissible with respect to
the operator $A$.

Our next aim will be derivation of a representation for the operator-valued function
$\mathfrak{F}(\lambda)$,  which connects it with the so-called characteristic function of a symmetric
operator $A$. The next subsection will be decoted to the definition of the characteristic function and some its properties.

\subsection{The characteristic function of a symmetric operator.}

Let $A$ be a closed symmetric operator in a Hilbert space $H$.
Choose an arbitrary number $\lambda\in \mathbb{R}_e$.
Notice that the operator $T = 0_{ \mathcal{N}_{ \overline{\lambda} } }$ is admissible with respect to $A$.
In fact, suppose that for some $h\in H$ holds
$(U_{\overline{\lambda}}\oplus T) h = h$. Let $h = h_1 + h_2$, $h_1\in \mathcal{M}_{\overline{\lambda}}$,
$h_2\in \mathcal{N}_{\overline{\lambda}}$, then
$$ (U_{\overline{\lambda}}\oplus T) h = U_{\overline{\lambda}} h_1 = h = h_1 + h_2, $$
and therefore
$$ \| h_1 \|_H^2 = \| U_{\overline{\lambda}} h_1 \|_H^2 = \| h_1 \|_H^2 + \| h_2 \|_H^2. $$
Consequently, we get $h_2=0$. Since the operator $U_{\overline{\lambda}}$ has no non-zero fixed
elements, then $h_1 = h =0$. By Theorem~\ref{t2_1_p2_1} we obtain that $T$ is admissible with respect to $A$.

\noindent
Consider a quasi-self-adjoint extension
$A_\lambda:= A_{ 0_{ \mathcal{N}_{\overline{\lambda}} }, \overline{\lambda} }$ of the operator $A$.
By the generalized Neumann formulas for $\lambda$ from $\mathbb{C}_+$ ($\mathbb{C}_-$)
the operator $A_\lambda$ is maximal closed dissipative (respectively accumulative)
extension of $A$. By Proposition~\ref{p6_1_p2_1} it follows that
$$ A_\lambda^* = A_{ 0_{ \mathcal{N}_{\overline{\lambda}} }, \overline{\lambda} }^* =
A_{ 0_{ \mathcal{N}_\lambda }, \lambda } = A_{\overline{\lambda}}. $$
Moreover, by~(\ref{f4_1_p2_1}),(\ref{f4_2_p2_1}) we may write:
\begin{equation}
\label{f7_2_p2_1}
D(A_\lambda) = D(A)\dotplus \mathcal{N}_{\overline{\lambda}}(A),
\end{equation}
and
\begin{equation}
\label{f7_3_p2_1}
A_\lambda ( f + \psi ) = Af + \lambda\psi,\qquad f\in D(A),\ \psi\in \mathcal{N}_{\overline{\lambda}}(A).
\end{equation}

Choose and fix an arbitrary number $\lambda_0\in \mathbb{R}_e$ and consider the half-plane
$\Pi_{\lambda_0}$ containing the point $\lambda_0$.
If we additionally suppose that $\lambda\in\Pi_{\lambda_0}$, for the quasi-self-adjoint extension $A_\lambda$
by the generalized Neumann formulas it corresponds a non-expanding operator
$C(\lambda)$ with $D(C(\lambda)) = \mathcal{N}_{\overline{\lambda_0}}(A)$,
$R(C(\lambda)) \subseteq \mathcal{N}_{\lambda_0}(A)$.
Namely, the operator $C(\lambda)$ is given by the following equality:
$$ C(\lambda) = C(\lambda; \lambda_0, A) = \left.(A_\lambda - \lambda_0 E_H) (A_\lambda - \overline{\lambda_0} E_H)
\right|_{ \mathcal{N}_{\overline{\lambda_0}}(A) },\qquad \lambda\in\Pi_{\lambda_0}. $$
This operator-valued function $C(\lambda)$ in the half-plane $\Pi_{\lambda_0}$
is said to be {\bf the characteristic function of the symmetric operator $A$}.
Observe that
\begin{equation}
\label{f7_4_p2_1}
A_\lambda = A_{C(\lambda), \overline{\lambda_0}}.
\end{equation}

\begin{prop}
\label{p7_1_p2_1}
Let $A$ be a closed symmetric operator in a Hilbert space $H$, and $\lambda,z\in \mathbb{R}_e$ be
arbitrary numbers: $\lambda\in\Pi_z$.
The space  $H$ can be represented as a direct sum:
$$ H = \mathcal{M}_\lambda(A) \dotplus \mathcal{N}_z(A). $$
The projection operator $\mathcal{P}_{\lambda,z}$ in $H$ on the subspace $\mathcal{N}_z(A)$
parallel to the subspace $\mathcal{M}_\lambda(A)$ (i.e. the operator which to arbitrary vector $h\in H$,
$h = h_1 + h_2$, $h_1\in \mathcal{M}_\lambda(A)$, $h_2\in \mathcal{N}_z(A)$, put into correspondence
a vector $\mathcal{P}_{\lambda,z} h = h_2$) has the following form:
\begin{equation}
\label{f7_5_p2_1}
\mathcal{P}_{\lambda,z} = \frac{ \lambda - \overline{z} }{ z - \overline{z} }
P^H_{\mathcal{N}_z(A)} (A_{ \overline{z} } - zE_H) (A_{ \overline{z} } - \lambda E_H)^{-1}.
\end{equation}
\end{prop}
{\bf Proof. } Since $A_{ \overline{z} }$ for $\overline{z}$ from $\mathbb{C}_+$ ($\mathbb{C}_-$)
is a maximal dissipative (respectively accumulative) extension of $A$, then
there exist the inverses $(A_{ \overline{z} } - \lambda E_H)^{-1}$ and
$(A_{ \overline{z} } - z E_H)^{-1}$, defined on the whole $H$ (see Theorem~\ref{t3_1_p2_1}).
The following operator:
$$ S_{\lambda,z} := (A_{ \overline{z} } - \lambda E_H) (A_{ \overline{z} } - z E_H)^{-1} =
E_H + (z-\lambda) (A_{ \overline{z} } - z E_H)^{-1}, $$
maps bijectively the whole space $H$ on the whole $H$, and
$$ S_{\lambda,z} := (A_{ \overline{z} } - z E_H) (A_{ \overline{z} } - \lambda E_H)^{-1}. $$
By~(\ref{f7_3_p2_1}) it follows that $(A_{\overline{z}} -zE_H)\psi = (\overline{z} - z)\psi$, for an arbitrary
$\psi\in \mathcal{N}_z(A)$. Therefore
$(A_{\overline{z}} -zE_H)^{-1} = \frac{1}{ \overline{z} - z } E_{\mathcal{N}_z(A)}$.
For an arbitrary $g\in \mathcal{N}_z(A)$ we may write:
\begin{equation}
\label{f7_6_p2_1}
S_{\lambda,z} g = g + (z-\lambda) (A_{ \overline{z} } - z E_H)^{-1} g =
g + \frac{z - \lambda}{ \overline{z} - z } g = \frac{ \overline{z} - \lambda }{ \overline{z} - z } g,
\end{equation}
to get
$$ S_{\lambda,z} \mathcal{N}_z(A) = \mathcal{N}_z(A). $$
For an arbitrary vector $f\in D(A)$ holds:
$$ S_{\lambda,z} (A-zE_H)f = (A-\lambda E_H)f, $$
and therefore
$$ S_{\lambda,z} \mathcal{M}_z(A) = \mathcal{M}_\lambda(A). $$
Choose an arbitrary element $h\in H$, and set $f = S_{\lambda,z}^{-1}h$. Let
$f = f_1 + f_2$, $f_1\in \mathcal{M}_z(A)$, $f_2\in \mathcal{N}_z(A)$, then
$$ h = S_{\lambda,z} f_1 + S_{\lambda,z} f_2\in \mathcal{M}_\lambda(A) + \mathcal{N}_z(A). $$
Therefore $H = \mathcal{M}_\lambda(A) + \mathcal{N}_z(A)$.
Suppose that
$v\in \mathcal{M}_\lambda(A) \cap \mathcal{N}_z(A)$. Then there exist elements
$u\in \mathcal{M}_z(A)$ and $w\in \mathcal{N}_z(A)$ such that
$S_{\lambda,z} u = v$ and $S_{\lambda,z} w = v$. Therefore $S_{\lambda,z} (u-w)=0$,
and by the invertibility of  $S_{\lambda,z}$ we get $u=w$. Since elements $u$ and $w$ are orthogonal,
then $u=w=0$ and $v=0$. Thus, the required decomposition of the space $H$ is proved.
Moreover, that for the element $h$ holds
$$ \mathcal{P}_{\lambda,z} h = S_{\lambda,z} f_2 = S_{\lambda,z} P^H_{\mathcal{N}_z(A)} f =
S_{\lambda,z} P^H_{\mathcal{N}_z(A)} S_{\lambda,z}^{-1}h =
\frac{ \overline{z} - \lambda }{ \overline{z} - z } P^H_{\mathcal{N}_z(A)} S_{\lambda,z}^{-1}h $$
$$ = \frac{ \overline{z} - \lambda }{ \overline{z} - z } P^H_{\mathcal{N}_z(A)}
(A_{ \overline{z} } - z E_H) (A_{ \overline{z} } - \lambda E_H)^{-1} h, $$
where we used relation~(\ref{f7_6_p2_1}). Thus, equality~(\ref{f7_5_p2_1}) is proved.
$\Box$

Continue our considerations started before the formulation of the last proposition. Define the following two
operator-valued functions:
$$ Q(\lambda) = \mathcal{P}_{\lambda,\lambda_0},\quad
K(\lambda) = \left. \mathcal{P}_{\lambda,\lambda_0} \right|_{ \mathcal{N}_{ \overline{\lambda_0} } },\qquad
\lambda\in\Pi_{\lambda_0}. $$
By representation~(\ref{f7_5_p2_1}) it follows that
$Q(\lambda)$ and $K(\lambda)$ are analytic in the half-plane $\Pi_{\lambda_0}$.
In fact, since the operator $A_{\overline{\lambda_0}}$ is dissipative or accumulative, then
estimate~(\ref{f3_4_p2_1}) holds. By Proposition~\ref{p1_1_p1_1}                           
we obtain that the operator-valued function $(A_{ \overline{z} } - \lambda E_H)^{-1}$
is analytic in $\Pi_{\lambda_0}$. Then the same is true for
$Q(\lambda)$ and $P(\lambda)$.

Let us check that {\it for each $\lambda\in\Pi_{\lambda_0}$ the operator $K(\lambda)$ is non-expanding}.
In fact, for arbitrary $\lambda\in\Pi_{\lambda_0}$ and $\psi\in \mathcal{N}_{\overline{\lambda_0}}$
we can assert that $\psi - K(\lambda)\psi$ belongs to $\mathcal{M}_\lambda(A)$, and therefore
$\psi - K(\lambda)\psi = (A - \lambda E_H) f$, where $f\in D(A)$.
Let $f = (U_{\lambda_0}(A) - E_H) h$, where $h\in \mathcal{M}_{\lambda_0}(A)$. Then
$A f = (\lambda_0 U_{\lambda_0}(A) - \overline{\lambda_0} E_H) h$.
Substituting expressions for $f$ and $Af$ in the expression for $\psi - K(\lambda)\psi$ we get
$$ \psi - K(\lambda)\psi =
(\lambda_0 U_{\lambda_0}(A) - \overline{\lambda_0} E_H) h - \lambda (U_{\lambda_0}(A) - E_H) h. $$
Therefore
$$ \psi + (\lambda - \lambda_0) U_{\lambda_0}(A) h = K(\lambda) \psi + (\lambda - \overline{\lambda_0}) h. $$
Using the orthogonality of the summands in the left-hand side and also in the right-hand side we get
$$ \| \psi \|_H^2 + |\lambda - \lambda_0|^2 \| U_{\lambda_0}(A) h \|_H^2 =
\| K(\lambda)\psi \|_H^2 + |\lambda - \overline{\lambda_0}|^2 \| h \|_H^2. $$
Therefore
$$ \| \psi \|_H^2 - \| K(\lambda)\psi \|_H^2 = \delta_{\lambda,\lambda_0} \| h \|_H^2, $$
где $\delta_{\lambda,\lambda_0} := |\lambda - \overline{\lambda_0}|^2 - |\lambda - \lambda_0|^2$.
It remains to notice that
$$ \delta_{\lambda,\lambda_0} = |\mathop{\rm Re}\nolimits \lambda - \mathop{\rm Re}\nolimits \overline{\lambda_0}|^2
+ |\mathop{\rm Im}\nolimits \lambda - \mathop{\rm Im}\nolimits \overline{\lambda_0}|^2
- |\mathop{\rm Re}\nolimits \lambda - \mathop{\rm Re}\nolimits \lambda_0|^2 $$
\begin{equation}
\label{f7_6_1_p2_1}
- |\mathop{\rm Im}\nolimits \lambda - \mathop{\rm Im}\nolimits \lambda_0|^2
= |\mathop{\rm Im}\nolimits \lambda + \mathop{\rm Im}\nolimits \lambda_0|^2 -
|\mathop{\rm Im}\nolimits \lambda - \mathop{\rm Im}\nolimits \lambda_0|^2 > 0,
\end{equation}
since $\mathop{\rm Im}\nolimits\lambda$ and $\mathop{\rm Im}\nolimits\lambda_0$ have the same sign.

Let us find a relation between the function $K(\lambda)$ and the characteristic function $C(\lambda)$.
Choose an arbitrary $\lambda\in\Pi_{\lambda_0}$ and an element $\psi\in N_{\overline{\lambda_0}}$.
By the Neumann formula~(\ref{f4_2_p2_1}) we may write:
$$ A_{C(\lambda),\overline{\lambda_0}} (-C(\lambda)\psi + \psi) =
-\overline{\lambda_0}C(\lambda)\psi + \lambda_0 \psi. $$
On the other hand, by~(\ref{f7_4_p2_1}),(\ref{f7_2_p2_1}) and~(\ref{f7_3_p2_1}) we get
$$ A_{C(\lambda),\overline{\lambda_0}} (-C(\lambda)\psi + \psi) =
A_\lambda (-C(\lambda)\psi + \psi) = A_\lambda (f + \psi_1) = Af + \lambda \psi_1, $$
where $f\in D(A)$ and $\psi_1\in \mathcal{N}_{\overline{\lambda}}(A)$ are some elements:
$$ -C(\lambda)\psi + \psi = f + \psi_1 \in D(A_\lambda). $$
Therefore
$$ -\overline{\lambda_0}C(\lambda)\psi + \lambda_0 \psi =
Af + \lambda \psi_1. $$
Subtracting fron th last equality the previpus one, multiplied by $\lambda$, we obtain that
$$ (\lambda_0 - \lambda)\psi - (\overline{\lambda_0} - \lambda) C(\lambda)\psi = (A-\lambda E_H) f. $$
Thus, the following inclusion holds:
\begin{equation}
\label{f7_7_p2_1}
(\lambda - \lambda_0)\psi - (\lambda - \overline{\lambda_0}) C(\lambda)\psi \in
\mathcal{M}_\lambda(A).
\end{equation}
Divide the both sides of the last equality by $\lambda - \overline{\lambda_0}$ and apply
the operator $\mathcal{P}_{\lambda,\lambda_0}$:
$$ \frac{ \lambda - \lambda_0 }{ \lambda - \overline{\lambda_0} } \mathcal{P}_{\lambda,\lambda_0}\psi -
C(\lambda)\psi = 0, $$
i.e. $C(\lambda)\psi =
\frac{ \lambda - \lambda_0 }{ \lambda - \overline{\lambda_0} } \mathcal{P}_{\lambda,\lambda_0}\psi$.
Therefore
$$ C(\lambda) = \frac{ \lambda - \lambda_0 }{ \lambda - \overline{\lambda_0} } K(\lambda),\qquad
\lambda\in\Pi_{\lambda_0}. $$
Since $K(\lambda)$ is analytic in $\Pi_{\lambda_0}$ and
non-expanding, then {\it the characteristic function $C(\lambda)$ is analytic and }
\begin{equation}
\label{f7_8_p2_1}
\| C(\lambda) \| \leq \left|
\frac{ \lambda - \lambda_0 }{ \lambda - \overline{\lambda_0} }
\right|,\qquad \lambda\in\Pi_{\lambda_0}.
\end{equation}

Now we shall obtain another formula for the characteristic function $C(\lambda)$.
Let $\lambda\in\Pi_{\lambda_0}$ be an arbitrary number.
Let us check that {\it elements $\psi\in \mathcal{N}_{ \overline{\lambda_0} }(A)$ and
$\widehat\psi\in \mathcal{N}_{\lambda_0}(A)$
are connected by the following relation
\begin{equation}
\label{f7_9_p2_1}
(\lambda - \lambda_0)\psi - (\lambda - \overline{\lambda_0}) \widehat\psi \in
\mathcal{M}_\lambda(A),
\end{equation}
if and only if
$$ \widehat\psi = C(\lambda)\psi. $$
}
The sufficiency of this assertion follows from relation~(\ref{f7_7_p2_1}).
Suppose that relation~(\ref{f7_9_p2_1}) holds. Apply the operator $\mathcal{P}_{\lambda,\lambda_0}$ to the both sides of this equality:
$$ (\lambda - \lambda_0) \mathcal{P}_{\lambda,\lambda_0} \psi -
(\lambda - \overline{\lambda_0}) \widehat\psi = 0, $$
and therefore
$$ \widehat\psi =
\frac{ \lambda - \lambda_0 }{ \lambda - \overline{\lambda_0} } \mathcal{P}_{\lambda,\lambda_0}\psi =
C(\lambda)\psi, $$
what we needed to prove.

\subsection{A connection of the operator-valued function $\mathfrak{F}(\lambda)$ with the characteristic function.}

Consider a closed symmetric operator $A$ в in a Hilbert space $H$.
Let $\widetilde A$ be a self-adjoint extension of  $A$ acting in
a Hilbert space $\widetilde H\supseteq H$.
Choose and fix an arbitrary number $\lambda_0$ from $\mathbb{R}_e$.
For the extension $\widetilde A$ there correspond operator-valued functions $\mathfrak{B}_\lambda =
\mathfrak{B}_\lambda(A,\widetilde A)$ and
$\mathfrak{F}(\lambda) =
\mathfrak{F}(\lambda;\lambda_0,A,\widetilde A)$ (see~(\ref{f6_10_1_p2_1}),(\ref{f6_14_1_p2_1})).

As it was done above, the space $H$ we represent as~(\ref{f5_1_p2_1}), where
$H_e = \widetilde H\ominus H$, and the operator $A$ can be identified with the operator
$A\oplus o_{H_e}$.

Consider an operator $\mathcal{A}$ of the form~(\ref{f5_2_p2_1}),
where $A_e$ is a symmetric operator in a Hilbert space $H_e$, such that
$A\oplus A_e\subseteq \widetilde A$. In particular, we can take as $A_e$ the operator $o_{H_e}$.
Thus, the operator $\widetilde A$ is a self-adjoint extension of the operator $\mathcal{A}$.

By the generalized Neumann formulas, for the self-adjoint extension $\widetilde A$ of
$\mathcal{A}$ there corresponds an isometric operator $T$ with the domain $D(T)=\mathcal{N}_{\lambda_0}(\mathcal{A})$
and the range $R(T)=\mathcal{N}_{\overline{\lambda_0}}(\mathcal{A})$, haning
block representation~(\ref{f5_6_p2_1}),
where $T_{11} = P_{\mathcal{N}_{\overline{\lambda_0}}(A)} T P_{\mathcal{N}_{\lambda_0}(A)}$,
$T_{12} = P_{\mathcal{N}_{\overline{\lambda_0}}(A)} T P_{\mathcal{N}_{\lambda_0}(A_e)}$,
$T_{21} = P_{\mathcal{N}_{\overline{\lambda_0}}(A_e)} T P_{\mathcal{N}_{\lambda_0}(A)}$,
$T_{22} = P_{\mathcal{N}_{\overline{\lambda_0}}(A_e)} T P_{\mathcal{N}_{\lambda_0}(A_e)}$.

\noindent
Denote $C_e(\lambda) := C(\lambda; \lambda_0, A_e)$, i.e. $C_e(\lambda)$ is a characteristic
function of the operator $A_e$ in $H_e$.

By~(\ref{f4_1_p2_1}) the domain of $\widetilde A$ consists of
elements $\widetilde f$ of the following form:
$$ \widetilde f = f_1 + f_2 + T (\psi_1 + \psi_2) - \psi_1 - \psi_2 $$
\begin{equation}
\label{f8_5_p2_1}
= [ f_1 + T_{11}\psi_1 + T_{12}\psi_2 - \psi_1 ] +
[ f_2 + T_{21}\psi_1 + T_{22}\psi_2 - \psi_2 ],
\end{equation}
where $f_1\in D(A)$, $f_2\in D(A_e)$, $\psi_1\in \mathcal{N}_{\lambda_0}(A)$, $\psi_2\in \mathcal{N}_{\lambda_0}(A_e)$.
Notice that an expression in the first square brackets belongs to $H$, and in the second square brackets belongs to $H_e$.
Formula~(\ref{f4_2_p2_1}) gives the following relation for the operator $\widetilde A$:
\begin{equation}
\label{f8_6_p2_1}
\widetilde A \widetilde f =
[ A f_1 + \lambda_0(T_{11}\psi_1 + T_{12}\psi_2) - \overline{\lambda_0} \psi_1 ] +
[ A_e f_2 + \lambda_0 (T_{21}\psi_1 + T_{22}\psi_2) - \overline{\lambda_0} \psi_2 ].
\end{equation}
Here also an expression in the first square brackets belongs to $H$, and in the second square brackets belongs to $H_e$.
Recall that the manifold $\mathfrak{L}_\lambda$ ($\lambda\in \mathbb{C}$) consists of those
elements $\widetilde f\in D(\widetilde A)$, for which $(\widetilde A - \lambda E_{\widetilde H})
\widetilde f$ belongs to $H$.
Subtracting from relation~(\ref{f8_6_p2_1}) relation~(\ref{f8_5_p2_1}), multiplied by $\lambda$,
we see that the element $\widetilde f$ of form~(\ref{f8_5_p2_1}) belongs to $\mathfrak{L}_\lambda$
if and only if
\begin{equation}
\label{f8_7_p2_1}
(A_e - \lambda E_{H_e}) f_2 + (\lambda_0 - \lambda) (T_{21}\psi_1 + T_{22}\psi_2) - (\overline{\lambda_0} -
\lambda) \psi_2 = 0.
\end{equation}

Suppose now that $\lambda\in\Pi_{\lambda_0}$.
If relation~(\ref{f8_7_p2_1}) is satisfied then  by the property of the characteristic function,
see the text with formula~(\ref{f7_9_p2_1}), we get
\begin{equation}
\label{f8_8_p2_1}
\psi_2 = C_e(\lambda) (T_{21}\psi_1 + T_{22}\psi_2).
\end{equation}
Moreover, by~(\ref{f8_7_p2_1}) we get
\begin{equation}
\label{f8_9_p2_1}
f_2 = (A_e - \lambda E_{H_e})^{-1} \left( (\lambda - \lambda_0) (T_{21}\psi_1 + T_{22}\psi_2) -
(\lambda - \overline{\lambda_0}) \psi_2 \right).
\end{equation}
Since $\lambda$ and $\lambda_0$ lie in the same half-plane, then the number
$\delta_{\lambda,\lambda_0} = |\lambda - \overline{\lambda_0}|^2 - |\lambda - \lambda_0|^2$ is positive,
see~(\ref{f7_6_1_p2_1}).
By the property of the characteristic function~(\ref{f7_8_p2_1}) we obtain that
$$ \| C_e(\lambda) \| \leq \left|
\frac{ \lambda - \lambda_0 }{ \lambda - \overline{\lambda_0} }
\right| < 1,\qquad \lambda\in\Pi_{\lambda_0}. $$
The operator $C_e(\lambda) T_{22}$ may be considered as an operator in a Hilbert space
$\mathcal{N}_{\lambda_0}(A_e)$, defined on the whole $\mathcal{N}_{\lambda_0}(A_e)$.
Since $T_{22}$ is non-expanding, then $C_e(\lambda) T_{22}$ is contractive
and there exists the inverse $(E_{\mathcal{N}_{\lambda_0}(A_e)} - C_e(\lambda) T_{22})^{-1}$,
defined on the whole $\mathcal{N}_{\lambda_0}(A_e)$.
By~(\ref{f8_8_p2_1}) we get
\begin{equation}
\label{f8_10_p2_1}
\psi_2 = (E_{\mathcal{N}_{\lambda_0}(A_e)} - C_e(\lambda) T_{22})^{-1} C_e(\lambda) T_{21} \psi_1.
\end{equation}

Conversely, if elements $f_1\in D(A)$ and $\psi_1\in \mathcal{N}_{\lambda_0}(A)$ we choose arbitrarily,
the element $\psi_2$ we define by~(\ref{f8_10_p2_1}), and $f_2$ we define by~(\ref{f8_9_p2_1})
(this is correct since~(\ref{f8_10_p2_1}) implies~(\ref{f8_8_p2_1}))
and therefore $\left( (\lambda - \lambda_0) (T_{21}\psi_1 + T_{22}\psi_2) -
(\lambda - \overline{\lambda_0}) \psi_2 \right)$ belongs to $\mathcal{M}_\lambda(A_e)$),
then the element $\widetilde f$ of form~(\ref{f8_5_p2_1}) belongs to $\mathfrak{L}_\lambda$.

Thus, the set $\mathfrak{L}_\lambda$ ($\lambda\in\Pi_{\lambda_0}$) consists of
elements $\widetilde f$ of form~(\ref{f8_5_p2_1}), where
$f_1\in D(A)$ and $\psi_1\in \mathcal{N}_{\lambda_0}(A)$ are arbitrary,
$\psi_2$ and $f_2$ are defined by~(\ref{f8_10_p2_1}),(\ref{f8_9_p2_1}).
Since the domain of $\mathfrak{B}_\lambda$ is $P^{\widetilde H}_H \mathfrak{L}_\lambda$,
and its action is defined by~(\ref{f6_10_1_p2_1}), then
$D(\mathfrak{B}_\lambda)$ consists of elements of form:
$$ f = f_1 + T_{11}\psi_1 + T_{12}\psi_2 - \psi_1, $$
where elements $f_1\in D(A)$ and $\psi_1\in \mathcal{N}_{\lambda_0}(A)$ are arbitrary,
and $\psi_2$ is defined by~(\ref{f8_10_p2_1}). Moreover, we have
$$ \mathfrak{B}_\lambda f =
A f_1 + \lambda_0(T_{11}\psi_1 + T_{12}\psi_2) - \overline{\lambda_0} \psi_1. $$
Substituting an expression for $\psi_2$ from~(\ref{f8_10_p2_1}) in these expressions we obtain that
$D(\mathfrak{B}_\lambda)$ consists of elements of the following form:
\begin{equation}
\label{f8_11_p2_1}
f = f_1 + \left(
T_{11} + T_{12} (E_{\mathcal{N}_{\lambda_0}(A_e)} - C_e(\lambda) T_{22})^{-1} C_e(\lambda) T_{21}
\right) \psi_1 - \psi_1,
\end{equation}
where elements $f_1\in D(A)$ and $\psi_1\in \mathcal{N}_{\lambda_0}(A)$ are arbitrary.
The action of the operator $\mathfrak{B}_\lambda$ has the following form:
\begin{equation}
\label{f8_12_p2_1}
\mathfrak{B}_\lambda f =
A f_1 + \lambda_0
\left(
T_{11} + T_{12} (E_{\mathcal{N}_{\lambda_0}(A_e)} - C_e(\lambda) T_{22})^{-1} C_e(\lambda) T_{21}
\right)
- \overline{\lambda_0} \psi_1.
\end{equation}

Using the obtained expression for the operator $\mathfrak{B}_\lambda$ ($\lambda\in\Pi_{\lambda_0}$),
we find expression for $\mathfrak{F}(\lambda) =
\mathfrak{F}(\lambda;\lambda_0,A,\widetilde A)$. Since
$\mathfrak{F}(\lambda)$ has the following form (see~(\ref{f6_14_1_p2_1})):
\begin{equation}
\label{f8_14_p2_1}
\mathfrak{F}(\lambda) = E_{\mathcal{N}_{\lambda_0}(A)} +
(\lambda_0 - \overline{\lambda_0})
\left. (\mathfrak{B}_\lambda - \lambda_0 E_H)^{-1} \right|_{\mathcal{N}_{\lambda_0}(A)},
\end{equation}
then we need to define the action of the operator
$(\mathfrak{B}_\lambda - \lambda_0 E_H)^{-1}$  on the set $\mathcal{N}_{\lambda_0}(A)$.
By~(\ref{f8_12_p2_1}), for an arbitrary element $f\in D(\mathfrak{B}_\lambda)$
of the form~(\ref{f8_11_p2_1}) we may write:
$$ (\mathfrak{B}_\lambda - \lambda_0 E_H) f =
(A - \lambda_0 E_H) f_1 + (\lambda_0 - \overline{\lambda_0}) \psi_1. $$
Applying the operator $(\mathfrak{B}_\lambda - \lambda_0 E_H)^{-1}$ to the both sides of the latter equality
we get
$$ (\mathfrak{B}_\lambda - \lambda_0 E_H)^{-1}
\left(
(A - \lambda_0 E_H) f_1 + (\lambda_0 - \overline{\lambda_0}) \psi_1
\right)
=
f $$
$$ =
f_1 + \left(
T_{11} + T_{12} (E_{\mathcal{N}_{\lambda_0}(A_e)} - C_e(\lambda) T_{22})^{-1} C_e(\lambda) T_{21}
\right) \psi_1 - \psi_1. $$
In particular, choosing $f_1 = 0$ we get
$$ (\mathfrak{B}_\lambda - \lambda_0 E_H)^{-1}
\psi_1
= \frac{1}{\lambda_0 - \overline{\lambda_0}}
\left(
T_{11} + T_{12} (E_{\mathcal{N}_{\lambda_0}(A_e)} - C_e(\lambda) T_{22})^{-1} C_e(\lambda) T_{21}
\right) \psi_1 $$
$$ - \frac{1}{\lambda_0 - \overline{\lambda_0}} \psi_1,\qquad \psi_1\in \mathcal{N}_{\lambda_0}(A). $$
Substituting the obtained relation for the operator
$(\mathfrak{B}_\lambda - \lambda_0 E_H)^{-1}|_{\mathcal{N}_{\lambda_0}(A)}$ into~(\ref{f8_14_p2_1}),
we obtain that
\begin{equation}
\label{f8_15_p2_1}
\mathfrak{F}(\lambda;\lambda_0,A,\widetilde A) =
T_{11} + T_{12} (E_{\mathcal{N}_{\lambda_0}(A_e)} - C_e(\lambda) T_{22})^{-1} C_e(\lambda) T_{21},\qquad
\lambda\in\Pi_{\lambda_0}.
\end{equation}
This representation will play an important role for the analytic description of the generalized resolvents of a symmetric
operator $A$.

\subsection{Boundary properties
of the operator-valued function $\mathfrak{F}(\lambda)$: the case of a densely defined
symmetric operator $A$.}

Consider an arbitrary closed symmetric operator $A$ in a Hilbert space $H$,
having a dense domain: $\overline{D(A)} = H$.
Let $\widetilde A$ be a self-adjoint extension of the operator $A$, acting in
a Hilbert space $\widetilde H\supseteq H$.
Let us check that
\begin{equation}
\label{f9_5_1_p2_1}
P^{\widetilde H}_H D(\widetilde A) \subseteq D(A^*),\quad
P^{\widetilde H}_H \widetilde A h = A^* P^{\widetilde H}_H h,\quad h\in D(\widetilde A).
\end{equation}
In fact, for arbitrary elements $h\in D(\widetilde A)$ and $f\in D(A)$ we write:
$$ (Af, P^{\widetilde H}_H h)_H = (Af,h)_{\widetilde H} = (\widetilde Af, h)_{\widetilde H}
= (f, \widetilde A h)_{\widetilde H} = (f, P^{\widetilde H}_H \widetilde A h)_H, $$
and the required relation follows.
Using~(\ref{f9_5_1_p2_1}) we may write:
\begin{equation}
\label{f9_6_p2_1}
\widetilde A h = A^* P^{\widetilde H}_H h + P^{\widetilde H}_{\widetilde H\ominus H} \widetilde A h,\qquad
h\in D(\widetilde A).
\end{equation}

Consider the sets $\widetilde{\mathfrak{L}}_\lambda$, $\mathfrak{L}_\lambda$ and
$\mathfrak{L}_\infty$, defined by equalities~(\ref{f6_3_p2_1})-(\ref{f6_5_p2_1}) for $\lambda\in \mathbb{C}$.
By the definition of $\widetilde{\mathfrak{L}}_\lambda$, for an arbitrary element
$g$ from $\widetilde{\mathfrak{L}}_\lambda$ it holds:
$0 = P^{\widetilde H}_{\widetilde H\ominus H} (\widetilde A - \lambda E_{\widetilde H}) g =
P^{\widetilde H}_{\widetilde H\ominus H} \widetilde A g - \lambda P^{\widetilde H}_{\widetilde H\ominus H} g$,
i.e.
\begin{equation}
\label{f9_6_1_p2_1}
P^{\widetilde H}_{\widetilde H\ominus H} \widetilde A g =
\lambda P^{\widetilde H}_{\widetilde H\ominus H} g,\qquad
g\in \widetilde{\mathfrak{L}}_\lambda,\quad \lambda\in \mathbb{C}.
\end{equation}
Therefore
\begin{equation}
\label{f9_7_p2_1}
\widetilde A g = A^* P^{\widetilde H}_H g +
\lambda P^{\widetilde H}_{\widetilde H\ominus H} g,\qquad
g\in \widetilde{\mathfrak{L}}_\lambda,\quad \lambda\in \mathbb{C}.
\end{equation}
It follows that
\begin{equation}
\label{f9_8_p2_1}
(\widetilde A - \lambda E_{\widetilde H}) g = (A^* - \lambda E_H) P^{\widetilde H}_H g,\qquad
g\in \widetilde{\mathfrak{L}}_\lambda,\quad \lambda\in \mathbb{C}.
\end{equation}

\begin{prop}
\label{p9_1_p2_1}
Let $A$ be a closed symmetric operator in a Hilbert space $H$, $\overline{D(A)} = H$, and
$\widetilde A$ be a self-adjoint extension of $A$, acting in a Hilbert space
 $\widetilde H\supseteq H$.
Let $\lambda_0\in \mathbb{R}_e$ be an arbitrary number.
For arbitrary complex numbers $\lambda_1$ and $\lambda_2$, $\lambda_1\not= \overline{\lambda_2}$,
and arbitrary elements $h_1\in \widetilde{\mathfrak{L}}_{\lambda_1}$,
$h_2\in \widetilde{\mathfrak{L}}_{\lambda_2}$, set
\begin{equation}
\label{f9_10_p2_1}
g_k := P^{\widetilde H}_H h_k = f_k + \varphi_k + \psi_k,\qquad k=1,2,
\end{equation}
where $f_k\in D(A)$, $\varphi_k\in \mathcal{N}_{ \overline{\lambda_0} }(A)$,
$\psi_k\in \mathcal{N}_{\lambda_0}(A)$.
Then
\begin{equation}
\label{f9_11_p2_1}
(h_1,h_2)_{\widetilde H} = (g_1,g_2)_H +
\frac{ \lambda_0 - \overline{\lambda_0} }{ \lambda_1 - \overline{\lambda_2} }
\left(
(\psi_1,\psi_2)_H - (\varphi_1,\varphi_2)_H
\right).
\end{equation}
\end{prop}
{\bf Proof. }
Notice that
$$ \frac{1}{ \lambda_1 - \overline{\lambda_2} }
\left[
(h_1, (\widetilde A - \lambda_2 E_{\widetilde H}) h_2)_{\widetilde H}
-
((\widetilde A - \lambda_1 E_{\widetilde H}) h_1, h_2)_{\widetilde H}
\right] $$
$$ = \frac{1}{ \lambda_1 - \overline{\lambda_2} }
\left[
(h_1, \widetilde A h_2)_{\widetilde H} - \overline{\lambda_2} (h_1, h_2)_{\widetilde H} -
(\widetilde A h_1, h_2)_{\widetilde H} + \lambda_1 (h_1, h_2)_{\widetilde H}
\right] $$
$$ = (h_1,h_2)_{\widetilde H}. $$
Using relation~(\ref{f9_8_p2_1}) to a transformation of the left-hand side of the last expression we get
$$ (h_1,h_2)_{\widetilde H} =
\frac{1}{ \lambda_1 - \overline{\lambda_2} }
\left[
(h_1, (A^* - \lambda_2 E_H) P^{\widetilde H}_H h_2)_{\widetilde H}
-
((A^* - \lambda_1 E_H) P^{\widetilde H}_H h_1, h_2)_{\widetilde H}
\right] $$
$$ =
\frac{1}{ \lambda_1 - \overline{\lambda_2} }
\left[
(g_1, (A^* - \lambda_2 E_H) g_2)_H
-
((A^* - \lambda_1 E_H) g_1, g_2)_H
\right] $$
$$ = (g_1,g_2)_H +
\frac{1}{ \lambda_1 - \overline{\lambda_2} }
\left[
(g_1, A^* g_2)_H
-
(A^* g_1, g_2)_H
\right]. $$
Substituting here the expression for $g_k$ from~(\ref{f9_10_p2_1}), and taking into account that
$A^* \varphi_k = \lambda_0 \varphi_k$, $A^* \psi_k = \overline{\lambda_0} \psi_k$, $k=1,2$,
we obtain the required equality~(\ref{f9_11_p2_1}).
$\Box$

Consider the operator-valued function $\mathfrak{B}_\lambda = \mathfrak{B}_\lambda(A,\widetilde A)$,
$\lambda\in \mathbb{R}_e$, and the operator
$\mathfrak{B}_\infty = \mathfrak{B}_\infty(A,\widetilde A)$ (see~(\ref{f6_10_0_p2_1}), (\ref{f6_10_1_p2_1}),
(\ref{f6_10_2_p2_1})).
As it was said above, the operatorsоператоры $\mathfrak{B}_\lambda$ ($\lambda\in \mathbb{R}_e$) and
$\mathfrak{B}_\infty$ are extensions of the operator $A$ in $H$. Therefore
\begin{equation}
\label{f9_50_p2_1}
\mathfrak{B}_\infty = A^*|_{\mathfrak{L}_\infty},\quad
\mathfrak{B}_\lambda = A^*|_{\mathfrak{L}_\lambda},\qquad
\lambda\in \mathbb{R}_e.
\end{equation}

Choose and fix an arbitrary number $\lambda_0\in \mathbb{R}_e$.
Consider the operator-valued function $\mathfrak{F}(\lambda) = \mathfrak{F}(\lambda; \lambda_0, A, \widetilde A)$,
$\lambda\in\Pi_{\lambda_0}$, defined by~(\ref{f6_14_1_p2_1}).
Choose arbitrary two points $\lambda',\lambda''\in\Pi_{\lambda_0}$, and arbitrary
elements $h'\in \widetilde{\mathfrak{L}}_{\lambda'}$, $h''\in \widetilde{\mathfrak{L}}_{\lambda''}$.
Then the elements $P^{\widetilde H}_H h'$ and $P^{\widetilde H}_H h''$ belong to
$D(\mathfrak{B}_{\lambda'})$ and $D(\mathfrak{B}_{\lambda''})$, respectively.
By~(\ref{f6_14_2_p2_1}) and the generalized Neumann formulas we may write:
$$ P^{\widetilde H}_H h' = f' + \mathfrak{F}(\lambda') \psi' - \psi',\quad
P^{\widetilde H}_H h'' = f'' + \mathfrak{F}(\lambda'') \psi'' - \psi'', $$
where
$f',f''\in D(A)$, $\psi',\psi''\in N_{\lambda_0}(A)$.
By Proposition~\ref{p9_1_p2_1} we get
$$ (h',h'')_{\widetilde H} = ( P^{\widetilde H}_H h', P^{\widetilde H}_H h'' )_H +
\frac{ \lambda_0 - \overline{\lambda_0} }{ \lambda' - \overline{\lambda''} }
\left(
(\psi',\psi'')_H - (\mathfrak{F}(\lambda') \psi', \mathfrak{F}(\lambda'') \psi'')_H
\right). $$
On the other hand, we may write:
$$ (h',h'')_{\widetilde H} = ( P^{\widetilde H}_H h' + P^{\widetilde H}_{\widetilde H\ominus H} h',
P^{\widetilde H}_H h'' + P^{\widetilde H}_{\widetilde H\ominus H} h'' )_{\widetilde H} $$
$$ = ( P^{\widetilde H}_H h', P^{\widetilde H}_H h'' )_H +
( P^{\widetilde H}_{\widetilde H\ominus H} h', P^{\widetilde H}_{\widetilde H\ominus H} h'' )_H. $$
Therefore
$$ ( P^{\widetilde H}_{\widetilde H\ominus H} h', P^{\widetilde H}_{\widetilde H\ominus H} h'' )_H =
\frac{ \lambda_0 - \overline{\lambda_0} }{ \lambda' - \overline{\lambda''} }
\left(
(\psi',\psi'')_H - (\mathfrak{F}(\lambda') \psi', \mathfrak{F}(\lambda'') \psi'')_H
\right). $$
Multiplying the both sides of the last equality by $\lambda' \overline{\lambda''}$
and taking into account~(\ref{f9_6_1_p2_1}) we get
\begin{equation}
\label{f9_51_p2_1}
( P^{\widetilde H}_{\widetilde H\ominus H} \widetilde A h', P^{\widetilde H}_{\widetilde H\ominus H}
\widetilde A h'' )_H =
(\lambda_0 - \overline{\lambda_0}) \frac{ \lambda' \overline{\lambda''}  }{ \lambda' - \overline{\lambda''} }
\left(
(\psi',\psi'')_H - (\mathfrak{F}(\lambda') \psi', \mathfrak{F}(\lambda'') \psi'')_H
\right).
\end{equation}
Equality~(\ref{f9_51_p2_1}) will be used in the sequel.
Also we shall need the following proposition.

\begin{prop}
\label{p9_2_p2_1}
If a sequence of elements of a Hilbert space is bounded and all its weakly convergent subsequences
have the same weak limit, then the sequence converges weakly.
\end{prop}
{\bf Proof. }
Let $\mathfrak{H} = \{ h_n \}_{n=1}^\infty$, $h_n\in H$,  be a bounded sequence of elements of a Hilbert space $H$.
Suppose to the contrary that all weakly convergent subsequences of
 $\mathfrak{H}$ converge to an element $h\in H$, but
the sequence $\mathfrak{H}$ does not converge to $h$.
In this case there exists an element $g\in H$ such that $(h_n,g)_H$ does not converge to $(h,g)_H$,
as $n\rightarrow\infty$. Therefore there exists a number $\varepsilon > 0$, and a subsequence
$\{ h_{n_k} \}_{k=1}^\infty$, 
such that
\begin{equation}
\label{f9_52_p2_1}
| (h_{n_k},g)_H - (h,g)_H | = | (h_{n_k} - h,g)_H | \geq \varepsilon,\qquad
k\in \mathbb{N}.
\end{equation}
Since the sequence $\{ h_{n_k} \}_{k=1}^\infty$ is bounded, it contains a weakly convergent
subsequence. This subsequence, by assumption, should converge weakly to to an element $h$.
But this is impossible according to~(\ref{f9_52_p2_1}).
The obtained contradiction completes the proof.
$\Box$

Let continue our considerations started before the statement of the last proposition.
Recall that the operator $\mathfrak{B}_\infty = \mathfrak{B}_\infty(A,\widetilde A)$
is a (not necessarily closed) symmetric extension of the operator $A$ in $H$.
By the Neumann formulas for it there corresponds an isometric operator
$\Phi_\infty = \Phi_\infty(\lambda_0; A,\widetilde A)$ with the domain $D(\Phi_\infty)\subseteq N_{\lambda_0}(A)$ and
the range $R(\Phi_\infty)\subseteq N_{ \overline{\lambda_0} }(A)$.
The operator $\Phi_\infty$ was already defined for the general case of a not necessarily densely defined
operator $A$ (below~(\ref{f6_14_2_p2_1})).

\begin{thm}
\label{t9_1_p2_1}
Let $A$ be a closed symmetric operator in a Hilbert space $H$, $\overline{D(A)} = H$, and
$\widetilde A$ be a self-adjoint extension of $A$, acting in a Hilbert space
 $\widetilde H\supseteq H$.
Let $\lambda_0\in \mathbb{R}_e$ be an arbitrary number,
$\{ \lambda_n \}_{n=1}^\infty$, $\lambda_n\in\Pi_{\lambda_0}$, be a number sequence,
tending to $\infty$,
$\{ \psi_n \}_{n=1}^\infty$, $\psi_n\in \mathcal{N}_{\lambda_0}(A)$, be a sequence,
weakly converging to an element $\psi$, and
\begin{equation}
\label{f9_52_1_p2_1}
\sup_{n\in \mathbb{N}} \left[
\frac{ |\lambda_n|^2 }{ | \mathop{\rm Im}\nolimits \lambda_n| }
\left(
\| \psi_n \|_H - \| \mathfrak{F}(\lambda_n; \lambda_0, A, \widetilde A) \psi_n \|_H
\right)
\right] < \infty.
\end{equation}
Then $\psi\in D(\Phi_\infty(\lambda_0; A,\widetilde A))$, and the sequence
$\{ \mathfrak{F}(\lambda_n) \psi_n \}_{n=1}^\infty$ converges weakly to
an element $\Phi_\infty \psi$.
Moreover, if the sequence $\{ \psi_n \}_{n=1}^\infty$ converges to $\psi$ in a strong sence, then
\begin{equation}
\label{f9_53_p2_1}
\lim_{n\to\infty} \mathfrak{F}(\lambda_n) \psi_n =
\lim_{n\to\infty} \mathfrak{F}(\lambda_n) \psi
= \Phi_\infty \psi.
\end{equation}
\end{thm}
{\bf Proof. }
Since elements $\lambda_n$ belong to $\Pi_{\lambda_0}$, then by relation~(\ref{f6_14_2_p2_1}) we get
$\mathfrak{B}_{\lambda_n}(A,\widetilde A) = A_{ \mathfrak{F}(\lambda_n; \lambda_0, A, \widetilde A) ,\lambda_0}$,
for every $n\in \mathbb{N}$.
Since $\psi_n\in \mathcal{N}_{\lambda_0}(A) = D(\mathfrak{F}(\lambda_n))$, then
by the generalized Neumann formulas
elements $\mathfrak{F}(\lambda_n)\psi_n - \psi_n \in D(\mathfrak{B}(\lambda_n)) = \mathfrak{L}_{\lambda_n}$,
$n\in \mathbb{N}$.
Consider a sequence of elements $g_n\in \widetilde{ \mathfrak{L} }_{\lambda_n}$, $n\in \mathbb{N}$,
such that
$$ P^{\widetilde H}_H g_n = \mathfrak{F}(\lambda_n)\psi_n - \psi_n,\qquad n\in  \mathbb{N}. $$
By~(\ref{f9_6_p2_1}) we may write:
$$ \widetilde A g_n = A^* ( \mathfrak{F}(\lambda_n)\psi_n - \psi_n ) +
P^{\widetilde H}_{\widetilde H\ominus H} \widetilde A g_n $$
$$ = \lambda_0 \mathfrak{F}(\lambda_n)\psi_n - \overline{\lambda_0} \psi_n +
P^{\widetilde H}_{\widetilde H\ominus H} \widetilde A g_n,\qquad n\in \mathbb{N}. $$
Therefore
\begin{equation}
\label{f9_54_p2_1}
P^{\widetilde H}_H \widetilde A g_n = \lambda_0 \mathfrak{F}(\lambda_n)\psi_n - \overline{\lambda_0} \psi_n,\qquad
n\in \mathbb{N}.
\end{equation}
Observe that
$$ \| P^{\widetilde H}_H \widetilde A g_n \|_H \leq
|\lambda_0| \| \mathfrak{F}(\lambda_n) \| \| \psi_n \|_H + |\lambda_0| \| \psi_n \|_H \leq K,\qquad
n\in \mathbb{N}, $$
where $K$ is a constant. Here we used the fact that $\mathfrak{F}(\lambda)$ is non-expanding and a weakly convergent sequence
is bounded.
Let us use formula~(\ref{f9_51_p2_1}), with $\lambda' = \lambda'' = \lambda_n$,
$h' = h'' = g_n$:
$$ ( P^{\widetilde H}_{\widetilde H\ominus H} \widetilde A g_n, P^{\widetilde H}_{\widetilde H\ominus H}
\widetilde A g_n )_H =
\mathop{\rm Im}\nolimits\lambda_0 \frac{ |\lambda_n|^2 }{ \mathop{\rm Im}\nolimits\lambda_n}
\left(
(\psi_n,\psi_n)_H - (\mathfrak{F}(\lambda_n) \psi_n, \mathfrak{F}(\lambda_n) \psi_n)_H
\right)
$$
$$ = \mathop{\rm Im}\nolimits\lambda_0 \frac{ |\lambda_n|^2 }{ \mathop{\rm Im}\nolimits\lambda_n}
\left(
\| \psi_n \|_H - \| \mathfrak{F}(\lambda_n) \psi_n \|_H
\right)
\left(
\| \psi_n \|_H + \| \mathfrak{F}(\lambda_n) \psi_n \|_H
\right) $$
$$ \leq |\mathop{\rm Im}\nolimits\lambda_0| \frac{ |\lambda_n|^2 }{ |\mathop{\rm Im}\nolimits\lambda_n| }
\left(
\| \psi_n \|_H - \| \mathfrak{F}(\lambda_n) \psi_n \|_H
\right)
2\| \psi_n \|_H. $$
From~(\ref{f9_52_1_p2_1}) it follows that $\| P^{\widetilde H}_{\widetilde H\ominus H} \widetilde A g_n \|_H \leq K_2$,
where $K_2$ is a constant. Thus, we conclude that the sequence
$\{ \widetilde A g_n \}_{n=1}^\infty$ is bounded.
By~(\ref{f9_6_1_p2_1}) the following equality holds:
$$ \frac{1}{ \lambda_n }  P^{\widetilde H}_{\widetilde H\ominus H} \widetilde A g_n =
 P^{\widetilde H}_{\widetilde H\ominus H} g_n,\qquad n\in \mathbb{N}, $$
and therefore
\begin{equation}
\label{f9_55_p2_1}
P^{\widetilde H}_{\widetilde H\ominus H} g_n \rightarrow 0,\qquad n\rightarrow\infty.
\end{equation}
Since a bounded set in a Hilbert space is weakly compact, then
from the sequence $\{ \widetilde A g_n \}_{n=1}^\infty$ we can select a weakly convergent
subsequence. Let
$\{ \widetilde A g_{n_k} \}_{k=1}^\infty$ be an arbitrary subsequence which converges weakly
to an element $h\in H$:
\begin{equation}
\label{f9_56_p2_1}
\widetilde A g_{n_k} \rightharpoonup h,\qquad k\rightarrow\infty.
\end{equation}
Then $P^{\widetilde H}_H \widetilde A g_{n_k} \rightharpoonup P^{\widetilde H}_H h$, as $k\rightarrow\infty$.
From~(\ref{f9_54_p2_1}) it follows that
\begin{equation}
\label{f9_57_p2_1}
\mathfrak{F}(\lambda_{n_k})\psi_{n_k}
\rightharpoonup
\frac{1}{\lambda_0} P^{\widetilde H}_H h + \frac{ \overline{\lambda_0} }{ \lambda_0 }\psi,\qquad
k\rightarrow\infty.
\end{equation}
Denote $\varphi :=
\frac{1}{\lambda_0} P^{\widetilde H}_H h + \frac{ \overline{\lambda_0} }{ \lambda_0 } \psi$.
By the weak completeness of a Hilbert space the element $\varphi$ belongs to
$\mathcal{N}_{ \overline{\lambda_0} }(A)$.
By the definition of elements $g_n$ we may write:
$$ g_n = \mathfrak{F}(\lambda_n)\psi_n - \psi_n
+ P^{\widetilde H}_{\widetilde H\ominus H} g_n,\qquad n\in \mathbb{N}. $$
By~(\ref{f9_55_p2_1}),(\ref{f9_57_p2_1}) we get
$$ g_{n_k} \rightharpoonup \varphi -\psi,\qquad k\rightarrow\infty. $$
Since $\widetilde A$ is closed, and the weak closeness of an operator is equivalent to its closeness,
then from~(\ref{f9_56_p2_1}) and the last relation it follows that
$\varphi -\psi\in D(\widetilde A)$ and
\begin{equation}
\label{f9_58_p2_1}
\widetilde A (\varphi -\psi) = h.
\end{equation}
Since elements $\varphi$ and $\psi$ belong to the space $H$,
then $\varphi -\psi\in D(\widetilde A)\cap H = \mathfrak{L}_\infty = D(\mathfrak{B}_\infty)$.
As it was mentioned before the statement of the theorem, by the Neumann formulas
for the operator $\mathfrak{B}_\infty$ there corresponds
an isometric operator $\Phi_\infty$. By~(\ref{f4_1_p2_1}) the following equality holds:
$$ D(\mathfrak{B}_\infty) = D(A) \dotplus (\Phi_\infty-E_H) D(\Phi_\infty). $$
Therefore $\varphi -\psi = f + \Phi_\infty u - u$, $f\in D(A)$, $u\in D(\Phi_\infty)$.
Since $\varphi\in\mathcal{N}_{ \overline{\lambda_0} }(A)$, $\psi\in \mathcal{N}_{ \lambda_0 }(A)$,
by the linear independence of the manifolds $D(A)$, $\mathcal{N}_{ \lambda_0 }(A)$ and
$\mathcal{N}_{ \overline{\lambda_0} }(A)$, we obtain that
$f=0$, $u=\psi$ и $\Phi_\infty u = \varphi$. Thus, the element $\psi$ belongs to
$D(\Phi_\infty)$ and $\Phi_\infty\psi = \varphi$.
Therefore relation~(\ref{f9_58_p2_1}) takes the following form:
$$ \widetilde A (\Phi_\infty\psi -\psi) = h. $$
Then relation~(\ref{f9_56_p2_1}) becomes
$$ \widetilde A g_{n_k} \rightharpoonup \widetilde A (\Phi_\infty\psi -\psi),\qquad k\rightarrow\infty. $$
Thus, an arbitrary weakly convergent subsequence of a sequence $\{ \widetilde A g_n \}_{n=1}^\infty$
has the same limit. By Proposition~\ref{p9_2_p2_1} we conclude that the sequence
$\{ \widetilde A g_n \}_{n=1}^\infty$ is a weakly convergent and its limit is equal to
$\widetilde A (\Phi_\infty\psi -\psi)$. Then
$$ P^{\widetilde H}_H \widetilde A g_n \rightharpoonup P^{\widetilde H}_H \widetilde A (\Phi_\infty\psi -\psi),\qquad
n\rightarrow\infty. $$
By relation~(\ref{f9_6_p2_1}) we may write:
$$ P^{\widetilde H}_H \widetilde A (\Phi_\infty\psi -\psi) =
A^* P^{\widetilde H}_H (\Phi_\infty\psi -\psi) = \lambda_0 \Phi_\infty\psi -
\overline{ \lambda_0 } \psi. $$
By~(\ref{f9_54_p2_1}) and two last relations we get
$$ \lambda_0 \mathfrak{F}(\lambda_n)\psi_n - \overline{\lambda_0} \psi_n
\rightharpoonup
\lambda_0 \Phi_\infty\psi -
\overline{ \lambda_0 } \psi,\qquad n\rightarrow\infty, $$
i.e.
$$ \mathfrak{F}(\lambda_n)\psi_n
\rightharpoonup \Phi_\infty\psi,\qquad n\rightarrow\infty, $$
Thus, the first assertion of the theorem is proved.
Suppose that the sequence $\{ \psi_n \}_{n=1}^\infty$ converges strongle to the element
$\psi$.

It is known that if a sequence $F = \{ f_n \}_{n=1}^\infty$, $f_n\in \mathcal{H}$,
of a Hilbert space $\mathcal{H}$ converges weakly to an element $f$ then
\begin{equation}
\label{f9_59_p2_1}
\| f \|_{\mathcal{H}} \leq \underline\lim_{n\to\infty} \| f_n \|_{ \mathcal{H} }.
\end{equation}
In fact, we may write
$$ \| f \|_{\mathcal{H}}^2 = \lim_{n\to\infty} |(f_n,f)_{\mathcal{H}}|. $$
If $\{ \| f_{n_k} \|_{\mathcal{H}} \}_{k=1}^\infty$ is a convergent
subsequence of a sequence $\{ \| f_n \|_{\mathcal{H}} \}_{n=1}^\infty$,
then
$$ \| f \|_{\mathcal{H}}^2 =
\lim_{k\to\infty} |(f_{n_k},f)_{\mathcal{H}}| \leq \lim_{k\to\infty} \| f_{n_k} \|_{\mathcal{H}}
\| f \|_{\mathcal{H}}. $$
From the last relation it follows the required estimate~(\ref{f9_59_p2_1}).

Applying estimate~(\ref{f9_59_p2_1}) to the sequence $\{ \mathfrak{F}(\lambda_n)\psi_n \}_{n=1}^\infty$ we obtain
an estimate:
$$ \| \Phi_\infty \psi \|_H \leq
\underline\lim_{n\to\infty} \| \mathfrak{F}(\lambda_n) \psi_n \|_H. $$
On the other hand, since $\mathfrak{F}(\lambda)$ is non-expanding,
then
$$ \overline{\lim}_{n\to\infty} \| \mathfrak{F}(\lambda_n) \psi_n \|_H \leq
\lim_{n\to\infty} \| \psi_n \|_H = \| \psi \|_H = \| \Phi_\infty \psi \|_H. $$
Therefore
$\lim_{n\to\infty} \| \mathfrak{F}(\lambda_n) \psi_n \|_H = \| \Phi_\infty \psi \|_H$,
and a sequence $\{ \mathfrak{F}(\lambda_n) \psi_n \}_{n=1}^\infty$ converges
strongly to an element $\Phi_\infty \psi$.
It remains to notice that
$$ \| \mathfrak{F}(\lambda_n) \psi - \Phi_\infty \psi \|_H =
\| \mathfrak{F}(\lambda_n) \psi - \mathfrak{F}(\lambda_n) \psi_n +
\mathfrak{F}(\lambda_n) \psi_n -\Phi_\infty \psi \|_H $$
$$ \leq \| \mathfrak{F}(\lambda_n) \psi - \mathfrak{F}(\lambda_n) \psi_n \|_H +
\| \mathfrak{F}(\lambda_n) \psi_n -\Phi_\infty \psi \|_H $$
$$ \leq \| \psi - \psi_n \|_H +
\| \mathfrak{F}(\lambda_n) \psi_n -\Phi_\infty \psi \|_H \rightarrow 0,\qquad n\rightarrow\infty. $$
Thus, relation~(\ref{f9_53_p2_1}) is proved.
$\Box$

\begin{cor}
\label{c9_1_p2_1}
In conditions of Theorem~\ref{t9_1_p2_1}, if one additionally assumes that the sequence
$\{ \lambda_n \}_{n=1}^\infty$ belongs to a set
$$ \Pi_{\lambda_0}^\varepsilon =
\{ \lambda\in\Pi_{\lambda_0}:\ \varepsilon < |\arg\lambda| < \pi - \varepsilon \},\qquad
0 < \varepsilon < \frac{\pi}{2}, $$
then condition~(\ref{f9_52_1_p2_1}) can be replaced by the following condition:
\begin{equation}
\label{f9_60_p2_1}
\sup_{n\in \mathbb{N}} \left[
|\lambda_n|
\left(
\| \psi_n \|_H - \| \mathfrak{F}(\lambda_n; \lambda_0, A, \widetilde A) \psi_n \|
\right)
\right] < \infty.
\end{equation}
\end{cor}
{\bf Proof. }
In fact, if the sequence $\{ \lambda_n \}_{n=1}^\infty$ belonfs to the set
$\Pi_{\lambda_0}^\varepsilon$, then
$$ \frac{ |\lambda_n| }{ |\mathop{\rm Im}\nolimits \lambda_n| } \leq M,\qquad n\in \mathbb{N}, $$
where $M$ is a constant, hot depending of $n$. Therefore from~(\ref{f9_60_p2_1})
it follows~(\ref{f9_52_1_p2_1}).
$\Box$

Notice that from relation~(\ref{f6_9_p2_1}) it follows that
\begin{equation}
\label{f9_61_p2_1}
\widetilde{\mathfrak{L}}_\lambda = (\widetilde A - \lambda E_{\widetilde H})^{-1} H ,\quad
\mathfrak{L}_\lambda = \mathbf{R}_\lambda(A) H,\qquad \lambda\in \mathbb{R}_e,
\end{equation}
where $\mathbf{R}_\lambda(A)$ is the generalized resolvent of  $A$, corresponding to the self-adjoint extension $\widetilde A$.

\begin{prop}
\label{p9_3_p2_1}
Let $A$ be a closed symmetric operator in a Hilbert space $H$, and
$\widetilde A$ be a self-adjoint extension of $A$, acting in
a Hilbert space $\widetilde H\supseteq H$.
Let $R_\lambda$ and $\mathbf{R}_\lambda$, $\lambda\in \mathbb{R}_e$, denote the resolvent of $\widetilde A$ and
the generalized resolvent of $A$, corresponding to the extension $\widetilde A$, respectively,
and $\lambda_0\in \mathbb{R}_e$, $0 < \varepsilon < \frac{\pi}{2}$ be arbitrary numbers.
Then the following relations hold:
$$ \lim_{\lambda\in\Pi_{\lambda_0}^\varepsilon,\ \lambda\to\infty} (-\lambda R_\lambda h) = h,\qquad
h\in \widetilde H; $$
$$ \lim_{\lambda\in\Pi_{\lambda_0}^\varepsilon,\ \lambda\to\infty} (-\lambda \mathbf{R}_\lambda g) = g,\qquad g\in H. $$
\end{prop}
{\bf Proof. }
Choose an arbitrary sequence $\{ \lambda_n \}_{n=1}^\infty$, $\lambda_n\in \Pi_{\lambda_0}^\varepsilon$,
tending to $\infty$, and an arbitrary element $h\in \widetilde H$. Let us check that
\begin{equation}
\label{f9_61_1_p2_1}
\lambda_n R_{\lambda_n} h \rightarrow - h,\qquad n\rightarrow\infty.
\end{equation}
Let $\lambda_n = \sigma_n + i\tau_n$, $\sigma_n,\tau_n\in \mathbb{R}$, $n\in \mathbb{N}$.
Notice that since numbers $\lambda_n$ belong to the set $\Pi_{\lambda_0}^\varepsilon$, it follows an estimate:
\begin{equation}
\label{f9_61_2_p2_1}
\left| \frac{ \sigma_n }{ \tau_n } \right| < M_1,\qquad n\in \mathbb{N},
\end{equation}
where $M_1$ is a constant, not depending on $n$.
Using the functional equation for the resolvent we write
$$ R_{\lambda_n} - R_{\lambda_0} = (\lambda_n - \lambda_0) R_{\lambda_n} R_{\lambda_0}, $$
\begin{equation}
\label{f9_61_3_p2_1}
\lambda_n R_{\lambda_n} R_{\lambda_0} = - R_{\lambda_0} + R_{\lambda_n}
+ \lambda_0 R_{\lambda_n} R_{\lambda_0},\qquad n\in \mathbb{N}.
\end{equation}
Since
$$ |\lambda_n| = |\tau_n| \left|
i + \frac{\sigma_n}{\tau_n}
\right| \leq
|\tau_n|
\left(
1 + \left| \frac{\sigma_n}{\tau_n} \right|
\right)
\leq |\tau_n| M_1, $$
then $\tau_n\rightarrow\infty$ при $n\rightarrow\infty$.
On the other hand,
by the functional equation for the resolvent
we may write:
$$ \| R_\lambda f \|_H^2 = (R_\lambda^* R_\lambda f,f)_H = (R_{\overline{\lambda}} R_\lambda f,f)_H
= \frac{1}{ \overline{\lambda} - \lambda } ((R_{\overline{\lambda}} - R_\lambda) f,f)_H $$
$$ = \frac{1}{ \lambda - \overline{\lambda} } ( R_\lambda f,f)_H - (R_\lambda^* f,f)_H)
= \frac{1}{ \mathop{\rm Im}\nolimits \lambda } \mathop{\rm Im}\nolimits (R_\lambda f,f)_H,\quad
\lambda\in \mathbb{R}_e,\ f\in \widetilde H. $$
Therefore
$$ \| R_{\lambda_n} f \|_H^2 \leq \frac{1}{ |\tau_n| }
| (R_{\lambda_n} f,f)_H | \leq
\frac{1}{ |\tau_n| }
\| R_{\lambda_n} f \|_H \| f \|_H; $$
$$ \| R_{\lambda_n} f \|_H \leq \frac{1}{ | \tau_n | } \| f \|_H,\qquad
n\in \mathbb{N}, f\in \widetilde H; $$
\begin{equation}
\label{f9_61_4_p2_1}
\| R_{\lambda_n} \| \leq \frac{1}{ | \tau_n | },\qquad
n\in \mathbb{N}.
\end{equation}
It follows that
$$ \| R_{\lambda_n} \|_H \rightarrow 0,\qquad n\rightarrow\infty. $$
By~(\ref{f9_61_3_p2_1}) we get
$$ u.-\lim_{n\to\infty} \lambda_n R_{\lambda_n} R_{\lambda_0} = - R_{\lambda_0}. $$
In particular, we obtain that
$$ \lambda_n R_{\lambda_n} u  = -u,\qquad u\in R_{\lambda_0}\widetilde H = D(\widetilde A). $$
By~(\ref{f9_61_4_p2_1}) and~(\ref{f9_61_2_p2_1}) the following inequality holds:
$$ \| \lambda_n R_{\lambda_n} \| \leq \left| \frac{ \sigma_n + i\tau_n }{ \tau_n } \right|
\leq
\left| \frac{ \sigma_n  }{ \tau_n } \right| + 1 < M_1 + 1,\qquad
n\in \mathbb{N}. $$
Since the sequence of linear operators $\{ \lambda_n R_{\lambda_n} \}_{n=1}^\infty$
converges (in the strong operator topology) on a dense set and the norms of the operators are uniformly bounded,
by the Banach-Steinhaus Theorem the sequence $\{ \lambda_n R_{\lambda_n} \}_{n=1}^\infty$
converges on the whole $\widetilde H$ to a continuous linear operator.
By the continuity, this limit operator coincides with $-E_{\widetilde H}$.
Thus, relation~(\ref{f9_61_1_p2_1}) is proved, and the required relations in the statement of the theorem follow.
$\Box$

\begin{thm}
\label{t9_2_p2_1}
Let $A$ be a closed symmetric operator in a Hilbert space $H$, $\overline{D(A)} = H$, and
$\widetilde A$ be a self-adjoint extension of $A$, acting in
a Hilbert space $\widetilde H\supseteq H$.
Let $\lambda_0\in \mathbb{R}_e$ be an arbitrary number,
$\psi$ be an arbitrary number from $D(\Phi_\infty(\lambda_0; A,\widetilde A))$, and a vector-valued function
$\psi(\lambda)$, $\lambda\in\Pi_{\lambda_0}$, with values in $\mathcal{N}_{\lambda_0}(A)$,
is defined by the following formula:
\begin{equation}
\label{f9_62_p2_1}
\psi(\lambda) = \frac{\lambda}{ \lambda_0 - \overline{\lambda_0} }
P^H_{ \mathcal{N}_{\lambda_0}(A) } (A^* - \lambda_0 E_H) \mathbf{R}_\lambda (\psi - \Phi_\infty\psi),
\end{equation}
where $\mathbf{R}_\lambda$ is a generalized resolvent of $A$, corresponding to the self-adjoint
extension $\widetilde A$.
Then
\begin{equation}
\label{f9_63_p2_1}
\lim_{\lambda\in\Pi_{\lambda_0}^\varepsilon,\ \lambda\to\infty} \psi(\lambda) = \psi,\quad
\lim_{\lambda\in\Pi_{\lambda_0}^\varepsilon,\ \lambda\to\infty}
\mathfrak{F}(\lambda; \lambda_0, A, \widetilde A) \psi(\lambda) = \Phi_\infty \psi,
\end{equation}
and there exists a finite limit
\begin{equation}
\label{f9_64_p2_1}
\lim_{\lambda,\zeta\in\Pi_{\lambda_0}^\varepsilon,\ \lambda,\zeta\to\infty}
\left\{
\frac{ \lambda \overline{\zeta} }{ \lambda - \overline{\zeta} }
\left[
( \psi(\lambda), \psi(\zeta) )_H - ( \mathfrak{F}(\lambda) \psi(\lambda), \mathfrak{F}(\zeta) \psi(\zeta) )_H
\right]
\right\}.
\end{equation}
Here $0 < \varepsilon < \frac{\pi}{2}$.
\end{thm}
{\bf Proof. }
Denote $g = \Phi_\infty(\lambda_0; A,\widetilde A) \psi - \psi \in D(\mathfrak{B}_\infty) =
\mathfrak{L}_\infty = D(\widetilde A)\cap H$, and by
$R_\lambda$ and $\mathbf{R}_\lambda$, $\lambda\in \mathbb{R}_e$, we denote the resolvent of $\widetilde A$ and
the generalized resolvent of $A$, corresponding to the extension $\widetilde A$, respectively.
Consider the following two vector-valued functions:
$$ \widetilde g (\lambda) = - \lambda R_\lambda g,\quad
g(\lambda) = P^{\widetilde H}_H \widetilde g(\lambda) = - \lambda \mathbf{R}_\lambda g,\qquad
\lambda\in\Pi_{\lambda_0}. $$
Then the function $\psi(\lambda)$ from the statement of the theorem takes the following form:
$$ \psi(\lambda) =
\frac{1}{ \lambda_0 - \overline{\lambda_0} }
P^H_{ \mathcal{N}_{\lambda_0}(A) } (A^* - \lambda_0 E_H) g(\lambda),\qquad \lambda\in\Pi_{\lambda_0}. $$
Notice that by~(\ref{f9_61_p2_1}), the element $g(\lambda)$ belongs to $\mathfrak{L}_\lambda$,
and therefore by~(\ref{f9_50_p2_1}) the definition of the function $\psi(\lambda)$ is correct.
By Proposition~\ref{p9_3_p2_1} we get:
\begin{equation}
\label{f9_64_1_p2_1}
\lim_{\lambda\in\Pi_{\lambda_0}^\varepsilon,\ \lambda\to\infty} \widetilde g(\lambda) = g,\qquad
\lim_{\lambda\in\Pi_{\lambda_0}^\varepsilon,\ \lambda\to\infty} g(\lambda) = g,\qquad
0 < \varepsilon < \frac{\pi}{2}.
\end{equation}
Moreover, by Proposition~\ref{p9_3_p2_1} with $h=\widetilde Ag$ we obtain that
\begin{equation}
\label{f9_65_p2_1}
\widetilde A g =
\lim_{\lambda\in\Pi_{\lambda_0}^\varepsilon,\ \lambda\to\infty} (-\lambda R_\lambda \widetilde A g) =
\lim_{\lambda\in\Pi_{\lambda_0}^\varepsilon,\ \lambda\to\infty} (-\lambda \widetilde A R_\lambda g) =
\lim_{\lambda\in\Pi_{\lambda_0}^\varepsilon,\ \lambda\to\infty} \widetilde A \widetilde g(\lambda).
\end{equation}
By~(\ref{f9_61_p2_1}) it follows that the element $\widetilde g(\lambda)$ belongs to
$\widetilde{\mathfrak{L}}_\lambda \subseteq D(\widetilde A)$.
By~(\ref{f9_5_1_p2_1}) we write
$$ g(\lambda) = P^{\widetilde H}_H \widetilde g(\lambda) \in D(A^*),\quad
P^{\widetilde H}_H \widetilde A \widetilde g(\lambda)
= A^* g(\lambda),\qquad \lambda\in\Pi_{\lambda_0}. $$
From the last expression and~(\ref{f9_65_p2_1}) we obtain that
$$ P^{\widetilde H}_H \widetilde A g =
\lim_{\lambda\in\Pi_{\lambda_0}^\varepsilon,\ \lambda\to\infty} P^{\widetilde H}_H \widetilde A \widetilde g(\lambda)
= \lim_{\lambda\in\Pi_{\lambda_0}^\varepsilon,\ \lambda\to\infty} A^* g(\lambda),\qquad
\lambda\in\Pi_{\lambda_0}, $$
or, taking into account property~(\ref{f9_5_1_p2_1}), we write:
\begin{equation}
\label{f9_66_p2_1}
\lim_{\lambda\in\Pi_{\lambda_0}^\varepsilon,\ \lambda\to\infty} A^* g(\lambda) = A^* g,\qquad
\lambda\in\Pi_{\lambda_0}.
\end{equation}
Passing to the limit in the above expression for the function $\psi(\lambda)$,
using~(\ref{f9_64_1_p2_1}) and~(\ref{f9_66_p2_1}),  we get:
$$ \lim_{\lambda\in\Pi_{\lambda_0}^\varepsilon,\ \lambda\to\infty} \psi(\lambda) =
\frac{1}{ \lambda_0 - \overline{\lambda_0} }
P^H_{ \mathcal{N}_{\lambda_0}(A) } (A^* - \lambda_0 E_H) g. $$
By the definition of the element $g$ the following equality holds:
$(A^* - \lambda_0 E_H) g = (\lambda_0 - \overline{\lambda_0}) \psi$.
Substituting it into
the previous relation we obtain the first relation in~(\ref{f9_63_p2_1}).

As it was already noticed, the element $g(\lambda)$ belongs to $\mathfrak{L}_\lambda =
D( \mathfrak{B}_\lambda(A,\widetilde A) )$, and therefore by the definition of the function
$\mathfrak{F}(\lambda; \lambda_0, A, \widetilde A)$ and the generalized Neumann formulas
we write:
$$ g(\lambda) = f(\lambda) + \mathfrak{F}(\lambda) \psi_1 (\lambda) - \psi_1 (\lambda),\qquad
\lambda\in\Pi_{\lambda_0}, $$
where $f(\lambda)\in D(A)$, $\psi_1(\lambda)\in \mathcal{N}_{\lambda_0}(A)$.
Using the above expression for $\psi(\lambda)$ we directly calculate:
$$ \psi(\lambda) = \frac{1}{ \lambda_0 - \overline{\lambda_0} }
P^H_{ \mathcal{N}_{\lambda_0}(A) } (A^* - \lambda_0 E_H) g(\lambda) = \psi_1(\lambda),\qquad
\lambda\in\Pi_{\lambda_0}. $$
Therefore
$$ g(\lambda) = f(\lambda) + \mathfrak{F}(\lambda) \psi (\lambda) - \psi (\lambda),\qquad
\lambda\in\Pi_{\lambda_0}. $$
Observe that
$$ \frac{1}{ \lambda_0 - \overline{\lambda_0} } P^H_{ \mathcal{N}_{\overline{\lambda_0}} }
( A^* - \overline{\lambda_0} E_H ) g(\lambda) = \mathfrak{F}(\lambda) \psi(\lambda),\qquad
\lambda\in\Pi_{\lambda_0}. $$
Consequently, by~(\ref{f9_64_1_p2_1}),(\ref{f9_66_p2_1}) and the definition of $g$,  we get:
$$ \lim_{\lambda\in\Pi_{\lambda_0}^\varepsilon,\ \lambda\to\infty}
\mathfrak{F}(\lambda) \psi(\lambda) =
\frac{1}{ \lambda_0 - \overline{\lambda_0} }
\lim_{\lambda\in\Pi_{\lambda_0}^\varepsilon,\ \lambda\to\infty}
P^H_{ \mathcal{N}_{\overline{\lambda_0}} }
( A^* - \overline{\lambda_0} E_H ) g(\lambda) $$
$$ = \frac{1}{ \lambda_0 - \overline{\lambda_0} }
P^H_{ \mathcal{N}_{\overline{\lambda_0}} } (A^* g - \overline{\lambda_0} g)
= \Phi_\infty \psi, $$
i.e. the second relation in~(\ref{f9_63_p2_1}) is satisfied, as well.

Let us use relation~(\ref{f9_51_p2_1}) for $\lambda' = \lambda$, $\lambda'' = \zeta$,
$h' = \widetilde g(\lambda)$, $h'' = \widetilde g(\zeta)$, $P^{\widetilde H}_H h' = g(\lambda)$,
$P^{\widetilde H}_H h'' = g(\zeta)$:
$$ ( P^{\widetilde H}_{\widetilde H\ominus H} \widetilde A \widetilde g(\lambda),
P^{\widetilde H}_{\widetilde H\ominus H} \widetilde A \widetilde g(\zeta) )_H =
(\lambda_0 - \overline{\lambda_0}) \frac{ \lambda \overline{\zeta}  }{ \lambda - \overline{\zeta} }
(
(\psi(\lambda),\psi(\zeta))_H $$
$$ - (\mathfrak{F}(\lambda) \psi(\lambda), \mathfrak{F}(\zeta) \psi(\zeta))_H
). $$
By relation~(\ref{f9_65_p2_1}), the left-hand side has a finite limit as
$\lambda,\zeta\rightarrow\infty$, $\lambda,\zeta\in\Pi_{\lambda_0}^\varepsilon$, and
the last assertion of the theorem follows.
$\Box$

\begin{prop}
\label{p9_4_p2_1}
Let $A$ be a closed symmetric operator in a Hilbert space $H$, $\overline{D(A)} = H$, and
$\widetilde A$ be a self-adjoint extension of $A$, acting in
a Hilbert space $\widetilde H\supseteq H$.
Suppose that the operator $A$ is a direct sum of two maximal symmetric
operators.
Let $\lambda_0\in \mathbb{R}_e$ be an arbitrary number,
$\psi$ be an arbitrary element from $D(\Phi_\infty(\lambda_0; A,\widetilde A))$.
Then there exists a finite limit
\begin{equation}
\label{f9_67_p2_1}
\lim_{\lambda,\zeta\in\Pi_{\lambda_0}^\varepsilon,\ \lambda,\zeta\to\infty}
\left\{
\frac{ \lambda \overline{\zeta} }{ \lambda - \overline{\zeta} }
\left[
( \psi, \psi )_H - ( \mathfrak{F}(\lambda; \lambda_0, A, \widetilde A) \psi,
\mathfrak{F}(\zeta; \lambda_0, A, \widetilde A) \psi )_H
\right]
\right\}.
\end{equation}
Here $0 < \varepsilon < \frac{\pi}{2}$.
\end{prop}
{\bf Proof. }
The operator $A$ has the following form: $A = A_1\oplus A_2$, where $A_k$ is a maximal symmetric operator
in a Hilbert space $H_k$, $k=1,2$, and $H = H_1\oplus H_2$.
Suppose at first that the domain of the regularity of $A_2$ is a half-plane
$\Pi_{\lambda_0}$, and the domain of the regularity of $A_1$ is $\Pi_{-\lambda_0}$.
Recalling a similar situation in~(\ref{f5_4_p2_1}) we may write:
$$ \mathcal{N}_\lambda(A) = \mathcal{N}_{\lambda}(A_1),\quad
\mathcal{N}_{\overline{\lambda}}(A) = \mathcal{N}_{\overline{\lambda}}(A_2),\qquad
\lambda\in\Pi_{\lambda_0}. $$
Calculate the function $\psi(\lambda)$, defined in~(\ref{f9_62_p2_1}).
Let $\psi'$ be an arbitrary element from $\mathcal{N}_{\lambda_0}(A)$, and $\lambda\in\Pi_{\lambda_0}$.
From~(\ref{f9_61_p2_1}) it follows that the element $\mathbf{R}_\lambda \psi'$ belongs to the manifold
$\mathfrak{L}_\lambda = D(\mathfrak{B}_\lambda(A,\widetilde A))$, which by~(\ref{f9_50_p2_1}) lies in $D(A^*)$.
Using~(\ref{f6_12_p2_1}) we write:
$$ (A^* - \lambda E_H) \mathbf{R}_\lambda \psi' = (\mathfrak{B}_\lambda - \lambda E_H) \mathbf{R}_\lambda \psi'
= (\mathfrak{B}_\lambda - \lambda E_H) (\mathfrak{B}_\lambda - \lambda E_H)^{-1} \psi' = \psi'. $$
Conversely, if an element $g$ belongs to $\mathfrak{L}_\lambda$ and
$(A^* - \lambda E_H) g = \psi'$, then $g = (\mathfrak{B}_\lambda - \lambda E_H)^{-1} \psi' =
\mathbf{R}_\lambda\psi'$.

Since $\mathbf{R}_\lambda \psi'\in \mathfrak{L}_\lambda = D(\mathfrak{B}_\lambda)$, then by definition
of the function $\mathfrak{F}(\lambda; \lambda_0, A, \widetilde A)$ and the generalized Neumann formulas we write:
$$ \mathbf{R}_\lambda \psi' = f + \mathfrak{F}(\lambda) \psi_1 - \psi_1, $$
where $f\in D(A)$, $\psi_1\in \mathcal{N}_{\lambda_0}(A)$.
On the other hand, consider an element
$$ g' = \frac{ \lambda - \lambda_0 }{ \lambda - \overline{\lambda_0} }
(A_2 - \lambda E_{H_2})^{-1} \mathfrak{F}(\lambda) \psi' + (\mathfrak{F}(\lambda) - E_H)
\frac{1}{ \lambda - \overline{\lambda_0} } \psi'. $$
Observe that $\mathfrak{F}(\lambda) \psi\in \mathcal{N}_{\overline{\lambda_0}}(A) =
\mathcal{N}_{\overline{\lambda}}(A_2)$, and the point $\lambda$ belongs to the domain of the regularity
of $A_2$, and therefore the first summand in the right-hand side is defined correctly.
By the generalized Neumann formulas the element $g'$ belongs to $D(\mathfrak{B}_\lambda) =
\mathfrak{L}_\lambda$.
Moreover, we have:
$$ (A^* - \lambda E_H) g' = \frac{ \lambda - \lambda_0 }{ \lambda - \overline{\lambda_0} }
(A-\lambda E_H)  (A_2 - \lambda E_{H_2})^{-1} \mathfrak{F}(\lambda) \psi' $$
$$ +
(\lambda_0 - \lambda)
\mathfrak{F}(\lambda)
\frac{1}{ \lambda - \overline{\lambda_0} } \psi'
- (\overline{\lambda_0} - \lambda)
\frac{1}{ \lambda - \overline{\lambda_0} } \psi' = \psi'. $$
By the above conciderations we conclude that $\mathbf{R}_\lambda \psi' = g'$.
Notice that
\begin{equation}
\label{f9_68_p2_1}
P^H_{\mathcal{N}_{\lambda_0}(A)} (A^* - \lambda_0 E_H) \mathbf{R}_\lambda \psi' =
P^H_{\mathcal{N}_{\lambda_0}(A)} (A^* - \lambda_0 E_H) g' =
\frac{ \lambda_0 - \overline{\lambda_0} }{ \lambda - \overline{\lambda_0} } \psi'.
\end{equation}
Consider an arbitrary element $\varphi\in \mathcal{N}_{\overline{\lambda_0}}(A) =
\mathcal{N}_{\overline{\lambda_0}}(A_2)\subseteq H_2$.
Since $\lambda$ belongs to the domain of the regularity of $A_2$, we may write:
$$ \mathbf{R}_\lambda \varphi = P^{\widetilde H}_H (\widetilde A - \lambda E_H)^{-1} \varphi =
(A_2 - \lambda E_H)^{-1} \varphi. $$
It follows that
$$ P^H_{\mathcal{N}_{\lambda_0}(A)} (A^* - \lambda_0 E_H) \mathbf{R}_\lambda \varphi =
P^H_{\mathcal{N}_{\lambda_0}(A)} (A^* - \lambda_0 E_H) (A_2 - \lambda E_H)^{-1} \varphi
$$
\begin{equation}
\label{f9_69_p2_1}
= P^H_{\mathcal{N}_{\lambda_0}(A)} (A - \lambda_0 E_H) (A - \lambda E_H)^{-1} \varphi = 0.
\end{equation}

Let $\psi$ be an arbitrary element from $D(\Phi_\infty(\lambda_0; A,\widetilde A))\subseteq
\mathcal{N}_{\lambda_0}(A)$, and $\lambda\in\Pi_{\lambda_0}$.
Using~(\ref{f9_68_p2_1}) with $\psi' = \psi$ and~(\ref{f9_69_p2_1}), calculate the function
$\psi(\lambda)$ from~(\ref{f9_62_p2_1}):
$$ \psi(\lambda) =
\frac{\lambda}{ \lambda_0 - \overline{\lambda_0} }
P^H_{ \mathcal{N}_{\lambda_0}(A) } (A^* - \lambda_0 E_H) \mathbf{R}_\lambda (\psi - \Phi_\infty\psi)
= \frac{\lambda}{ \lambda - \overline{\lambda_0} } \psi,\qquad \lambda\in\Pi_{\lambda_0}. $$
By Theorem~\ref{t9_2_p2_1} we conclude that there exists a limit~(\ref{f9_67_p2_1}).

In the case when the domain of the regularity of $A_2$ is a half-plane
$\Pi_{-\lambda_0}$, and the domain of the regularity of $A_1$ is $\Pi_{\lambda_0}$, we can reassign the operators
 $A_1$ and $A_2$.

In the case when the domain of the regularity of the operators $A_1$ and $A_2$ is the same half-plane
($\Pi_{\lambda_0}$ or $\Pi_{-\lambda_0}$), the operator $A$ is itself a maximal symmetric
operator.
In this case we set $H_1 = H$, $H_2 = \{ 0 \}$, $A_1 = A$ and $A_2 = 0_{H_2}$.
For the self-adjoint operator $A_2$ the both half-planes $\Pi_{\lambda_0}$ and
$\Pi_{-\lambda_0}$ are the domains of the regularity and we can apply the proved part
of the theorem.
$\Box$

\begin{thm}
\label{t9_3_p2_1}
Let $\lambda_0\in \mathbb{R}_e$ be an arbitrary number and $F(\lambda)$ be an analytic in the half-plane
$\Pi_{\lambda_0}$ operator-valued function, which values are linear non-expanding
operators with the domain  $D(F(\lambda)) = N_1$ and
the range $R(F(\lambda)) \subseteq N_2$, where $N_1$ and $N_2$ are some Hilbert spaces.
Let $\{ \lambda_n \}_{n=1}^\infty$, $\lambda_n\in\Pi_{\lambda_0}$, be a sequence of numbers,
converging to $\infty$, $\{ g_n \}_{n=1}^\infty$, $g_n\in N_1$, be a sequence of elements,
weakly converging to an element $g$, and
\begin{equation}
\label{f9_70_p2_1}
\sup_{n\in \mathbb{N}} \left[
\frac{ |\lambda_n|^2 }{ | \mathop{\rm Im}\nolimits \lambda_n| }
\left(
\| g_n \|_{N_1} - \| F(\lambda_n) g_n \|_{N_2}
\right)
\right] < \infty.
\end{equation}
Then there exists a finite limit
\begin{equation}
\label{f9_71_p2_1}
\lim_{\lambda,\zeta\in\Pi_{\lambda_0}^\varepsilon,\ \lambda,\zeta\to\infty}
\left\{
\frac{ \lambda \overline{\zeta} }{ \lambda - \overline{\zeta} }
\left[
( g, g )_{N_1} - ( F(\lambda) g, F(\zeta)g )_{N_2}
\right]
\right\}.
\end{equation}
Here $0 < \varepsilon < \frac{\pi}{2}$.
\end{thm}
{\bf Proof. }
Consider the operator $A$ of the following form: $A = A_1\oplus A_2$, where $A_k$ is a maximal densely defined
symmetric operator
in a Hilbert space $H_k$, $k=1,2$, and define $H = H_1\oplus H_2$.
Choose the operators $A_1$ and $A_2$ such that the domain of the regularity of the operator $A_2$ will be
$\Pi_{\lambda_0}$, the domain of the regularity of $A_1$ will be $\Pi_{-\lambda_0}$,
and it holds:
$$ \dim \mathcal{N}_{\lambda_0}(A_1) = \dim N_1,\quad \dim \mathcal{N}_{ \overline{\lambda_0} }(A_2) = \dim N_2. $$
The required operators $A_1$ and $A_2$ is not hard to construct by using the Neumann formulas.
As in~(\ref{f5_4_p2_1}) and in the proof of the previous proposition we may write:
$$ \mathcal{N}_\lambda(A) = \mathcal{N}_{\lambda}(A_1) = \dim N_1,\quad
\mathcal{N}_{\overline{\lambda}}(A) = \mathcal{N}_{\overline{\lambda}}(A_2) = \dim N_2,\qquad
\lambda\in\Pi_{\lambda_0}. $$
Consider arbitrary isometric operators $U$ and $W$, mapping respectively $N_1$ on
$\mathcal{N}_{\lambda_0}(A_1)$, and $N_2$ on $\mathcal{N}_{ \overline{\lambda_0} }(A_2)$.
Consider the following operator-valued function:
$$ F_1(\lambda) := W F(\lambda) U^{-1},\qquad \lambda\in\Pi_{\lambda_0}. $$
The function $F_1(\lambda)$ is analytic in the half-plane $\Pi_{\lambda_0}$, and its values are
linear non-expanding operators from
$\mathcal{N}_{\lambda_0}(A)$ into $\mathcal{N}_{ \overline{\lambda_0} }(A)$.
Since $A$ is densely defined, we can apply the Shtraus formula~(\ref{f5_2_p1_1}) 
Namely, the following formula
$$ \mathbf R_\lambda = \left( A_{F_1(\lambda)} - \lambda E_H \right)^{-1},\qquad  \lambda\in \Pi_{\lambda_0}, $$
defines a generalized resolvent $\mathbf R_\lambda$ of $A$.
It follows that
$$  A_{F_1(\lambda)}  = \mathbf R_\lambda^{-1} + \lambda E_H,\qquad \lambda\in \Pi_{\lambda_0}. $$
The generalized resolvent $\mathbf R_\lambda$ is generated by a (not necessarly unique)
self-adjoint extension $\widetilde A$ in a Hilbert space $\widetilde H\supseteq H$.
Consider the operator-valued function $\mathfrak{B}_\lambda (A,\widetilde A)$, $\lambda\in \mathbb{R}_e$.
Comparing the above expression for~$A_{F_1(\lambda)}$ with fornula~(\ref{f6_12_1_p2_1})
we conclude that
$$  \mathfrak{B}_\lambda =  A_{F_1(\lambda)},\qquad \lambda\in \Pi_{\lambda_0}. $$
Recalling the definition of the function $\mathfrak{F}(\lambda; \lambda_0, A, \widetilde A)$
we see that
$$ \mathfrak{F}(\lambda; \lambda_0, A, \widetilde A) = F_1(\lambda),\qquad \lambda\in \Pi_{\lambda_0}. $$
Denote
$$ \psi_n := U g_n,\quad \psi := Ug,\qquad n\in \mathbb{N}. $$
By the condition of the theorem we get that the sequence $\psi_n$ weakly converges to the element $\psi$.
From~(\ref{f9_70_p2_1}) it follows that
$$ \sup_{n\in \mathbb{N}} \left[
\frac{ |\lambda_n|^2 }{ | \mathop{\rm Im}\nolimits \lambda_n| }
\left(
\| \psi_n \|_H - \| \mathfrak{F}(\lambda; \lambda_0, A, \widetilde A) \psi_n \|_H
\right)
\right] < \infty. $$
By Theorem~\ref{t9_3_p2_1} we conclude that $\psi\in D(\Phi_\infty(\lambda_0; A,\widetilde A))$.
We may apply Proposition~\ref{p9_4_p2_1}, and~(\ref{f9_67_p2_1}) follows. From that relation it follows the required
relation~(\ref{f9_71_p2_1}).
$\Box$

Now we can state a theorem connecting the operator $\Phi_\infty(\lambda_0; A,\widetilde A)$
with the operator-valued function $\mathfrak{F}(\lambda; \lambda_0, A, \widetilde A)$.

\begin{thm}
\label{t9_4_p2_1}
Let $A$ be a closed symmetric operator in a Hilbert space $H$, $\overline{D(A)} = H$, and
$\widetilde A$ be a self-adjoint extension of $A$, acting in
a Hilbert space $\widetilde H\supseteq H$.
Let $\lambda_0\in \mathbb{R}_e$, $0 < \varepsilon < \frac{\pi}{2}$, be arbitrary numbers.
Then the following relations hold:
$$ D( \Phi_\infty(\lambda_0; A,\widetilde A) ) =
\left\{
\psi\in \mathcal{N}_{\lambda_0}(A):\ \right. $$
\begin{equation}
\label{f9_72_p2_1}
\left.
\underline{\lim}_{ \lambda\in\Pi_{\lambda_0}^\varepsilon, \lambda\to\infty }
\left[
|\lambda| ( \| \psi \|_H - \| \mathfrak{F}(\lambda; \lambda_0, A, \widetilde A) \psi \|_H )
\right] < +\infty
\right\},
\end{equation}
\begin{equation}
\label{f9_73_p2_1}
\Phi_\infty(\lambda_0; A,\widetilde A) \psi =
\lim_{ \lambda\in\Pi_{\lambda_0}^\varepsilon, \lambda\to\infty }
\mathfrak{F}(\lambda; \lambda_0, A, \widetilde A) \psi,\qquad \psi\in D(\Phi_\infty(\lambda_0; A,\widetilde A)).
\end{equation}
\end{thm}
{\bf Proof. } Suppose that $\psi\in D( \Phi_\infty(\lambda_0; A,\widetilde A) )
\subseteq \mathcal{N}_{\lambda_0}(A)$.
Check the inequality in~(\ref{f9_72_p2_1}).
If $\psi = 0$, then the inequality is obvious.
Assume that $\psi\not= 0$.
Consider an arbitrary sequence $\{ \lambda_n \}_{n=1}^\infty$, $\lambda_n\in
\Pi_{\lambda_0}^\varepsilon$, converging to $\infty$. By~(\ref{f9_62_p2_1})
we define the function $\psi(\lambda)$, $\lambda\in\Pi_{\lambda_0}$.
Set
$$ g_n := \psi(\lambda_n)\in \mathcal{N}_{\lambda_0}(A),\qquad n\in \mathbb{N}. $$
By Theorem~\ref{t9_2_p2_1} the following relation hold:
$$
\lim_{n\to\infty} g_n = \psi,\quad
\lim_{n\to\infty}
\mathfrak{F}(\lambda_n; \lambda_0, A, \widetilde A) g_n = \Phi_\infty \psi, $$
and there exists a finite limit
$$
\lim_{ n\to\infty}
\left\{
\frac{ |\lambda_n|^2 }{ \lambda_n - \overline{\lambda_n} }
\left[
( g_n, g_n )_H - ( \mathfrak{F}(\lambda_n) g_n, \mathfrak{F}(\lambda_n) g_n )_H
\right]
\right\} $$
$$ = \frac{1}{2i}
\lim_{ n\to\infty}
\left\{
\frac{ |\lambda_n|^2 }{ \mathop{\rm Im}\nolimits \lambda_n }
\left[
\|  g_n \|_H - \| \mathfrak{F}(\lambda_n) g_n \|_H
\right]
\left[
\|  g_n \|_H + \| \mathfrak{F}(\lambda_n) g_n \|_H
\right]
\right\}. $$
Notice that
$$ \|  g_n \|_H + \| \mathfrak{F}(\lambda_n) g_n \|_H \rightarrow
\|  \psi \|_H + \| \Phi_\infty \psi \|_H
= 2 \|  \psi \|_H,\qquad n\rightarrow\infty, $$
since $\Phi_\infty$ is isometric.
Since $\psi\not= 0$, it follows that there exists a finite limit
$$
\lim_{ n\to\infty}
\left\{
\frac{ |\lambda_n|^2 }{ \mathop{\rm Im}\nolimits \lambda_n }
\left[
\|  g_n \|_H - \| \mathfrak{F}(\lambda_n) g_n \|_H
\right]
\right\}.
$$
We can apply Theorem~\ref{t9_3_p2_1} to the sequences $\{ \lambda_n \}_{n=1}^\infty$, $\{ g_n \}_{n=1}^\infty$, which converges to $g$,
and the function $\mathfrak{F}(\lambda; \lambda_0, A, \widetilde A)$, mapping
$\mathcal{N}_{\lambda_0}(A)$ in $\mathcal{N}_{ \overline{\lambda_0} }(A)$,
and obtain that there exists a finite limit
$$
\lim_{  \lambda\in\Pi_{\lambda_0}^\varepsilon, \lambda\to\infty }
\left\{
\frac{ |\lambda|^2 }{ \mathop{\rm Im}\nolimits \lambda }
\left[
\|  \psi \|_H^2 - \| \mathfrak{F}(\lambda) \psi \|_H^2
\right]
\right\}
$$
\begin{equation}
\label{f9_74_p2_1}
\lim_{  \lambda\in\Pi_{\lambda_0}^\varepsilon, \lambda\to\infty }
\left\{
\frac{ |\lambda|^2 }{ \mathop{\rm Im}\nolimits \lambda }
\left[
\|  \psi \|_H - \| \mathfrak{F}(\lambda) \psi \|_H
\right]
\left[
\|  \psi \|_H + \| \mathfrak{F}(\lambda) \psi \|_H
\right]
\right\}.
\end{equation}
Since the operator $\mathfrak{F}(\lambda)$ is non-expanding, then
$$ \| \mathfrak{F}(\lambda_n) \psi - \Phi_\infty \psi \|_H =
\| \mathfrak{F}(\lambda_n) \psi - \mathfrak{F}(\lambda_n) g_n + \mathfrak{F}(\lambda_n) g_n
- \Phi_\infty \psi \|_H $$
$$ \leq \| \| \mathfrak{F}(\lambda_n) \psi - \mathfrak{F}(\lambda_n) g_n \|_H +
\| \mathfrak{F}(\lambda_n) g_n - \Phi_\infty \psi \|_H $$
$$ \leq \| \| \psi - g_n \|_H +
\| \mathfrak{F}(\lambda_n) g_n - \Phi_\infty \psi \|_H\rightarrow 0,\qquad n\rightarrow\infty. $$
Since the sequence $\{ \lambda_n \}_{n=1}^\infty$ was an arbitrary sequence from
$\Pi_{\lambda_0}^\varepsilon$, converging to $\infty$, then
$$ \lim_{\lambda\in \Pi_{\lambda_0}^\varepsilon,\ \lambda\to\infty} \mathfrak{F}(\lambda) \psi =
\Phi_\infty \psi. $$
Thus, relation~(\ref{f9_73_p2_1}) is proved.
We may write:
$$ \lim_{\lambda\in \Pi_{\lambda_0}^\varepsilon,\ \lambda\to\infty} (\| \mathfrak{F}(\lambda) \psi \|_H +
\| \psi \|_H) = \| \Phi_\infty \psi \|_H + \| \psi \|_H = 2 \| \psi \|_H. $$
From the latter relation and~(\ref{f9_74_p2_1}) it follows that there exists a finite limit:
$$ \lim_{  \lambda\in\Pi_{\lambda_0}^\varepsilon, \lambda\to\infty }
\left\{
\frac{ |\lambda|^2 }{ \mathop{\rm Im}\nolimits \lambda }
\left[
\|  \psi \|_H - \| \mathfrak{F}(\lambda) \psi \|_H
\right]
\right\} =: l(\psi) . $$
Therefore there exists a finite limit
$$ \lim_{  \lambda\in\Pi_{\lambda_0}^\varepsilon, \lambda\to\infty }
\left\{
\frac{ |\lambda|^2 }{ | \mathop{\rm Im}\nolimits \lambda | }
\left[
\|  \psi \|_H - \| \mathfrak{F}(\lambda) \psi \|_H
\right]
\right\} = |l(\psi)|. $$
In particular, there exists a finite limit
$$ \lim_{n\to\infty}
\left\{
\frac{ |\lambda_n|^2 }{ |\mathop{\rm Im}\nolimits \lambda_n| }
\left[
\|  \psi \|_H - \| \mathfrak{F}(\lambda_n) \psi \|_H
\right]
\right\} = |l(\psi)|. $$
Now suppose that the sequence $\{ \lambda_n \}_{n=1}^\infty$ is such that there exists a (finite or infinite) limit
$$ \lim_{n\to\infty}
\left\{
|\lambda_n|
\left[
\|  \psi \|_H - \| \mathfrak{F}(\lambda_n) \psi \|_H
\right]
\right\}. $$
Since $\left| \frac{ \mathop{\rm Im}\nolimits \lambda_n }{ \lambda_n } \right| \leq 1$, $n\in \mathbb{N}$,
then
$$ \lim_{n\to\infty}
\left\{
|\lambda_n|
\left[
\|  \psi \|_H - \| \mathfrak{F}(\lambda_n) \psi \|_H
\right]
\right\} $$
$$ =
\lim_{n\to\infty}
\left\{
\frac{ |\lambda_n|^2 }{ | \mathop{\rm Im}\nolimits \lambda | }
\left[
\|  \psi \|_H - \| \mathfrak{F}(\lambda_n) \psi \|_H
\right]
\right\}
\left| \frac{ \mathop{\rm Im}\nolimits \lambda_n }{ \lambda_n } \right| \leq |l(\psi)|. $$
The required estimate in~(\ref{f9_72_p2_1}) follows.

Conversely, suppose that $\psi\in \mathcal{N}_{\lambda_0}(A)$ and
$$ \underline{\lim}_{ \lambda\in\Pi_{\lambda_0}^\varepsilon, \lambda\to\infty }
\left[
|\lambda| ( \| \psi \|_H - \| \mathfrak{F}(\lambda; \lambda_0, A, \widetilde A) \psi \|_H )
\right] < +\infty. $$
This means that there exists a sequence of numbers $\{ \lambda_n \}_{n=1}^\infty$,
$\lambda_n\in \Pi_{\lambda_0}^\varepsilon$, converging to $\infty$, and such that
there exists a finite limit
$$ \lim_{n\to\infty}
\left\{
|\lambda_n|
\left[
\|  \psi \|_H - \| \mathfrak{F}(\lambda_n) \psi \|_H
\right]
\right\}. $$
In particular, this means that
$$ \sup_{n\in \mathbb{N}}
\left\{
|\lambda_n|
\left[
\|  \psi \|_H - \| \mathfrak{F}(\lambda_n) \psi \|_H
\right]
\right\} < \infty. $$
By Corollary~\ref{c9_1_p2_1} with $\psi_n = \psi$, $n\in \mathbb{N}$, we get
$\psi\in D( \Phi_\infty(\lambda_0; A,\widetilde A) )$.
$\Box$

\subsection{A connection of the operator $\Phi_\infty$ and the operator-valued function $\mathfrak{F}(\lambda)$.}

Our aim here will be to prove the following theorem.

\begin{thm}
\label{t10_1_p2_1}
Theorem~\ref{t9_4_p2_1} remains valid, if one removes the condition $\overline{D(A)} = H$.
\end{thm}
{\bf Proof. }
In the case $\overline{D(A)} = H$ the theorem is already proved.
Suppose now that $\overline{D(A)} \not= H$.
For the extension $\widetilde A$ there correspond an operator
$\mathfrak{B}_\infty = \mathfrak{B}_\infty(A,\widetilde A)$ and
operator-valued functions $\mathfrak{B}_\lambda =
\mathfrak{B}_\lambda(A,\widetilde A)$ and
$\mathfrak{F}(\lambda) =
\mathfrak{F}(\lambda;\lambda_0,A,\widetilde A)$ (see~(\ref{f6_10_2_p2_1}), (\ref{f6_10_0_p2_1}),
(\ref{f6_10_1_p2_1}), (\ref{f6_14_1_p2_1})).
Recall that for the operator $\mathfrak{B}_\infty$  there corresponds
an admissible with respect to $A$ isometric operator
$\Phi_\infty = \Phi_\infty(\lambda_0; A,\widetilde A)$ with the domain
$D(\Phi_\infty)\subseteq N_{\lambda_0}(A)$ and
the range $R(\Phi_\infty)\subseteq N_{ \overline{\lambda_0} }(A)$ (see a text below~(\ref{f6_14_2_p2_1})).
Observe that in the case $\overline{D(A)} = H$ the operator $\Phi_\infty$ was intensively used
in the previous subsection.

As before, the space $H$ decompose as in~(\ref{f5_1_p2_1}), where
$H_e = \widetilde H\ominus H$. The operator $A$ may be identified with an operator
$\mathcal{A} := A\oplus A_e$, $A_e = o_{H_e}$.
Thus, the operator $\widetilde A$ is a self-adjoint extension of $\mathcal{A}$.
By the generalized Neumann formulas for the self-adjoint extension $\widetilde A$ of
$\mathcal{A}$ there corresponds an admissible with respect to $A$
isometric operator $T$ with the domain $D(T)=\mathcal{N}_{\lambda_0}(\mathcal{A})$
and the range $R(T)=\mathcal{N}_{\overline{\lambda_0}}(\mathcal{A})$, having
the block representation~(\ref{f5_6_p2_1}),
where $T_{11} = P_{\mathcal{N}_{\overline{\lambda_0}}(A)} T P_{\mathcal{N}_{\lambda_0}(A)}$,
$T_{12} = P_{\mathcal{N}_{\overline{\lambda_0}}(A)} T P_{\mathcal{N}_{\lambda_0}(A_e)}$,
$T_{21} = P_{\mathcal{N}_{\overline{\lambda_0}}(A_e)} T P_{\mathcal{N}_{\lambda_0}(A)}$,
$T_{22} = P_{\mathcal{N}_{\overline{\lambda_0}}(A_e)} T P_{\mathcal{N}_{\lambda_0}(A_e)}$.
Moreover, we have $\widetilde A = \mathcal{A}_T$.
Notice that the operator $\mathfrak{B} = \mathfrak{B}(\lambda_0;\mathcal{A},T)$,
defined in~(\ref{f5_22_p2_1}), coincides with $\mathfrak{B}_\infty$.
By Theorem~\ref{t5_3_p2_1} we have:
$\mathfrak{B}(\lambda_0;\mathcal{A},T) = A_{\Phi(\lambda_0;\mathcal{A},T)}$.
Since the Neumann formulas define a one-to-one correspondence, then
$$ \Phi(\lambda_0;\mathcal{A},T) = \Phi_\infty(\lambda_0; A,\widetilde A). $$

Consider an arbitrary closed symmetric operator $A_1$ in a Hilbert
space $H_1$, $\overline{ D(A_1) } = H_1$, which has the same defect numbers as
 $A$.
Consider arbitrary isometric operators $U$ and $W$, mapping respectively
$\mathcal{N}_{\lambda_0}(A_1)$ on
$\mathcal{N}_{\lambda_0}(A)$, and
$\mathcal{N}_{ \overline{\lambda_0} }(A_1)$
on $\mathcal{N}_{ \overline{\lambda_0} }(A)$.
Then the following operator
$$ V := \left(
\begin{array}{cc} W^{-1} T_{11} U & W^{-1} T_{12} \\
T_{21} U & T_{22}\end{array}
\right) $$
$$ =
\left(
\begin{array}{cc} W^{-1} & 0 \\
0 & E\end{array}
\right)
\left(
\begin{array}{cc} T_{11}  & T_{12} \\
T_{21} & T_{22}\end{array}
\right)
\left(
\begin{array}{cc} U & 0 \\
0 & E \end{array}
\right), $$
maps isometrically all the subspace $\mathcal{N}_{\lambda_0}(A_1)\oplus \mathcal{N}_{\lambda_0}(A_e)$
on all the subspace $\mathcal{N}_{ \overline{\lambda_0} }(A_1)\oplus \mathcal{N}_{ \overline{\lambda_0} }(A_e)$.
By the first assertion of Theorem~\ref{t5_0_p2_1}, the operator $T_{22}$ is adnissible with respect to $A_e$.
By the second assertion of Theorem~\ref{t5_0_p2_1} we get that
$V$ is admissible with respect to $\mathcal{A}_1 := A_1\oplus A_e$, acting
in a Hilbert space $\widetilde H_1 := H_1\oplus H_e$.
By the generalized Neumann formulas for the operator $V$ there corresponds a self-adjoint operator
$\widetilde A_1\supseteq \mathcal{A}_1$ in a Hilbert space $\widetilde H_1$:
$\widetilde A_1 = (\mathcal{A}_1)_V$.

Consider the operator-valued functions
$\mathfrak{B}_\lambda' = \mathfrak{B}_\lambda(A_1,\widetilde A_1)$ и
$\mathfrak{F}'(\lambda) =
\mathfrak{F}(\lambda;\lambda_0,A_1,\widetilde A_1)$.
By~(\ref{f8_15_p2_1}) we may write:
$$ \mathfrak{F}(\lambda;\lambda_0,A_1,\widetilde A_1) =
V_{11} + V_{12} (E_{\mathcal{N}_{\lambda_0}(A_e)} - C_e(\lambda) V_{22})^{-1} C_e(\lambda) V_{21},\qquad
\lambda\in\Pi_{\lambda_0}, $$
where $C_e(\lambda) = C(\lambda; \lambda_0, A_e)$ is the characteristic function
of the operator $A_e$ in $H_e$. Then
$$ \mathfrak{F}'(\lambda) =
W^{-1} T_{11} U + W^{-1} T_{12} (E_{\mathcal{N}_{\lambda_0}(A_e)} - C_e(\lambda) T_{22})^{-1} C_e(\lambda)
T_{21} U = W^{-1} \mathfrak{F}(\lambda) U, $$
\begin{equation}
\label{f10_4_p2_1}
\lambda\in\Pi_{\lambda_0},
\end{equation}
what follows from~(\ref{f8_15_p2_1}).

Consider the following operators
$\mathfrak{B}'_\infty = \mathfrak{B}_\infty(A_1,\widetilde A_1)$ and
$\Phi'_\infty = \Phi_\infty(\lambda_0; A_1,\widetilde A_1)$.
Also consider the operator $\mathfrak{B}' = \mathfrak{B}(\lambda_0;\mathcal{A}_1,V)$,
defined in~(\ref{f5_22_p2_1}). It coincides with the operator $\mathfrak{B}'_\infty$.
By Theorem~\ref{t5_3_p2_1} we get
$\mathfrak{B}' = A_{\Phi(\lambda_0;\mathcal{A}_1,V)}$.
Since the Neumann formulas define a one-to-one correspondence then
$$ \Phi' := \Phi(\lambda_0;\mathcal{A},T) = \Phi'_\infty. $$
By Remark~\ref{r5_1_p2_1} and formula~(\ref{f5_10_p2_1}) we may write:
\begin{equation}
\label{f10_5_p2_1}
\Phi' \psi_1' = V_{11} \psi_1' + V_{12} (X_z(A_e) - V_{22})^{-1} V_{21} \psi_1',\qquad \psi_1'\in D(\Phi'),
\end{equation}
and
$$ \Phi \psi_1 = T_{11} \psi_1 + T_{12} (X_z(A_e) - T_{22})^{-1} T_{21} \psi_1,\qquad \psi_1\in D(\Phi). $$
By the definition of the operator $\Phi$, $D(\Phi)$ consists of elements of
$\psi_1\in \mathcal{N}_{\lambda_0}(A)$ such that there exists $\psi_2\in D(X_{\lambda_0}(A_e))$:
$$ T_{21}\psi_1 + T_{22}\psi_2 = X_{\lambda_0}(A_e)\psi_2, $$
and $D(\Phi')$ consists of elements
$\psi_1'\in \mathcal{N}_{\lambda_0}(A_1)$ such that there exists $\psi_2'\in D(X_{\lambda_0}(A_e))$:
$$ V_{21}\psi_1' + V_{22}\psi_2' =
T_{21} U \psi_1' + T_{22}\psi_2' =
X_{\lambda_0}(A_e)\psi_2'. $$
Therefore
\begin{equation}
\label{f10_5_1_p2_1}
D(\Phi) = U D(\Phi').
\end{equation}
By~(\ref{f10_5_p2_1}) we get
$$ \Phi' \psi_1' = W^{-1} T_{11} U \psi_1' +
W^{-1} T_{12} (X_z(A_e) - T_{22})^{-1} T_{21} U \psi_1' =
W^{-1} \Phi U \psi_1', $$
$$ \psi_1'\in D(\Phi'). $$
Therefore
\begin{equation}
\label{f10_6_p2_1}
\Phi' = W^{-1} \Phi U.
\end{equation}
By Theorem~\ref{t9_4_p2_1} for the operator $A_1$ in a Hilbert space $H_1$
and its self-adjoint extension $\widetilde A_1$ in $\widetilde H_1$ we get
$$ D( \Phi_\infty' ) =
\left\{
\psi'\in \mathcal{N}_{\lambda_0}(A_1):\
\underline{\lim}_{ \lambda\in\Pi_{\lambda_0}^\varepsilon, \lambda\to\infty }
\left[
|\lambda| ( \| \psi' \|_{H_1} - \| \mathfrak{F}'(\lambda) \psi' \|_H )
\right] < +\infty
\right\}, $$
$$
\Phi'_\infty \psi' =
\lim_{ \lambda\in\Pi_{\lambda_0}^\varepsilon, \lambda\to\infty }
\mathfrak{F}'(\lambda) \psi',\qquad \psi'\in D(\Phi_\infty'). $$
By relations~(\ref{f10_5_1_p2_1}), (\ref{f10_6_p2_1}) and (\ref{f10_4_p2_1})
we obtain the required relations~(\ref{f9_72_p2_1}) and~(\ref{f9_73_p2_1}).
$\Box$

\subsection{Shtraus's formula for the generalized resolvents of symmetric operator.}

Consider a closed symmetric operator $A$ in a Hilbert space $H$.
Choose and fix an arbitrary point $\lambda_0\in \mathbb{R}_e$.
A function $F(\lambda)\in \mathcal{S}(\Pi_{\lambda_0}; \mathcal{N}_{\lambda_0}(A),
\mathcal{N}_{ \overline{\lambda_0} }(A))$ is said to be  {\it $\lambda_0$-admissible (admissible) with respect to the operator $A$},
if the validity of
\begin{equation}
\label{f11_1_p2_1}
\lim_{\lambda\in\Pi_{\lambda_0}^\varepsilon,\ \lambda\to\infty} F(\lambda) \psi = X_{\lambda_0}\psi,
\end{equation}
\begin{equation}
\label{f11_2_p2_1}
\underline{\lim}_{\lambda\in\Pi_{\lambda_0}^\varepsilon,\ \lambda\to\infty}
\left[
|\lambda| (\| \psi \|_H - \| F(\lambda) \psi \|_H)
\right] < +\infty,
\end{equation}
for some $\varepsilon$: $0<\varepsilon <\frac{\pi}{2}$,
implies $\psi = 0$.

A set of all operator-valued functions $F(\lambda)\in \mathcal{S}(\Pi_{\lambda_0}; \mathcal{N}_{\lambda_0}(A),
\mathcal{N}_{ \overline{\lambda_0} }(A))$, which are $\lambda_0$-admissible with respect to the operator $A$,
we shall denote by
$$ \mathcal{S}_{a; \lambda_0} (\Pi_{\lambda_0}; \mathcal{N}_{\lambda_0}(A),
\mathcal{N}_{ \overline{\lambda_0} }(A)) = \mathcal{S}_{a} (\Pi_{\lambda_0}; \mathcal{N}_{\lambda_0}(A),
\mathcal{N}_{ \overline{\lambda_0} }(A)). $$

In the case $\overline{D(A)} = H$, we have $D(X_{\lambda_0}) = \{ 0 \}$.
Therefore in this case an arbitrary function from
$\mathcal{S}(\Pi_{\lambda_0}; \mathcal{N}_{\lambda_0}(A),
\mathcal{N}_{ \overline{\lambda_0} }(A))$ is admissible with respect to $A$.
Thus, {\it if $\overline{D(A)} = H$, then
$\mathcal{S}_{a;\lambda_0}(\Pi_{\lambda_0}; \mathcal{N}_{\lambda_0}(A),
\mathcal{N}_{ \overline{\lambda_0} }(A)) =
\mathcal{S}(\Pi_{\lambda_0}; \mathcal{N}_{\lambda_0}(A),
\mathcal{N}_{ \overline{\lambda_0} }(A))$}.

\begin{prop}
\label{p11_1_p2_1}
Let $A$ be a closed symmetric operator in a Hilbert space $H$,  and $\lambda_0\in \mathbb{R}_e$ be an arbitrary number.
If $F(\lambda)\in \mathcal{S}(\Pi_{\lambda_0}; \mathcal{N}_{\lambda_0}(A),
\mathcal{N}_{ \overline{\lambda_0} }(A))$ is
$\lambda_0$-admissible with respect to $A$, then
$F(\zeta)$ is a $\lambda_0$-admissible operator with respect to $A$,
for all $\zeta$ from $\Pi_{\lambda_0}$.
\end{prop}
{\bf Proof. }
Choose an arbitrary point $\zeta$ from $\Pi_{\lambda_0}$. Suppose that for an element
$\psi\in D(X_{\lambda_0}(A))$ we have the equality:
$$ F(\zeta) \psi = X_{\lambda_0}(A) \psi. $$
Since the forbidden ooperator is isometric, we get
$$ \| F(\zeta) \psi \|_H = \| \psi \|_H. $$
On the other hand, since the operator $F(\lambda)$ is non-expanding, we may write:
$$ \| F(\lambda) \psi \|_H \leq \| \psi \|_H,\qquad \lambda\in\Pi_{\lambda_0}. $$
By the maximum principle for analytic vector-valued functions we get
$$ F(\lambda) \psi = X_{\lambda_0}(A) \psi,\qquad \lambda\in\Pi_{\lambda_0}. $$
Therefore
$\| F(\lambda) \psi \|_H = \| \psi \|_H$, $\lambda\in\Pi_{\lambda_0}$, and
relations~(\ref{f11_1_p2_1}),(\ref{f11_2_p2_1}) hold.
Since $F(\lambda)$ is admissible with respect to $A$, then $\psi = 0$.
Thus, the operator $F(\zeta)$ is admissible with respect to $A$.
$\Box$

\begin{thm}
\label{t11_1_p2_1}
Let $A$ be a closed symmetric operator in a Hilbert space $H$, and $\lambda_0\in \mathbb{R}_e$ be an arbitrary
point.
An arbitrary generalized resolvent $\mathbf{R}_{s;\lambda}$ of the operator $A$ has the following form:
\begin{equation}
\label{f11_3_p2_1}
\mathbf R_{s;\lambda} = \left\{ \begin{array}{cc}
\left( A_{F(\lambda)} - \lambda E_H \right)^{-1}, & \lambda\in \Pi_{\lambda_0}\\
\left( A_{F^*(\overline{\lambda})} - \lambda E_H \right)^{-1}, &
\overline{\lambda}\in \Pi_{\lambda_0}
\end{array}
\right.,
\end{equation}
where $F(\lambda)$ is a function from $\mathcal{S}_{a;\lambda_0}(\Pi_{\lambda_0}; \mathcal{N}_{\lambda_0},
\mathcal{N}_{\overline{\lambda_0}})$.
Conversely, an arbitrary function $F(\lambda)\in \mathcal{S}_{a;\lambda_0}(\Pi_{\lambda_0}; \mathcal{N}_{\lambda_0},
\mathcal{N}_{\overline{\lambda_0}})$
defines by relation~(\ref{f11_3_p2_1}) a generalized resolvent
$\mathbf{R}_{s;\lambda}$ of the operator $A$.
Moreover, for different functions from
$\mathcal{S}_{a;\lambda_0}(\Pi_{\lambda_0}; \mathcal{N}_{\lambda_0},
\mathcal{N}_{\overline{\lambda_0}})$
there correspond different generalized resolvents of the operator $A$.
\end{thm}
{\bf Proof. }
Let $A$ be a closed symmetric operator in a Hilbert space $H$, and
$\lambda_0\in \mathbb{R}_e$ be a fixed number.
Consider an arbitrary generalized resolvent $\mathbf{R}_{s;\lambda} (A)$ of $A$.
It is generated by a self-adjoint extension $\widetilde A$ of $A$ in a Hilbert space
 $\widetilde H\supseteq H$.
Consider the function $\mathfrak{B}_\lambda(A,\widetilde A)$
(see~(\ref{f6_10_0_p2_1}),(\ref{f6_10_1_p2_1})) and the function
$\mathfrak{F}(\lambda; \lambda_0, A, \widetilde A)$ (see~(\ref{f6_14_1_p2_1})).
From relations~(\ref{f6_12_p2_1}) and~(\ref{f6_15_p2_1}) it follows~(\ref{f11_3_p2_1}),
where $F(\lambda)=\mathfrak{F}(\lambda)$ is a function from $\mathcal{S}(\Pi_{\lambda_0}; \mathcal{N}_{\lambda_0},
\mathcal{N}_{\overline{\lambda_0}})$,
which values are admissible with respect to $A$ operators.
Let us show that the function $F(\lambda)$ is admissible with respect to $A$.
Suppose that for an element $\psi\in D(X_{\lambda_0}(A))\subseteq \mathcal{N}_{\lambda_0}(A)$ hold~(\ref{f11_1_p2_1}) and~(\ref{f11_2_p2_1}).
By Theorem~\ref{t10_1_p2_1} we get
$\psi\in D(\Phi_\infty(\lambda_0; A,\widetilde A))$, and
$$ \Phi_\infty(\lambda_0; A,\widetilde A) \psi = X_{\lambda_0}(A) \psi. $$
Since $\Phi_\infty$ is admissible with respect to $A$, then $\psi = 0$.
Therefore the function $F(\lambda)$ is admissible with respect to $A$.

Conversely, let $F(\lambda)$ be an arbitrary function from
$\mathcal{S}_{a;\lambda_0}(\Pi_{\lambda_0}; \mathcal{N}_{\lambda_0}, \mathcal{N}_{\overline{\lambda_0}})$.
Consider an arbitrary closed symmetric operator $A_1$ in a Hilbert space $H_1$, $\overline{ D(A_1) } = H_1$,
which has the same defect numbers as  $A$.
Consider arbitrary isometric operators $U$ and $W$, whic map respectively
$\mathcal{N}_{\lambda_0}(A_1)$ on
$\mathcal{N}_{\lambda_0}(A)$, and
$\mathcal{N}_{ \overline{\lambda_0} }(A_1)$
on $\mathcal{N}_{ \overline{\lambda_0} }(A)$.
Set
\begin{equation}
\label{f11_3_1_p2_1}
F_1(\lambda) = W^{-1} F(\lambda) U,\qquad \lambda\in\Pi_{\lambda_0}.
\end{equation}
Еру агтсешщт $F_1(\lambda)$ иудщтпы ещ $\mathcal{S}(\Pi_{\lambda_0}; \mathcal{N}_{\lambda_0}(A_1),
\mathcal{N}_{\overline{\lambda_0}}(A_1))$.
Since in the case of the densely defined operators the theorem was proved, we can assert that $F_1(\lambda)$
generates a generalized resolvent
$\mathbf{R}_\lambda(A_1)$ of the operator $A_1$ in $H_1$:
\begin{equation}
\label{f11_4_p2_1}
\mathbf R_{\lambda}(A_1) = \left\{ \begin{array}{cc}
\left( (A_1)_{F_1(\lambda)} - \lambda E_H \right)^{-1}, & \lambda\in \Pi_{\lambda_0}\\
\left( (A_1)_{F^*(\overline{\lambda})} - \lambda E_H \right)^{-1}, &
\overline{\lambda}\in \Pi_{\lambda_0}
\end{array}
\right..
\end{equation}
The generalized resolvent $\mathbf{R}_\lambda(A_1)$ is generated by a self-adjoint extension
$\widetilde A_1$ in a Hilbert space $\widetilde H_1\supseteq H_1$.
Set
$$ H_e := \widetilde H_1\ominus H_1. $$
The operator $A_1$ may be considered as an operator $\mathcal{A}_1 := A_1\oplus A_e$ in a Hilbert space
$\widetilde H_1 = H_1\oplus H_e$, where $A_e = o_{H_e}$.

By the generalized Neumann formulas for the self-adjoint extension $\widetilde A_1$ of
$\mathcal{A}_1$ there correspond an admissible with respect to $A_1$
isometric operator $T$ with the domain $D(T)=\mathcal{N}_{\lambda_0}(\mathcal{A}_1)$
and the range $R(T)=\mathcal{N}_{\overline{\lambda_0}}(\mathcal{A}_1)$, having
the block representation~(\ref{f5_6_p2_1}),
where $T_{11} = P_{\mathcal{N}_{\overline{\lambda_0}}(A_1)} T P_{\mathcal{N}_{\lambda_0}(A_1)}$,
$T_{12} = P_{\mathcal{N}_{\overline{\lambda_0}}(A_1)} T P_{\mathcal{N}_{\lambda_0}(A_e)}$,
$T_{21} = P_{\mathcal{N}_{\overline{\lambda_0}}(A_e)} T P_{\mathcal{N}_{\lambda_0}(A_1)}$,
$T_{22} = P_{\mathcal{N}_{\overline{\lambda_0}}(A_e)} T P_{\mathcal{N}_{\lambda_0}(A_e)}$.
Moreover, we have $\widetilde A_1 = (\mathcal{A}_1)_T$.

Notice that the operator $\mathfrak{B} = \mathfrak{B}(\lambda_0;\mathcal{A}_1,T)$,
defined by~(\ref{f5_22_p2_1}), coincides with the operator $\mathfrak{B}_\infty(A_1,\widetilde A_1)$.
By Theorem~\ref{t5_3_p2_1} we get:
$\mathfrak{B}(\lambda_0;\mathcal{A}_1,T) = A_{\Phi(\lambda_0;\mathcal{A}_1,T)}$.
Since the Neumann formulas establish a one-to-one correspondence, then
$$ \Phi' := \Phi(\lambda_0;\mathcal{A}_1,T) = \Phi_\infty(\lambda_0; A_1,\widetilde A_1). $$
For the extension $\widetilde A_1$ there correspond
functions $\mathfrak{B}'_\lambda =
\mathfrak{B}_\lambda(A_1,\widetilde A_1)$ and
$\mathfrak{F}'(\lambda) =
\mathfrak{F}(\lambda;\lambda_0,A_1,\widetilde A_1)$.
By~(\ref{f8_15_p2_1}) we obtain the following equality:
\begin{equation}
\label{f11_4_0_p2_1}
\mathfrak{F}'(\lambda) =
T_{11} + T_{12} (E_{\mathcal{N}_{\lambda_0}(A_e)} - C_e(\lambda) T_{22})^{-1} C_e(\lambda) T_{21},\qquad
\lambda\in\Pi_{\lambda_0},
\end{equation}
where $C_e(\lambda) = C(\lambda; \lambda_0, A_e)$ is the characteristic
function of the operator $A_e$ in $H_e$.

Since $T$ is admissible with respect to $\mathcal{A}_1$, by Theorem~\ref{t5_0_p2_1}
the operator $T_{22}$ is admissible with respect to $A_e$. By Remark~\ref{r5_1_p2_1} the operator
$\Phi' = \Phi(\lambda_0;\mathcal{A}_1,T)$ admits the following representation:
\begin{equation}
\label{f11_4_1_p2_1}
\Phi' \psi_1' = T_{11} \psi_1' + T_{12} (X_{\lambda_0}(A_e) - T_{22})^{-1} T_{21} \psi_1',\qquad
\psi_1'\in D(\Phi').
\end{equation}
By Theorem~\ref{t10_1_p2_1} we get
\begin{equation}
\label{f11_4_2_p2_1}
D( \Phi' ) =
\left\{
\psi'\in \mathcal{N}_{\lambda_0}(A_1):\
\underline{\lim}_{ \lambda\in\Pi_{\lambda_0}^\varepsilon, \lambda\to\infty }
\left[
|\lambda| ( \| \psi' \|_{H_1} - \| \mathfrak{F}'(\lambda) \psi' \|_{H_1} )
\right] < +\infty
\right\},
\end{equation}
\begin{equation}
\label{f11_4_3_p2_1}
\Phi' \psi' =
\lim_{ \lambda\in\Pi_{\lambda_0}^\varepsilon, \lambda\to\infty }
\mathfrak{F}'(\lambda) \psi',\qquad \psi'\in D(\Phi'),
\end{equation}
where $0 < \varepsilon < \frac{\pi}{2}$.

Comparing relations~(\ref{f11_4_p2_1}) and~(\ref{f6_12_p2_1}), taking into account~(\ref{f6_14_2_p2_1}) we see that
\begin{equation}
\label{f11_4_4_p2_1}
\mathfrak{F}'(\lambda) =  F_1(\lambda),\qquad \lambda\in\Pi_{\lambda_0}.
\end{equation}

Consider the following operator
$$ V = \left(
\begin{array}{cc} V_{11} & V_{12} \\
V_{21} & V_{22}\end{array}
\right)
:= \left(
\begin{array}{cc} W T_{11} U^{-1} & W T_{12} \\
T_{21} U^{-1} & T_{22}\end{array}
\right) $$
$$ =
\left(
\begin{array}{cc} W & 0 \\
0 & E\end{array}
\right)
\left(
\begin{array}{cc} T_{11}  & T_{12} \\
T_{21} & T_{22}\end{array}
\right)
\left(
\begin{array}{cc} U^{-1} & 0 \\
0 & E \end{array}
\right), $$
whic maps isometrically all the subspace $\mathcal{N}_{\lambda_0}(A)\oplus \mathcal{N}_{\lambda_0}(A_e)$
on the whole subspace $\mathcal{N}_{ \overline{\lambda_0} }(A)\oplus \mathcal{N}_{ \overline{\lambda_0} }(A_e)$.

Сщтышвук фт щзукфещк $\mathcal{A} = A\oplus A_e$ шт ф Ршдиуке ызфсу
$\widetilde H := H\oplus H_e$.
As it was already noticed, the operator $T_{22}$ is admissible with respect to $A_e$.
By Remark~\ref{r5_1_p2_1} we obtain that the operator
$\Phi = \Phi(\lambda_0;\mathcal{A},V)$ admits the followig representation:
$$ \Phi \psi_1 = V_{11} \psi_1 + V_{12} (X_{\lambda_0}(A_e) - V_{22})^{-1} V_{21} \psi_1,\qquad
\psi_1\in D(\Phi). $$

By the definition of the operator $\Phi$, $D(\Phi)$ consists of elements
$\psi_1\in \mathcal{N}_{\lambda_0}(A)$ such that there exists $\psi_2\in D(X_{\lambda_0}(A_e))$:
$$ V_{21}\psi_1 + V_{22}\psi_2 = X_{\lambda_0}(A_e)\psi_2, $$
and $D(\Phi')$ consists of elements
$\psi_1'\in \mathcal{N}_{\lambda_0}(A_1)$ such that there exists $\psi_2'\in D(X_{\lambda_0}(A_e))$:
$$ T_{21}\psi_1' + T_{22}\psi_2' =
V_{21} U \psi_1' + T_{22}\psi_2' =
X_{\lambda_0}(A_e)\psi_2'. $$
It follows that
\begin{equation}
\label{f11_5_p2_1}
D(\Phi) = U D(\Phi').
\end{equation}
From~(\ref{f11_4_1_p2_1}) it follows that
$$ \Phi' \psi_1' = T_{11} \psi_1' + T_{12} (X_{\lambda_0}(A_e) - T_{22})^{-1} T_{21} \psi_1' $$
$$ = W^{-1} V_{11} U \psi_1' +
W^{-1} V_{12} (X_{\lambda_0}(A_e) - V_{22})^{-1} V_{21} U \psi_1' $$
$$ = W^{-1} \Phi U \psi_1',\qquad
\psi_1'\in D(\Phi'). $$
Therefore
\begin{equation}
\label{f11_6_p2_1}
\Phi' = W^{-1} \Phi U.
\end{equation}
From~(\ref{f11_4_2_p2_1}),(\ref{f11_4_3_p2_1}) it follows that
\begin{equation}
\label{f11_7_p2_1}
D( \Phi ) =
\left\{
\psi\in \mathcal{N}_{\lambda_0}(A):\
\underline{\lim}_{ \lambda\in\Pi_{\lambda_0}^\varepsilon, \lambda\to\infty }
\left[
|\lambda| ( \| \psi \|_{H} - \| F(\lambda) \psi \|_{H} )
\right] < +\infty
\right\},
\end{equation}
\begin{equation}
\label{f11_8_p2_1}
\Phi \psi =
\lim_{ \lambda\in\Pi_{\lambda_0}^\varepsilon, \lambda\to\infty }
F(\lambda) \psi,\qquad \psi\in D(\Phi),
\end{equation}
where $0 < \varepsilon < \frac{\pi}{2}$.

Let us check that the operator $\Phi$ is admissible with respect to $A$.
Consider an arbitrary element $\psi\in D(\Phi)\cap D(X_{\lambda_0}(A))$ such that
$$ \Phi \psi = X_{\lambda_0}(A) \psi. $$
Using~(\ref{f11_8_p2_1}) we get
$$ \lim_{ \lambda\in\Pi_{\lambda_0}^\varepsilon, \lambda\to\infty }
F(\lambda) \psi = X_{\lambda_0}(A) \psi. $$
From~(\ref{f11_7_p2_1}) it follows that
$$ \underline{\lim}_{ \lambda\in\Pi_{\lambda_0}^\varepsilon, \lambda\to\infty }
\left[
|\lambda| ( \| \psi \|_{H} - \| F(\lambda) \psi \|_{H} )
\right] < +\infty. $$
Since the function $F(\lambda)$ is $\lambda_0$-admissible with respect to $A$,
then $\psi = 0$.
Therefore $\Phi$ is admissible with respect to $A$.
Moreover, $V_{22} = T_{22}$ is admissible with respect to $A_e$. By Theorem~\ref{t5_1_p2_1}
the operator $V$ is admissible with respect to $\mathcal{A} = A\oplus A_e$.
By the generalized Neumann formulas for the operator $V$ there corresponds a self-adjoint extension
$\widetilde A := \mathcal{A}_V$ оператора $A$.
Denote by
$$ \mathbf{R}_\lambda = P^{\widetilde H}_H ( \widetilde A - \lambda E_{\widetilde H} )^{-1},\qquad
\lambda\in \mathbb{R}_e, $$
the generalized resolvent of $A$, which corresponds to the self-adjoint extension $\widetilde A$.

From~(\ref{f11_3_1_p2_1}),(\ref{f11_4_0_p2_1}) and~(\ref{f11_4_4_p2_1}) it follows that
$$ F(\lambda) = W F_1(\lambda) U^{-1} =
W
(T_{11} + T_{12} (E_{\mathcal{N}_{\lambda_0}(A_e)} - C_e(\lambda) T_{22})^{-1} C_e(\lambda) T_{21})
U^{-1} $$
$$ = V_{11} + V_{12} (E_{\mathcal{N}_{\lambda_0}(A_e)} - C_e(\lambda) V_{22})^{-1} C_e(\lambda) V_{21},\qquad
\lambda\in\Pi_{\lambda_0}. $$
By~(\ref{f8_15_p2_1}) we conclude that
$$ F(\lambda) = \mathfrak{F}(\lambda;\lambda_0,A,\widetilde A),\qquad
\lambda\in\Pi_{\lambda_0}. $$
From~(\ref{f6_12_p2_1}) and~(\ref{f6_15_p2_1}) it follows~(\ref{f11_3_p2_1})
for the constructed generalized resolvent $\mathbf{R}_\lambda$ and
for the given function $F(\lambda)$.

Let us check the last assertion of the theorem.
Suppose that two functions $F_1(\lambda),F_1(\lambda)\in
\mathcal{S}_{a;\lambda_0}(\Pi_{\lambda_0}; \mathcal{N}_{\lambda_0},
\mathcal{N}_{\overline{\lambda_0}})$ generate by~(\ref{f11_3_p2_1})
the same generalized resolvent of $A$. Then
$$ A_{F_1(\lambda)} = A_{F_2(\lambda)},\qquad \lambda\in\Pi_{\lambda_0}. $$
By the generalized Neumann formulas we obtain that
$$ F_1(\lambda) = F_2(\lambda),\qquad \lambda\in\Pi_{\lambda_0}. $$

$\Box$


\section{Generalized resolvents of isometric and symmetric operators with gaps in their spectrum.}\label{chapter3}

\subsection{Spectral functions of an isometric operator having a constant value on an arc of the circle.}\label{section3_1}

In the investigation of interpolation problems we shall use spectral functions of the operator related to a problem.
There will appear problems with spectral functions, which are constant on the prescribed
arcs of the unit circle. Such spectral functions
and the corresponding generalized resolvents will be studied in this subsection.

\begin{prop}
\label{p3_1_p3_1}
Let $V$ be a closed isometric operator in a Hilbert space $H$, and $\mathbf{F}(\delta)$,
$\delta\in \mathfrak{B}(\mathbb{T})$, be its spectral measure.
The following two conditions are equivalent:

\begin{itemize}
\item[(i)]  $\mathbf{F}(\Delta) = 0$, for an open arc $\Delta$ of the unit circle $\mathbb{T}$;
\item[(ii)] The generalized resolvent $\mathbf{R}_z(V)$, corresponding to the spectral measure
$\mathbf{F}(\delta)$, admits analytic continuation on the set
$\mathbb{D}\cup\mathbb{D}_e\cup \overline{\Delta}$,
where $\overline{\Delta} = \{ z\in \mathbb{C}:\ \overline{z}\in \Delta \}$, for an open открытой arc
$\Delta$ of $\mathbb{T}$.
\end{itemize}
\end{prop}
{\bf Proof. }
(i)$\Rightarrow$(ii).
In this case relation~(\ref{f1_1_p1_1}) 
takes the following form:
\begin{equation}
\label{ff3_1}
(\mathbf R_z h,g)_H = \int_{\mathbb{T}\backslash\Delta} \frac{1}{1-z\zeta} d(\mathbf{F}(\cdot) h,g)_H,\quad
\forall h,g\in H.
\end{equation}
Choose an arbitrary number $z_0\in \overline{\Delta}$. Since the function $\frac{1}{1-z_0\zeta}$
is continuous and bounded on $\mathbb{T}\backslash\Delta$, then there exists an integral
$$ I_{z_0}(h,g) := \int_{\mathbb{T}\backslash\Delta} \frac{1}{1-z_0\zeta} d(\mathbf{F}(\cdot) h,g)_H. $$
Тогда
$$ | (\mathbf R_z h,h)_H - I_{z_0}(h,h) | =
|z-z_0| \left|
\int_{\mathbb{T}\backslash\Delta} \frac{\zeta}{(1-z\zeta)(1-z_0\zeta)} d(\mathbf{F}(\cdot) h,h)_H
\right|
$$
$$ \leq
|z-z_0| \int_{\mathbb{T}\backslash\Delta} \frac{|\zeta|}{|1-z\zeta||1-z_0\zeta|} d(\mathbf{F}(\cdot) h,h)_H,\quad
z\in \mathbb{T}_e. $$
There exists a neighborhood $U(z_0)$ of $z_0$ such that $|z-\overline{\zeta}|\geq M_1 > 0$, $\forall\zeta
\in \mathbb{T}\backslash\Delta$, $\forall z\in U(z_0)$.
Thus, the integral in the last relation is bounded  in the neighborhood $U(z_0)$.
Therefore we get
$$ (\mathbf R_z h,h)_H \rightarrow I_{z_0}(h,h),\quad z\in\mathbb{T}_e,\ z\rightarrow z_0,\quad \forall h\in H. $$
Using the properties of sesquilinear forms we get
$$ (\mathbf R_z h,g)_H \rightarrow I_{z_0}(h,g),\quad z\in\mathbb{T}_e,\ z\rightarrow z_0,\quad \forall h,g\in H. $$
Set
$$ \mathbf R_{\widetilde z} :=
w.-\lim_{z\in \mathbb{T}_e,\ z\to \widetilde z} \mathbf R_z,\quad \forall \widetilde z\in\overline{\Delta}, $$
where the limit is understood in a sense of the weak operator topology.
We may write
$$ \left(
\frac{1}{z-z_0}(\mathbf R_z -  \mathbf R_{z_0})h,h
\right)_H =
\int_{\mathbb{T}\backslash\Delta} \frac{\zeta}{(1-z\zeta)(1-z_0\zeta)} d(\mathbf{F}(\cdot) h,h)_H, $$
$$ z\in U(z_0),\ h\in H. $$
The function under the integral sign is bounded in $U(z_0)$, and it tends to
$\frac{\zeta}{(1-z_0\zeta)^2}$. By the Lebesgue theorem on a limit we obtain that
$$ \lim_{z\to z_0} \left(
\frac{1}{z-z_0}(\mathbf R_z -  \mathbf R_{z_0})h,h
\right)_H =
\int_{\mathbb{T}\backslash\Delta} \frac{\zeta}{(1-z_0\zeta)^2} d(\mathbf{F}(\cdot) h,h)_H; $$
and therefore
$$ \lim_{z\to z_0} \left(
\frac{1}{z-z_0}(\mathbf R_z -  \mathbf R_{z_0})h,g
\right)_H =
\int_{\mathbb{T}\backslash\Delta} \frac{\zeta}{(1-z_0\zeta)^2} d(\mathbf{F}(\cdot) h,g)_H, $$
for $h,g\in H$.
Thus, there exists the derivative of the function $\mathbf R_z$ at $z=z_0$.

\noindent
(ii)$\Rightarrow$(i).
Choose an arbitrary element $h\in H$, and consider the following function $\sigma_h(t) := (\mathbf{F}_t h,h)_H$,
$t\in [0,2\pi)$, where $\mathbf{F}_t$ is a left-continuous spectral function of $V$,
corresponding to a spectral measure $\mathbf{F}(\delta)$. Consider the following function:
$$ f_h(z) = \frac{1}{2} \int_0^{2\pi} \frac{1+e^{it} z}{1-e^{it} z} d\sigma_h(t)
= \int_0^{2\pi} \frac{1}{1-e^{it} z} d\sigma_h(t) - \frac{1}{2}
\int_0^{2\pi} d\sigma_h(t) $$
\begin{equation}
\label{ff3_1_1}
= \int_0^{2\pi} \frac{1}{1-e^{it} z} d\sigma_h(t) - \frac{1}{2} \| h \|^2_H
= (\mathbf{R}_z h,h)_H - \frac{1}{2} \| h \|^2_H.
\end{equation}
Choose arbitrary numbers $t_1,t_2$, $0\leq t_1<t_2\leq 2\pi$, such that
\begin{equation}
\label{ff3_2}
l=l(t_1,t_2)= \{ z=e^{it}:\ t_1\leq t\leq t_2 \} \subset \Delta.
\end{equation}
We assume that $t_1$ and $t_2$ are points of the continuity of the function $\mathbf{F}_t$.
By the inversion formula we may write:
$$ \sigma_h(t_2) - \sigma_h(t_1) = \lim_{r\to 1-0} \int_{t_1}^{t_2} {\mathrm Re}
\left\{ f_h(re^{-i\tau}) \right\} d\tau. $$
Observe that
$$ {\mathrm Re} \left\{ f_h(re^{-i\tau}) \right\} =
{\mathrm Re} \left\{ (\mathbf{R}_{re^{-i\tau}} h,h)_H \right\} - \frac{1}{2} \| h \|^2_H $$
\begin{equation}
\label{ff3_2_1}
= \frac{1}{2} \left(
((\mathbf{R}_{re^{-i\tau}} + \mathbf{R}_{re^{-i\tau}}^*) h,h)_H
\right) - \frac{1}{2} \| h \|^2_H,\quad t_1\leq \tau \leq t_2.
\end{equation}
By~(\ref{ff3_2}) we obtain that $e^{-i\tau}$ belongs to $\overline\Delta$, for $t_1\leq \tau \leq t_2$.
Therefore
\begin{equation}
\label{ff3_2_2}
\lim_{r\to 1-0} ((\mathbf{R}_{re^{-i\tau}} + \mathbf{R}_{re^{-i\tau}}^*) h,h)_H =
((\mathbf{R}_{e^{-i\tau}} + \mathbf{R}_{e^{-i\tau}}^*) h,h)_H.
\end{equation}
By Theorem~\ref{t1_1_p1_1} 
the generalized resolvents of an isometric operator have the following property:
\begin{equation}
\label{ff3_3}
\mathbf{R}_z^* = E_H - \mathbf{R}_{\frac{1}{ \overline{z} }},\quad z\in \mathbb{T}_e.
\end{equation}
Passing to a limit in~(\ref{ff3_3}) as $z$ tends to $e^{-i\tau}$, we get
\begin{equation}
\label{ff3_4}
\mathbf{R}_{e^{-i\tau}}^* = E_H - \mathbf{R}_{e^{-i\tau}},\quad t_1\leq \tau \leq t_2.
\end{equation}
By~(\ref{ff3_2_1}),(\ref{ff3_2_2}) and~(\ref{ff3_4}) we get
\begin{equation}
\label{ff3_5}
\lim_{r\to 1-0} {\mathrm Re} \left\{ f_h(re^{-i\tau}) \right\} = 0,\quad t_1\leq \tau \leq t_2.
\end{equation}
Consider the followinf sector:
$$ L(t_1,t_2) = \{ z=re^{-it}:\ t_1\leq t\leq t_2,\ 0\leq r \leq 1 \}. $$
The generalized resolvent is analytic in each point of the closed sector $L(t_1,t_2)$.
Consequently, the function
${\mathrm Re} (\mathbf{R}_z h,h)$ is continuous and bounded in $L(t_1,t_2)$.
By the Lebesgue theorem on a limit we conclude that $\sigma_h(t_1)=\sigma_h(t_2)$.
If $1\notin\Delta$ the required result follows easily.
In the case $1\in\Delta$, we may write $\Delta = \Delta_1\cup\{ 1 \}\cup\Delta_2$,
where the open arcs $\Delta_1$ and $\Delta_2$ do not contain $1$. Therefore $\sigma_h(t)$
is constant in intervals, corresponding to $\Delta_1$ and $\Delta_2$.
Suppose that there exists a non-zero jump of $\sigma_h(t)$ at $t=0$.
By~(\ref{f1_1_p1_1}) 
we may write:
$$ (\mathbf{R}_z h,h)_H  =
\int_0^{2\pi} \frac{1}{1-e^{it} z} d\sigma_h(t) =
\int_0^{2\pi} \frac{1}{1-e^{it} z} d\widehat\sigma_h(t) + \frac{1}{1-z} a,\quad a>0, $$
where $\widehat\sigma_h(t) = \sigma_h(t) + \sigma_h(+0) - \sigma_h(0)$, $t\in [0,2\pi]$.
In a neighborhood of $1$
the left-hand side and the first summand in the right-hand side are bounded.
We obtained a contradiction.
$\Box$

\begin{thm}
\label{t3_1_p3_1}
Let $V$ be a closed isometric operator in a Hilbert space $H$, and $\mathbf{R}_z(V)$ be an arbitrary
generalized resolvent of the operator $V$.
Let $\{ \lambda_k \}_{k=1}^\infty$ be a sequence of numbers from $\mathbb{D}$ such that
$\lambda_k\rightarrow \widehat\lambda$, as $k\rightarrow\infty$; $\widehat\lambda\in \mathbb{T}$.
Suppose that for a number $z_0\in \mathbb{D}\backslash\{ 0 \}$,
the function $C(\lambda;z_0)$, corresponding to $\mathbf{R}_z(V)$ in Inin's formula,
satisfies the following relation:
\begin{equation}
\label{ff3_6}
\exists u.-\lim_{k\to\infty} C(\lambda_k;z_0) =: C(\widehat\lambda;z_0).
\end{equation}
Then for an arbitrary $z_0'\in \mathbb{D}\backslash\{ 0 \}$, the function $C(\lambda;z_0')$, corresponding to
the generalized resolvent
$\mathbf{R}_z(V)$ by Inin's formula~(\ref{f1_7}) 
satisfies the following relation:
\begin{equation}
\label{ff3_7}
\exists u.-\lim_{k\to\infty} C(\lambda_k;z_0') =: C(\widehat\lambda;z_0').
\end{equation}
Moreover, $C(\widehat\lambda;z_0')$ is a linear non-expanding operator, which maps the whole
$N_{z_0'}$ into
$N_{\frac{1}{ \overline{z_0'} }}$, and the corresponding orthogonal extension
$V_{C(\widehat\lambda;z_0');z_0'}$ does not depend on the choice of the point $z_0'\in \mathbb{D}\backslash\{ 0 \}$.
\end{thm}
{\bf Proof. }
Suppose that relation~(\ref{ff3_6}) holds for a point $z_0\in \mathbb{D}\backslash\{ 0 \}$.
Choose an arbitrary point $z_0'\in \mathbb{D}\backslash\{ 0 \}$.
Comparing the Inin formula for $z_0$ and for $z_0'$, we see that
\begin{equation}
\label{ff3_8}
V_{C(\lambda;z_0);z_0} = V_{C(\lambda;z_0');z_0'},\qquad \lambda\in \mathbb{D}.
\end{equation}
By~(\ref{f1_6_1}) 
we may write
\begin{equation}
\label{ff3_9}
V_{C(\lambda;z_0);z_0} = \frac{1}{z_0} E_H +
\frac{|z_0|^2-1}{z_0}
\left(
E_H + z_0 V^+_{z_0;C(\lambda;z_0)}
\right)^{-1},\qquad \lambda\in \mathbb{D}.
\end{equation}
Substituting in~(\ref{ff3_8}) similar relations for $z_0$ and $z_0'$, and multiplying by $z_0 z_0'$
we get
$$ z_0' E_H + z_0'(|z_0|^2-1)
\left(
E_H + z_0 V^+_{z_0;C(\lambda;z_0)}
\right)^{-1} $$
$$ =
z_0 E_H + z_0(|z_0'|^2-1)
\left(
E_H + z_0' V^+_{z_0';C(\lambda;z_0')}
\right)^{-1}. $$
Then
$$ \left(
E_H + z_0' V^+_{z_0';C(\lambda;z_0')}
\right)^{-1} = \frac{1}{z_0(|z_0'|^2-1)}
\left(
(z_0'-z_0) E_H \right.
$$
$$ \left.
+
z_0'(|z_0|^2-1)
\left(
E_H + z_0 V^+_{z_0;C(\lambda;z_0)}
\right)^{-1}
\right) $$
\begin{equation}
\label{ff3_10}
=
\frac{1-z_0'\overline{z_0}}{1-|z_0'|^2}
\left(
E_H + \frac{z_0-z_0'}{1-z_0'\overline{z_0}} V^+_{z_0;C(\lambda;z_0)}
\right)
\left(
E_H + z_0 V^+_{z_0;C(\lambda;z_0)}
\right)^{-1},\ \lambda\in \mathbb{D}.
\end{equation}
\begin{lem}
\label{l3_1_p3_1}
Let $z_0,z_0'\in \mathbb{D}$. Then
\begin{equation}
\label{ff3_11}
\left|
\frac{z_0-z_0'}{1-z_0'\overline{z_0}}
\right|<1.
\end{equation}
\end{lem}
{\bf Proof. }
Consider a linear fractional transformation: $w=w(u)=\frac{z_0-u}{1-\overline{z_0}u}$.
If $|u|=1$, then  $|1-\overline{z_0}u|=|u(\overline{u}-\overline{z_0})|=|u-z_0|$.
Moreover, we have $w(z_0)=0$. Therefore $w$ maps $\mathbb{D}$ on $\mathbb{D}$.
$\Box$

Using~(\ref{ff3_10}) and (\ref{ff3_11}) we may write:
$$ E_H + z_0' V^+_{z_0';C(\lambda;z_0')} $$
\begin{equation}
\label{ff3_12}
=
\frac{1-|z_0'|^2}{1-z_0'\overline{z_0}}
\left( E_H + z_0 V^+_{z_0;C(\lambda;z_0)} \right)
\left(
E_H + \frac{z_0-z_0'}{1-z_0'\overline{z_0}} V^+_{z_0;C(\lambda;z_0)}
\right)^{-1},\quad \lambda\in \mathbb{D}.
\end{equation}
Using~(\ref{f1_5}) 
we write:
$$ V^+_{z_0;C(\lambda;z_0)} = V_{z_0} \oplus C(\lambda;z_0),\quad \lambda\in \mathbb{D}. $$
By the statement of the theorem we easily obtain that $C(\widehat\lambda;z_0)$ is a non-expanding operator and
\begin{equation}
\label{ff3_13}
\exists u.-\lim_{k\to\infty} V^+_{z_0;C(\lambda_k;z_0)} = V_{z_0} \oplus C(\widehat\lambda;z_0)
= V^+_{z_0;C(\widehat\lambda;z_0)}.
\end{equation}
We may write
$$ \left\|
\left(
E_H + \frac{z_0-z_0'}{1-z_0'\overline{z_0}} V^+_{z_0;C(\lambda_k;z_0)}
\right)^{-1}
-
\left(
E_H + \frac{z_0-z_0'}{1-z_0'\overline{z_0}} V^+_{z_0;C(\widehat\lambda;z_0)}
\right)^{-1}
\right\|
$$
$$ \leq
\left\|
\left(
E_H + \frac{z_0-z_0'}{1-z_0'\overline{z_0}} V^+_{z_0;C(\lambda_k;z_0)}
\right)^{-1}
\right\| $$
\begin{equation}
\label{ff3_14}
* \left|
\frac{z_0-z_0'}{1-z_0'\overline{z_0}} \right|
\left\|
V^+_{z_0;C(\lambda_k;z_0)} - V^+_{z_0;C(\widehat\lambda;z_0)}
\right\|
\left\|
\left(
E_H + \frac{z_0-z_0'}{1-z_0'\overline{z_0}} V^+_{z_0;C(\widehat\lambda;z_0)}
\right)^{-1}
\right\|.
\end{equation}
Since
$$ \left| \frac{z_0-z_0'}{1-z_0'\overline{z_0}} \right|
\left\| V^+_{z_0;C(\lambda_k;z_0)} \right\| \leq \delta < 1, $$
then
$$ \left\|
\left(
E_H + \frac{z_0-z_0'}{1-z_0'\overline{z_0}} V^+_{z_0;C(\lambda_k;z_0)}
\right) h
\right\|
\geq \left| \| h \| - \left\|
\frac{z_0-z_0'}{1-z_0'\overline{z_0}} V^+_{z_0;C(\lambda_k;z_0)}
\right|
\right| $$
$$ \geq (1-\delta)\| h \|; $$
$$ \left\|
\left(
E_H + \frac{z_0-z_0'}{1-z_0'\overline{z_0}} V^+_{z_0;C(\lambda_k;z_0)}
\right)^{-1} \right\| \leq \frac{1}{1-\delta}. $$
Passing to a limit in~(\ref{ff3_14}) we get
\begin{equation}
\label{ff3_15}
u.-\lim_{k\to\infty} \left( E_H + \frac{z_0-z_0'}{1-z_0'\overline{z_0}} V^+_{z_0;C(\lambda_k;z_0)}
\right)^{-1} =
\left( E_H + \frac{z_0-z_0'}{1-z_0'\overline{z_0}} V^+_{z_0;C(\widehat\lambda;z_0)}
\right)^{-1}.
\end{equation}
By~(\ref{ff3_15}),(\ref{ff3_13}),(\ref{ff3_12}) we obtain that there exists a limit
\begin{equation}
\label{ff3_16}
u.-\lim_{k\to\infty} V^+_{z_0';C(\lambda_k;z_0')}
= u.-\lim_{k\to\infty} V_{z_0'} \oplus C(\lambda_k;z_0')
=: V',
\end{equation}
such that
$$ E_H + z_0' V' $$
\begin{equation}
\label{ff3_17}
=
\frac{1-|z_0'|^2}{1-z_0'\overline{z_0}}
\left( E_H + z_0 V^+_{z_0;C(\widehat\lambda;z_0)} \right)
\left(
E_H + \frac{z_0-z_0'}{1-z_0'\overline{z_0}} V^+_{z_0;C(\widehat\lambda;z_0)}
\right)^{-1}.
\end{equation}
Relation~(\ref{ff3_16}) shows that $V'|_{M_{z_0'}} = V_{z_0'}$. Set
$$ C(\widehat\lambda;z_0') = V'|_{N_{z_0'}}. $$
By~(\ref{ff3_16}) we obtain that $C(\widehat\lambda;z_0')$ is a linear non-expanding
operator, which maps $N_{z_0'}$ in $N_{\frac{1}{ \overline{z_0'} }}$. Thus,
we get
\begin{equation}
\label{ff3_18}
V' = V_{z_0'}\oplus C(\widehat\lambda;z_0') = V^+_{z_0';C(\widehat\lambda;z_0')}.
\end{equation}
From~(\ref{ff3_16}),(\ref{ff3_18}) it easily follows that
\begin{equation}
\label{ff3_18_1}
u.-\lim_{k\to\infty} C(\lambda_k;z_0')
= C(\widehat\lambda;z_0'),
\end{equation}
and relation~(\ref{ff3_7}) is proved.

By~(\ref{ff3_17}),(\ref{ff3_18}) we get
$$ \left( E_H + z_0' V^+_{z_0';C(\widehat\lambda;z_0')} \right)^{-1} $$
\begin{equation}
\label{ff3_19}
= \frac{1-z_0'\overline{z_0}}{1-|z_0'|^2}
\left(
E_H + \frac{z_0-z_0'}{1-z_0'\overline{z_0}} V^+_{z_0;C(\widehat\lambda;z_0)}
\right)
\left( E_H + z_0 V^+_{z_0;C(\widehat\lambda;z_0)} \right)^{-1};
\end{equation}
$$ (1-|z_0'|^2) \left( E_H + z_0' V^+_{z_0';C(\widehat\lambda;z_0')} \right)^{-1} $$
$$ = \left(
(1-z_0'\overline{z_0}) E_H + (z_0-z_0') V^+_{z_0;C(\widehat\lambda;z_0)}
\right)
\left( E_H + z_0 V^+_{z_0;C(\widehat\lambda;z_0)} \right)^{-1};
$$
Subtracting $E_H$ from the both sides of the last relation and dividing by $-z_0'$ we get:
$$ \frac{1}{z_0'}E_H + \frac{|z_0'|^2-1}{z_0'} \left( E_H + z_0' V^+_{z_0';C(\widehat\lambda;z_0')} \right)^{-1} $$
$$ = -\frac{1}{z_0'}\left(
(1-z_0'\overline{z_0}) E_H + (z_0-z_0') V^+_{z_0;C(\widehat\lambda;z_0)}
- ( E_H + z_0 V^+_{z_0;C(\widehat\lambda;z_0)} )
\right) $$
$$*
\left( E_H + z_0 V^+_{z_0;C(\widehat\lambda;z_0)} \right)^{-1}
$$
$$ = \left( \overline{z_0} E_H + V^+_{z_0;C(\widehat\lambda;z_0)} \right)
\left( E_H + z_0 V^+_{z_0;C(\widehat\lambda;z_0)} \right)^{-1}. $$
From~(\ref{f1_6}) 
and~(\ref{f1_6_1}) 
we get
\begin{equation}
\label{ff3_20}
V_{C(\widehat\lambda;z_0');z_0'} = V_{C(\widehat\lambda;z_0);z_0},\qquad \forall z_0'\in \mathbb{D}.
\end{equation}
$\Box$

\begin{thm}
\label{t3_2_p3_1}
Let $V$ be a closed isometric operator in a Hilbert space $H$, and
$z_0\in \mathbb{D}\backslash\{ 0 \}$ be a fixed point.
Let $\mathbf{R}_z = \mathbf{R}_z(V)$ be an arbitrary generalized resolvent of the operator $V$,
and $C(\lambda;z_0)\in \mathcal{S}(N_{z_0};N_{\frac{1}{ \overline{z_0} }})$
corresponds to $\mathbf{R}_z(V)$ by Inin's formula~(\ref{f1_7}). 
The generalized resolvent $\mathbf{R}_z(V)$ admits an analytic continuation on a set
$\mathbb{D}\cup\mathbb{D}_e\cup \Delta$,
where $\Delta$ is an open arc of $\mathbb{T}$,
if and only if the following conditions hold:

\begin{itemize}
\item[1)] The function $C(\lambda;z_0)$ admits a continuation on a set $\mathbb{D}\cup\Delta$,
which is continuous in the uniform operator topology;

\item[2)] The continued function $C(\lambda;z_0)$ maps isometriclly $N_{z_0}$ on the whole
$N_{ \frac{1}{ \overline{z_0} } }$, for all points $\lambda\in\Delta$;

\item[3)] The operator $(E_H - \lambda V_{C(\lambda;z_0);z_0})^{-1}$ exists and it is defined on
the whole $H$, for all points $\lambda\in\Delta$.
\end{itemize}
\end{thm}
{\bf Proof. }{\it Necessity.}
Choose an arbitrary point $\widehat\lambda\in\Delta$.
Let $z\in \mathbb{D}\backslash\{ 0 \}$ be an arbitrary point, and
$C(\lambda;z)\in \mathcal{S}(N_{z};N_{\frac{1}{ \overline{z} }})$
correspond to the generalized resolvent $\mathbf{R}_z(V)$ by the Inin formula.
Using the Inin formula we may write:
$$ E_H - z V_{C(\lambda;z);z} = \frac{z}{\lambda} (E_H - \lambda V_{C(\lambda;z);z}) +
\left( 1-\frac{z}{\lambda}\right) E_H $$
$$ = \frac{z}{\lambda} \mathbf{R}_\lambda^{-1}
+ \left( 1-\frac{z}{\lambda}\right) E_H =
\frac{z}{\lambda} \left[
E_H + \left( \frac{\lambda}{z} - 1 \right) \mathbf{R}_\lambda
\right] \mathbf{R}_\lambda^{-1},\ \forall\lambda\in \mathbb{D}\backslash\{ 0 \}.
$$
Therefore the following operator
$$ \left[
E_H + \left( \frac{\lambda}{z} - 1 \right) \mathbf{R}_\lambda
\right]
=
\frac{\lambda}{z} ( E_H - z V_{C(\lambda;z);z} ) \mathbf{R}_\lambda, $$
has a bounded inverse, defined on the whole $H$:
\begin{equation}
\label{ff3_20_1}
\left[
E_H + \left( \frac{\lambda}{z} - 1 \right) \mathbf{R}_\lambda
\right]^{-1}
=
\frac{z}{\lambda} \mathbf{R}_\lambda^{-1} ( E_H - z V_{C(\lambda;z);z} )^{-1},\quad
\lambda \in \mathbb{D}\backslash\{ 0 \}.
\end{equation}
Then
\begin{equation}
\label{ff3_21}
(E_H - z V_{C(\lambda;z);z})^{-1} = \frac{\lambda}{z}
\mathbf{R}_\lambda
\left[
E_H + \left( \frac{\lambda}{z} - 1 \right) \mathbf{R}_\lambda
\right]^{-1},\quad \lambda \in \mathbb{D}\backslash\{ 0 \}.
\end{equation}
Choose an arbitrary $\delta$: $0<\delta<1$. Assume that $z\in \mathbb{D}\backslash\{ 0 \}$
satisfies the following additional condition:
\begin{equation}
\label{ff3_22}
\left| \frac{\widehat\lambda}{z} - 1 \right| \| \mathbf{R}_{\widehat\lambda} \| < \delta.
\end{equation}
Let us check that such points exist. If $\| \mathbf{R}_{\widehat\lambda} \|=0$, then it is obvious.
In the opposite case we seek $z$ in the following form: $z=\varepsilon\widehat\lambda$,
$0<\varepsilon <1$. Then condition~(\ref{ff3_22}) will take the following form:
$$ \left| \frac{1}{\varepsilon} - 1 \right| <
\frac{\delta}{ \| \mathbf{R}_{\widehat\lambda} \| }; $$
or $\varepsilon > \frac{1}{ 1 +
\frac{\delta}{ \| \mathbf{R}_{\widehat\lambda} \| } }$.

\noindent
In this case, there exists an inverse $\left[
E_H + \left( \frac{\widehat\lambda}{z} - 1 \right) \mathbf{R}_{\widehat\lambda}
\right]^{-1}$, which is bounded and defined on the whole $H$.
Moreover, by the continuity the following inequality holds:
\begin{equation}
\label{ff3_22_1}
\left| \frac{\lambda}{z} - 1 \right| \| \mathbf{R}_{\lambda} \| < \delta,
\end{equation}
in an open neighborhood $U(\widehat\lambda)$
of $\widehat\lambda$. Therefore there exists the inverse
\begin{equation}
\label{ff3_23}
\left[
E_H + \left( \frac{\lambda}{z} - 1 \right) \mathbf{R}_{\lambda}
\right]^{-1},\quad \forall\lambda\in U(\widehat\lambda),
\end{equation}
which is bounded and defined on the whole $H$.
We may write
$$ \left\|
\left( E_H + \left( \frac{\lambda}{z} - 1 \right) \mathbf{R}_{\lambda} \right) h
\right\| \geq
\left|
\| h \| - \left| \frac{\lambda}{z} - 1 \right|
\| \mathbf{R}_{\lambda} h \|
\right| $$
$$ \geq (1-\delta) \| h \|,\qquad h\in H. $$
Therefore
\begin{equation}
\label{ff3_24}
\left\|
\left[
E_H + \left( \frac{\lambda}{z} - 1 \right) \mathbf{R}_{\lambda}
\right]^{-1}
\right\| \leq \frac{1}{1-\delta},\quad
\lambda\in U(\widehat\lambda).
\end{equation}
Choose an arbitrary sequence $\{ \lambda_k \}_{k=1}^\infty$
of points from $\mathbb{D}$ such that $\lambda_k\rightarrow\widehat\lambda$, as $k\rightarrow\infty$.
There exists a number $k_0\in \mathbb{N}$ such that $\lambda_k\in U(\widehat\lambda)\cap \mathbb{D}$,
$k\geq k_0$.
We write
$$ \left\|
\left[
E_H + \left( \frac{\lambda_k}{z} - 1 \right) \mathbf{R}_{\lambda_k}
\right]^{-1}
-
\left[
E_H + \left( \frac{\widehat\lambda}{z} - 1 \right) \mathbf{R}_{\widehat\lambda}
\right]^{-1} \right\| $$
$$ \leq
\left \| \left[
E_H + \left( \frac{\lambda_k}{z} - 1 \right) \mathbf{R}_{\lambda_k}
\right]^{-1}
\right\|
\left\|
\left( \frac{\lambda_k}{z} - 1 \right) \mathbf{R}_{\lambda_k} -
\left( \frac{\widehat\lambda}{z} - 1 \right) \mathbf{R}_{\widehat\lambda}
\right\| $$
$$ *
\left\|
\left[
E_H + \left( \frac{\widehat\lambda}{z} - 1 \right) \mathbf{R}_{\widehat\lambda}
\right]^{-1}
\right\|. $$
The first factor on the right in the last equality is uniformly bounded by~(\ref{ff3_24}).
Thus, we get
\begin{equation}
\label{ff3_25}
u.-\lim_{\lambda\in \mathbb{D},\ \lambda\to\widehat\lambda}
\left[
E_H + \left( \frac{\lambda}{z} - 1 \right) \mathbf{R}_{\lambda}
\right]^{-1} =
\left[
E_H + \left( \frac{\widehat\lambda}{z} - 1 \right) \mathbf{R}_{\widehat\lambda}
\right]^{-1}.
\end{equation}
From relations~(\ref{ff3_21}),(\ref{ff3_25}) it follows that
\begin{equation}
\label{ff3_26}
u.-\lim_{\lambda\in \mathbb{D},\ \lambda\to\widehat\lambda}
(E_H - z V_{C(\lambda;z);z})^{-1} = \frac{\widehat\lambda}{z}
\mathbf{R}_{\widehat\lambda}
\left[
E_H + \left( \frac{\widehat\lambda}{z} - 1 \right) \mathbf{R}_{\widehat\lambda}
\right]^{-1}.
\end{equation}
Thus, there exists the following limit:
$$ u.-\lim_{\lambda\in \mathbb{D},\ \lambda\to\widehat\lambda}
V^+_{z;C(\lambda;z)} =
u.-\lim_{\lambda\in \mathbb{D},\ \lambda\to\widehat\lambda}
\left(
-\frac{1}{z} E_H + \frac{1-|z|^2}{z} (E_{H} - z V_{C(\lambda;z);z})^{-1}
\right) $$
\begin{equation}
\label{ff3_27}
=
-\frac{1}{z} E_H +
\frac{1-|z|^2}{z^2}
\widehat\lambda
\mathbf{R}_{\widehat\lambda}
\left[
E_H + \left( \frac{\widehat\lambda}{z} - 1 \right) \mathbf{R}_{\widehat\lambda}
\right]^{-1}
=: V_z',
\end{equation}
where we used~(\ref{f1_6_2}). 
Observe that
$V^+_{z;C(\lambda;z)} = V_z\oplus C(\lambda;z)$.
Set
$C(\widehat\lambda;z) = V_z'|_{N_z}$. By~(\ref{ff3_27}) we conclude that
$C(\widehat\lambda;z)$ is a linear non-expanding operator which maps
$N_z$ into $N_{\frac{1}{ \overline{z} }}$. Moreover, we have $V_z'|_{M_z} = V_z$, and therefore
$$ V_z' = V_z\oplus C(\widehat\lambda;z) = V^+_{z;C(\widehat\lambda;z)}. $$
Using~(\ref{ff3_27}) we easily obtain that
\begin{equation}
\label{ff3_28}
u.-\lim_{\lambda\in \mathbb{D},\ \lambda\to\widehat\lambda}
C(\lambda;z) = C(\widehat\lambda;z).
\end{equation}
By~\ref{t3_1_p3_1} we conclude that the last relation holds for $z_0$, where
$C(\widehat\lambda;z_0)$ is a linear non-expanding operator which maps $N_{z_0}$
into
$N_{\frac{1}{ \overline{z_0} }}$, and
$V_{C(\widehat\lambda;z_0);z_0} = V_{C(\widehat\lambda;z);z}$.

\noindent
Comparing relation~(\ref{ff3_27}) for $V_z' = V^+_{z;C(\widehat\lambda;z)}$ and
relation~(\ref{f1_6_2}) 
we get
\begin{equation}
\label{ff3_28_1}
\frac{\widehat\lambda}{z}
\mathbf{R}_{\widehat\lambda}
\left[
E_H + \left( \frac{\widehat\lambda}{z} - 1 \right) \mathbf{R}_{\widehat\lambda}
\right]^{-1} = \left( E_H - z V_{C(\widehat\lambda;z);z} \right)^{-1},
\end{equation}
for the prescribed choice of $z$.

Thus, we continued by the continuity the function $C(\lambda;z_0)$ on a set
$\mathbb{D}\cup\Delta$. Let us check that this continuation is
continuous in the uniform operator topology.
It remains to check that for an arbitrary $\widehat\lambda\in\Delta$ it holds
\begin{equation}
\label{ff3_29}
u.-\lim_{\lambda\in \mathbb{D}\cup\Delta,\ \lambda\to\widehat\lambda}
C(\lambda;z_0) = C(\widehat\lambda;z_0).
\end{equation}
Choose a number $z\in \mathbb{D}\backslash\{ 0 \}$, which satisfies~(\ref{ff3_22}) and construct a neighborhood
$U(\widehat\lambda)$ as it was done before.
For an arbitrary $\lambda\in U(\widehat\lambda)$ we may write:
$$ \left\|
\left[
E_H + \left( \frac{\lambda}{z} - 1 \right) \mathbf{R}_{\lambda}
\right]^{-1}
-
\left[
E_H + \left( \frac{\widehat\lambda}{z} - 1 \right) \mathbf{R}_{\widehat\lambda}
\right]^{-1}
\right\| $$
$$ \leq
\left\|
\left[
E_H + \left( \frac{\lambda}{z} - 1 \right) \mathbf{R}_{\lambda}
\right]^{-1}
\right\|
\left\|
\left( \frac{\widehat\lambda}{z} - 1 \right) \mathbf{R}_{\widehat\lambda} -
\left( \frac{\lambda}{z} - 1 \right) \mathbf{R}_{\lambda}
\right\| $$
$$ *
\left\|
\left[
E_H + \left( \frac{\widehat\lambda}{z} - 1 \right) \mathbf{R}_{\widehat\lambda}
\right]^{-1}
\right\|. $$
Using~(\ref{ff3_24}) we get
\begin{equation}
\label{ff3_30}
u.-\lim_{\lambda\to\widehat\lambda}
\left[
E_H + \left( \frac{\lambda}{z} - 1 \right) \mathbf{R}_{\lambda}
\right]^{-1}
=
\left[
E_H + \left( \frac{\widehat\lambda}{z} - 1 \right) \mathbf{R}_{\widehat\lambda}
\right]^{-1}.
\end{equation}
By~(\ref{ff3_28_1}) we get
\begin{equation}
\label{ff3_31}
\frac{\lambda}{z}
\mathbf{R}_{\lambda}
\left[
E_H + \left( \frac{\lambda}{z} - 1 \right) \mathbf{R}_{\lambda}
\right]^{-1} = \left( E_H - z V_{C(\lambda;z);z} \right)^{-1},\quad \forall\lambda\in
(U(\widehat\lambda)\cap\overline{\mathbb{D}})\backslash\{ 0 \},
\end{equation}
for the prescribed choice of $z$.
In fact, for an arbitrary $\widetilde\lambda\in\Delta\cap U(\widehat\lambda)$,
there exists a neighborhood $\widetilde U(\widetilde\lambda)\subset U(\widehat\lambda)$, where
inequality~(\ref{ff3_22_1}) is satisfied for the same choice of $z$.
Repeating the arguments below~(\ref{ff3_22_1}) with $\widetilde\lambda$ instead of $\widehat\lambda$,
we obtain that что~(\ref{ff3_31}) holds for $\widetilde\lambda$.
For points inside $\mathbb{D}$ we can apply~(\ref{ff3_21}).

From~(\ref{ff3_30}),(\ref{ff3_31}) it follows that
\begin{equation}
\label{ff3_32}
u.-\lim_{\lambda\in \mathbb{D}\cup\Delta,\ \lambda\to\widehat\lambda}
\left( E_H - z V_{C(\lambda;z);z} \right)^{-1}
=
\left( E_H - z V_{C(\widehat\lambda;z);z} \right)^{-1}.
\end{equation}
Since it was proved that $V_{C(\widehat\lambda;z);z}$ does not depend on the choice of
$z\in \mathbb{D}\backslash\{ 0 \}$ (and for $z\in \mathbb{D}$ this follows from the Inin formula),
then the last relation holds for all $z\in \mathbb{D}\backslash\{ 0 \}$.

Using relation~(\ref{f1_6_2}) 
we get
\begin{equation}
\label{ff3_33}
u.-\lim_{\lambda\in \mathbb{D}\cup\Delta,\ \lambda\to\widehat\lambda}
V^+_{z;C(\lambda;z)}
= V^+_{z;C(\widehat\lambda;z)},\quad \forall z\in \mathbb{D}\backslash\{ 0 \},
\end{equation}
and therefore relation (\ref{ff3_29}) holds.
Thus, condition~1) from the statement of the theorem is satisfied.

From~(\ref{ff3_28_1}) it is seen that
\begin{equation}
\label{ff3_34}
\mathbf{R}_{\widehat\lambda}
=
\frac{z}{\widehat\lambda}
\left( E_H - z V_{C(\widehat\lambda;z);z} \right)^{-1}
\left[
E_H + \left( \frac{\widehat\lambda}{z} - 1 \right) \mathbf{R}_{\widehat\lambda}
\right],
\end{equation}
for the above prescribed choice of $z$. Therefore the operator
$\mathbf{R}_{\widehat\lambda}$ has a bounded inverse, defined on the whole space $H$. Then
$$ E_H = \frac{z}{\widehat\lambda}
\left( E_H - z V_{C(\widehat\lambda;z);z} \right)^{-1}
\left[
E_H + \left( \frac{\widehat\lambda}{z} - 1 \right) \mathbf{R}_{\widehat\lambda}
\right]
\mathbf{R}_{\widehat\lambda}^{-1}; $$
$$ E_H - z V_{C(\widehat\lambda;z);z} =
\frac{z}{\widehat\lambda}
\left[
E_H + \left( \frac{\widehat\lambda}{z} - 1 \right) \mathbf{R}_{\widehat\lambda}
\right]
\mathbf{R}_{\widehat\lambda}^{-1} $$
$$ = \frac{z}{\widehat\lambda}
\mathbf{R}_{\widehat\lambda}^{-1} + \left( 1 - \frac{z}{\widehat\lambda} \right) E_H. $$
From the latter relation we get
\begin{equation}
\label{ff3_35}
\mathbf{R}_{\widehat\lambda} =
\left( E_H - \widehat\lambda V_{C(\widehat\lambda;z);z} \right)^{-1}.
\end{equation}
Since $V_{C(\widehat\lambda;z);z}$ does not depend on the choice of $z$,
then the last relation holds for all $z\in \mathbb{D}\backslash\{ 0 \}$.
Consequently, condition~3) from the statement of the theorem holds.

Using property~(\ref{ff3_3}), and passing to the limit as $z\rightarrow\widehat\lambda$, we get
\begin{equation}
\label{ff3_36}
\mathbf{R}_{\widehat\lambda}^* = E_H - \mathbf{R}_{\widehat\lambda}.
\end{equation}
On the other hand, from~(\ref{ff3_35}) we get:
\begin{equation}
\label{ff3_37}
\mathbf{R}_{\widehat\lambda}^* = \left( E_H - \overline{\widehat\lambda}
V_{C(\widehat\lambda;z);z}^* \right)^{-1}.
\end{equation}
From~(\ref{ff3_35})-(\ref{ff3_37}) it is seen that
\begin{equation}
\label{ff3_38}
E_H = \left( E_H - \overline{\widehat\lambda}
V_{C(\widehat\lambda;z);z}^* \right)^{-1}
+
\left( E_H - \widehat\lambda V_{C(\widehat\lambda;z);z} \right)^{-1},\quad
\forall z\in \mathbb{D}\backslash\{ 0 \}.
\end{equation}
Multiplying the both sides of the last relation by $(E_H - \overline{\widehat\lambda}
V_{C(\widehat\lambda;z);z}^*)$ from the left and by
$(E_H - \widehat\lambda V_{C(\widehat\lambda;z);z})$ from the right we get
$$ (E_H - \overline{\widehat\lambda}
V_{C(\widehat\lambda;z);z}^*)(E_H - \widehat\lambda V_{C(\widehat\lambda;z);z}) =
E_H - \widehat\lambda V_{C(\widehat\lambda;z);z}
+
E_H - \overline{\widehat\lambda}
V_{C(\widehat\lambda;z);z}^*. $$
After the multiplication in the left-hand side and simplifying the expression we obtain that
$$ V_{C(\widehat\lambda;z);z}^* V_{C(\widehat\lambda;z);z} = E_H,\qquad \forall z\in \mathbb{D}\backslash\{ 0 \}. $$
On the other hand, by~(\ref{ff3_38}) we may write:
$$ \left( E_H - \overline{\widehat\lambda}
V_{C(\widehat\lambda;z);z}^* \right)^{-1}
= E_H - \left( E_H - \widehat\lambda V_{C(\widehat\lambda;z);z} \right)^{-1} $$
$$ = -\widehat\lambda V_{C(\widehat\lambda;z);z}
\left( E_H - \widehat\lambda V_{C(\widehat\lambda;z);z} \right)^{-1}; $$
$$ V_{C(\widehat\lambda;z);z} = -\frac{1}{ \widehat\lambda }
\left( E_H - \overline{\widehat\lambda}
V_{C(\widehat\lambda;z);z}^* \right)^{-1}
\left( E_H - \widehat\lambda V_{C(\widehat\lambda;z);z} \right). $$
Since the operator $\left( E_H - \widehat\lambda V_{C(\widehat\lambda;z);z} \right)^{-1}$
is defined on the whole $H$ and it is bounded, then we conclude that $R(V_{C(\widehat\lambda;z);z}) = H$.
Thus, the operator $V_{C(\widehat\lambda;z);z}$ is unitary in $H$.
Therefore the corresponding operator $V^+_{z;C(\widehat\lambda;z)} = V_z\oplus
C(\widehat\lambda;z)$ is unitary, as well.
In particular, this means that the operator $C(\widehat\lambda;z)$ is isometric and it maps
$N_z$ on the whole $N_{ \frac{1}{\overline{z}} }$. Since $z$ is an arbitrary point from
$\mathbb{D}\backslash\{ 0 \}$, we obtain that condition~2) from the statement of the theorem is satisfied.

{\it Sufficiency. }
Let conditions~1)-3) be satisfied.
Choose an arbitrary point $\widehat\lambda\in\Delta$, and an arbitrary sequence
$\{ \lambda_k \}_{k=1}^\infty$
points from $\mathbb{D}\cup\Delta$ such that $\lambda_k\rightarrow\widehat\lambda$, as $k\rightarrow\infty$.
Using~(\ref{f1_6_1}) 
we write:
$$ E_H - \lambda_k V_{C(\lambda_k;z_0);z_0} =
\left( 1-\frac{1}{z_0} \right) E_H - \frac{|z_0|^2-1}{z_0} (E_{H} + z_0 V^+_{z_0;C(\lambda_k;z_0)})^{-1} $$
$$ = (1-\lambda_k \overline{z_0})
\left[
E_H + \frac{z_0-\lambda_k}{1-\lambda_k \overline{z_0}} V^+_{z_0;C(\lambda_k;z_0)}
\right]
(E_{H} + z_0 V^+_{z_0;C(\lambda_k;z_0)})^{-1}; $$
$$ E_H - \widehat\lambda V_{C(\widehat\lambda;z_0);z_0} = $$
\begin{equation}
\label{ff3_39}
= (1-\widehat\lambda \overline{z_0})
\left[
E_H + \frac{z_0-\widehat\lambda}{1-\widehat\lambda \overline{z_0}} V^+_{z_0;C(\widehat\lambda;z_0)}
\right]
(E_{H} + z_0 V^+_{z_0;C(\widehat\lambda;z_0)})^{-1}.
\end{equation}
Using~(\ref{ff3_39}) write
$$
E_H + \frac{z_0-\widehat\lambda}{1-\widehat\lambda \overline{z_0}} V^+_{z_0;C(\widehat\lambda;z_0)}
=
\frac{1}{1-\widehat\lambda \overline{z_0}}
\left( E_H - \widehat\lambda V_{C(\widehat\lambda;z_0);z_0} \right)
\left( E_{H} + z_0 V^+_{z_0;C(\widehat\lambda;z_0)} \right). $$
Since the operator $V_{C(\widehat\lambda;z_0);z_0}$ is closed, by condition~3) it follows that
there exists the inverse
$( E_H - \widehat\lambda V_{C(\widehat\lambda;z_0);z_0} )^{-1}$, which is defined on the whole $H$
and bounded.
Therefore there exists
$[ E_H + \frac{z_0-\widehat\lambda}{1-\widehat\lambda \overline{z_0}} V^+_{z_0;C(\widehat\lambda;z_0)} ]^{-1}$,
which is bounded and defined on the whole $H$.
By~(\ref{ff3_39}) we get
$$ \left( E_H - \widehat\lambda V_{C(\widehat\lambda;z_0);z_0} \right)^{-1} $$
\begin{equation}
\label{ff3_41}
=
\frac{1}{1-\widehat\lambda \overline{z_0}}
(E_{H} + z_0 V^+_{z_0;C(\widehat\lambda;z_0)})
\left[
E_H + \frac{z_0-\widehat\lambda}{1-\widehat\lambda \overline{z_0}} V^+_{z_0;C(\widehat\lambda;z_0)}
\right]^{-1}.
\end{equation}
For points $\lambda_k$, which belong to $\Delta$, we can apply the same argument. For points $\lambda_k$
from $\mathbb{D}$ we may apply Lemma~\ref{l3_1_p3_1}. We shall obtain an analogous representation:
$$ \left( E_H - \lambda_k V_{C(\lambda_k;z_0);z_0} \right)^{-1} $$
\begin{equation}
\label{ff3_41_1}
= \frac{1}{1-\lambda_k \overline{z_0}}
(E_{H} + z_0 V^+_{z_0;C(\lambda_k;z_0)})
\left[
E_H + \frac{z_0-\lambda_k}{1-\lambda_k \overline{z_0}} V^+_{z_0;C(\lambda_k;z_0)}
\right]^{-1},\quad k\in \mathbb{N}.
\end{equation}

By condition~1) we get
$$ \left\|
V^+_{z_0;C(\lambda_k;z_0)} - V^+_{z_0;C(\widehat\lambda;z_0)}
\right\|
=
\sup_{h\in H,\ \| h \| = 1} \left\|
(C(\lambda_k;z_0) - C(\widehat\lambda;z_0)) P_{N_{z_0}} h
\right\| $$
$$ \leq
\left\|
C(\lambda_k;z_0) - C(\widehat\lambda;z_0)
\right\| \rightarrow 0,\quad k\rightarrow\infty; $$
\begin{equation}
\label{ff3_42}
u.-\lim_{k\rightarrow\infty} V^+_{z_0;C(\lambda_k;z_0)} = V^+_{z_0;C(\widehat\lambda;z_0)}.
\end{equation}
Let us check that there exists an open neighborhood $U_1(\widehat\lambda)$ of  $\widehat\lambda$, and a number
$K>0$ such that
\begin{equation}
\label{ff3_43}
\left\|
\left[
E_H + \frac{z_0-\lambda_k}{1-\lambda_k \overline{z_0}} V^+_{z_0;C(\lambda_k;z_0)}
\right]^{-1}
\right\| \leq K,\quad
\forall\lambda_k:\ \lambda_k\in U_1(\widehat\lambda).
\end{equation}
The last condition is equivalent to the following condition:
\begin{equation}
\label{ff3_44}
\left\|
\left[
E_H + \frac{z_0-\lambda_k}{1-\lambda_k \overline{z_0}} V^+_{z_0;C(\lambda_k;z_0)}
\right] g
\right\| \geq \frac{1}{K} \| g \|,\quad \forall g\in H,\
\forall\lambda_k:\ \lambda_k\in U_1(\widehat\lambda).
\end{equation}
Suppose to the contrary that condition~(\ref{ff3_44}) is not satisfied.
Choose an arbitrary sequence of open discs
$U^n(\widehat\lambda)$ with centres at $\widehat\lambda$
and radii $\frac{1}{n}$; and set $K_n = n$, $n\in \mathbb{N}$.
Then for each $n\in \mathbb{N}$, there exists an element $g_n\in H$, and
$\lambda_{k_n}\in U^n(\widehat\lambda)$,
$k_n\in \mathbb{N}$, such that:
\begin{equation}
\label{ff3_45}
\left\|
\left[
E_H + \frac{z_0-\lambda_{k_n}}{1-\lambda_{k_n} \overline{z_0}} V^+_{z_0;C(\lambda_{k_n};z_0)}
\right] g_n
\right\| < \frac{1}{n} \| g_n \|.
\end{equation}
It is clear that $g_n$ are all non-zero. Denote $\widehat g_n = \frac{g_n}{\| g_n \|_H}$, $n\in \mathbb{N}$.
Then
\begin{equation}
\label{ff3_46}
\left\|
\left[
E_H + \frac{z_0-\lambda_{k_n}}{1-\lambda_{k_n} \overline{z_0}} V^+_{z_0;C(\lambda_{k_n};z_0)}
\right] \widehat g_n
\right\| < \frac{1}{n}.
\end{equation}
Since $|\lambda_{k_n} - \widehat\lambda| < \frac{1}{n}$, then
$\lim_{n\to\infty} \lambda_{k_n} = \widehat\lambda$.
We may write
$$ \frac{1}{n} >
\left\|
\left[
E_H + \frac{z_0-\widehat\lambda}{1-\widehat\lambda \overline{z_0}} V^+_{z_0;C(\widehat\lambda;z_0)}
\right] \widehat g_n \right. $$
$$\left. +
\left(
\frac{z_0-\lambda_{k_n}}{1-\lambda_{k_n} \overline{z_0}} V^+_{z_0;C(\lambda_{k_n};z_0)}
-
\frac{z_0-\widehat\lambda}{1-\widehat\lambda \overline{z_0}} V^+_{z_0;C(\widehat\lambda;z_0)}
\right)
\widehat g_n
\right\| $$
$$ \geq \left|
\left\|
\left[
E_H + \frac{z_0-\widehat\lambda}{1-\widehat\lambda \overline{z_0}} V^+_{z_0;C(\widehat\lambda;z_0)}
\right] \widehat g_n \right\| \right. $$
$$ \left. -
\left\|
\frac{z_0-\lambda_{k_n}}{1-\lambda_{k_n} \overline{z_0}} V^+_{z_0;C(\lambda_{k_n};z_0)}
-
\frac{z_0-\widehat\lambda}{1-\widehat\lambda \overline{z_0}} V^+_{z_0;C(\widehat\lambda;z_0)}
\right\|
\right|
$$
\begin{equation}
\label{ff3_47}
\geq
L
-
\left\|
\frac{z_0-\lambda_{k_n}}{1-\lambda_{k_n} \overline{z_0}} V^+_{z_0;C(\lambda_{k_n};z_0)}
-
\frac{z_0-\widehat\lambda}{1-\widehat\lambda \overline{z_0}} V^+_{z_0;C(\widehat\lambda;z_0)}
\right\|,\ L>0,
\end{equation}
for sufficiently large $n$,
since the operator $[E_H + \frac{z_0-\widehat\lambda}{1-\widehat\lambda \overline{z_0}} V^+_{z_0;C(\widehat\lambda;z_0)}]$
has a bounded inverse, defined on the whole $H$, and
the norm in the right-hand side tends to zero. Passing to the limit in~(\ref{ff3_47}) as
$n\rightarrow\infty$, we obtain a contradiction.

\noindent
Thus, there exists an open neighborhood $U_1(\widehat\lambda)$ of $\widehat\lambda$, and a number
$K>0$ such that inequality~(\ref{ff3_43}) holds.
We may write:
$$ \left\|
\left[
E_H + \frac{z_0-\lambda_k}{1-\lambda_k \overline{z_0}} V^+_{z_0;C(\lambda_k;z_0)}
\right]^{-1}
-
\left[
E_H + \frac{z_0-\widehat\lambda}{1-\widehat\lambda \overline{z_0}} V^+_{z_0;C(\widehat\lambda;z_0)}
\right]^{-1}
\right\| $$
$$ \leq
\left\|
\left[
E_H + \frac{z_0-\lambda_k}{1-\lambda_k \overline{z_0}} V^+_{z_0;C(\lambda_k;z_0)}
\right]^{-1}
\right\|
\left\|
\frac{z_0-\widehat\lambda}{1-\widehat\lambda \overline{z_0}} V^+_{z_0;C(\widehat\lambda;z_0)}
-
\frac{z_0-\lambda_k}{1-\lambda_k \overline{z_0}} V^+_{z_0;C(\lambda_k;z_0)}
\right\|
$$
$$ * \left\|
\left[
E_H + \frac{z_0-\widehat\lambda}{1-\widehat\lambda \overline{z_0}} V^+_{z_0;C(\widehat\lambda;z_0)}
\right]^{-1}
\right\|. $$
Using~(\ref{ff3_43}) we conclude that
\begin{equation}
\label{ff3_48}
u.-\lim_{k\rightarrow\infty} \left[
E_H + \frac{z_0-\lambda_k}{1-\lambda_k \overline{z_0}} V^+_{z_0;C(\lambda_k;z_0)}
\right]^{-1}
=
\left[
E_H + \frac{z_0-\widehat\lambda}{1-\widehat\lambda \overline{z_0}} V^+_{z_0;C(\widehat\lambda;z_0)}
\right]^{-1}.
\end{equation}
By~(\ref{ff3_41}),(\ref{ff3_41_1}),(\ref{ff3_42}),(\ref{ff3_48}) we conclude that
\begin{equation}
\label{ff3_49}
u.-\lim_{k\rightarrow\infty}
\left( E_H - \lambda_k V_{C(\lambda_k;z_0);z_0} \right)^{-1}
=
\left( E_H - \widehat\lambda V_{C(\widehat\lambda;z_0);z_0} \right)^{-1},
\end{equation}
and therefore
\begin{equation}
\label{ff3_50}
u.-\lim_{\lambda\in \mathbb{D}\cup\Delta,\ \lambda\to\widehat\lambda}
\left( E_H - \lambda V_{C(\lambda;z_0);z_0} \right)^{-1}
=
\left( E_H - \widehat\lambda V_{C(\widehat\lambda;z_0);z_0} \right)^{-1}.
\end{equation}
By the Inin formula for $\lambda\in \mathbb{D}$ we have:
$\left( E_H - \lambda V_{C(\lambda;z_0);z_0} \right)^{-1} = \mathbf{R}_\lambda$.
Thus, relation~(\ref{ff3_50}) shows that
the function $\mathbf{R}_\lambda$,
$\lambda\in \mathbb{D}$, admits a continuation on a set $\mathbb{D}\cup\Delta$, and this continuation
is continuous in the uniform operator topology.

Choose an arbitrary element $h\in H$ and consider the following analytic function:
\begin{equation}
\label{ff3_50_1}
f(\lambda) = f_h(\lambda) = (\mathbf{R}_\lambda h,h),\qquad \lambda\in \mathbb{D}.
\end{equation}
The function $f(\lambda)$ admits a continuous continuation on $\mathbb{D}\cup \Delta$, which has
the following form:
$$ f(\lambda) = \left( \left( E_H - \lambda V_{C(\lambda;z_0);z_0} \right)^{-1} h,h \right),\quad
\lambda\in \mathbb{D}\cup \Delta. $$
Let us check that
\begin{equation}
\label{ff3_51}
\left( E_H - \lambda V_{C(\lambda;z_0);z_0} \right)^{-1}
= E_H -
\left( E_H - \overline{\lambda} V^*_{C(\lambda;z_0);z_0} \right)^{-1},\quad
\forall\lambda\in\Delta.
\end{equation}
Choose an arbitrary number $\lambda\in\Delta$.
By condition~2) we conclude that $V_{C(\lambda;z_0);z_0}$ is unitary.
Then
\begin{equation}
\label{ff3_52}
(E_H - \overline{\lambda}
V_{C(\lambda;z);z}^*)(E_H - \lambda V_{C(\lambda;z);z}) =
E_H - \lambda V_{C(\lambda;z);z}
+
E_H - \overline{\lambda}
V_{C(\lambda;z);z}^*.
\end{equation}
To verify the last relation, it is enough to perform the multiplication in the left-hand side and simplify the obtained expression.
Multiplying relation~(\ref{ff3_52}) from the left by
$\left( E_H - \overline{\lambda} V^*_{C(\lambda;z_0);z_0} \right)^{-1}$
and from the right by $\left( E_H - \lambda V_{C(\lambda;z_0);z_0} \right)^{-1}$ we easily get~(\ref{ff3_51}).

We may write:
$$ \overline{f(\lambda)} = \left( h,\left( E_H - \lambda V_{C(\lambda;z_0);z_0} \right)^{-1} h \right) =
\left( \left( E_H - \overline{\lambda} V^*_{C(\lambda;z_0);z_0} \right)^{-1} h,h \right) $$
$$ = (h,h) - \left(\left( E_H - \lambda V_{C(\lambda;z_0);z_0} \right)^{-1} h,h \right) =
(h,h) - f(\lambda);$$
\begin{equation}
\label{ff3_53}
{\mathrm Re} f(\lambda) = \frac{1}{2} (h,h),\qquad  \lambda\in \Delta.
\end{equation}
Denote $g(\lambda) = if(\lambda) - \frac{i}{2} (h,h)$, $\lambda\in \mathbb{D}\cup\Delta$.
Then ${\mathrm Im} g(\lambda)=0$.
Thus, by the Schwarc principle, $g(\lambda)$
has an analytic continuation
$\widetilde g(\lambda) = \widetilde g_h(\lambda)$
on a set $\mathbb{D}\cup\Delta\cup \mathbb{D}_e$.
Moreover, the following relation holds:
\begin{equation}
\label{ff3_54}
\widetilde g(\lambda) = \overline{ \widetilde g\left( \frac{1}{ \overline{\lambda} } \right) },\qquad
\lambda\in \mathbb{D}_e.
\end{equation}
Then the following function
$$ \widetilde f(\lambda) = \widetilde f_h(\lambda) := \frac{1}{i} \widetilde g(\lambda) + \frac{1}{2}(h,h),\quad
\lambda\in \mathbb{D}\cup\Delta\cup \mathbb{D}_e, $$
is an analytic continuation of $f(\lambda)$. Using~(\ref{ff3_54}) we get
\begin{equation}
\label{ff3_55}
\widetilde f(\lambda) = - \overline{ \widetilde f\left( \frac{1}{ \overline{\lambda} } \right) }
+ (h,h),\qquad
\lambda\in \mathbb{D}_e.
\end{equation}
By~(\ref{ff3_3}) we may write
$$ \widetilde f(\lambda) = - \overline{ \left( \mathbf{R}_{\frac{1}{ \overline{\lambda} }} h,h \right)_H }
+ (h,h)_H =
- \left( h, ( E_H - \mathbf{R}^*_{\lambda} ) h \right)_H
+ (h,h)_H
$$
\begin{equation}
\label{ff3_56}
= (\mathbf{R}_\lambda h,h)_H,\qquad \lambda\in \mathbb{D}_e.
\end{equation}
Set
$$ R_\lambda (h,g) = \frac{1}{4} \left( \widetilde f_{h+g}(\lambda) - \widetilde f_{h-g}(\lambda) +
i \widetilde f_{h+ig}(\lambda) - i \widetilde f_{h-ig}(\lambda) \right), $$
\begin{equation}
\label{ff3_57}
h,g\in H,\quad \lambda\in \mathbb{D}\cup\Delta\cup \mathbb{D}_e.
\end{equation}
Notice that $R_\lambda (h,g)$ is an analytic function of $\lambda$ in $\mathbb{D}\cup\Delta\cup \mathbb{D}_e$.
From~(\ref{ff3_50_1}),(\ref{ff3_56}) we conclude that
\begin{equation}
\label{ff3_58}
R_\lambda (h,g) = (\mathbf{R}_\lambda h,g)_H,\quad
h,g\in H,\quad \lambda\in \mathbb{D}\cup \mathbb{D}_e.
\end{equation}
From~(\ref{ff3_50}) it is seen that
$$ R_{\widehat\lambda} (h,g) =
\lim_{\lambda\in \mathbb{D},\ \lambda\to\widehat\lambda}
(\mathbf{R}_\lambda h,g)_H $$
$$ =
\lim_{\lambda\in \mathbb{D},\ \lambda\to\widehat\lambda}
\left( \left( E_H - \lambda V_{C(\lambda;z_0);z_0} \right)^{-1} h, g\right)_H $$
\begin{equation}
\label{ff3_59}
=
\left( \left( E_H - \widehat\lambda V_{C(\widehat\lambda;z_0);z_0} \right)^{-1} h, g\right)_H,\quad
h,g\in H,\quad \widehat\lambda\in\Delta.
\end{equation}
Therefore the following operator-valued function
$$ T_\lambda = \left\{
\begin{array}{cc} \mathbf{R}_\lambda, & \lambda\in \mathbb{D}\cup \mathbb{D}_e \\
\left( E_H - \lambda V_{C(\lambda;z_0);z_0} \right)^{-1} & \lambda\in\Delta
\end{array}\right., $$
is a continuation of $\mathbf{R}_\lambda$, and it is analytic with respect to the weak operator topology,
and therefore with respect to the uniform operator topology, as well.
$\Box$

\begin{cor}
\label{c3_1_p3_1}
Theorem~\ref{t3_2_p3_1} remains valid for the choice $z_0=0$.
\end{cor}
{\bf Proof. }
Let $V$ be a closed isometric operator in a Hilbert space $H$, and $\mathbf{R}_z(V)$
be an arbitrary generalized resolvent of $V$.
Let $F(\lambda) = C(\lambda;0)\in
\mathcal{S}(N_{0};N_{\infty})$
correspond to $\mathbf{R}_z(V)$ by Inin's formula~(\ref{f1_7}) 
for $z_0=0$, whic in this case becomes the Chumakin formula~(\ref{f2_1_p1_1}). 
Consider an arbitrary open arc $\Delta$ of $\mathbb{T}$.

Choose an arbitrary point $z_0\in \mathbb{D}\backslash\{ 0 \}$.
Consider the following isometric operator:
$$ \mathbf{V} = (V+\overline{z_0}E_H)( E_H + z_0 V )^{-1},\quad D(\mathbf{V})=(E_H+z_0 V)D(V). $$
Then
$$ V = (\mathbf{V}-\overline{z_0}E_H)(E_H - z_0 \mathbf{V})^{-1} = \mathbf{V}_{z_0}. $$
Recall that the generalized resolvents $\mathbf{V}$ and $\mathbf{V}_{z_0}$ are related by~(\ref{f2_2}) 
and this correspondence is one-to-one.
Let $\mathbf{R}_z(\mathbf{V})$ be the generalized resolvent which by~(\ref{f2_2}) 
corresponds to the generalized resolvent $\mathbf{R}_z(\mathbf{V}_{z_0})=\mathbf{R}_z(V)$.

From~(\ref{f2_2}) 
weee that $\mathbf{R}_{\widetilde t}(\mathbf{V}_{z_0})$ has a limit as
$\widetilde t\rightarrow \widetilde{t}_0\in\Delta$, if and only if
$\mathbf{R}_{\widetilde u}(\mathbf{V})$ has a limit as
$\widetilde u\rightarrow \widetilde{u}_0\in\Delta_1$, where
$$ \Delta_1 = \left\{ \widetilde u:\ \widetilde u = \frac{\widetilde t + z_0}
{1+\overline{z_0} \widetilde t},\ \widetilde t\in\Delta \right\}. $$
Thus, $\mathbf{R}_{\widetilde t}(\mathbf{V}_{z_0})$
admits a continuation by the continuity on
$\mathbb{T}_e\cup\Delta$, if and only if $\mathbf{R}_{\widetilde u}(\mathbf{V})$
admits a continuation by the continuity on
$\mathbb{T}_e\cup\Delta_1$.
The limit values are connected by~(\ref{f2_2}), as well. 
By~(\ref{f2_2}) 
we see that the continuation $\mathbf{R}_{\widetilde t}(\mathbf{V}_{z_0})$ is analytic
if and only if
when the continuation
$\mathbf{R}_{\widetilde u}(\mathbf{V})$ is analytic.
Thus, $\mathbf{R}_{\widetilde t}(V) = \mathbf{R}_{\widetilde t}(\mathbf{V}_{z_0})$
admits an analytic continuation on
$\mathbb{T}_e\cup\Delta$, if and only if  $\mathbf{R}_{\widetilde u}(\mathbf{V})$
admits an analytic continuation on $\mathbb{T}_e\cup\Delta_1$.

By Theorem~\ref{t3_2_p3_1}, $\mathbf{R}_{\widetilde u}(\mathbf{V})$
admits an analytic continuation on $\mathbb{T}_e\cup\Delta_1$ if and only if

\begin{itemize}
\item[1)] $C(\lambda;z_0)$ has a continuation on a set $\mathbb{D}\cup\Delta_1$, which
is continuous in the uniform operator topology;

\item[2)] The continued function $C(\lambda;z_0)$ maps isometrically $N_{z_0}$ on the whole
$N_{ \frac{1}{ \overline{z_0} } }$, for all points $\lambda\in\Delta_1$;

\item[3)] The operator $(E_H - \lambda \mathbf{V}_{C(\lambda;z_0);z_0})^{-1}$ exists and it is defined on the whole
$H$, for all points $\lambda\in\Delta_1$,
\end{itemize}
where $C(\lambda;z_0)\in \mathcal{S}(N_{z_0};N_{\frac{1}{ \overline{z_0} }})$
corresponds to $\mathbf{R}_z(\mathbf{V})$ by Inin's formula.
Recall that $C(\lambda;z_0)$ is connected with $F(\widetilde t)$ in the following way:
$$ C(\widetilde u) = F\left( \frac{\widetilde u - z_0}{1-\overline{z_0} \widetilde u} \right),\quad
u\in \mathbb{T}_e. $$
Using this relation we easily obtain that condition~1) is equivalent to the following condition:

1') $F(\widetilde t)$ has a continuation on a set $\mathbb{D}\cup\Delta$,
which is continuous in the uniform operator topology;

and condition~2) is equivalent to the following condition:

2') The continued function $F(\widetilde t)$ maps isometrically $N_0$ on the whole
$N_{\infty}$, for all $\widetilde t \in\Delta$.

From~(\ref{ff3_39}) it follows that the operator
$(E_H - \lambda \mathbf{V}_{C(\lambda;z_0);z_0})^{-1}$ exists and it is defined on the whole
$H$,  for all $\lambda\in\Delta_1$, if and only if the operator
$$ \left[ E_H - \frac{\lambda - z_0}{1-\lambda\overline{z_0}} \left( \mathbf{V}_{z_0}\oplus C(\lambda;z_0) \right)
\right]^{-1} =
\left[ E_H - \widetilde t \left( \mathbf{V}_{z_0}\oplus F(\widetilde t) \right)
\right]^{-1}, $$
exists and it is defined on the whole $H$,  for all $\widetilde t\in\Delta$.
$\Box$

\subsection{Some decompositions of a Hilbert space in a direct sum of subspaces.}\label{section3_2}

In the sequel we shall use some decompositions of a Hilbert space related to
a given isometric operator.
\begin{thm}
\label{t4_1_p3_1}
Let $V$ be a closed isometric operator in a Hilbert space $H$.
Let $\zeta\in \mathbb{T}$,  and $\zeta^{-1}$ is a point of the regular type of $V$.
Then
$$ D(V) \dotplus N_\zeta(V) = H; $$
$$ R(V) \dotplus N_\zeta(V) = H. $$
\end{thm}
{\bf Proof. }
At first, assume that $\zeta = 1$.
Check that
\begin{equation}
\label{ff4_3}
\| Vf + g \|_H = \| f+g \|_H,\qquad f\in D(V),\ g\in N_1.
\end{equation}
In fact, we may write:
$$ \| Vf + g \|_H^2 = (Vf,Vf)_H + (Vf,g)_H + (g,Vf)_H + (g,g)_H $$
$$ = \| f \|_H^2 + (Vf,g)_H + (g,Vf)_H + \|g\|_H^2. $$
Since $0 = ((E-V)f, g)_H = (f,g)_H - (Vf,g)_H$, we get
\begin{equation}
\label{ff4_4}
(Vf,g)=(f,g),\qquad f\in D(V),\ g\in N_1,
\end{equation}
and therefore
$$ \| Vf + g \|_H^2 =
\| f \|_H^2 + (f,g)_H + (g,f)_H + \|g\|_H^2 = \| f+g \|_H^2. $$
Consider the following operator:
$$ U(f+g) = Vf+g,\qquad f\in D(V),\ g\in N_1. $$
Let us check that this operator is well-defined on $D(U)= D(V) + N_1$.
Suppose that an element $h\in D(U)$ admits two representations:
$$ h = f_1+g_1 = f_2+g_2,\qquad f_1,f_2\in D(V),\ g_1,g_2\in N_1. $$
By~(\ref{ff4_3}) we may write:
$$ \| Vf_1+g_1 - (Vf_2+g_2) \|_H^2 =  \| V(f_1-f_2) + (g_1-g_2) \|_H^2 =
\| f_1-f_2 + g_1-g_2 \|_H^2 = 0. $$
Thus, $U$ is well-defined.
Moreover, it is easily seen that $U$ is linear. Using~(\ref{ff4_4}) we write:
$$ (U(f+g),U(\widetilde f + \widetilde g)) =
(Vf+g,V\widetilde f + \widetilde g) = (Vf,V\widetilde f) +
(Vf,\widetilde g) + (g,V\widetilde f) + (g,\widetilde g) $$
$$ = (f,\widetilde f) + (f,\widetilde g) + (g,\widetilde f) + (g,\widetilde g) =
(f+g,\widetilde f + \widetilde g), $$
for arbitrary $f,\widetilde f\in D(V)$, $g,\widetilde g\in N_1$.
Therefore $U$ is isometric.

Suppose that there exists an element $h\in H$, $h\not= 0$, $h\in D(V)\cap N_1$. Then
$$ 0 = U 0 = U( h + (-h)) = Vh - h = (V-E_H)h, $$
what contradicts to the fact that the point $1$ is a point of the regular type of $V$.
Therefore
\begin{equation}
\label{ff4_6}
D(V)\cap N_1 = \{ 0 \}.
\end{equation}
Observe that we {\it a priori} do not know, if $D(U)$ is a closed manifold.
Consider the following operator:
$$ W = U|_S, $$
where the manifold $S$ is given by the following equality:
$$ S = \left\{ h\in D(U):\ h\perp N_1 \right\} = D(U)\cap M_1. $$
Thus, $W$ is an isometric operator with the domain $D(W) = D(U)\cap M_1$.
Choose an arbitrary element $g\in D(U)$. Let
$g = g_{M_1} + g_{N_1}$, $g_{M_1}\in M_1$, $g_{N_1}\in N_1\subseteq D(U)$.
Then $g_{M_1} = P^H_{M_1} g \in D(U)$, $g_{M_1}\perp N_1$. Therefore $g_{M_1}\in D(W)$;
$$ P^H_{M_1} D(U) \subseteq D(W). $$
On the other hand, choose an arbitrary element $h\in D(W)$. Then $h\in D(U)\cap M_1$, and therefore
$h = P^H_{M_1} h\in P^H_{M_1} D(U)$. Consequently, we get
$$ D(W) = P^H_{M_1} D(U) = P^H_{M_1} (D(V)+N_1) = P^H_{M_1}D(V) \subseteq M_1. $$
Choose an arbitrary element $f\in D(V)$. Let $f= f_{M_1} + f_{N_1}$, $f_{M_1}\in M_1$, $f_{N_1}\in N_1$.
Then $f-f_{N_1}\in D(U)$, and $f-f_{N_1}\perp N_1$, i.e. $f-f_{N_1}\in D(W)$.
We may write
$$ (W-E_H)(f-f_{N_1}) = (U-E_H)(f-f_{N_1}) = U(f-f_{N_1}) - f + f_{N_1} = Vf - f; $$
\begin{equation}
\label{ff4_9}
(W-E_H) D(W) \supseteq (V-E_H)D(V) = M_1.
\end{equation}
On the other hand, choose an arbitrary element $w\in D(W)$, $w = w_{D(V)} + w_{N_1}$, $w_{D(V)}\in D(V)$,
$w_{N_1}\in N_1$.
Since $w\perp N_1$, then
\begin{equation}
\label{ff4_9_1}
w = P^H_{M_1} w = P^H_{M_1} w_{D(V)}.
\end{equation}
We may write
$$ (W-E_H) w = Uw - w = Vw_{D(V)} + w_{N_1} -w_{D(V)} - w_{N_1} = (V-E_H) w_{D(V)}; $$
$$ (W-E_H) D(W) \subseteq (V-E_H) D(V) = M_1. $$
From the latter relation and relation~(\ref{ff4_9}) it follows that
\begin{equation}
\label{ff4_10}
(W-E_H) D(W) = (V-E_H)D(V) = M_1.
\end{equation}
Moreover, if $(W-E_H) w = 0$, then $(V-E_H) w_{D(V)}=0$; and therefore $w_{D(V)}=0$. By~(\ref{ff4_9_1})
it implies $w=0$. Consequently, there exists the inverse $(W-E_H)^{-1}$. Using~(\ref{ff4_10})
we get
$$ D(W) = (W-E_H)^{-1} M_1. $$
Since $D(W)\subseteq M_1$, then by~(\ref{ff4_10}) we may write:
$$ W w = (W-E_H)w + w\in M_1; $$
\begin{equation}
\label{ff4_12}
D(W)\subseteq M_1,\quad W D(W) \subseteq M_1.
\end{equation}
Consider the closure $\overline{W}$ of $W$ with the domain $\overline{D(W)}$.
By~(\ref{ff4_12}) we see that
$$ D(\overline{W})\subseteq M_1,\quad \overline{W} D(\overline{W}) \subseteq M_1. $$
Then
$$ (\overline{W} - E_H) D(\overline{W}) \subseteq M_1. $$
On the other hand, by~(\ref{ff4_10}) we get
$$ (\overline{W}  - E_H) D(\overline{W}) \supseteq (\overline{W}  - E_H) D(W) =
(W  - E_H) D(W) = M_1. $$
We conclude that
\begin{equation}
\label{ff4_14}
(\overline{W} - E_H) D(\overline{W}) = M_1.
\end{equation}
Let us check that there exists the inverse $(\overline{W} - E_H)^{-1}$.
Suppose to the contrary that there exists an element $h\in D(\overline{W})$, $h\not=0$,
such that $(\overline{W} - E_H)h = 0$.
By the definition of the closure there exists a sequence of elements $h_n\in D(W)$, $n\in \mathbb{N}$,
tending to $h$, and
$W h_n\rightarrow \overline{W}h$, as $n\rightarrow\infty$.
Then $(W-E_H)h_n \rightarrow \overline{W}h - h = 0$, as $n\rightarrow\infty$.
Let $h_n = h_{1;n} + h_{2;n}$, where $h_{1;n}\in D(V)$, $h_{2;n}\in N_1$, $n\in \mathbb{N}$.
Then $(W-E_H)h_n = U h_n - h_n = Vh_{1;n} + h_{2;n} - h_{1;n} - h_{2;n} = (V-E_H) h_{1;n}$.
Therefore $(V-E_H) h_{1;n} \rightarrow 0$, as $n\rightarrow\infty$.
Since $(V-E_H)$ has a bounded inverse, then $h_{1;n} \rightarrow 0$, as $n\rightarrow\infty$.
Then $h_n = P^H_{M_1} h_n = P^H_{M_1} h_{1;n} \rightarrow 0$, as $n\rightarrow\infty$.
Consequently, we get $h=0$, what contradicts to our assumption.

Thus, there exists the inverse $(\overline{W} - E_H)^{-1}\supseteq (W - E_H)^{-1}$.
By~(\ref{ff4_14}),(\ref{ff4_10}), we get
$$ (\overline{W} - E_H)^{-1} = (W - E_H)^{-1}, $$
and therefore
$$ \overline{W} = W. $$
hus, $W$ may be considered as a closed isometric operator in a Hilbert space $M_1$.
The operator $(W - E_H)^{-1}$ is closed and it is defined on the whole $M_1$.
Therefore $(W - E_H)^{-1}$ is bounded. This means that the point $1$
is a regular point of $W$.
Therefore $W$ is a unitary operator in $M_1$.
In particular, this implies that $D(W)=R(W)=M_1$.

By the definition of $W$ we get $D(W) = M_1\subseteq D(U)$, and $U|_{M_1} = W$.
On the other hand, we have $U|_{N_1} = E_{N_1}$. Therefore $D(U)=H$ and
$$ U = W\oplus E_{N_1}. $$
So, $U$ is a unitary operator. Therefore $D(U)=R(U)=H$, what implies
\begin{equation}
\label{ff4_15}
D(V)+N_1=H,\quad R(V)+N_1=H.
\end{equation}
The first sum is direct by~(\ref{ff4_6}).
Suppose that $h\in R(V)\cap N_1$.  Then $h = Vf$, $f\in D(V)$, and we may write:
$$ 0 = Vf + (-h) = U(f+(-h)). $$
Since $U$ is unitary, we get $f=h=Vf$, $(V-E_H)f=0$, and therefore $f=0$, and $h=0$.
Thus, the second sum in~(\ref{ff4_15}) is direct, as well.
So, we proved the theorem in the case $\zeta = 1$.

In the general case we can apply the proved part of the theorem to the operator
$\widehat V := \zeta V$.
$\Box$

\begin{cor}
\label{c4_1_p3_1}
In the conditions of Theorem~\ref{t4_1_p3_1} the following decompositions hold:
$$ \overline{ (H\ominus D(V)) \dotplus M_\zeta } = H; $$
$$ \overline{ (H\ominus R(V)) \dotplus M_\zeta } = H. $$
\end{cor}
{\bf Proof. }
The proof is based on the following lemma.
\begin{lem}
\label{l4_1_p3_1}
Let $M_1$ and $M_2$ be two subspaces of a Hilbert space $H$ such that
$M_1\cap M_2 = \{ 0 \}$, and
$$ M_1 \dotplus M_2 = H. $$
Then
$$ \overline{ (H\ominus M_1) \dotplus (H\ominus M_2) } = H. $$
\end{lem}
{\bf Proof. }
Suppose that an element $h\in H$ is such that $h\in ((H\ominus M_1) \cap (H\ominus M_2))$. Then
$h\perp M_1$, $h\perp M_2$, and therefore $h\perp (M_1+M_2)$, $h\perp H$, $h=0$.

Suppose that an element $g\in H$ satisfies the following condition: $g\perp ((H\ominus M_1) \dotplus (H\ominus M_2))$.
Then
$g\in M_1$, $g\in M_2$, and therefore $g=0$.
$\Box$

Applying the last lemma with $M_1 = D(V)$, $M_2 = N_\zeta$, and $M_1 = R(V)$, $M_2 = N_\zeta$,
we complete the proof of the corollary.
$\Box$

\subsection{Isometric operators with lacunas in a spectrum.}\label{section3_3}

Let $V$ be a closed isometric operator in a Hilbert space $H$.
An open arc $\Delta$ of $\mathbb{T}$ is said to be {\bf a gap in the spectrum
of the isometric operator $V$}, if all points of $\Delta$ are points of the regular type of $V$. 
Above we considered conditions when for a prescribed open arc
$\Delta\subseteq \mathbb{T}$ and a spectrak measure
$\mathbf{F}(\delta)$ there holds: $\mathbf{F}(\Delta) = 0$. We shall see later that a necessary condition
of the existence at least one such spectral function is that $\Delta$ is a gap in the spectrum of $V$.

\noindent
Theorem~\ref{t3_2_p3_1} provides conditions, extracting parameters $C(\lambda;z_0)$ in the Inin formula,
which generate generalized resolvents (and therefore  the corresponding spectral functions),
which have an analytic continuation on $\mathbb{D}\cup\mathbb{D}_e\cup \overline{\Delta}$
(we can apply this theorem with $\overline{\Delta}$ instesd of $\Delta$).
By Proposition~\ref{p3_1_p3_1} for these generalized resolvents there corresond
spectral functions $\mathbf{F}(\delta)$ such that: $\mathbf{F}(\Delta) = 0$.
Moreover,  the extracted class of parameters $C(\lambda;z_0)$ by condition~3) of Theorem~\ref{t3_2_p3_1}, depends on
the operator $V$. Our aim here will be to find instead of condition~3) other conditions, which will not use $V$.

\begin{lem}
\label{l5_1_p3_1}
Let $V$ be a closed isometric operator in a Hilbert space $H$, and $\zeta\in \mathbb{T}$.
Then
\begin{equation}
\label{ff5_1}
V P^H_{M_0(V)} f = \zeta^{-1} P^H_{M_\infty(V)} f,\qquad \forall f\in N_\zeta(V),
\end{equation}
and therefore
\begin{equation}
\label{ff5_2}
\left\| P^H_{M_0(V)} f \right\| =
\left\| P^H_{M_\infty(V)} f \right\|,\qquad \forall f\in N_\zeta(V);
\end{equation}
\begin{equation}
\label{ff5_3}
\left\| P^H_{N_0(V)} f \right\| =
\left\| P^H_{N_\infty(V)} f \right\|,\qquad \forall f\in N_\zeta(V).
\end{equation}
\end{lem}
{\bf Proof. } Choose an arbitrary element $f\in N_\zeta$.
For an arbitrary element $u\in D(V)=M_0$ we may write:
$$ \left( \zeta^{-1} f - V P^H_{D(V)} f, Vu \right)_H
= \zeta^{-1} (f,Vu)_H - \left( P^H_{D(V)} f, u \right)_H $$
$$ = (f,\zeta V u)_H - (f,u)_H = \left( f, (\zeta V - E_H) u \right)_H = 0. $$
Therefore $( \zeta^{-1} f - V P^H_{M_0} f )\perp M_\infty$.
Applying $P^H_{M_\infty}$ to this element we get~(\ref{ff5_1}).
Relation~(\ref{ff5_2}) is obvious, since $V$ is isometric. Then it easily follows~(\ref{ff5_3}).
$\Box$

\begin{lem}
\label{l5_2_p3_1}
Let $V$ be a closed isometric operator in a Hilbert space $H$, and
$C$ be a linear bounded operator in $H$, $D(C) = N_0(V)$, $R(C)\subseteq N_\infty(V)$.
Let $\zeta\in \mathbb{T}$ be such number that $\zeta^{-1}$ be an eigenvalue
of the operator $V^+_{0;C} = V\oplus C$.
If $f\in H$, $f\not= 0$, is an eigenvector of $V^+_{0;C}$, corresponding to $\zeta^{-1}$, then
$f\in N_\zeta(V)$, and
\begin{equation}
\label{ff5_4}
C P^H_{N_0(V)} f = \zeta^{-1} P^H_{N_\infty(V)} f.
\end{equation}
\end{lem}
{\bf Proof. }
LetП $f$ be an eigenvalue of $V^+_{0;C}$, corresponding to the eigenvalue $\zeta^{-1}\in \mathbb{T}$:
$$ (V\oplus C) f = V P^H_{M_0} f + C P^H_{N_0} f = \zeta^{-1} \left( P^H_{M_\infty} f + P^H_{N_\infty} f
\right). $$
By the orthogonality of summands, the last relation is equivalent to the following two conditions:
\begin{equation}
\label{ff5_5}
V P^H_{M_0} f = \zeta^{-1} P^H_{M_\infty} f;
\end{equation}
\begin{equation}
\label{ff5_6}
C P^H_{N_0} f = \zeta^{-1} P^H_{N_\infty} f,
\end{equation}
and~(\ref{ff5_4}) follows.
Relation~(\ref{ff5_5}) implies
$P^H_{M_\infty} (\zeta^{-1} f - V P^H_{M_0} f) = 0$;
$( \zeta^{-1} f - V P^H_{M_0} f )\perp M_\infty$.
For an arbitrary element $u\in D(V)$ we may write:
$$ 0 = \left( \zeta^{-1} f - V P^H_{M_0} f , Vu \right)_H
= \zeta^{-1} (f,Vu)_H -
\left( P^H_{M_0} f , u \right)_H $$
$$ = (f,\zeta Vu)_H - (f,u)_H
= (f, (\zeta V-E_H)u )_H. $$
Therefore $f\in N_\zeta$.
$\Box$

For an arbitrary number $\zeta\in \mathbb{T}$, define an operator $W_\zeta$ in the following way:
\begin{equation}
\label{ff5_7}
W_\zeta P^H_{N_0} f = \zeta^{-1} P^H_{N_\infty} f,\qquad f\in N_\zeta,
\end{equation}
with the domain $D(W_\zeta) = P^H_{N_0} N_\zeta$.
Let us check that this definition is correct.
If an element $g\in D(W_\zeta)$ admits two representations:
$g = P^H_{N_0} f_1 = P^H_{N_0} f_2$, $f_1,f_2\in N_\zeta$, then
$P^H_{N_0} (f_1-f_2) = 0$. By~(\ref{ff5_3}) this implies
$P^H_{N_\infty} (f_1-f_2) = 0$, and therefore the definition is correct.

\noindent
Observe that $W_\zeta$ is linear and
$$ \left\| W_\zeta P^H_{N_0} f \right\|_H =
\left\| P^H_{N_\infty} f \right\|_H =
\left\| P^H_{N_0} f \right\|_H. $$
Thus, the operator $W_\zeta$ is isometric. Notice that $R(W_\zeta) = P^H_{N_\infty} N_\zeta$.

\noindent
Set
$$ S = P^H_{N_0}|_{N_\zeta},\quad Q = P^H_{N_\infty}|_{N_\zeta}. $$
In what follows, we shall assume that
{\it $\zeta^{-1}$ is a point of the regular type of $V$}.
Let us check that in this case the operators $S$ and $Q$ are invertible.
Suppose to the contrary that there exists an element
$f\in N_\zeta$, $f\not=0$: $Sf = P^H_{N_0} f = 0$.
Then $f = P^H_{M_0} f \not= 0$. By Theorem~\ref{t4_1_p3_1} we get
$f\in M_0\cap N_\zeta = \{ 0 \}$. We come to a contradiction.

\noindent
In a similar manner, suppose that there exists an element $g\in N_\zeta$, $g\not=0$: $Qf = P^H_{N_\infty} f = 0$.
Then $g = P^H_{M_\infty} g \not= 0$. By Theorem~\ref{t4_1_p3_1} we obtain that
$g\in M_\infty\cap N_\zeta = \{ 0 \}$. We come to a contradiction, as well.

\noindent
By Theorem~\ref{t4_1_p3_1} we get
$$ P^H_{N_0} N_\zeta = N_0;\quad P^H_{N_\infty} N_\zeta = N_\infty. $$
Thus, the operators $S^{-1}$ and $Q^{-1}$ are closed and they are defined on subspaces $N_0$ and $N_\infty$,
respectively.
Therefore $S^{-1}$ and $Q^{-1}$ are bounded.
From~(\ref{ff5_7}) we see that $D(W_\zeta) = N_0$, $R(W_\zeta) = N_\infty$, and
\begin{equation}
\label{ff5_8_1}
W_\zeta = \zeta^{-1} Q S^{-1}.
\end{equation}

\begin{thm}
\label{t5_1_p3_1}
Let $V$  be a closed isometric operator in a Hilbert space  $H$, and
$C$ be a linear bounded operator in $H$, $D(C) = N_0(V)$, $R(C)\subseteq N_\infty(V)$.
Let $\zeta\in \mathbb{T}$, and $\zeta^{-1}$ be a point of the regular type of $V$.
The point $\zeta^{-1}$ is an eigenvalue of $V^+_{0;C} = V\oplus C$,
if and only if the following condition holds:
\begin{equation}
\label{ff5_9}
(C - W_\zeta) g = 0,\qquad g\in N_0(V),\qquad  g\not=0.
\end{equation}
\end{thm}
{\bf Proof. }
{\it Necessity. }
Since $\zeta^{-1}$ is an eigenvalue of $V^+_{0;C} = V\oplus C$,
then by~\ref{l5_2_p3_1} we obtain that there exists
$f\in N_\zeta$, $f\not= 0$, such that
\begin{equation}
\label{ff5_10}
C P^H_{N_0} f = \zeta^{-1} P^H_{N_\infty} f.
\end{equation}
Comparing the last relation with the definition of $W_\zeta$ we see that $C P^H_{N_0} f =
W_\zeta P^H_{N_0} f$. Set $g=P^H_{N_0} f = Sf$. Since $S$ is invertible, then $g\not= 0$.

{\it Sufficiency. } From~(\ref{ff5_9}) we get~(\ref{ff5_10}) with $f:= S^{-1}g$.
By Lemma~\ref{l5_1_p3_1} we see that relations~(\ref{ff5_5}),(\ref{ff5_6}) hold.
This is equivalent, as we have seen before~(\ref{ff5_5}), that
$\zeta^{-1}$ is an eigenvalue of $V^+_{0;C} = V\oplus C$,
corresponding to eigenvector $f$.
$\Box$

\begin{thm}
\label{t5_2_p3_1}
Let $V$   be a closed isometric operator in a Hilbert space   $H$, and
$C$ be a linear bounded operator in $H$, $D(C) = N_0(V)$, $R(C)\subseteq N_\infty(V)$.
Let $\zeta\in \mathbb{T}$, and $\zeta^{-1}$ is a point of the regular type of $V$.
The following relation
\begin{equation}
\label{ff5_11}
R\left(
V^+_{0;C} - \zeta^{-1} E_H
\right) = H,
\end{equation}
holds if and only if the following relations hold:
\begin{equation}
\label{ff5_12}
\left(
C - W_\zeta
\right) N_0(V) = N_\infty(V);
\end{equation}
\begin{equation}
\label{ff5_12_1}
P^H_{M_\infty} M_\zeta(V) = M_\infty(V).
\end{equation}
\end{thm}
{\bf Proof. }
{\it Necessity. }
Choose an arbitrary element $h\in N_\infty$. By~(\ref{ff5_11}) there exists an element $x\in H$ such that
\begin{equation}
\label{ff5_14}
\left( V^+_{0;C} - \zeta^{-1} E_H \right) x = (V\oplus C) x - \zeta^{-1} x = h.
\end{equation}
For an arbitrary $u\in D(V)$ we may write:
$$ (x,(E_H-\zeta V)u)_H = (x,(V^{-1} - \zeta E_H) Vu)_H =
(x,((V^+_{0;C})^* - \zeta E_H) Vu)_H $$
$$ = ((V^+_{0;C} - \zeta^{-1}E_H)x,Vu)_H = (h,Vu)_H = 0, $$
and therefore $x\in N_\zeta$. Set $g=Sx\in N_0$.
Using~(\ref{ff5_8_1}) write:
$$ (C-W_\zeta) g = CSx - W_\zeta Sx = CSx - \zeta^{-1} Q x. $$
Since $h\in N_\infty$, we apply $P^H_{N_\infty}$ to equality~(\ref{ff5_14}) and get
$$ C P^H_{N_0} x - \zeta^{-1} P^H_{N_\infty} x = h; $$
$$ C S x - \zeta^{-1} Q x = h. $$
Therefore
$$ (C-W_\zeta) g = h, $$
and relation~(\ref{ff5_12}) holds.

Choose an arbitrary element $\widehat h\in M_\infty$. By~(\ref{ff5_11})
there exists $\widehat x\in H$ such that
$$ \left( V^+_{0;C} - \zeta^{-1} E_H \right) \widehat x = (V\oplus C) \widehat x - \zeta^{-1} \widehat x
= \widehat h. $$
The last equality is equivalent to the following two equalities, obtained by an application
of projectors $P^H_{M_\infty}$ and $P^H_{N_\infty}$:
\begin{equation}
\label{ff5_16}
V P^H_{M_0} \widehat x - \zeta^{-1} P^H_{M_\infty} \widehat x = \widehat h;
\end{equation}
\begin{equation}
\label{ff5_17}
C P^H_{N_0} \widehat x - \zeta^{-1} P^H_{N_\infty} \widehat x = 0.
\end{equation}
By Theorem~\ref{t4_1_p3_1} we may write:
$$ \widehat x = x_{D(V)} + x_{N_\zeta},\qquad x_{D(V)}\in D(V),\ x_{N_\zeta}\in N_\zeta. $$
Substituting this decomposition in~(\ref{ff5_17}) we get
$$ C P^H_{N_0} x_{N_\zeta} - \zeta^{-1} P^H_{N_\infty} x_{N_\zeta} - \zeta^{-1} P^H_{N_\infty} x_{D(V)} = 0; $$
$$ (C - W_\zeta) P^H_{N_0} x_{N_\zeta} = \zeta^{-1} P^H_{N_\infty} x_{D(V)}. $$
On the other hand, substituting the decomposition in~(\ref{ff5_16}) we obtain that
$$ V x_{D(V)} + V P^H_{M_0} x_{N_\zeta} - \zeta^{-1} P^H_{M_\infty} x_{D(V)} - \zeta^{-1}
P^H_{M_\infty} x_{N_\zeta} = \widehat h; $$
$$ V x_{D(V)}  - \zeta^{-1} P^H_{M_\infty} x_{D(V)} = \widehat h, $$
where we used Lemma~\ref{l5_1_p3_1}. Then
$$ P^H_{M_\infty} ( V - \zeta^{-1} E_H ) x_{D(V)} = \widehat h, $$
and relation~(\ref{ff5_12_1}) follows directly.

{\it Sufficiency. }
Choose an arbitrary element $h\in H$, $h = h_1 + h_2$, $h_1\in M_\infty$, $h_2\in N_\infty$.
By~(\ref{ff5_12}) there exists an element $g\in N_0$ such that
$$ \left( C - W_\zeta \right) g = Cg - W_\zeta g = h_2. $$
Set $x = S^{-1} g\in N_\zeta$. Then
$$ CS x - \zeta^{-1} Q x = h_2; $$
$$ P^H_{N_\infty} (V\oplus C) P^H_{N_0} x  - \zeta^{-1} P^H_{N_\infty} x = h_2; $$
$$ P^H_{N_\infty} \left( V^+_{0;C} x - \zeta^{-1} x \right) = h_2. $$
By Lemma~\ref{l5_1_p3_1} we may write:
$$ P^H_{M_\infty} \left( V^+_{0;C} x - \zeta^{-1} x \right) =
V P^H_{M_0} x - \zeta^{-1} P^H_{M_\infty} x = 0. $$
Therefore
\begin{equation}
\label{ff5_20}
\left( V^+_{0;C} - \zeta^{-1} E_H \right) x = h_2.
\end{equation}
By~(\ref{ff5_12_1}) there exists an element $w\in M_\zeta$ such that
$$ P^H_{M_\infty} w = h_1. $$
Let
$$ w = ( V - \zeta^{-1} E_H ) \widetilde x_{D(V)},\quad \widetilde x_{D(V)}\in D(V). $$
Then
\begin{equation}
\label{ff5_21}
V \widetilde x_{D(V)}  - \zeta^{-1} P^H_{M_\infty} \widetilde x_{D(V)} = h_1.
\end{equation}
By~(\ref{ff5_12}) there exists $r\in N_0$ such that
$$ \left( C - W_\zeta \right) r = \zeta^{-1} P^H_{N_\infty} \widetilde x_{D(V)}. $$
Set
$\widetilde x_{N_\zeta} := S^{-1} r\in N_\zeta$. Then
$$ (C-W_\zeta) P^H_{N_0} \widetilde x_{N_\zeta} =
\zeta^{-1} P^H_{N_\infty} \widetilde x_{D(V)}; $$
\begin{equation}
\label{ff5_22}
C P^H_{N_0} \widetilde x_{N_\zeta} - \zeta^{-1} P^H_{N_\infty} \widetilde x_{N_\zeta}
- \zeta^{-1} P^H_{N_\infty} \widetilde x_{D(V)} = 0.
\end{equation}
Set $\widetilde x = \widetilde x_{D(V)} + \widetilde x_{N_\zeta}$.
By~(\ref{ff5_22}) we get
\begin{equation}
\label{ff5_23}
C P^H_{N_0} \widetilde x - \zeta^{-1} P^H_{N_\infty} \widetilde x = 0.
\end{equation}
Using relation~(\ref{ff5_21}) and Lemma~\ref{l5_1_p3_1} we write:
$$ h_1 = V \widetilde x_{D(V)}  - \zeta^{-1} P^H_{M_\infty} \widetilde x_{D(V)} +
V P^H_{M_0} \widetilde x_{N_\zeta} - \zeta^{-1} P^H_{M_\infty} \widetilde x_{N_\zeta} $$
\begin{equation}
\label{ff5_24}
= V P^H_{M_0} \widetilde x - \zeta^{-1} P^H_{M_\infty} \widetilde x.
\end{equation}
Adding relations~(\ref{ff5_23}) and~(\ref{ff5_24}) we get
\begin{equation}
\label{ff5_25}
(V\oplus C) \widetilde x - \zeta^{-1} \widetilde  x
= \left( V^+_{0;C} - \zeta^{-1} E_H \right) \widetilde  x
= h_1.
\end{equation}
Adding relations~(\ref{ff5_20}) and~(\ref{ff5_25}) we conclude that relation~(\ref{ff5_11})
holds.
$\Box$

\begin{thm}
\label{t5_3_p3_1}
Let $V$    be a closed isometric operator in a Hilbert space $H$, and
$\Delta$ be an open arc of $\mathbb{T}$ such that
$\zeta^{-1}$ is a point of the regular type of $V$, $\forall\zeta\in\Delta$.
Suppose that the following condition holds:
\begin{equation}
\label{ff5_26}
P^H_{M_\infty(V)} M_\zeta(V) = M_\infty(V),\qquad \forall\zeta\in\Delta.
\end{equation}
Let $\mathbf{R}_z = \mathbf{R}_z(V)$  be an arbitrary generalized resolvent of $V$,
and $C(\lambda;0)\in \mathcal{S}(\mathbb{D}; N_{0},N_\infty)$
corresponds to $\mathbf{R}_z(V)$ by Inin's formula~(\ref{f1_7}). 
The operator-valued function $\mathbf{R}_z(V)$ has an analytic continuation
on a set $\mathbb{D}\cup\mathbb{D}_e\cup \Delta$ if and only if
the following conditions hold:

\begin{itemize}
\item[1)] $C(\lambda;0)$  admits a continuation on a set $\mathbb{D}\cup\Delta$ and this continuation
is continuous in the unform operator topology;

\item[2)] The continued function $C(\lambda;0)$ maps isometrically $N_{0}(V)$ on the whole
$N_\infty(V)$, for all $\lambda\in\Delta$;

\item[3)] The operator $C(\lambda;0)-W_\lambda$ is invertible for all $\lambda\in\Delta$, and
\begin{equation}
\label{ff5_27}
(C(\lambda;0)-W_\lambda) N_0(V) = N_\infty(V),\qquad \forall\lambda\in\Delta.
\end{equation}
\end{itemize}
\end{thm}
{\bf Proof. }
{\it Necessity.}
Suppose that $\mathbf{R}_z(V)$ has an analytic continuation
on a set $\mathbb{D}\cup\mathbb{D}_e\cup \Delta$.
By Corollary~\ref{c3_1_p3_1} we conclude that conditions~1) and 2) hold, and the operator
 $(E_H - \lambda V_{C(\lambda;0);0})^{-1} = -\frac{1}{\lambda}
( (V\oplus C(\lambda;0)) - \frac{1}{\lambda} E_H )^{-1}$ exists and it is defined on the whole
$H$, for all $\lambda\in\Delta$.
By Theorem~\ref{t5_1_p3_1} we obtain that the operator
$C(\lambda;0)-W_\lambda$ is invertible for all $\lambda\in\Delta$.
By Theorem~\ref{t5_2_p3_1} we obtain that relation~(\ref{ff5_27}) holds.

{\it Sufficiency. } Suppose that conditions~1)-3) hold.
By Theorem~\ref{t5_2_p3_1} we obtain that $R( (V\oplus C(\lambda;0)) - \frac{1}{\lambda} E_H ) =
R(E_H - \lambda V_{C(\lambda;0);0}) = H$.
By Theorem~\ref{t5_1_p3_1} we see that the operator $(V\oplus C(\lambda;0)) - \frac{1}{\lambda} E_H$
is invertible. By Corollary~\ref{c3_1_p3_1} we obtain that  $\mathbf{R}_z(V)$ admits an analytic continuation
on a set  $\mathbb{D}\cup\mathbb{D}_e\cup \Delta$.
$\Box$

\begin{rmr}
\label{r5_1_p3_1}
By Corollary~\ref{c3_1_p3_1}, if $\mathbf{R}_z(V)$ has an analytic continuation
on a set  $\mathbb{D}\cup\mathbb{D}_e\cup \Delta$, then
$( (V\oplus C(\lambda;0)) - \frac{1}{\lambda} E_H )^{-1}$ exists and it is bounded.
Therefore the point $\lambda^{-1}$, $\lambda\in\Delta$, are points of the regular type of $V$.

\noindent
On the other hand, by Theorem~\ref{t5_2_p3_1}, in this case condition~(\ref{ff5_26}) holds.
Thus, by Proposition~\ref{p3_1_p3_1} condition~(\ref{ff5_26}) and a condition that points
$\lambda^{-1}$, $\lambda\in\Delta$, are points of the regular type of  $V$,
are {\it necessary} for the existence of a spectral function $\mathbf{F}$ of $V$ such that
$\mathbf{F}(\overline{\Delta}) = 0$. Consequently, these conditions
do not imply on the generality of Theorem~\ref{t5_3_p3_1}.
\end{rmr}

Now we shall obtain an analogous result but the corresponding conditions will be
put on the parameter $C(\lambda;z_0)$ of Inin's formula for an arbitrary $z_0\in \mathbb{D}$.
Before to do that, we shall prove the followng simple proposition:
\begin{prop}
\label{p2_1_p3_1}
Let $V$  be a closed isometric operator in a Hilbert space $H$, and $z_0\in \mathbb{D}$ be a fixed number.
For an arbitrary $\zeta\in \mathbb{T}$
the following two conditions are equivalent:
\begin{itemize}
\item[(i)]  $\zeta^{-1}\in \rho_r(V)$;
\item[(ii)] $\frac{1-\zeta\overline{z_0}}{\zeta-z_0}\in \rho_r(V_{z_0})$.
\end{itemize}
\end{prop}
{\bf Proof.}
$(i)\Rightarrow(ii)$.
We may write:
$$ V_{z_0} - \frac{1-\zeta\overline{z_0}}{\zeta-z_0} E_H =
(V- \overline{z_0} E_H)(E_H - z_0 V)^{-1} - \frac{1-\zeta\overline{z_0}}{\zeta-z_0} (E_H - z_0 V)(E_H - z_0 V)^{-1} $$
$$ = \frac{\zeta(1-|z_0|^2)}{\zeta - z_0} (V-\frac{1}{\zeta} E_H)(E_H - z_0 V)^{-1}. $$
The operator in the right-hand side has a bounded inverse, defined on $(V-\zeta^{-1} E_H)D(V)$.

\noindent
$(ii)\Rightarrow(i)$. We write:
$$ V - \frac{1}{\zeta} E_H =
(V_{z_0} + \overline{z_0}E_H) (E_H + z_0 V_{z_0})^{-1} - \frac{1}{\zeta} (E_H + z_0 V_{z_0})
(E_H + z_0 V_{z_0})^{-1} $$
$$ = \frac{\zeta - z_0}{\zeta} \left(
V_{z_0} - \frac{1-\zeta\overline{z_0}}{\zeta-z_0} E_H
\right)
(E_H + z_0 V_{z_0})^{-1}, $$
and therefore the operator on the right has a bounded inverse, defined on 
$(V_{z_0} - \frac{1-\zeta\overline{z_0}}{\zeta-z_0} E_H)D(V_{z_0})$.
$\Box$

\begin{thm}
\label{t5_4_p3_1}
Let $V$ be a closed isometric operator in a Hilbert space $H$, and
$\Delta$ be an open arc of $\mathbb{T}$ such that
$\zeta^{-1}$ is a point of the regular type of $V$,
$\forall\zeta\in\Delta$. Let $z_0\in \mathbb{D}$ be an arbitrary fixed point, and
suppose that the following condition holds:
\begin{equation}
\label{ff5_28}
P^H_{M_{ \frac{1}{ \overline{z_0} } }(V)} M_\frac{z_0 + \zeta}{1+\zeta \overline{z_0}}(V) =
M_{ \frac{1}{ \overline{z_0} } }(V),\qquad \forall\zeta\in\Delta.
\end{equation}
Let $\mathbf{R}_z = \mathbf{R}_z(V)$ be an arbitrary generalized resolvent of  $V$, and
$C(\lambda;z_0)\in \mathcal{S}(\mathbb{D}; N_{z_0},N_{ \frac{1}{ \overline{z_0} } })$
corresponds to $\mathbf{R}_z(V)$ by the Inin formula~(\ref{f1_7}). 
The operator-valued function $\mathbf{R}_z(V)$ has an analytic continuation
on a set $\mathbb{D}\cup\mathbb{D}_e\cup \Delta$ if and only if
the following conditions hold:

\begin{itemize}
\item[1)] $C(\lambda;z_0)$ admits a continuation on a set $\mathbb{D}\cup\Delta$, which
is continuous in the uniform operator topology;

\item[2)] The continued function $C(\lambda;z_0)$ maps isometrically $N_{z_0}(V)$ on the whole
$N_{ \frac{1}{ \overline{z_0} } }(V)$, for all $\lambda\in\Delta$;

\item[3)] The operator $C(\lambda;z_0)- W_{\lambda;z_0}$ is invertible for all $\lambda\in\Delta$, and
$$ (C(\lambda;z_0)- W_{\lambda;z_0}) N_{z_0}(V) = N_{ \frac{1}{ \overline{z_0} } }(V),\qquad
\forall\lambda\in\Delta. $$
Here the operator-valued function $W_{\lambda;z_0}$ is defined by the following equality:
$$ W_{\lambda;z_0} P^H_{N_{z_0}(V)} f = \frac{1-\overline{z_0}\lambda}{ \lambda - z_0}
P^H_{ N_{ \frac{1}{ \overline{z_0} } }(V) } f,\qquad f\in N_\lambda(V),\ \lambda\in \mathbb{T}. $$
\end{itemize}
\end{thm}
{\bf Proof. } At first, we notice that in the case~$z_0=0$ this theorem coincides with
Theorem~\ref{t5_3_p3_1}.
Thus, we can now assume that $z_0\in \mathbb{D}\backslash\{ 0 \}$.

Suppose that $\mathbf{R}_z(V)$ admits an analytic continuation
on a set
$\mathbb{D}\cup\mathbb{D}_e\cup \Delta$.
Recall that the generalized resolvent $\mathbf{R}_z(V)$
is connected with the generalized resolvent $\mathbf{R}_z(V_{z_0})$
of $V_{z_0}$ by~(\ref{f2_2}) 
and this correspondence is one-to-one.
Therefore $\mathbf{R}_z(V_{z_0})$ admits an analytic continuation
on a set $\mathbb{T}_e\cup \Delta_1$, where
$$ \Delta_1 = \left\{ \widetilde t:\ \widetilde t = \frac{\lambda - z_0}{1-\overline{z_0}\lambda},\
\lambda\in\Delta \right\}. $$
By Proposition~\ref{p2_1_p3_1}, points $\widetilde t^{-1}$, $\widetilde t\in\Delta_1$,
are points of the regular type of the operator $V_{z_0}$.
Moreover, relation~(\ref{ff5_26}), written for $V_{z_0}$ with $\zeta\in\Delta_1$,
coincides with relation~(\ref{ff5_28}).
We can apply Theorem~\ref{t5_3_p3_1} to $V_{z_0}$ and an open arc $\Delta_1$.
Then if we rewrite conditions~1)-3) of that theorem in terms of $C(\lambda;z_0)$, using
the one-to-one correspondence between $C(\lambda;z_0)$ for $V$, and $C(\lambda;0)$ for $V_{z_0}$,
we easily obtain conditions~1)-3) of the theorem.

On the other hand, let conditions~1)-3) be satisfied.
Then conditions of Theorem~\ref{t5_3_p3_1} for $V_{z_0}$ hold.
Therefore the following function
$\mathbf{R}_z(V_{z_0})$ admits an analytic continuation on $\mathbb{T}_e\cup \Delta_1$.
Consequently, the function $\mathbf{R}_z(V)$
admits an analytic continuation on $\mathbb{D}\cup\mathbb{D}_e\cup \Delta$.
$\Box$

\subsection{Spectral functions of a symmetric operator having a constant value on an open
real interval.}\label{section3_4}

Now we shall consider symmetric (not necessarily densely defined) operators in a Hilbert space
and obtain simlar results for them as for the isometric operators above.
The next proposition is an analogue of Proposition~\ref{p3_1_p3_1}.

\begin{prop}
\label{p7_1_p3_1}
Let $A$ be a closed symmetric operator in a Hilbert space $H$, and $\mathbf{E}(\delta)$,
$\delta\in \mathfrak{B}(\mathbb{R})$, be its spectral measure.
The following two conditions are equivalent:

\begin{itemize}
\item[(i)]  $\mathbf{E}(\Delta) = 0$, for an open (finite or infinite)
interval $\Delta\subseteq \mathbb{R}$;
\item[(ii)] The generalized resolvent $\mathbf{R}_z(A)$, corresponding to the spectral measure
$\mathbf{E}(\delta)$, admits an analytic continuation on a set
$\mathbb{R}_e\cup\Delta$, for an open (finite or infinite) real interval
$\Delta\subseteq \mathbb{R}$.
\end{itemize}
\end{prop}
{\bf Proof. }
At first, suppose that the interval $\Delta$ is finite.

\noindent
(i)$\Rightarrow$(ii).
In this case relation~(\ref{f3_1_p1_1}) 
has the following form:
\begin{equation}
\label{f7_1_p3_1}
(\mathbf R_\lambda h,g)_H = \int_{\mathbb{R}} \frac{1}{t-\lambda} d(\mathbf{E}(\cdot) h,g)_H =
\int_{\mathbb{R}} \frac{1}{t - \lambda} d(\mathbf{E}_t h,g)_H,\quad \forall h,g\in H.
\end{equation}
Choose an arbitrary number $z_0\in \Delta$. Since the function $\frac{1}{t-z_0}$
is continuous and bounded on $\mathbb{R}\backslash\Delta$, then there exists an integral
$$ I_{z_0}(h,g) := \int_{\mathbb{R}\backslash\Delta} \frac{1}{t-z_0} d(\mathbf{E}(\cdot) h,g)_H. $$
Then
$$ | (\mathbf R_z h,h)_H - I_{z_0}(h,h) | =
|z-z_0| \left|
\int_{\mathbb{R}\backslash\Delta} \frac{1}{(t-z)(t-z_0)} d(\mathbf{E}(\cdot) h,h)_H
\right|
$$
$$ \leq
|z-z_0| \int_{\mathbb{T}\backslash\Delta} \frac{1}{|t-z||t-z_0|} d(\mathbf{E}(\cdot) h,h)_H,\quad
z\in \mathbb{R}_e. $$
There exists a neighborhood $U(z_0)$ of the point $z_0$ such that $|z-t|\geq M_1 > 0$, $\forall t
\in \mathbb{R}\backslash\Delta$, $\forall z\in U(z_0)$.
Therefore the integral in the latter relation is bounded in the neighborhood $U(z_0)$.
We obtain that
$$ (\mathbf R_z h,h)_H \rightarrow I_{z_0}(h,h),\quad z\in\mathbb{R}_e,\ z\rightarrow z_0,\quad \forall h\in H. $$
Using the properties of sesquilinear forms we conclude that
$$ (\mathbf R_z h,g)_H \rightarrow I_{z_0}(h,g),\quad z\in\mathbb{R}_e,\ z\rightarrow z_0,\quad \forall h,g\in H. $$
Set
$$ \mathbf R_{\widetilde z} :=
w.-\lim_{z\in \mathbb{R}_e,\ z\to \widetilde z} \mathbf R_z,\quad \forall \widetilde z\in \Delta. $$
We may write
$$ \left(
\frac{1}{z-z_0}(\mathbf R_z -  \mathbf R_{z_0})h,h
\right)_H =
\int_{\mathbb{R}\backslash\Delta} \frac{1}{(t-z)(t-z_0)} d(\mathbf{E}(\cdot) h,h)_H, $$
$$ z\in U(z_0),\ h\in H. $$
The function under the integral sign is bounded in $U(z_0)$, and it tends to
$\frac{1}{(t-z_0)^2}$. By the Lebesgue theorem we get
$$ \lim_{z\to z_0} \left(
\frac{1}{z-z_0}(\mathbf R_z -  \mathbf R_{z_0})h,h
\right)_H =
\int_{\mathbb{R}\backslash\Delta} \frac{1}{(t-z_0)^2} d(\mathbf{E}(\cdot) h,h)_H; $$
and therefore
$$ \lim_{z\to z_0} \left(
\frac{1}{z-z_0}(\mathbf R_z -  \mathbf R_{z_0})h,g
\right)_H =
\int_{\mathbb{R}\backslash\Delta} \frac{1}{(t-z_0)^2} d(\mathbf{E}(\cdot) h,g)_H, $$
for $h,g\in H$.
Consequently, there exists the derivative of $\mathbf R_z$ at $z=z_0$.

\noindent
(ii)$\Rightarrow$(i).
Choose an arbitrary element $h\in H$, and consider the following function $\sigma_h(t) := (\mathbf{E}_t h,h)_H$,
$t\in \mathbb{R}$, where $\mathbf{E}_t$ is a left-continuous spectral function of $A$,
corresponding to the spectral measure $\mathbf{E}(\delta)$.
Consider an arbitrary interval $[t_1,t_2]\subset\Delta$.
Assume that $t_1$ and $t_2$ are points of the continuity of the function $\mathbf{E}_t$.
By the Stieltjes-Perron inversion formula цу ьфн цкшеу:
$$ \sigma_h(t_2) - \sigma_h(t_1) = \lim_{y\to +0} \int_{t_1}^{t_2}
{\mathrm Im}
\left\{ (\mathbf{R}_{x+iy} h,h)_H \right\} dx. $$
If $\{ z_n \}_{n=1}^\infty$, $z_n\in \mathbb{R}_e$, is an arbitrary sequence, tending to
some real point $x\in [t_1,t_2]$, then
$$ \overline z_n \rightarrow \overline{x} = x,\qquad n\rightarrow\infty. $$
By the analyticity of the generalized resolvent in a neighborhood of $x$, we may write:
$$ \mathbf{R}_{z_k} \rightarrow \mathbf{R}_x,\quad
\mathbf{R}_{z_k}^* \rightarrow \mathbf{R}_x^*,\qquad k\rightarrow\infty. $$
By property~(\ref{f3_1_1_p1_1}) 
we get
$$ \mathbf{R}_{z_k}^* = \mathbf{R}_{\overline{z_k}} \rightarrow \mathbf{R}_x,\qquad k\rightarrow\infty. $$
Therefore
$$ \mathbf{R}_x^* = \mathbf{R}_x,\qquad x\in [t_1,t_2]. $$
Using the last equality we write:
$$ {\mathrm Im} \left\{ (\mathbf{R}_{x+iy} h,h)_H \right\} \rightarrow
{\mathrm Im} \left\{ (\mathbf{R}_{x} h,h)_H \right\} = 0,\qquad y\rightarrow +0,\quad x\in [t_1,t_2].. $$
Since the generalized resolvent is analytic in $[t_1,t_2]$, then the function
${\mathrm Im}
\left\{ (\mathbf{R}_{x+iy} h,h)_H \right\}$ is continuous and bounded on $[t_1,t_2]$.
By the Lebesgue theorem we get
$$ \sigma_h(t_1) = \sigma_h(t_2). $$
Since points were arbitrary it follows that the function
$\sigma_h(t) = (\mathbf{E}_t h,h)_H$ is constant on $\Delta$.
Since $h$ was arbitrary,using properties of sesquilinear forms we obtain that
$E_t$ is constant.

Consider the case of an infinite $\Delta$. In this case we can
represent  $\Delta$ as a countable union (not necessarily disjunct)
finite open real intervals. Applying the proved part of the proposition for the finite intervals
we easily obtain the required statements. We shall also need to use
the $\sigma$-additivity of the orthogonal spectral mesures.
$\Box$

The following auxilliary peoposition is an analog of Proposition~\ref{p2_1_p3_1}.
\begin{prop}
\label{p7_2_p3_1}
Let $A$  be a closed symmetric operator in a Hilbert space  $H$, and $\lambda_0\in \mathbb{R}_e$ be a
 fixed number.
For an arbitrary $\lambda\in \mathbb{R}$
the following two conditions are equivalent:
\begin{itemize}
\item[(i)]  $\lambda\in \rho_r(A)$;
\item[(ii)] $\frac{ \lambda - \overline{\lambda_0} }{ \lambda - \lambda_0 }\in \rho_r(U_{\lambda_0}(A))$,
где $U_{\lambda_0}(A)$ - преобразование Кэли оператора $A$.
\end{itemize}
\end{prop}
{\bf Proof.}
$(i)\Rightarrow(ii)$.
We may write:
$$ U_{\lambda_0}(A) - \frac{ \lambda - \overline{\lambda_0} }{ \lambda - \lambda_0 } E_H =
E_H + (\lambda_0 - \overline{\lambda_0}) (A - \lambda_0 E_H)^{-1} -
\frac{ \lambda - \overline{\lambda_0} }{ \lambda - \lambda_0 } E_H $$
$$ = -\frac{ \lambda_0 - \overline{\lambda_0} }{ \lambda - \lambda_0 }
(A - \lambda E_H) (A - \lambda_0 E_H)^{-1}. $$
The operator on the right-hand side has the bounded inverse:
$$ -\frac{ \lambda - \lambda_0 }{ \lambda_0 - \overline{\lambda_0} }
(A - \lambda_0 E_H)(A - \lambda E_H)^{-1} =
-\frac{ \lambda - \lambda_0 }{ \lambda_0 - \overline{\lambda_0} } \left(
E_H + (\lambda - \lambda_0) (A - \lambda E_H)^{-1}
\right), $$
defined on $\mathcal{M}_\lambda(A)$.

\noindent
$(ii)\Rightarrow(i)$. The following relation holds:
$$ A - \lambda E_H = \lambda_0 E_H + (\lambda_0 - \overline{\lambda_0}) (U_{\lambda_0}(A) - E_H)^{-1}
- \lambda E_H $$
$$ = (\lambda_0 - \lambda)
( U_{\lambda_0}(A) - \frac{ \lambda - \overline{\lambda_0} }{ \lambda - \lambda_0 } E_H )
(U_{\lambda_0}(A) - E_H)^{-1}. $$
The operator on the right has the bounded inverse:
$$ \frac{1}{ \lambda_0 - \lambda }
(U_{\lambda_0}(A) - E_H)
( U_{\lambda_0}(A) - \frac{ \lambda - \overline{\lambda_0} }{ \lambda - \lambda_0 } E_H )^{-1}, $$
defined on
$( U_{\lambda_0}(A) - \frac{ \lambda - \overline{\lambda_0} }{ \lambda - \lambda_0 } E_H )
D(U_{\lambda_0}(A))$.
$\Box$

We shall need one more auxilliary proposition which is an addition to Theorem~\ref{t3_1_p1_1} 
\begin{prop}
\label{p7_3_p3_1}
Let $A$ be a closed symmetric operator in a Hilbert space $H$, and
$z\in \mathbb{R}_e$ be a fixed point.
Let $\Delta\subseteq \mathbb{R}$ be a (finite or infinite) open interval,
and $\mathbf{R}_{s;\lambda}(A)$ be a generalized resolvent of $A$.
The following two conditions are equivalent:

\begin{itemize}
\item[(i)]
The generalized resolvent $\mathbf{R}_{s;\lambda}(A)$ admits an analytic continuation on
a set $\mathbb{R}_e\cup\Delta$;

\item[(ii)]
The generalized resolvent $\mathbf{R}_{u;\zeta}(U_z(A))$, corresponding to the generalized resolvent
$\mathbf{R}_{s;\lambda}(A)$ by the one-to-one correspondence~(\ref{f3_4_p1_1}) 
from Theorem~\ref{t3_1_p1_1}, 
has an analytic continuation on $\mathbb{T}_e\cup\widetilde\Delta$, where
\begin{equation}
\label{f7_1_0_p3_1}
\widetilde\Delta =
\left\{
\zeta\in \mathbb{T}:\ \zeta = \frac{\lambda - z}{ \lambda - \overline{z} },\ \lambda\in\Delta
\right\}.
\end{equation}
\end{itemize}

\noindent
If condition~(i) is satisfied, the the following statements are true:
\begin{equation}
\label{f7_1_1_p3_1}
\mbox{the interval $\Delta$ consists of points of the regular type of $A$};
\end{equation}
\begin{equation}
\label{f7_1_2_p3_1}
P^H_{ \mathcal{M}_{ \overline{z} } (A) } \mathcal{M}_\lambda(A) = \mathcal{M}_{ \overline{z} } (A),\qquad
\forall\lambda\in\Delta,
\end{equation}
Moreover, condition~(\ref{f7_1_1_p3_1}) is equivalent to the following condition:
$$ \mbox{ the interval $\widetilde\Delta$ consists of points $\zeta$ such that $\zeta^{-1}$ is a point } $$
\begin{equation}
\label{f7_1_3_p3_1}
\mbox{
of the regular type of $U_z(A)$};
\end{equation}
and condition~(\ref{f7_1_2_p3_1}) is equivalent to the following condition:
\begin{equation}
\label{f7_1_4_p3_1}
P^H_{ M_\infty (U_z(A)) } M_\zeta (U_z(A)) = M_\infty (U_z(A)),\qquad \forall\zeta\in\widetilde\Delta.
\end{equation}
\end{prop}
{\bf Proof.}
(i)$\Rightarrow$(ii).
Notice that the linear fractional transformation
$\zeta = \frac{\lambda - z}{ \lambda - \overline{z} }$ maps the real line on the unit circle.
Relation~(\ref{f3_4_p1_1}) 
in terms of $\zeta$ reads
\begin{equation}
\label{f7_2_p3_1}
\mathbf{R}_{u;\zeta} (U_z)
= \frac{1}{1-\zeta} E_H +
(z-\overline{z}) \frac{\zeta}{ (1-\zeta)^2 }
\mathbf{R}_{s; \frac{ z - \overline{z}\zeta }{ 1-\zeta} } (A),\qquad
\zeta\in \mathbb{T}_e\backslash\{ 0 \}.
\end{equation}
Choose an arbitrary point $\widehat\zeta\in\widetilde\Delta$. Set
$\widehat\lambda = \frac{ z - \overline{z} \widehat\zeta }{ 1- \widehat\zeta } \in\Delta$.
Consider an arbitrary sequence $\{ \zeta_k \}_{k=1}^\infty$, $\zeta_k\in
\mathbb{T}_e\backslash\{ 0 \}$, tending to $\widehat\zeta$. The sequence
$\{ \lambda_k \}_{k=1}^\infty$, $\lambda_k := \frac{ z - \overline{z}\zeta_k }{ 1-\zeta_k }\in
\mathbb{R}_e\backslash \{ z,\overline{z} \}$,
tends to $\widehat\lambda$. Writing relation~(\ref{f7_2_p3_1}) for elements of this sequence
and passing to the limit we conclude that the generalized resolvent $\mathbf{R}_{u;\zeta} (U_z)$
can be continued by the continuity on $\widetilde\Delta$.
For the limit values relation~(\ref{f7_2_p3_1}) holds.
From the analyticity of thr right-hand side of this relation it follows the analyticity of
the continued generalized resolvent $\mathbf{R}_{u;\zeta} (U_z)$.

(ii)$\Rightarrow$(i). Let us express from~(\ref{f3_4_p1_1}) 
the generalized resolvent of $A$:
\begin{equation}
\label{f7_3_p3_1}
\mathbf{R}_{s;\lambda} (A) =
\frac{z-\overline{z}}{(\lambda - \overline{z})(\lambda - z)}
\left(
\mathbf{R}_{u;\frac{\lambda - z}{\lambda - \overline{z}}} (U_z)
- \frac{\lambda - \overline{z}}{z-\overline{z}} E_H
\right),\qquad
\lambda\in \mathbb{R}_e\backslash\{ z,\overline{z} \}.
\end{equation}
Proceeding in a similar manner as inthe proof of the previous assertion we show that
$\mathbf{R}_{s;\lambda} (A)$ admits a continuation by the continuity on $\Delta$,
and this continuation is analytic.

Let condition~(i) be satisfied. By the proved part condition~(ii) holds, as well.
As it was noticed in Remark~\ref{r5_1_p3_1}, in this case there hold
conditions~(\ref{f7_1_3_p3_1}) and~(\ref{f7_1_4_p3_1}). It is enough to check the equivalence
of relations~(\ref{f7_1_1_p3_1}) and~(\ref{f7_1_3_p3_1}), and also of relations~(\ref{f7_1_2_p3_1})
and~(\ref{f7_1_4_p3_1}). The first equivalence follows from the definition of  $\widetilde\Delta$
and by Proposition~\ref{p7_2_p3_1}.
The second equivalence follows from relation:
$$ M_\infty (U_z(A)) = \mathcal{M}_{ \overline{z} }(A),\quad
M_\zeta (U_z(A)) = \mathcal{M}_{ \frac{ z-\overline{z}\zeta  }{ 1-\zeta  } } (A),\qquad
\zeta\in \widetilde\Delta. $$
$\Box$

\begin{prop}
\label{p7_4_p3_1}
Let $A$ be a closed symmetric operator in a Hilbert space $H$, and
$z\in \mathbb{R}_e$ be an arbitrary fixed point.
Let $\Delta\subseteq \mathbb{R}$ be a (finite or infinite) open interval,
and conditions~(\ref{f7_1_1_p3_1}) and~(\ref{f7_1_2_p3_1}) hold.
Consider an arbitrary generalized resolvent $\mathbf{R}_{s;\lambda}(A)$ of $A$.
Let $\mathbf{R}_{u;\zeta}(U_z(A))$ be generalized resolvent, corresponding to
$\mathbf{R}_{s;\lambda}(A)$ by relation~(\ref{f3_4_p1_1}) 
from Theorem~\ref{t3_1_p1_1}, 
and $C(\lambda;0)\in \mathcal{S}( \mathbb{D}; N_{0}(U_z(A)),N_\infty(U_z(A)) )$
corresponds to  $\mathbf{R}_{u;\zeta}(U_z(A))$ by Inin's formula~(\ref{f1_7}). 

The generalized resolvent $\mathbf{R}_{s;\lambda}(A)$ admits an analytic continuation on
$\mathbb{R}_e\cup\Delta$ if and only if the following conditions hold:

\begin{itemize}
\item[1)] $C(\lambda;0)$ admits a continuation on $\mathbb{D}\cup\widetilde\Delta$ and this continuation
is continuous in the uniform operator topology;

\item[2)] The continued function $C(\lambda;0)$ maps isometrically $N_{0}(U_z(A))$ on the whole
$N_\infty(U_z(A))$, for all $\lambda\in\widetilde\Delta$;

\item[3)] The operator $C(\lambda;0)-W_\lambda$ is invertible for all $\lambda\in\widetilde\Delta$, and
\begin{equation}
\label{f7_5_p3_1}
(C(\lambda;0)-W_\lambda) N_0(U_z(A)) = N_\infty(U_z(A)),\qquad \forall\lambda\in\widetilde\Delta,
\end{equation}
where the operator-valued function $W_\lambda$ is defined for $U_z(A)$ by~(\ref{ff5_7}).
\end{itemize}
Here $\widetilde\Delta$ is from~(\ref{f7_1_0_p3_1}).
\end{prop}
{\bf Proof.}
By Proposition~\ref{p7_3_p3_1} a possibility of an analytic continuation of $\mathbf{R}_{s;\lambda}(A)$  on
 $\mathbb{R}_e\cup\Delta$ is equivalent to a possibility of an analytic continuation of
$\mathbf{R}_{u;\zeta}(U_z(A))$ on $\mathbb{T}_e\cup\widetilde\Delta$.
Using the equivalence of conditions~(\ref{f7_1_1_p3_1}),(\ref{f7_1_2_p3_1})
and conditions~(\ref{f7_1_3_p3_1}),(\ref{f7_1_4_p3_1}) we can apply Theorem~\ref{t5_3_p3_1}
for  $U_z(A)$ and $\widetilde\Delta$, and the required result follows.
$\Box$

Consider an arbitrary symmetric operator $A$ in a Hilbert space $H$.
Fix an arbitrary number $z\in \mathbb{R}_e$. Consider an arbitrary generalized resolvent
$\mathbf{R}_{s;\lambda} (A)$ of $A$ and the corresponding to it by~(\ref{f7_3_p3_1})
the generalized resolvent $\mathbf{R}_{u;\zeta} (U_z)$ of $U_z(A)$.
The linear fractional transformation
$\zeta = \frac{\lambda - z}{ \lambda - \overline{z} }$ maps
$\Pi_z$ on$\mathbb{D}$. Restrict relation~(\ref{f7_3_p3_1})
on $\Pi_z$:
$$ \mathbf{R}_{s;\lambda} (A) =
\frac{z-\overline{z}}{(\lambda - \overline{z})(\lambda - z)}
\left(
\mathbf{R}_{u;\frac{\lambda - z}{\lambda - \overline{z}}} (U_z)
- \frac{\lambda - \overline{z}}{z-\overline{z}} E_H
\right),\qquad
\lambda\in \Pi_z\backslash\{ z \}. $$
For the generalized resolvent $\mathbf{R}_{u;\frac{\lambda - z}{\lambda - \overline{z}}} (U_z)$ we can use
the Chumakin formula, and for the generalized resolvent $\mathbf{R}_{s;\lambda} (A)$ we can apply the Shtraus formula:
$$ \left( A_{F(\lambda),z} - \lambda E_H \right)^{-1} $$
$$ =
\frac{z-\overline{z}}{(\lambda - \overline{z})(\lambda - z)}
\left(
\left[ E_H - \frac{\lambda - z}{\lambda - \overline{z}} (U_z\oplus \Phi_{\frac{\lambda - z}{\lambda - \overline{z}}}
) \right]^{-1}
- \frac{\lambda - \overline{z}}{z-\overline{z}} E_H
\right), $$
$$ \lambda\in \Pi_z\backslash\{ z \}. $$
where
$F(\lambda)$ is a function from $\mathcal{S}_{a;z}(\Pi_{z}; \mathcal{N}_{z}(A),
\mathcal{N}_{\overline{z}}(A))$, а
$\Phi_\zeta$ is a function from
$\mathcal{S}(\mathbb{D};N_0(U_z(A)),N_\infty(U_z(A))) =
\mathcal{S}(\mathbb{D};\mathcal{N}_{z}(A),\mathcal{N}_{\overline{z}}(A))$.
Observe that $\Phi_\zeta = C(\zeta;0)$, where $C(\zeta;0)$ is a parameter from the Inin formula, corresponding to
$\mathbf{R}_{u;\frac{\lambda - z}{\lambda - \overline{z}}} (U_z)$.
After some transformations of the right-hand side we get
$$ \left( A_{F(\lambda),z} - \lambda E_H \right)^{-1} $$
\begin{equation}
\label{f7_6_p3_1}
=
\frac{1}{ \lambda - \overline{z} }
\left(
(U_z\oplus \Phi_{\frac{\lambda - z}{\lambda - \overline{z}}}) - E_H
\right)
\left[ E_H - \frac{\lambda - z}{\lambda - \overline{z}} (U_z\oplus \Phi_{\frac{\lambda - z}{\lambda - \overline{z}}}
) \right]^{-1},\quad  \lambda\in \Pi_z\backslash\{ z \}.
\end{equation}
If for some $h\in H$, $h\not= 0$, we had
$\left(
(U_z\oplus \Phi_{\frac{\lambda - z}{\lambda - \overline{z}}}) - E_H
\right) h = 0$,
then we we would set
$g = \left[ E_H - \frac{\lambda - z}{\lambda - \overline{z}} (U_z\oplus \Phi_{\frac{\lambda - z}{\lambda - \overline{z}}}
) \right] h \not=0$, and from~(\ref{f7_6_p3_1}) we would get
$$ \left( A_{F(\lambda),z} - \lambda E_H \right)^{-1} g = 0, $$
what is impossible. Therefore there exists the inverse
$\left(
(U_z\oplus \Phi_{\frac{\lambda - z}{\lambda - \overline{z}}}) - E_H
\right)^{-1}$.
Then
$$ A_{F(\lambda),z} - \lambda E_H  $$
$$ = (\lambda - \overline{z})
\left[ E_H - \frac{\lambda - z}{\lambda - \overline{z}} (U_z\oplus \Phi_{\frac{\lambda - z}{\lambda - \overline{z}}}
) \right]
\left(
(U_z\oplus \Phi_{\frac{\lambda - z}{\lambda - \overline{z}}}) - E_H
\right)^{-1},\quad  \lambda\in \Pi_z\backslash\{ z \}. $$
After elementary transformations we get
$$ A_{F(\lambda),z} = zE_H + (z-\overline{z}) \left(
(U_z\oplus \Phi_{\frac{\lambda - z}{\lambda - \overline{z}}}) - E_H
\right)^{-1},\qquad \lambda\in \Pi_z\backslash\{ z \}. $$
Recalling the definition of the quasi-self-adjoint extension $A_{F(\lambda),z}$ we write:
$$ \left(
(U_z\oplus F(\lambda)) - E_H
\right)^{-1} =
\left(
(U_z\oplus \Phi_{\frac{\lambda - z}{\lambda - \overline{z}}}) - E_H
\right)^{-1}, $$
$$ U_z\oplus F(\lambda) = U_z\oplus \Phi_{\frac{\lambda - z}{\lambda - \overline{z}}},\qquad
\lambda\in \Pi_z\backslash\{ z \}. $$
Therefore
\begin{equation}
\label{f7_7_p3_1}
F(\lambda) = \Phi_{\frac{\lambda - z}{\lambda - \overline{z}}} =
C\left(\frac{\lambda - z}{\lambda - \overline{z}};0 \right),\qquad
\lambda\in \Pi_z,
\end{equation}
where the equality for $\lambda = z$ follows from the continuity of
$F(\lambda)$ and $\Phi_\zeta$.

\begin{thm}
\label{t7_1_p3_1}
Let $A$ be a closed symmetric operator in a Hilbert space $H$, and
$z\in \mathbb{R}_e$ be an arbitrary fixed point.
Let $\Delta\subseteq \mathbb{R}$ be a (finite or infinite) open interval,
and conditions~(\ref{f7_1_1_p3_1}) and~(\ref{f7_1_2_p3_1}) hold.
Consider an arbitrary  generalized resolvent $\mathbf{R}_{s;\lambda}(A)$ of $A$.
Let $F(\lambda)\in \mathcal{S}_{a;z}(\Pi_{z}; \mathcal{N}_{z}(A),
\mathcal{N}_{\overline{z}}(A))$
corresponds to  $\mathbf{R}_{s;\lambda}(A)$ by the Shtraus formula~(\ref{f11_3_p2_1}). 
The generalized resolvent $\mathbf{R}_{s;\lambda}(A)$ admits an analytic continuation on
 $\mathbb{R}_e\cup\Delta$ if and only if the following conditions hold:

\begin{itemize}
\item[1)] $F(\lambda)$ admits a continuation on$\mathbb{R}_e\cup\Delta$ and this continuation
is continuous in the uniform operator topology;

\item[2)] The continued function $F(\lambda)$ maps isometrically $\mathcal{N}_{z}(A)$ on the whole
$\mathcal{N}_{\overline{z}}(A)$, for all $\lambda\in\Delta$;

\item[3)] The operator $F(\lambda)-\mathcal{W}_\lambda$ is invertible for all $\lambda\in\Delta$, and
\begin{equation}
\label{f7_8_p3_1}
(F(\lambda)-\mathcal{W}_\lambda) \mathcal{N}_{z}(A) = \mathcal{N}_{\overline{z}}(A),\qquad
\forall\lambda\in\Delta,
\end{equation}
where $\mathcal{W}_\lambda = W_{\frac{\lambda - z}{ \lambda - \overline{z} }}$,
$\lambda\in\Delta$, and $W_\zeta$ is defined for $U_z(A)$ by~(\ref{ff5_7}).
\end{itemize}
\end{thm}
{\bf Proof.}
For the generalized resolvent $\mathbf{R}_{s;\lambda}(A)$ by~(\ref{f7_3_p3_1}) it corresponds
a generalized resolvent of the Cayley transformation $\mathbf{R}_{u;\zeta}(U_z(A))$.
Let $C(\zeta;0)\in \mathcal{S}( \mathbb{D}; N_{0}(U_z(A)),N_\infty(U_z(A)) )$ be a parameter from
Inin's formula, corresponding to $\mathbf{R}_{u;\zeta}(U_z(A))$.
Consider an arc $\widetilde\Delta$ from~(\ref{f7_1_0_p3_1}).

\noindent
{\it Necessity. }
Let the generalized resolvent $\mathbf{R}_{s;\lambda}(A)$ admits an analytic continuation on
 $\mathbb{R}_e\cup\Delta$.
By~\ref{p7_4_p3_1} for the parameter
$C(\zeta;0)$ it holds:
\begin{itemize}
\item[(a)] $C(\zeta;0)$ admits a continuation on $\mathbb{D}\cup\widetilde\Delta$ and this continuation
is continuous in the uniform operator topology;

\item[(b)] The continued function $C(\zeta;0)$ maps isometrically $N_{0}(U_z(A))=\mathcal{N}_z(A)$ on the whole
$N_\infty(U_z(A))=\mathcal{N}_{\overline{z}}(A)$, for all $\zeta\in\widetilde\Delta$;

\item[(c)] The operator $C(\zeta;0)-W_\zeta$ is invertible for all $\zeta\in\widetilde\Delta$, and
\begin{equation}
\label{f7_9_p3_1}
(C(\zeta;0)-W_\zeta) N_0(U_z(A)) = N_\infty(U_z(A)),\qquad \forall\zeta\in\widetilde\Delta,
\end{equation}
where  $W_\zeta$ is defined for $U_z(A)$ by~(\ref{ff5_7}).
\end{itemize}
Define a function $F_+(\lambda)$ in the following way:
\begin{equation}
\label{f7_10_p3_1}
F_+(\lambda) =
C\left(\frac{\lambda - z}{\lambda - \overline{z}};0 \right),\qquad
\lambda\in \Pi_z\cup\Delta.
\end{equation}
By~(\ref{f7_7_p3_1}) this function is an extension of the parameter $F(\lambda)$,
corresponding to the generalized resolvent by the Shtraus formula. From condition~(a) it follows that this continuation
is continuous  in the uniform operator topology.
From condition~(b) it follows that $F_+(\lambda)$ maps isometrically $\mathcal{N}_{z}(A)$ on the whole
$\mathcal{N}_{\overline{z}}(A)$, for all $\lambda\in\Delta$.
Consider a function $\mathcal{W}_\lambda = W_{\frac{\lambda - z}{ \lambda - \overline{z} }}$,
$\lambda\in\Delta$. From condition~(c) it follows that
$F_+(\lambda;0)-\mathcal{W}_\lambda$ is invertible for all $\lambda\in\Delta$, and
$$ (F_+(\lambda)-\mathcal{W}_\lambda) \mathcal{N}_{z}(A) = \mathcal{N}_{\overline{z}}(A),\qquad
\forall\lambda\in\Delta. $$

\noindent
{\it Sufficiency. } Suppose that for the parameter $F(\lambda)$, corresponding to the generalized resolvent
$\mathbf{R}_{s;\lambda}(A)$ by the Shtraus formula, conditions~1)-3) hold.
Define a function $C_+(\zeta;0)$ in the following way:
\begin{equation}
\label{f7_11_p3_1}
C_+(\zeta;0) = F\left( \frac{ z - \overline{z}\zeta }{ 1-\zeta } \right),\qquad
\zeta\in \mathbb{D}\cup\widetilde\Delta.
\end{equation}
By~(\ref{f7_7_p3_1}) this function is
a continuation of the parameter $C(\zeta;0)$. From condition~1) it follows that this continuation
is continuous in the uniform operator topology.
From condition~2) it follows that $C_+(\zeta;0)$ maps isometrically $N_{0}(U_z(A))=\mathcal{N}_z(A)$ on the whole
$N_\infty(U_z(A))=\mathcal{N}_{\overline{z}}(A)$, for all $\zeta\in\widetilde\Delta$.
Finally, condition~3) means that
$C_+(\zeta;0)-W_\zeta$ is invertible for all $\zeta\in\widetilde\Delta$, and
$$ (C_+(\zeta;0)-W_\zeta) N_0(U_z(A)) = N_\infty(U_z(A)),\qquad \forall\zeta\in\widetilde\Delta. $$
By Proposition~\ref{p7_4_p3_1} we conclude that
the generalized resolvent $\mathbf{R}_{s;\lambda}(A)$ admits an analytic continuation on
$\mathbb{R}_e\cup\Delta$.
$\Box$

Notice that by~(\ref{ff5_7}) the function $\mathcal{W}_\lambda$ fron the formulation of the last theorem
has the following form:
$$ \mathcal{W}_\lambda P^H_{\mathcal{N}_z(A)} f =
\frac{\lambda - \overline{z}}{\lambda - z} P^H_{\mathcal{N}_{\overline{z}}(A)} f,\qquad
f\in \mathcal{N}_\lambda(A),\ \lambda\in\Delta. $$
Moreover, we emphasize that conditions~(\ref{f7_1_1_p3_1}) and~(\ref{f7_1_2_p3_1}) are
necessary for the existence at least one generalized resolvent of $A$,
which admits an analytic continuation on $\mathbb{R}_e\cup\Delta$.
Consequently, these conditions does not imply on the generality of a description such generalized resolvents in
Theorem~\ref{t7_1_p3_1}.


\section{Formal credits.}

\textbf{ Section~\ref{chapter1}. } \S \ref{section1_1}.
Results of this subsection, except for Propositions, belong to Chumakin, see~\cite{cit_3000_Chumakin}, \cite{cit_4000_Chumakin}.
Propositions~\ref{p1_1_p1_1}-\ref{p1_3_p1_1} belong to Shtraus, see~\cite[footnote on page 83]{cit_900_Shtraus},
\cite[Lemma 1.1, Lemma 1.2]{cit_970_Shtraus}.

\noindent
\S \ref{section1_2}.
Chumakin's formula appeared in papers of Chumakin~\cite{cit_3000_Chumakin} and \cite{cit_4000_Chumakin}.

\noindent
\S \ref{section1_3}.
Inin's formula was obtained in the paper~\cite[Theorem]{cit_5000_Inin}. Our proof differs from the original
proof and it is based directly on the use of Chumakin's formula for the generalized resolvents and
a method which is similar to the metod of Chumakin in the paper~\cite{cit_3500_Chumakin}.

\noindent
\S \ref{section1_4}.
The definition of the generalized resolvent of a symmetric, not necessarily densely defined operator in a Hilbert space
belongs to Shtraus and he established its basic properties, see~\cite{cit_900_Shtraus}, \cite{cit_1500_Shtraus}.
A connection between the generalized resolvents  of a symmetric operator and the generalized resolvent its Cayley transformation 
(Theorem~\ref{t3_1_p1_1})
ina somewhat less general form was established by Chumakin in~\cite{cit_3500_Chumakin}.

\noindent
\S \ref{section1_5}. The Shtraus formula was obtained in~\cite{cit_900_Shtraus}. Our proof differs form the original one
and we follow the idea of Chumakin, proposed in~\cite{cit_3500_Chumakin}.

\noindent
\textbf{ Section~\ref{chapter2}. } This section is based on results of Shtraus from papers~\cite{cit_900_Shtraus},
\cite{cit_950_Shtraus}, \cite{cit_955_Shtraus}, \cite{cit_970_Shtraus}, \cite{cit_1500_Shtraus}, see also references
in these papers.

\noindent
\textbf{ Section~\ref{chapter3}. } \S \ref{section3_1}.
Proposition~\ref{p3_1_p3_1} is an analog of Theorem~3.1 (A),(B) in a paper of McKelvey~\cite{cit_5500_McKelvey},
see also~\cite[Lemma 1.1]{cit_5700_Varlamova_Luks}.
Theorem~\ref{t3_1_p3_1} is an analog of Theorem~2.1 (A) in a paper of McKelvey~\cite{cit_5500_McKelvey}.
Theorem~\ref{t3_2_p3_1} is an analog of Theorem~2.1 (B) in a paper of McKelvey~\cite{cit_5500_McKelvey}.

\noindent
\S \ref{section3_2}. Theorem~\ref{t4_1_p3_1} appeared in~\cite[Lemma 1]{cit_5950_Ryabtseva}.
However her proof was based on a lemma of Shmulyan~\cite[Lemma 3]{cit_5000_Shmulyan}.
As far as we know, no correct proof of this lemma was published.
An attempt to proof the lemma of Shmulyan was performed by L.A.~Shtraus in~\cite[Lemma]{cit_5900_Shtraus_LA}.
Unfortunately, the proof was not complete. We used the idea of L.A.~Shtraus for the proof
a weaker result: Theorem~\ref{t4_1_p3_1}. 

\noindent
\S \ref{section3_3}.
Results of this subsection mostly belong to Ryabtseva, in some cases with our filling of gaps in proofs, corrections and a generalization.
For Lemmas~\ref{l5_1_p3_1}, \ref{l5_2_p3_1} see Lemma~2 and its corollary, and also Lemma~3 in~\cite{cit_5950_Ryabtseva}.
Theorem~\ref{t5_1_p3_1} was obtained in~\cite[Theorem 1]{cit_5950_Ryabtseva}.
Theorem~\ref{t5_2_p3_1} is a corrected version of Theorem~2 in~\cite{cit_5950_Ryabtseva}.
Theorem~\ref{t5_3_p3_1} is a little corrected version of Theorem~4 in~\cite{cit_5950_Ryabtseva}.
For Proposition~\ref{p2_1_p3_1} see~\cite[p. 34]{cit_5000_Inin}.

\noindent
\S \ref{section3_4}.
Proposition~\ref{p7_1_p3_1} is contained in Theorem~3.1 (A),(B) in a paper of McKelvey~\cite{cit_5500_McKelvey},
see also~\cite[Lemma 1.1]{cit_5700_Varlamova_Luks}.
Theorem~\ref{t7_1_p3_1} is close to results of Varlamova-Luks, see~\cite[Theorem]{cit_5900_Varlamova_Luks},
\cite[Theorem~2]{cit_5800_Varlamova_Luks}, \cite[Theorem 3.1]{cit_5700_Varlamova_Luks}.

}
\end{document}